%% file: ArticleVersion.tex
\documentclass[a4paper, 11pt]{report}

\usepackage[T1]{fontenc}
\usepackage[utf8]{inputenc}
\usepackage[main=english, french]{babel}															 

\usepackage{modello}
\usepackage[a4paper,width=160mm,top=25mm,bottom=25mm,bindingoffset=6mm]{geometry}  

\usepackage[autostyle]{csquotes}
\usepackage[style=alphabetic, backend=biber, backref=true]{biblatex}
\DeclareFieldFormat*{title}{\mkbibemph{#1}}         

\newcommand{\Crys}{\infcatname{Crys}}

\hypersetup{bookmarksnumbered=true, colorlinks=true}

\title{On Derived D-Modules and their several definitions}
\author{Carlo Buccisano\\\vspace{3mm} PhD Thesis\\ \vspace{5mm} Université Paris-Saclay -- France}
\date{October 10, 2025}

\begin{document}

	\maketitle

	\begin{abstract}
		The so called theory of derived D-modules is an extension of classical D-modules to derived algebraic geometry, which uses the derived information of the base scheme. \\
		We prove that the three different definitions of derived D-modules, given by Beraldo, Nuiten and Toën-Vezzosi, on a (nice) derived scheme yield equivalent symmetric monoidal $\infty$-categories. We deduce this as a corollary of more general statements about Chevalley-Eilenberg cohomology of dg-Lie algebroids, proving a conjecture by E. Pavia, and about the relation between representations of a dg-Lie algebroid and some class of ind-coherent sheaves on the associated formal moduli problem, which can be of independent interest.	
	\end{abstract}
	
	\pagestyle{plain}
	{
		\hypersetup{linkcolor=black}
		\tableofcontents
	}

	\newpage

	\pagestyle{fancy}
	\fancyhf{}
    \fancyhead[L]{\bfseries \thepage}
    \fancyhead[C]{\bfseries }
    \fancyhead[R]{\bfseries \nouppercase{\leftmark}}

	\setlength{\headheight}{13.6pt}
	
	\selectlanguage{english}

\subfile{introduction}

   \chapter{Preliminaries}
       \label{section:preliminaries}

	   In this chapter, we review some definitions and constructions that will be used later on. 
	   Most of the content is well known, but we decided to record it anyway for the sake of readability of the manuscript. 
	   We advise the expert reader to skip this part on a first read and only come back to it when needed.

	   \subfile{morita_categorical_stuff}

	   \subfile{RecallsDAG}

	   \subfile{IndCoh}

	   \subfile{ProCoherentSheaves}

       \subfile{filtered}

\subfile{filtered_first_definition}

	\chapter{Derived D-modules and dg-Lie algebroids}
		\label{section:second_definition_derived_DModules}

		In this chapter, we will introduce another way of defining derived D-modules, namely by giving a ``derived'' definition of the ring $\D_{ X }$.
		As mentioned in the introduction, in the classical case (where $X$ is a smooth scheme) one has the equivalence of (sheaves of) associative algebras 
		\[
			\D_{ X } \simeq U(\mathcal{T}_{ X }),
		\]
		where the right-hand side is the universal enveloping algebra of the tangent Lie algebroid (that is, the tangent sheaf equipped with the classical Lie bracket of derivations).
		The correct homotopical generalization of the tangent sheaf is the tangent \emph{complex} $\T_{ X }$, defined as the (derived) $\O_{ X }$-linear dual of the cotangent complex.\footnote{Since we always work with locally finitely presented derived schemes, which have dualizable cotangent complex, we don't need to worry about eventual loss of information by the dualization procedure.}
		On a derived affine scheme $\Spec A$, we can interpret $\T_{ A }$ as the $A$-module of (derived) derivations from $A$ to $A$; as in the classical case, this inherits a (dg-)Lie bracket over the base field $\C$.
		\begin{remark}
			\label{remark:lie_bracket_tangent_complex}
			The fact claimed above is well known. 
			Here is one way to see this explicitly: let $A$ be a cofibrant (so in particular quasi-free) connective cdga. Then the dg-$A$-module $\L_{ A }$ is computed as (dg) K{\"a}hler differentials of $A$ and it is free as a graded $A$-module. Therefore, its $A$-linear dual (for which there is no need to derive, and which is first computed as a graded dual and then differential is added) is exactly 
				\[
					\T_{ A } \simeq \mathrm{Der}^{ \dg }_{ \C }(A, A) \subset \mathrm{End}^{ \dg }_{ \C }(A, A)
				\]
				and inherits the Lie bracket from the right-hand side (which is an associative dg-algebra).

		\end{remark}

		There is also a generalization of Lie algebroids to derived algebraic geometry, the so called dg-Lie algebroids, see \cite{Nuiten:KoszulLieAlgbd}. There is a version of the universal enveloping (dg-)algebra, behaving as one would expect.
		Therefore the definition above can be read in a derived setting, giving a second definition of derived left D-modules. \\
		We will then review the third possible definition of derived D-modules introduced in \cite{Ber:DerivedDMod}, and then relate it to representations of dg-Lie algebroids.
	
	    \subfile{lie_algebroids}

		\subfile{BeraldoVsNuiten}

	\chapter{Filtered Chevalley-Eilenberg cohomology of dg-Lie algebroids}
		\label{section:main_result}

		In this chapter, we will prove our main result about the Chevalley-Eilenberg cohomology of a dg-Lie algebroid. Here is its statement.

		\begin{thm}
			Let $A \in \CAlg^{ \cn }_{ \C }$ and consider $\g \in \LieAlgbd_{ A }$ with an underlying perfect $A$-module.
			Then the complete Hodge-filtered Chevalley-Eilenberg $\infty$-functor induces a symmetric monoidal equivalence of $\infty$-categories
				\[
					\CE^{\cpl}(\g, -)\colon \LMod_{ U(\g) } \stackrel{\sim}{\longrightarrow} \Mod^{ \const }_{ \CE^{ \cpl }(\g, A) }.
				\]
		\end{thm}	

		\subfile{construction}

		\subfile{proof}

		\printbibliography[heading=bibintoc]
\end{document}

%% file: introduction.tex
\chapter*{Introduction}
	\addcontentsline{toc}{chapter}{Introduction}

	\section*{Context}

		\subsection*{Derived algebraic geometry}

			The general domain of this thesis is derived algebraic geometry, over a characteristic $0$ field. This is a subject which applies techniques coming from homotopy theory and higher category theory to algebraic geometry. The foundations of the subject have been laid out in \cite{TV:HAG1}, \cite{TV:HAG2} in the beginning of 2000's and, more recently, in \cite{GR:DAG1}, \cite{GR:DAG2} and \cite{Lurie:SAG}.

			One could say that derived algebraic geometry's purpose is to handle ``singular'' situations in classical algebraic geometry, for example non-transverse (and Tor-dependent) intersections or quotient stacks $[X/G]$ by non free actions (although, by now, the number of applications of the subject are way more than just this one). The philosophy adopted by derived algebraic geometry to better tackle those questions is, in a way, conceptually similar to the passage from algebraic varieties to schemes in algebraic geometry. 
			In fact, it is now commonly accepted that schemes are a better framework to work with than classical algebraic varieties, even if one is only interested in the latter ones. One of the many motivating examples comes again from intersection and multiplicities, which are lost when working only with algebraic varieties: it suffices to intersect a parabola with the tangent line at its vertex (which is a non-transverse Tor-independent intersection). Working with schemes, we accept a more abstract definition of the geometric objects of interest (and shift our perspective from the topological space to the functions on it) to have a way better behaved theory.
			The shift to derived algebraic geometry is similar: by redefining, in a further abstract way, what are the geometric objects of interests (and therefore how intersections are performed and tangent spaces are replaced), we obtain a better behaved theory, with cleaner results. To be more specific, we replace commutative rings (affine schemes) with connective $\E_{ \infty }$-ring spectra (derived affine schemes) and classical étale-stacks in groupoids by étale-hypersheaves valued in spaces. These objects and definitions are purely $\infty$-categorical, adding a further layer of homotopical complexity.

			Aside from intersections and (co)tangent complexes, another important area in which derived algebraic geometry sheds further light is (derived) deformation theory, where the use of simplicial rings and cotangent complexes allows us to gain a better understanding of several higher cohomology groups, which lacked a geometric interpretation beforehand.
			For a more extensive motivation and history of derived algebraic geometry we redirect the reader to \cite{Toen:DAG}.

		\subsection*{Algebraic D-modules}

	    	Classical algebraic D-modules have been a well-established theory since the 1970s, mainly born after the ideas of Sato, Bernstein, and then Kashiwara, Deligne, Schapira (and many others). 
			Some classical references are \cite{Bor:AlgDMod}, \cite{Hotta:GreenBookDMod}.
			The initial idea was to tackle analytic questions, like solving systems of linear partial differential equations, on an algebraic variety. Due to the local nature of solutions, the main idea was then to use sheaf theory, as it had been developed few years before in algebraic geometry by Grothendieck. This idea revealed itself very fruitful, yielding celebrated results like the Riemann-Hilbert correspondences (which links D-modules on a smooth complex variety, a purely algebraic data, with perverse sheaves on the analytification of such variety, a purely topological data) and the proof of Kazhdan-Lusztig conjecture. More recently, let us also mention the geometric Langlands conjecture (and its proof) which, on an algebraic curve $X$ over $\C$, proves that (some subcategories) of (twisted) D-modules on the moduli stack of $G$-bundles on $X$ and of ind-coherent sheaves on the moduli stack of $\check{G}$-local systems of $X$ are equivalent. See the series of papers \cite{Gaitsgory:GLC1}, \cite{Gaistgory:GLC2}, \cite{Gaistgory:GLC3}, \cite{Gaitsgory:GLC4} and \cite{Gaitsgory:GLC5}.\\ 
			Here's a brief rundown of how the theory of D-modules is built, in algebraic geometry.
			Given a smooth scheme $X$ over $\C$, one can define, using Zariski-local coordinates, a sheaf of associative algebras $\D_{ X }$, the so called (sheaf of) \emph{differential operators on $X$}. It is an associative ring which contains $\O_{ X }$ (but in which $\O_{ X }$ is not central). Left (right) modules over it will be left (right) D-modules on $X$, having an underlying quasi-coherent $\O_{ X }$-module. One proves that left D-modules can be pulled back along smooth morphisms of smooth schemes (this operation corresponds to the $*$-pullback of the underlying quasi-coherent sheaves) and right D-modules can be pushed forward (although this operation is, in general, only well defined at the level of derived categories and corresponds to the $*$-pushforward of the underlying \emph{ind-coherent} sheaves). Left D-modules admit a symmetric monoidal structure via relative tensor product over $\O_{ X }$ and they are equivalent to right D-modules by tensoring with $\omega_{ X } \coloneqq \Lambda^{ \mathrm{top} }\Omega^{ 1 }_{ X }$. Finally, the Kashiwara lemma proves that pushforward of D-modules along a closed immersion of smooth schemes is fully faithful, and thus allows one to extend the definition of D-modules to any scheme of finite type (by locally choosing a closed embedding into a smooth scheme, considering D-modules with a support condition on the underlying sheaf and then glueing together these categories).

	\section*{What is in this thesis?}
		Let us now discuss the content and contributions of this thesis. We need to start with a story.

		\subsection*{A brief tale about D-modules in derived algebraic geometry}	

			The driving question at the core of this work is:
			\begin{quote}
				\textbf{Question 0:} What are D-modules in derived algebraic geometry?
			\end{quote}
			There are several possible answers (to different interpretations) to the above question. 
			One way of understanding the question above is to find how to generalize the theory of D-modules on schemes (whose construction is very ``hands-on'' and given by explicit formulas) to derived algebraic geometry (where explicit formulas are usually not homotopically sound). Let us remark that one preliminary question, entirely in the realm of classical algebraic geometry, can also be how to define D-modules on a stack in groupoids (for one example of such construction see \cite[Chapter 7]{BeilDrin:QuantHitchin}).\\
			The main idea that allows us to carry classical D-modules to derived algebraic geometry (that is, to basically any prestack) is the identification between (left) D-modules on $X$ and crystals of quasi-coherent sheaves on $X$, where $X$ is a smooth scheme. This is an idea that goes back to Grothendieck; for a modern account, using $\infty$-categorical language, see \cite{GR:CrystalsDMod}. A crystal of quasi-coherent sheaves on $X$ is an object in the $\infty$-category obtained via the totalization of the following cosimplicial diagram
			\begin{diag}
				\Qcoh(X) \ar[r, shift left] \ar[r, shift right] & \Qcoh((X \times X)^{\wedge}_{\Delta}) \ar[r, shift left] \ar[r] \ar[r, shift right] & \Qcoh((X \times X \times X)^{\wedge}_{\Delta}) \dots
			\end{diag}
			It corresponds to the datum of a quasi-coherent sheaf $\mathcal{F}$ on $X$ equipped with equivalences $f^{ * }(\mathcal{F}) \stackrel{\sim}{\to} g^{ * }(\mathcal{F})$ for every pair of infinitesimally close points $f,g\colon \Spec R \to X$, and higher homotopies (realizing higher cocycle ``conditions''). Reformulating this definition via formal completions along diagonals, it is an easy computation to link this additional structure on $\mathcal{F}$ with a left action of $\mathcal{D}_{ X }$ on it, as explained in the notes \cite{Lurie:DMod}.
			Therefore, one possible answer to the above question can be to define the $\infty$-category of left D-modules on any (derived) prestack $X\colon \CAlg^{ \cn }_{ \C } \to \S$ as the totalization of the cosimplicial diagram drawn above, since all the terms there can be defined in derived algebraic geometry. It is proved in \cite[Proposition 3.4.3]{GR:CrystalsDMod} that the limit above yields\footnote{Technically, this is only true for bounded $X$ or classically formally smooth $X$, and the ``correct'' result uses ind-coherent sheaves, but we will simplify the story here.}
			\[
				\Qcoh(X_{ \dR }),
			\]
			the $\infty$-category of quasi-coherent sheaves on the de Rham prestack of $X$.
			Let us recall that the de Rham prestack of $X$ is defined as the functor
			\[
				X_{ \dR }\colon \CAlg^{ \cn }_{ \C } \to \S, \qquad A \mapsto X( H^{ 0 }(A)^{ \red } ),
			\]
			see \cref{defn:de_rham_prestack}. We have a first answer.
			\begin{quote}
				\textbf{Answer 0:} Given $X$ a (derived) prestack, the $\infty$-category of left D-modules on $X$ is $\Qcoh(X_{ \dR })$.
			\end{quote}
			With this interpretation, most of the formal constructions (like inverse and direct images, equivalence between left and right D-modules, Kashiwara lemma) go through: see \cite{GR:CrystalsDMod} and \cite[Chapter 4]{GR:DAG2}. As it is evident from the definition of de Rham prestack associated with $X$, one sees that the derived structure on $X$ is irrelevant when computing $X_{\dR}$ (because we only evaluate $X$ on classical algebras). That is, the functor
			\[
				(-)_{ \dR }\colon \PreStk \to \PreStk, \qquad X \mapsto X_{ \dR }
			\]
			can be factored as 
			\[
				X \mapsto t_{ 0 }(X) \mapsto X_{ \dR } \simeq (t_{ 0 }(X))_{ \dR }
			\]
			where $t_{ 0 }$ forgets the derived structure of $X$. This means that passing from $X$ to $X_{ \dR }$ ``forgets'' the derived structure present on $X$, see \cref{prop:de_rham_prestack_is_classical}. This leads us to the following question.
			\begin{quote}
				\textbf{Question 1:} Can one define D-modules in derived algebraic geometry in such a way that derived information is remembered? This is what we want to refer to as ``derived D-modules''.
			\end{quote}
			This question has already been answered as well. In fact, there are several, conjecturally equivalent, definitions of derived D-modules present in the literature: one uses derived foliations, introduced in \cite{TV:AlgFoliations2} (see \cref{defn:toen_vezzosi_derived_d_modules}), one uses dg-Lie algebroids and their representations, introduced in \cite{Nuiten:KoszulLieAlgbd} (although derived D-modules are never explicitly mentioned there) and the third one comes from \cite{Ber:DerivedDMod}, where ind-coherent sheaves and de Rham stack are used (see \cref{defn:beraldo_bounded_case} for the bounded case). 
			To get an intuitive idea, let us take a step back and work on a classical smooth scheme $X \to \Spec(\C)$. We can define left D-modules on $X$ in two equivalent ways:
			\begin{enumerate}
				\item there exists a sheaf of (filtered) associative algebras $\D_{ X }$ on $X$, obtained as universal enveloping algebra (equipped with the PBW filtration) of the tangent Lie algebroid $\mathcal{T}_{ X }$ (which is the classical tangent sheaf equipped with the standard Lie bracket between vector fields), i.e.\ $\D_{ X } = U(\mathcal{T}_{ X })$; left modules over $\D_{ X }$ are then representations of $\mathcal{T}_{ X }$;
				\item another way of giving a left $\D_{ X }$ action on a quasicoherent sheaf $M$ on $X$ is to give a $\C_{ X }$-linear \emph{flat} connection $\nabla\colon M \to M \otimes_{ \O_{ X } } \Omega^{ 1 }_{ X }$; without spelling out this definition in details, let us observe that it is equivalent to considering 
					\[
						M \otimes_{ \O_{ X } }^{ \gr } \Omega^{ * }_{ X } 
					\]
					as a graded module over $\Omega^{ * }_{ X }$, which is the algebraic de Rham complex of $X$ without differential (and equipped with wedge product of differential forms, making it a commutative graded algebra). Then, the datum of a flat connection $\nabla$ on $M$, is exactly the datum of a $\C_{ X }$-linear differential $\partial$ on the graded object $M \otimes_{ \O_{ X } }^{ \gr } \Omega^{ * }_{ X }$ that makes the former a dg-module over $(\Omega^{ * }_{ X }, d_{ \dR }, \wedge )$, i.e.\ the actual cdga of algebraic de Rham forms of $X$ with wedge product. That is, giving a left $\D_{ X }$-action on $M$ amounts to fill the dashed arrows of 
					\begin{diag}
						M \ar[r, dashed] & M \otimes_{ \O_{ X } } \Omega^{ 1 }_{ X } \ar[r, dashed] & M \otimes_{ \O_{ X } } \Omega^{ 2 }_{ X } \ar[r, dashed] & \dots
					\end{diag}
					so to get a complex ``compatible'' with the de Rham differential (which basically forces all the other arrows to be defined starting from the first one). This complex is sometimes referred to as the ``de Rham complex of the left D-module $M$''.
			\end{enumerate}
			This equivalence is classical, and used everywhere; let us just observe how we cannot just ``dualize'' $\nabla$ due to the fact that it is not $\O_{ X }$-linear (although it is clear, by hand, how to prove the above equivalence). \\
			Let us now explain the idea behind the three different definitions of derived D-modules.
			\begin{enumerate}[label=(\roman*)]
				\item Trying to generalize approach 1 to derived algebraic geometry, one wants to replace the tangent sheaf $\mathcal{T}_{ X }$ with the tangent \emph{complex} $\T_{ X }$. Using the theory of dg-Lie algebroids developed in \cite{Nuiten:HomotopicalLieAlgbd} and \cite{Nuiten:KoszulLieAlgbd} (see also \cite{Vezzosi:LieAlgebroidsModel}), it is possible to consider $\T_{ X }$ as a dg-Lie algebroid over $X$, take its universal enveloping dg-algebra $U(\T_{ X })^{ \PBW }$ and consider left modules over it (all in a homotopy-coherent way). We explain this in \cref{subsection:dg_lie_algebroids}.

				\item To generalize approach 2, one needs to generalize the classical algebraic de Rham complex in a derived friendly way. Clearly, replacing the cotangent sheaf $\Omega^{ 1 }_{ X }$ with the cotangent \emph{complex} $\L_{ X }$ is the first correct move. One needs then to consider the correct homotopical generalization of chain complexes (which are defined by explicit equalities like $d^{ 2 } = 0$), for which one can equivalently consider complete filtered objects, mixed graded objects or coherent chain complexes. Then, using these technical tools, one can formulate a derived version of $2$, as explained in the article \cite{TV:AlgFoliations2} (see \cref{section:first_definition_derived_DModules} for a more detailed review).
					
				\item The third definition of derived D-modules, introduced in \cite{Ber:DerivedDMod}, is at the level of sheaves and can be thought as a modification of the previously introduced notion of crystals on $X$. On a bounded stack $X$, it is given by the full subcategory of $\IndCoh(X_{ \dR })$ (which corresponds to ``classical'' right D-modules on $X$) spanned by those objects whose underlying ind-coherent sheaf on $X$ lies in the essential image of the embedding $\Qcoh(X) \hookrightarrow \IndCoh(X)$, see \cref{defn:beraldo_bounded_case}. It is, in particular, this $\infty$-category is monadic over $\Qcoh(X)$. In the general, unbounded case, the monad must be modified (via the process of filtered renormalization described in \cite[§4]{Ber:DerivedDMod}). The idea behind these technical definitions is to try to obtain ``quasi-coherent representations'' of $\T_{ X }$ by modifying the ind-coherent counterparts $\T_{ X }^{ \IC }$ and $U(\T_{ X }^{ \IC })$ (arising from the abstract theory of ``ind-coherent Lie algebroids'' introduced in \cite[Chapter 9]{GR:DAG2}).

			\end{enumerate}
		
			Let us remark how the derived structure is important in all these three definitions: for example, in the first two definitions (co)tangent complexes are used (instead of (co)tangent sheaves). The third one also remembers derived information, thanks to the presence $\Qcoh(X)$.
			Our speculation then arrives to the following question.

			\begin{quote}
				\textbf{Question 2:} What is the relation between these three possible definitions of derived D-modules? 
			\end{quote}

			One contribution of this thesis is proving the following conjecture.

			\begin{quote}
				\textbf{Conjecture:} These three definitions yield equivalent symmetric monoidal $\infty$-categories on a derived affine scheme $X$ satisfying suitable finiteness conditions.
			\end{quote}
	
		Let us conclude by motivating why it is interesting to study derived D-modules. For example, as explained in \cite{Ber:DerivedDMod}, they naturally appear in geometric Langlands when one computes the Drinfeld center of the monoidal $\infty$-category $\mathbb{H}(Y)$ for a quasi-smooth stack $Y$, defined as right modules over the dg-algebra of Hochschild cochains on $Y$. It is proved there that one obtains the $\infty$-category of \emph{derived} D-modules on the unbounded derived loop stack $LY \coloneqq Y \times_{ Y \times Y } Y$.
		Another application comes from \cite[§7.1]{TV:AlgFoliations2}, where a version of the Grothendieck-Riemann-Roch formula is proved for derived D-modules (generalizing the classical version for D-modules over smooth schemes).
		Finally, we believe that in trying to study an analogue of the geometric Langlands conjecture (now theorem) for surfaces (or higher dimensional varieties), derived D-modules will naturally arise, since $\mathrm{Bun}_{ G }(X)$ (the moduli stack of $G$-bundles on $X$) is not smooth if $\dim X \geq 2$ (so that one needs to consider its derived version $\R\IntHom(X, BG)$).

	\subsection*{Thesis outline and contributions}
		Let us describe more specifically what is written in this thesis and what are its mathematical contributions.\\ 
		We start off with \cref{section:preliminaries} where we recall various preliminary constructions that will be useful later on. We first review the notions of enriched hom objects in presentable $\infty$-categories (and their induced ``Morita adjunctions''), some finiteness conditions on algebras and modules and some general definitions of stacks and prestacks (like the de Rham prestack and formal completions) in derived algebraic geometry. We then review pro-coherent sheaves, ind-coherent sheaves and their functoriality and state their equivalence by Serre duality, so as to be able to use them interchangeably later on. We conclude this preliminary chapter with a short review of graded and filtered objects and their main properties.\\
		All of its mathematical content is already known, but we decided to write it anyway to make this document more self contained. Some statements that we proved by ourselves and we don't know a written account of, are \cref{lemma:symm_monoidal_functor_mod_procoh} (which gives a procedure to build symmetric monoidal colimit preserving exact functor out of pro-coherent modules), \cref{lemma:graded_modules_weight_0} and \cref{coroll:filtered_complete_modules_gr_0} (which characterize filtered and graded modules whose associated graded is concentrated in weight $0$).

		In \cref{section:first_definition_derived_DModules} we review the relevant definitions of \cite{TV:AlgFoliations2} and introduce their definition of derived D-modules, which we call derived crystals in \cref{defn:toen_vezzosi_derived_d_modules}. Aside from the definitions, the main result we prove is \cref{coroll:constant_modules_de_rham_conservative}, which proves that, for a derived affine scheme $\Spec A$ locally of finite presentation, the $\infty$-category of derived crystal is monadic over $\Mod_{ A }$ (let us remark that we do not invoke \cite[Theorem 3.2.1]{TV:AlgFoliations2}). 

		We then continue with \cref{section:second_definition_derived_DModules}, where we start by reviewing the theory of dg-Lie algebroids and their representations by \cite{Nuiten:KoszulLieAlgbd}. In \cref{subsubsection:representations_dg_lie_algebroids}, we give a different proof of the PBW theorem for cofibrant dg-Lie algebroids, which differs from the original one (see \cref{remark:joost_pbw}).
		We introduce the definition of derived D-modules given by Beraldo in \cref{subsubsection:derived_DModules_IndCoh_bounded_case}, where we only consider the case of bounded derived schemes. 
		After a brief review of formal moduli problems under an affine scheme, we prove our first original result (see \cref{thm:symm_monoidal_representation_lie_algbd}).
		\begin{thm}
			Let $A$ be a bounded connective cdga over $\C$ and $\g$ be a dg-Lie algebroid over $A$. There exists a \emph{symmetric monoidal} equivalence 
			\[
				\LMod_{ U(\g) } \stackrel{\sim}{\longrightarrow} \ProCoh_{ 0 }(\MC(\g)/A).
			\]
		\end{thm}
		The right-hand side is a symmetric monoidal full subcategory of pro-coherent sheaves on $\MC(\g)$ (see \cref{defn:procoh_with_underlying_qcoh}), which is the formal moduli problem under $\Spec A$ associated to $\g$ via the Koszul duality result of \cref{thm:koszul_duality_lie_algbd}.
		This proposition, at the level of underlying $\infty$-categories, is already known and proved in \cite[Theorem 7.1]{Nuiten:KoszulLieAlgbd}: our contribution is its enhancement to a symmetric monoidal equivalence.
		As a corollary, we obtain the following (see \cref{coroll:nuiten_is_beraldo_monoidal}).
		\begin{coroll}
			Let $A$ be a bounded connective cdga over $\C$, $X \coloneqq \Spec A$. We have symmetric monoidal equivalences
			\[
				\LMod_{ U(\T_{ A }) } \simeq \ProCoh_{ 0 }(X_{ \dR }/A) \simeq \IndCoh(X_{ \dR }) \times_{ \IndCoh(X) } \Qcoh(X).
			\]
		\end{coroll}
		This proves that, on a bounded derived affine scheme, derived D-modules as defined by Beraldo and representations of the tangent dg-Lie algebroid are equivalent (compatibly with their natural tensor product structures). 
		\begin{quote}
			\textbf{Partial answer:} The definitions of derived D-modules (ii) and (iii) are equivalent over a bounded derived affine schemes, compatibly with their respective symmetric monoidal structure.
		\end{quote}

		In the final \cref{section:main_result}, we relate derived crystals with representations of dg-Lie algebroids. In \cref{subsection:chevalley_eilenberg}, using the PBW filtration on the universal enveloping algebra of a dg-Lie algebroid $\g$ over $A$, we construct the $\infty$-functor of complete Hodge-filtered Chevalley-Eilenberg cohomology, see \cref{defn:mixed_graded_CE}. It is given by the formula
		\[
			\CE^{ \cpl }(\g, -)\colon \LMod_{ U(\g) } \to \Mod_{ \C }^{ \cpl }, \qquad M \mapsto \IntHom_{ U(\g)^{ \PBW } }^{ \cpl }(A, M),
		\]
		where we use the $\Mod_{ \C }^{ \cpl }$-hom object of the $\infty$-category of complete filtered $\g$-representations $\LMod_{ U(\g)^{ \PBW } }^{ \cpl }$ (where $A$ and $M$ are endowed with the positive filtration).
		It is a filtered enhancement of the classical Chevalley-Eilenberg cohomology functor, as we prove in \cref{prop:CE_lifts_classical}. Moreover, it generalizes the construction of the graded mixed Chevalley-Eilenberg cohomology complex of a dg-Lie algebra introduced in \cite{Pavia:MixedCE}.
		After having proved that, by monoidal properties, we can lift this functor to 
		\[
			\CE^{ \cpl }(\g, -)\colon \LMod_{ U(\g) } \to \Mod_{ \CE^{ \cpl }(\g, A) }(\Mod_{ \C }^{ \cpl }),
		\]
		we arrive to \cref{subsection:main_proof}, where we then prove that it is fully faithful and identify its essential image with the \emph{constant} modules over the (complete filtered) commutative algebra $\CE^{ \cpl }(\g, A)$ (defined in \cref{defn:constant_graded_mixed_modules}). Precisely, we prove the following theorem (see \cref{thm:main_thm}), giving a positive answer to \cite[Conjecture 4.2.3]{Pavia:MixedCE} when specialized to the dg-Lie algebra case.
		\begin{thm}
			If $\g$ has an underlying perfect $A$-module, the Hodge-completed Chevalley-Eilenberg cohomology functor induces a symmetric monoidal equivalence of $\infty$-categories
			\[
				\CE^{ \cpl }(\g, -)\colon \LMod_{ U(\g) } \stackrel{\sim}{\longrightarrow} \Mod^{ \const }_{ \CE^{ \cpl }(\g, A) }.
			\]
		\end{thm}
		This is a result about cohomology of dg-Lie algebroids, that can be of independent interest.
		Plugging in $\g = \T_{ A }$ (assuming that $A$ is locally of finite presentation, so that $\L_{ A }$ is perfect) we obtain the following (see \cref{coroll:equivalence_derived_d_mod}).
		\begin{coroll}
			Let $A$ be a locally finitely presented and bounded cdga over $\C$. There is a symmetric monoidal equivalence
		\[
			\LMod_{ U(\T_{ A }) } \simeq \Mod^{ \const }_{ \CE^{ \cpl }(\T_{ A }, A) },
		\]
			where the right-hand side is the symmetric monoidal $\infty$-category of derived crystals over $\Spec A$.
		\end{coroll}
		We deduce that all three different definitions of derived D-modules are equivalent in \cref{coroll:equivalence_derived_d_mod}. 

		\begin{quote}
			\textbf{Final answer:} The three different $\infty$-categories of derived D-modules are equivalent, compatibly with their symmetric monoidal structures, on a bounded derived affine scheme which is locally finitely presented.
		\end{quote}

	\section*{Future research directions}
		Let us now give some research paths that we hope to pursue in the future.\\
		The first natural step to pursue after the results proved in this work is to prove that the three definitions of derived D-modules coincide over \emph{unbounded} derived schemes and stacks. In particular, we would like to prove a generalization to the unbounded case of \cref{thm:symm_monoidal_representation_lie_algbd}. Some complications arise for an unbounded $A$: for example all arguments using double linear duality (like \cref{prop:almost_perfect_duality_bounded}) break down and there is no more a Koszul duality equivalence between Lie algebroids over $A$ and formal moduli problems under $\Spec A$.\\
		To be able to globalize our results we need to understand the functoriality of derived D-modules, in particular pullbacks (which are naturally defined in \cite{TV:AlgFoliations2} and \cite{Ber:DerivedDMod}). We claim that the symmetric monoidal equivalences given here are compatible with the natural operations of pullbacks on all sides. \\
		One more complete and satisfactory answer would be to prove that all these versions of derived D-modules are $3$-functor formalisms (as defined in \cite{Scholze:SixFunctor}) and the equivalences we gave are equivalences of $3$-functor formalisms. \\
		We also wish to investigate the relationship between derived D-modules and sheaves (quasi-coherent and ind-coherent) on the \emph{filtered de Rham stack}, as defined in \cite[Definition 2.3.5]{Bhatt:PrismaticFGauges}. If these notions are equivalent, studying the functoriality of the passage to filtered de Rham stack, we would automatically have a $3$-functor formalism as described above.\\
		Finally, we believe that our results can be applied to generalize to the derived setting the results of \cite{Chen:FilteredDModules}, which relates different notions of Koszul duality (both in the linear, filtered and unfiltered setting) between $\D$-modules and $\Omega$-modules over a smooth classical stack.

	\section*{Notations and conventions} 
		Throughout this thesis, we use $(\infty, 1)$-categories, that we will denote by $\infty$-categories (and even by ``categories'', sometimes, as opposed to $1$-categories). We try to be model independent, although most of the references we use (\cite{Lurie:HTT} and \cite{Lurie:HA}) adopt quasi-categories. We try to follow the notations of Lurie:  $\S$ is the $\infty$-category of spaces, $\Cat$ is the (large) $\infty$-category of (small) $\infty$-categories, $\CAT$ is the (very large) $\infty$-category of (large) $\infty$-categories and $\Pr^{ L }$ is the (very large) $\infty$-category of presentable $\infty$-categories and colimit-preserving functors between them. We implicitly fix universes to deal with set theoretic issues, as is done in \cite{Lurie:HTT}. 
		All the objects and operations are to be understood $\infty$-categorically; that is, when we write $\Mod$ or $\LMod$ we mean modules (or left modules) in the sense of \cite[Chapter 4]{Lurie:HA}.
		We will refer to $\E_{ \infty }$-algebra objects as ``commutative algebras'', denoted by $\CAlg$, and to $\E_{ 1 }$-algebra objects as ``associative algebras'', denoted by $\Alg$. We always assume to work over $\C$, i.e.\ if not specified our algebras are taken in the symmetric monoidal $\infty$-category $\Mod_{ \C }$.
		If we need to consider abelian categories of (left) modules over a discrete ring $R$, we will use the notation $\Mod^{ \heartsuit }_{ R }$ (since they can be seen as heart of the canonical t-structure on the unbounded derived category $\mathcal{D}(R)$).
		We will sometimes work with model categories as well, for which we will use the superscript ``dg''. For example, $\Mod_{ \C }^{ \dg }$ denotes the $1$-category of dg-modules over $\C$ which we usually endow with the projective model structure. 
		Moreover, grading conventions are always cohomological: our differentials raise by $1$ the degree, so connective objects (denoted by the superscript ``cn'') will be those concentrated in degrees $(-\infty, 0]$.

		The terminology we use for derived algebraic geometry comes partially from \cite{GR:DAG1}. In particular, we denote by $\dAff_{ \C }$ the $\infty$-category of derived affine schemes over $\C$ (defined as the opposite of $\CAlg_{ \C }^{ \cn }$). We drop the derived prefix when talking about (étale) stacks, denoted by $\Stk$, and prestacks, denoted by $\PreStk$: we always considered derived (pre)stacks that are valued in spaces. We always work over $\Spec \C$, even though sometimes we might omit the subscript $\C$ from the notations. We consider $0$ to be a natural number.

%% file: morita_categorical_stuff.tex
\section{Enriched hom objects and adjunctions}
	\label{subsection:categorical_higher_algebra}

	In this section, we will introduce enriched hom objects in \cref{lemma:internal_hom_categories}, enriched adjunctions and give a criterion to be able to lift a standard adjunction to an enriched one in \cref{lemma:enriched_adjunctions}.
	We will restrict ourselves to the presentable realm (or stable or $\C$-linear, as explained below); for a more general context, see \cite[§4.2.1]{Lurie:HA}.\\
	Consider the (very large) $\infty$-category $\Pr^{ L }$ of (large) presentable $\infty$-categories and colimit preserving functors between them, endowed with the symmetric monoidal structure $\otimes$ known as ``Lurie tensor product'' (see \cite[Prop 4.8.1.15]{Lurie:HA}). 
	There, the symmetric monoidal structure on $\Sp$ is introduced and it is proved that there exists a symmetric monoidal localization\footnote{I.e.\ there exists a symmetric monoidal left adjoint of the embedding above.}
	\[
		\Mod_{ \Sp }(\Pr^{ L }) \simeq \Pr^{ L, st } \hookrightarrow \Pr^{ L },
	\]
	which means that being stable (i.e.\ having an $\Sp$-tensoring in $\Pr^{ L }$) is a property and not a structure.
	In what follows, $\Pr^{ L }$ can be replaced with $\Pr^{ L, st }$ and everything goes through. Moreover, everything goes through also in $\Mod_{ \Mod_{ \C }}(\Pr^{ L })$, which is denoted by $\mathrm{DGCat}_{ \mathrm{cont} }$ in \cite{GR:DAG1} (this time thanks to rigidity of the symmetric monoidal category $\Mod_{ \C}$).\\
	Throughout this section, let $\infcatname{C} \in \CAlg(\Pr^{ L })$, that is, a presentably symmetric monoidal $\infty$-category. 
	The basic object of interest is the so called ``enriched hom'' functor, introduced in the following lemma.

	\begin{lemma}
	   \label{lemma:internal_hom_categories}
	   Let $\infcatname{M} \in \Mod_{ \infcatname{C} }(\Pr^{ L})$, i.e.\ a presentable $\infty$-category (left) tensored over $\infcatname{C}$ (denote by $\odot$ this action). Then $\infcatname{M}$ is enriched over $\infcatname{C}$, meaning there exists an \emph{enriched hom functor}
	   \[
		   \IntHom_{ \infcatname{M} }(-,-)\colon \infcatname{M}^{ \op } \times \infcatname{M} \to \infcatname{C}
		\]
		which preserves limits in each variable and such that for each $A \in \infcatname{C}$ and $M, N \in \infcatname{M}$ we have natural isomorphisms of spaces
		\[
		   \Map_{ \infcatname{M} }(A \odot M, N) \simeq \Map_{ \infcatname{C} }(A, \IntHom_{ \infcatname{M} }(M, N)).
		\]
		\end{lemma}
	\begin{proof}
	   	Consider the action functor
		\[
			- \odot - \colon \infcatname{C} \otimes \infcatname{M} \to \infcatname{C}
		\]
		in $\Pr^{ L }$. Following the chain of equivalences in $\Pr^{ L }$ (using that $\Fun^{ L }(-, -)$ is the internal hom of $\Pr^{ L}$, see \cite[Remark 4.8.1.18]{Lurie:HA}) we have
		\[
			\Fun^{ L }(\infcatname{C} \otimes \infcatname{M}, \infcatname{M}) \simeq \Fun^{ L }(\infcatname{M}, \Fun^{ L }(\infcatname{C}, \infcatname{M})).
		\]
		Recall that ``passing to right adjoints'' defines an equivalence of $\infty$-categories between $\Fun^{ L }(\infcatname{A}, \infcatname{B})$ and $\Fun^{ R }(\infcatname{B}, \infcatname{A})^{ \op }$, as proved in \cite[Prop 5.2.6.2]{Lurie:HTT}.\\
		Passing to right adjoints and then forgetting (faithfully) to $\CAT$, we can write
		\begin{gather*}
			\Fun^{ L }(\infcatname{M}, \Fun^{ R }(\infcatname{M}, \infcatname{C})^{ \op }) \to \Fun^{ L }(\infcatname{M}, \Fun(\infcatname{M}^{ \op }, \infcatname{C}^{ \op })) \to \Fun(\infcatname{M}, \Fun(\infcatname{M}^{ \op }, \infcatname{C}^{ \op })) \simeq \\
			\simeq \Fun(\infcatname{M}^{ \op} \times \infcatname{M}, \infcatname{C})^{ \op }
		\end{gather*}
		where we used, similarly, that $(-)^{ \op }$ is a $2$-functor (that is, it acts on natural transformation in a natural way, which this time is covariant).
		Re-assembling everything and considering only the maximal $\infty$-groupoid underlying (because we just care about points now), we built a map
		\[
				\Fun^{ L }(\infcatname{C} \otimes \infcatname{M}, \infcatname{M})^{ \simeq } \to \Fun(\infcatname{M}^{ \op } \times \infcatname{M}, \infcatname{C})^{ \simeq }
			\]
			which sends the action to a functor that we call $\IntHom_{ \infcatname{M} }(-, -)$.
			The universal property is now verified by construction, and also the usual limit and colimit preservation properties of internal homs are verified (again, by following the faithful arrows of the construction).
	 \end{proof}
	 The next lemma explains that $\infcatname{C}$-linear functors can be considered as $\infcatname{C}$-enriched functors.
		   \begin{lemma}
			   \label{lemma:C_linear_functors_are_enriched}
			   Let $\infcatname{C} \in \CAlg(\Pr^{ L })$ and $\infcatname{M}, \infcatname{N} \in \Mod_{ \infcatname{C} }(\Pr^{ L })$ and let 
			   \[
				   F\colon \infcatname{M} \to \infcatname{N}
			   \]
			   be a $\infcatname{C}$-linear functor (that is, a morphism in $\Mod_{ \infcatname{C} }(\CAT)$).
			   Then $F$ induces a natural transformation
			   \[
				   \IntHom_{ \infcatname{M} }(-, -) \to \IntHom_{ \infcatname{N} }(F(-), F(-))
			   \]
			   of functors $\infcatname{M}^{ \op } \times \infcatname{M} \to \infcatname{C}$.
		   \end{lemma}
		   \begin{proof}
			   By adjunction, it suffices to give natural morphisms in $\infcatname{N}$
			   \[
				   \IntHom_{ \infcatname{M} }(-, \bullet) \odot F(-) \dashrightarrow F(\bullet)
			   \]
			   where $-, \bullet \in \infcatname{M}$.
				By $\infcatname{C}$-linearity of $F$ we have
				\[
					\IntHom_{ \infcatname{M} }(-, \bullet) \odot F(-) \simeq F\left( \IntHom_{ \infcatname{M} }(-, \bullet) \odot - \right) \stackrel{F(\ev_{ (-) })}{\longrightarrow} F(\bullet)
				\]
				where we then applied $F$ to the natural evaluation maps $\IntHom_{ \infcatname{M} }(-, \bullet) \odot (-) \to \bullet$ in $\infcatname{M}$.
		   \end{proof}
		   Let us now consider adjunctions between $\infcatname{C}$-modules where the left adjoint is $\infcatname{C}$-linear. The following lemma proves that they can be considered \emph{$\infcatname{C}$-enriched} adjunctions; that is, the usual property of adjunctions for mapping spaces is true at the level of $\infcatname{C}$-mapping objects.

		   \begin{lemma}
			   \label{lemma:enriched_adjunctions}
			   Let $\infcatname{C} \in \CAlg(\Pr^{ L })$ and $\infcatname{M}, \infcatname{N} \in \Mod_{ \infcatname{C} }(\Pr^{ L})$. Consider an adjunction 
			   \[
			   	   \adjunction{F}{\infcatname{M}}{\infcatname{N}}{G}
			   \]
			   such that $F$ is $\infcatname{C}$-linear (that is, $F \in \Mod_{ \infcatname{C} }(\Pr^{ L})^{ \Delta^{ 1 } }$).
			   Then it is a $\infcatname{C}$-enriched adjunction, meaning the familiar isomorphism of spaces comes from applying $\Map_{ \infcatname{C} }(1_{ \infcatname{C} }, -)$ to natural equivalences in $\infcatname{C}$
			   \[
				   \IntHom_{ \infcatname{N} }(F(M), N) \simeq \IntHom_{ \infcatname{M} }(M, G(N)).
			   \]
		   \end{lemma}
		   \begin{proof}
			   By \cref{lemma:C_linear_functors_are_enriched}, we can build natural morphisms in $\infcatname{C}$
			   \[
				   \phi\colon \IntHom_{ \infcatname{M} }(M, G(N)) \to \IntHom_{ \infcatname{N} }(F(M), FG(N)) \to \IntHom_{ \infcatname{N} }(F(M), N)
			   \]
			   where we ``composed'' by the counit in $N$ in the last arrow.
			   We can check that it is an isomorphism by using the Yoneda lemma (in $\infcatname{C}^{ \op }$); for an arbitrary $C \in \infcatname{C}$, consider the diagram
			   \begin{diag}
			   	  \Map_{\infcatname{C}}(C, \IntHom_{\infcatname{M}}(M, G(N))) \ar[r, "\sim"] \ar[d] & \Map_{\infcatname{M}}(C \odot M, G(N)) \ar[d, "F"] \\
				  \Map_{\infcatname{C}}(C, \IntHom_{\infcatname{N}}(F(M), FG(N))) \ar[d] \ar[r, "\sim"] & \Map_{\infcatname{N}}(C \odot F(M), FG(N)) \ar[d, "\mathrm{counit}_N"] \\
				  \Map_{\infcatname{C}}(C, \IntHom_{\infcatname{N}}(F(M), N)) \ar[r, "\sim"] & \Map_{\infcatname{N}}(C \odot F(M), N)
			   \end{diag}	
			   where the left column is obtained using the map $\phi$ described above and the line isomorphisms hold by tensor-hom adjunction. Notice how the $\infcatname{C}$-linearity of $F$ is crucially used in the top rightmost arrow.
			   One verifies, by unravelling all the definitions, that this diagram commutes.
			   Finally, composing the two functors in the right column we get, by definition, the adjunction isomorphism
			   \[
				   \Map_{ \infcatname{M} }(C \odot M, G(N)) \simeq \Map_{ \infcatname{N} }(F(C \odot M), N) \simeq \Map_{ \infcatname{N} }(C \odot F(M), N).
			   \]
			   This means that also the left composite must be an isomorphism of spaces, which implies, by Yoneda (varying $C$), that the initial morphism $\phi$ is an equivalence in $\infcatname{C}$.
		   \end{proof}
			
		   \begin{remark}
			   \label{remark:right_adjoint_lax_C_linear}
			   Let $\infcatname{C}$ be as above and $(F, G)$ be an adjunction between $\infcatname{C}$-module categories, such that $F$ is $\infcatname{C}$-linear.  This does not imply that $G$ is $\infcatname{C}$-linear; in fact it will only be \emph{lax} $\infcatname{C}$-linear in general.
			   This comes from considering the (first) commutative diagram accounting for the $\infcatname{C}$-linearity of $F$
			   \begin{diag}
			   	   \infcatname{C} \times \infcatname{M} \ar[r, "\id \times F"] \ar[d, "- \odot -"] & \infcatname{C} \times \infcatname{N} \ar[d, "- \odot -"] \\
				   \infcatname{M} \ar[r, "F"] & \infcatname{N},
			   \end{diag}
			   which, when passing to horizontal right adjoints, induces a Beck-Chevalley natural transformation
			   \[
				   C \odot G(N) \to G(C \odot N)
			   \]
			   for all $C \in \infcatname{C}$ and $N \in \infcatname{N}$.
		   \end{remark}

		   \begin{remark}
		   		\label{remark:adjoint_functor_rigid_module_categories}
			In the same setting as \cref{remark:right_adjoint_lax_C_linear}, suppose now that $\infcatname{C}$ is, additionally, locally rigid as in \cite[Definition D.7.4.1]{Lurie:SAG}. Then, if $G$ preserves colimits, it is automatically $\infcatname{C}$-linear (and not only lax), thanks to \cite[Remark D.7.4.4]{Lurie:SAG}.
		   \end{remark}

		   \begin{lemma}
		   		\label{lemma:enriched_right_adjoint}
				Let $\infcatname{C} \in \CAlg(\Pr^{ L })$ and $\infcatname{M}, \infcatname{N} \in \Mod_{ \infcatname{C} }(\Pr^{ L })$. Consider an adjunction
				\[
					\adjunction{F}{\infcatname{M}}{\infcatname{N}}{G}
				\]
				such that $G$ is $\infcatname{C}$-linear. Then there are natural morphisms in $\infcatname{C}$
				\[
					\IntHom_{ \infcatname{N} }(F(M), N) \longrightarrow \IntHom_{ \infcatname{M} }(M, G(N))
				\]
				which becomes the adjunction equivalences (of spaces) after applying $\Map_{ \infcatname{C} }(1_{ \infcatname{C} }, -)$.
		   \end{lemma}
		   \begin{proof}
		   		The natural maps are obtained as composition
				\[
					\phi\colon \IntHom_{ \infcatname{N} }(F(M), N) \longrightarrow \IntHom_{ \infcatname{M} }(GF(M), G(N)) \longrightarrow \IntHom_{ \infcatname{M} }(M, G(N))
				\]
				where the first arrow uses $\infcatname{C}$-linearity of $G$ (see \cref{lemma:C_linear_functors_are_enriched}) and the second arrow is precomposition with the unit morphism $M \to GF(M)$.
				As in the last part of the proof of \cref{lemma:enriched_adjunctions}, it is easily verified that $\Map_{ \infcatname{C} }(1_{ \infcatname{C} }, \phi)$ induces the adjunction equivalences.
		   \end{proof}

		   \begin{example}
			   \label{example:canonical_enrichment_lmod}
			   Most of our examples of $\infcatname{C}$-module categories will be given by $\LMod_{ A }(\infcatname{C})$ for $A \in \Alg(\infcatname{C})$ with the canonical $\infcatname{C}$-module structure coming from \cite[§4.3.2]{Lurie:HA}.
		   \end{example}

		   In the particular case in which $\infcatname{M} = \Mod_{ B }(\infcatname{C})$ for $B \in \CAlg(\infcatname{C})$, the canonical $\infcatname{C}$-enrichment of $\infcatname{M}$ comes from a symmetric monoidal left adjoint functor $\infcatname{C} \to \infcatname{M}$, which implies that the symmetric monoidal structure on $\Mod_{ B }(\infcatname{C})$ (relative tensor product over $B$) is compatible with the enrichment. In other words,
		   \[
			   \Mod_{ B }(\infcatname{C}) \in \CAlg(\Mod_{ \infcatname{C} }(\Pr^{ L})) \simeq \CAlg(\Pr^{ L})_{ \infcatname{C}/- }
		   \]
		   where the last equivalence holds thanks to \cite[§3.4.1]{Lurie:HA}.
		   This is the content of the following lemma.
		   
		   \begin{lemma}
			   \label{lemma:enrichment_mod_A}
				Let $\infcatname{C} \in \CAlg(\Pr^{ L})$ and $B \in \CAlg(\infcatname{C})$. The canonical $\infcatname{C}$-module structure on $\Mod_{ B }(\infcatname{C})$ upgrades to a commutative $\infcatname{C}$-algebra structure (in $\Pr^{ L}$) on $\Mod_{ B }(\infcatname{C})$, which is induced by the symmetric monoidal left adjoint ``free $B$-module'' functor 
				\[
					B \otimes -\colon \infcatname{C} \to \Mod_{ B }(\infcatname{C}).
				\]
			\end{lemma}
			\begin{proof}
				The free $B$-module functor exists, is left adjoint to the forgetful functor and is symmetric monoidal by \cite[Theorem 4.5.3.1]{Lurie:HA}. The $\infcatname{C}$-action on $\Mod_{ B }(\infcatname{C})$ is given by
				\[
					- \odot -\colon \infcatname{C} \times \Mod_{ B }(\infcatname{C}) \to \Mod_{ B }(\infcatname{C}),\qquad (V, M) \mapsto V \odot M = V \otimes M
				\]
				where the $B$-module structure on $V \otimes M$ is given by $B$ acting on $M$, that is, 
				\[
					B \otimes (V \otimes M) \simeq V \otimes (B \otimes M) \to V \otimes M = V \odot M.
				\]
				This is equivalent to the $B$-module obtained by $(B \otimes V) \otimes_{ B } M$ by elementary properties of relative tensor product in a presentably symmetric monoidal $\infty$-category.
			\end{proof}

			As a corollary, since $\infcatname{C}$ is canonically a module over itself (using $\id\colon \infcatname{C} \to \infcatname{C}$), we obtain that it is enriched over itself, with internal hom denoted by $\IntHom(-, -)$; in other words, $\infcatname{C}$ is \emph{closed} in the sense of \cite[Definition 4.1.1.15]{Lurie:HA}.

		   \begin{coroll}
			   \label{coroll:tensor_hom_enriched}
			   Let $\infcatname{C} \in \CAlg(\Pr^{ L })$ and $\infcatname{M} \in \Mod_{ \infcatname{C} }(\Pr^{ L })$. The natural equivalences of spaces in \cref{lemma:internal_hom_categories} come from applying $\Map_{ \infcatname{C} }(1, -)$ to the corresponding (natural) equivalences in $\infcatname{C}$ \[
				   \IntHom_{ \infcatname{M} }(A \odot M, N) \simeq \IntHom(A, \IntHom_{ \infcatname{M} }(M, N)).
		   	  \]
			  In other words, the adjunction
			  \[
				  \adjunction{- \odot M}{\infcatname{C}}{\infcatname{M}}{\IntHom_{ \infcatname{M} }(M, -) }		
			  \]
			  is $\infcatname{C}$-enriched for all $M \in \infcatname{M}$.
	       \end{coroll}
		   \begin{proof}
			   Let us first observe that $- \odot M\colon \infcatname{C} \to \infcatname{M}$ is a $\infcatname{C}$-linear functor, by standard properties of the $\infcatname{C}$-action on $\infcatname{M}$. It suffices then to apply the general \cref{lemma:enriched_adjunctions}.
		   \end{proof}

		  \begin{lemma}
			  \label{lemma:algebra_endomorphism}
			  Let $\infcatname{C}$ and $\infcatname{M}$ be as above and $X \in \infcatname{M}$.
			  Then $\IntHom(X, X) \in \infcatname{C}$ has a canonical associative algebra structure.
		  \end{lemma}
		  \begin{proof}
			   Well known: see \cite[§4.7.1]{Lurie:HA}. The idea is that the product map
			   \[
			   	   \IntHom(X, X) \otimes \IntHom(X, X) \dashrightarrow \IntHom(X, X)
			   \]
			   is adjoint to
			   \[
			   	  (\IntHom(X, X) \otimes \IntHom(X,X)) \odot X \simeq \IntHom(X, X) \odot (\IntHom(X, X) \odot X) \to \IntHom(X, X) \odot X \to X
			   \]
			   which is simply applying twice the ``evaluation'' map. The unit is the identity of $X$.
		  \end{proof}

		  We proved that all internal/enriched endomorphism objects are endowed with a canonical algebra structure. We will now prove that each object $M$ is canonically a left module for the algebra of its endomorphisms (which is the final algebra in this sense, see \cite[§4.7.1]{Lurie:HA}).

		  \begin{lemma} 
				\label{lemma:standard_bimodule_structure}
				Let $\infcatname{C}$ and $\infcatname{M}$ be as above and $M \in \infcatname{M}$. Then there exists a canonical (dashed) lift (along $\Lambda^{ 2 }_{ 2 } \hookrightarrow \Delta^{ 2 }$) of the (solid) diagram
					\begin{diag}
						& \LMod_{\IntHom_{\infcatname{M}}(M, M)}(\infcatname{M}) \ar[d] \\ 
						\Delta^0 \ar[r, "M"] \ar[ur, dashed] & \infcatname{M}.
					\end{diag}

			\end{lemma}
			\begin{proof}
				This is proved in \cite[§4.7.1]{Lurie:HA}. 
            \end{proof}

		  Sometimes we will wonder whether some functors between closed symmetric monoidal categories ``commute'' with the internal hom objects, especially with the algebra of endomorphisms.
		  \begin{lemma}
			  \label{lemma:adjunction_and_enriched_homs_projection_formula}
			  Let $\infcatname{C}, \infcatname{D} \in \CAlg(\Pr^{ L })$ and consider an adjunction
			  \[
				  \adjunction{F}{\infcatname{C}}{\infcatname{D}}{G}.
			  \]
			  Assume that:
			  \begin{enumerate}
				  \item $G$ is lax symmetric monoidal\footnote{Or, equivalently by \cite[Corollary 3.4.8]{HHLN:LaxMonoidalAdjunctions}, that $F$ is oplax symmetric monoidal.};
				  \item The ``projection formula'' natural morphisms, obtained as the composition
					  \[
						  F(A \otimes G(X)) \longrightarrow F(A) \otimes FG(X) \stackrel{F(A) \otimes \mathrm{counit}_{ X }}{\longrightarrow} F(A) \otimes X
					  \]
					  are equivalences in $\infcatname{D}$, for each $A \in \infcatname{C}$ and $X \in \infcatname{D}$.
			  \end{enumerate}
			  Then, using $\infcatname{C}$-enriched hom objects, there are natural isomorphisms
			  \[
					\epsilon_{ X, Y }\colon	G(\IntHom_{ \infcatname{D}}(X, Y))  \stackrel{\sim}{\longrightarrow} \IntHom_{ \infcatname{C} }(G(X), G(Y))		  
			  \]
			  in $\infcatname{C}$, for each $X, Y \in \infcatname{D}$. 
			  Moreover, if the adjunction above upgrades to the algebra level
			  \[
				  \adjunction{F}{\Alg(\infcatname{C})}{\Alg(\infcatname{D})}{G},
			  \]
			  then the isomorphisms $\epsilon_{ X, X }$ respect the ``endomorphism'' algebra structures, for each $X \in \infcatname{D}$.
		  \end{lemma}
		  \begin{proof}
			  First of all, let us observe that both $\infcatname{C}$ and $\infcatname{D}$ admit internal hom functors by \cref{lemma:internal_hom_categories}. The lax symmetric monoidal structure on $G$ induces a canonical op-lax symmetric monoidal structure on $F$ (easy to see by hand or as a corollary of more refined statements like \cite[Corollary 3.4.8]{HHLN:LaxMonoidalAdjunctions}).
			  Let us first build a map and then prove it is an equivalence. By the tensor-hom adjunction, it suffices to give
			  \[
				  G(\IntHom_{ \infcatname{D} }(X, Y)) \otimes G(X) \to G\left( \IntHom_{ \infcatname{D} }(X, Y) \otimes X \right) \to G(Y)
			  \]
			  where the first arrow comes from the lax monoidal structure of $G$, and the second arrow is simply $G$ of the evaluation in $X$. We built a map
			  \[
				  \epsilon_{ X, Y }\colon G(\IntHom_{ \infcatname{D} }(X, Y)) \to \IntHom_{ \infcatname{C} }(G(X), G(Y))
			  \]
			  in $\infcatname{C}$. To prove it's an equivalence it suffices to apply the Yoneda lemma (in $\infcatname{C}^{ \op }$) so choose an arbitrary $A \in \infcatname{C}$ and consider
			  \[
				  \Map_{ \infcatname{C} }(A, G(\IntHom_{ \infcatname{D} }(X, Y))) \simeq \Map_{ \infcatname{D} }(F(A), \IntHom_{ \infcatname{D} }(X, Y)) \simeq \Map_{ \infcatname{D} }(F(A) \otimes X, Y).
			  \]
			  Precomposing with the equivalence given by the projection formula, we can continue writing
			  \begin{gather*}
				  \Map_{ \infcatname{D} }(F(A) \otimes X, Y) \simeq \Map_{ \infcatname{D} }(F(A \otimes G(X)), Y) \simeq \Map_{ \infcatname{C} }(A \otimes G(X), G(Y)) \simeq \\
				  \simeq \Map_{ \infcatname{C} }(A, \IntHom_{ \infcatname{C} }(G(X), G(Y)).
			  \end{gather*}
			  Unravelling all adjunctions and lax monoidal structure maps, one can tediously verify that this isomorphism of spaces implies that $\Map_{ \infcatname{C} }(A, \epsilon_{ X,Y })$ is an isomorphism as well.
			  Varying $A \in \infcatname{C}$ we conclude, by Yoneda, that $\phi$ is an isomorphism in $\infcatname{C}$ for each $X, Y \in \infcatname{D}$.\\
			  Let us now consider the case $X = Y$, so that $\IntHom_{ \infcatname{D} }(X, X)$ and $\IntHom_{ \infcatname{C} }(G(X), G(X))$ are both algebras by \cref{lemma:algebra_endomorphism}. 
			  We want to verify that the map $\epsilon_{ X,X }$ respects this algebra structure.
			  We know by assumption that the adjunction lifts at the algebra level
			  \[
				  \adjunction{F}{\Alg(\infcatname{C})}{\Alg(\infcatname{D})}{G}.
			  \]
			  The fact that $G$ is lax symmetric monoidal implies that $G(\IntHom_{ \infcatname{D} }(X, X))$ has an induced algebra structure in $\infcatname{C}$.
			  Let us now observe that giving an algebra map
			  \[
				  G(\IntHom_{ \infcatname{D} }(X, X)) \to \IntHom_{ \infcatname{C} }(G(X), G(X))
			  \]
			  is equivalent, by the universal property in \cite[§4.7.1]{Lurie:HA}, to giving $G(X) \in \infcatname{C}$ a left action of $G(\IntHom_{ \infcatname{D} }(X, X))$.
			  By \cref{lemma:standard_bimodule_structure}, $X$ has a canonical left $\IntHom_{ \infcatname{D} }(X, X)$-module structure and, since lax symmetric monoidal functors preserve modules, we obtain an induced left $G(\IntHom_{ \infcatname{D} }(X, X))$-module structure on $G(X) \in \infcatname{C}$.
			  Finally, since forgetting the algebra structure is a conservative functor, we can verify that the above map of algebras is an isomorphism just as a map of $\infcatname{C}$.
			  This follows from a Yoneda type of argument exactly like above.
		  \end{proof}

		  Let us single out an easier condition (which is by no way necessary) to check for the above lemma to hold.

		  \begin{coroll}
			  \label{coroll:adjunction_enriched_hom_easy_version}
			   Let $\infcatname{C}, \infcatname{D} \in \CAlg(\Pr^{ L})$ and consider an adjunction
			  \[
				  \adjunction{F}{\infcatname{C}}{\infcatname{D}}{G}.
			  \]
			  Suppose $F$ is symmetric monoidal and $G$ is fully faithful.
			   Then, using enriched homs, we have natural isomorphisms
			  \[
				  \IntHom_{ \infcatname{C} }(G(X), G(Y)) \stackrel{\sim}{\longleftarrow} G(\IntHom_{ \infcatname{D}}(X, Y))
			  \]
			  in $\infcatname{C}$, for each $X, Y \in \infcatname{D}$. Moreover, if $X = Y$, then this isomorphism respects the ``endomorphism'' algebra structures.
		  \end{coroll}
		  \begin{proof}
			  It suffices to observe that the lax symmetric monoidal structure on $G$ is canonical and that the projection formula is satisfied, by symmetric monoidality of $F$ and by fully faithfulness of $G$ (which gives that the counit is an equivalence). Moreover, the adjunction $(F, G)$ lifts at the algebra level by \cite[Corollary 7.3.2.12]{Lurie:HA}. We then apply \cref{lemma:adjunction_and_enriched_homs_projection_formula}.
		  \end{proof}
		
	   \section{Morita adjunctions}
	   		\label{subsection:morita_adjunctions}
			
			In this section, we will introduce the so-called ``Morita adjunctions''. The starting point is the well known $\infty$-categorical version of the \emph{Eilenberg-Watts} theorem.

			\begin{thm}
				\label{thm:eilenberg_watts}
				Let $\infcatname{C} \in \CAlg(\Pr^{ L })$, $\infcatname{M} \in \Mod_{ \infcatname{C} }(\Pr^{ L })$ and $X, Y \in \Alg(\infcatname{C})$.
				There is an equivalence of $\infty$-categories
				\[
					\Fun^{ L }_{ \infcatname{C} }(\LMod_{ Y }(\infcatname{C}), \LMod_{ X }(\infcatname{M})) \simeq (X, Y)-\BiMod(\infcatname{M}).
				\]
			\end{thm}
			\begin{proof}
				This is proved in \cite[Theorem 4.8.4.1]{Lurie:HA}: choose $\mathcal{K}$ to denote the collection of all (small) simplicial sets, observe that all colimit preservation assumptions are satisfied (since we work in $\Pr^{ L }$) and notice that we denote by $\Fun^{ L }_{ \infcatname{C} }$ what Lurie denotes by $\mathrm{LinFun}_{ \infcatname{C} }^{ \mathcal{K} }$.
				The functor $\theta$ inducing the equivalence is the composition
				\begin{gather*}
					\Fun^{ L }_{ \infcatname{C} }(\LMod_{ Y }(\infcatname{C}), \LMod_{ X }(\infcatname{M})) \to \Fun_{ \infcatname{C} }(\LMod_{ Y }(\infcatname{C}), \LMod_{ X }(\infcatname{M})) \to \\
					\to \Fun_{ \infcatname{C} }(\RMod_{ Y }(\LMod_{ Y }(\infcatname{C})), \RMod_{ Y }(\LMod_{ X }(\infcatname{M}))) \to \RMod_{ Y }(\LMod_{ X }(\infcatname{M})) \simeq \\
					\simeq (X, Y)-\BiMod(\infcatname{M}) 
				\end{gather*}
				where we used $\infcatname{C}$-linearity to consider the induced functor on right $Y$-modules and then evaluation at the bimodule $Y$.
				The inverse of $\theta$ is given, on objects, by
				\[
					(X, Y)-\BiMod(\infcatname{M}) \ni Q \mapsto Q \otimes_{ Y } -\colon \LMod_{ Y }(\infcatname{C}) \to \LMod_{ X }(\infcatname{M}),
				\]
				where we used the same notation for the relative tensor product in $\infcatname{C}$ and the action of $\infcatname{C}$ on $\infcatname{M}$; see \cite[Remark 4.8.4.9]{Lurie:HA}.
			\end{proof}

			\begin{remark}
				For the case of $\infcatname{C} = \infcatname{M} = \Sp$, one has the more classical \cite[Proposition 7.1.2.4]{Lurie:HA}.
			\end{remark}
			
			We can now introduce another version of internal hom objects, which doesn't always come directly from \cref{lemma:internal_hom_categories} (and when it does, there is no ambiguity, see \cref{prop:how_to_compute_morita_homs} ).

			\begin{defn}
				\label{defn:morita_hom_objects}

				Let $\infcatname{C} \in \CAlg(\Pr^{ L })$, $\infcatname{M} \in \Mod_{ \infcatname{C} }(\Pr^{ L })$ and $X, Y \in \Alg(\infcatname{C})$. 
				For each $Q \in (X, Y)-\BiMod(\infcatname{M})$, corresponding to a left adjoint functor $F_{ Q }\coloneqq Q \otimes_{ Y } -$, we can consider the functor 
				\[
					\IntHom_{ X }^{ Y }(Q, -)\colon \LMod_{ X }(\infcatname{M}) \to \LMod_{ Y }(\infcatname{C}),
				\]
				defined as the right adjoint of $F_{ Q }$. We will omit the superscript $Y$ if it is clear from the context.
			\end{defn}

			\begin{prop}
				\label{prop:how_to_compute_morita_homs}
				Let $\infcatname{C} \in \CAlg(\Pr^{ L })$, $\infcatname{M} \in \Mod_{ \infcatname{C} }(\Pr^{ L })$, $X, Y \in \Alg(\infcatname{C})$ and $Q \in (X, Y)-\BiMod(\infcatname{M})$.
				There exists a natural transformation making the following diagram commutative 
				\begin{diag}
					& \LMod_Y(\infcatname{C}) \ar[d, "\oblv"] \\
					\LMod_X(\infcatname{M}) \ar[ur, "\IntHom^Y_X(Q{,} -)"] \ar[r, "\IntHom_X(Q{,} -)"] & \infcatname{C}.
				\end{diag}
			\end{prop}
			\begin{proof}
				Forgetting the right module structure along the algebra map $1_{ \infcatname{C} } \to Y$ (the unit of $Y$) gives us the diagram
				\begin{diag}
					\Fun^L_{\infcatname{C}}(\LMod_Y(\infcatname{C}), \LMod_X(\infcatname{M})) \ar[r, "\simeq"] \ar[d, "- \circ (Y \otimes -)"] & (X, Y)-\BiMod(\infcatname{M}) \ar[d, "\oblv_Y"] \\
					\Fun^L_{\infcatname{C}}(\infcatname{C}, \LMod_X(\infcatname{M})) \ar[r, "\simeq"] & \LMod_X(\infcatname{M}),
				\end{diag}
				which commutes by the explicit form of the equivalence in \cref{thm:eilenberg_watts} and by the fact that 
				\[
					Q \otimes_{ Y } (Y \otimes -) \simeq Q \otimes -
				\]
				where we recall that $Y \otimes -\colon \infcatname{C} \to \LMod_{ Y }(\infcatname{C})$ is the free left $Y$-module functor, left adjoint to $\oblv$.
				Passing to right adjoint functors in the above formula one obtains the desired equivalence.
			\end{proof}
			In particular, this proposition tells us that all ``Morita'' hom objects from \cref{defn:morita_hom_objects} can be computed (once we forget the $Y$-action) as enriched hom functors introduced in \cref{lemma:internal_hom_categories} (because for $Y = 1_{ \infcatname{C} }$ we can lift the $\infcatname{C}$-action on $\infcatname{M}$ to an action on $\LMod_{ X }(\infcatname{M})$ by \cite[§4.3.2]{Lurie:HA} and then pass to right adjoints). 
			In a nutshell, this tells us that enriched hom objects from a right $Y$-module inherit a canonical left $Y$-action by ``pre-multiplication''.\\
			Let us conclude this subsection with a final remark.
			\begin{remark}
				\label{remark:colimit_in_morita_hom}

				It is easily seen that Morita hom objects satisfy the usual properties of hom objects, for example sending colimits in the first variable to limits.
				In fact, a right adjoint of the functor $\IntHom^{ Y }_{ X }(\colim_{ i } Q_{ i }, -)$ is given by 
				\[
					(\colim_{ i } Q_{ i }) \otimes_{ Y } - \simeq \colim_{ i } (Q_{ i } \otimes_{ Y } -)
				\]
				and we can then conclude by observing that $\lim_{ i } \IntHom^{ Y }_{ X }(Q_{ i }, -)$ is the right adjoint of the right-hand side (where we used $\Fun^{ L }(\infcatname{A}, \infcatname{B}) \simeq \Fun^{ R }(\infcatname{B}, \infcatname{A})^{ \op }$).
			\end{remark}

%% file: RecallsDAG.tex
 \section{Some notions of Derived Algebraic Geometry}
	\label{subsection:quick_recalls_dag}

	Let us now recall some definitions of interest from derived algebraic geometry, among which some finiteness conditions on algebras and modules, the (algebraic) cotangent complex and the de Rham prestack.

	\subsection{Truncated and almost compact objects}
		\label{subsubsection:truncated_almost_compact_objects}
		
		This subsection is a brief overview of truncated and almost compact objects in a general $\infty$-category $\infcatname{C}$, as introduced in \cite[§5.5.6]{Lurie:HTT}. Its main result is \cref{prop:truncated_connective_algebras}, which states that truncatedness of an algebra (in a pre-stable $\infty$-category) is checked at the level of underlying object. The main interest of this subsection, for us, is to better understand the standard definitions of almost perfect modules (\cref{defn:almost_perfect_modules}) and of almost finitely presented connected algebras (\cref{defn:almost_finite_presentation}).
		We start by recalling the definitions of truncated/connective spaces. We will follow \cite{Lurie:HTT}.

		\begin{defn}
			\label{defn:truncated_spaces}
			Let $X \in \S$. We say that
			\begin{itemize}
				\item $X$ is \emph{(-2)-truncated} if it is weakly contractible\footnote{Meaning all its homotopy groups are trivial. Note that $\emptyset$ is automatically excluded.};
				\item $X$ is \emph{(-1)-truncated} if it is either empty or weakly contractible;
				\item $X$ is \emph{$n$-truncated}, for $n \geq 0$, if for each $x \in \pi_{ 0 }(X)$ and $i > n$ we have $\pi_{ i }(X, x) = 0$\footnote{$\emptyset$ satisfies this condition.}.
			\end{itemize}
			A map of spaces $f\colon X \to Y$ is \emph{$n$-truncated} if all of its (homotopy) fibers are $n$-truncated.
		\end{defn}
		Observe how every $n$-truncated space/map is automatically $(n+1)$-truncated, for each $n \geq -2$. Using mapping spaces, we can define truncated objects in any $\infty$-category, see \cite[§5.5.6]{Lurie:HTT}.
		\begin{defn}[Truncated object]
			\label{defn:truncated_objects}
			Let $\infcatname{C}$ be an $\infty$-category and $X \in \infcatname{C}$. We say that $X$ is \emph{$n$-truncated}, for $n \geq -2$, if $\Map_{ \infcatname{C} }(E, X)$ is an $n$-truncated space for each $E \in \infcatname{C}$.
			A morphism $f\colon X \to Y$ in $\infcatname{C}$ is \emph{$n$-truncated} if for each $E \in \infcatname{C}$ the induced morphism of spaces
			\[
				\Map_{ \infcatname{C} }(E, X) \to \Map_{ \infcatname{C} }(E, Y)
			\]
			is $n$-truncated.
			The $n$-truncated objects of $\infcatname{C}$ span a full subcategory denoted by $\infcatname{C}^{ \textrm{$n$-truncated }}$. 
		\end{defn}
		Thanks to \cite[Proposition 5.5.6.5]{Lurie:HTT}, the embedding $i\colon \infcatname{C}^{ \textrm{$n$-truncated} } \hookrightarrow \infcatname{C}$ preserves all existing limits in $\infcatname{C}$, so that if $\infcatname{C}$ is presentable we have an adjunction
		\[
			\adjunction{\mathfrak{t}_{ \leq n }}{\infcatname{C}}{\infcatname{C}^{ \textrm{$n$-truncated} }}{i}.
		\]
		The next proposition gives a simple criterion for a functor to preserve truncatedness.
		\begin{prop}[{\cite[Proposition 5.5.6.16]{Lurie:HTT}}]
			\label{prop:left_exact_functor_preserve_truncatedness}
			Let $\infcatname{C}, \infcatname{D}$ be $\infty$-categories that admit finite limits, and let $F\colon \infcatname{C} \to \infcatname{D}$ be a left exact functor between them.
			Then $F$ preserve $k$-truncated morphisms for all $k \geq -2$. In particular we have a commutative diagram
			\begin{diag}
				\infcatname{C} \ar[r, "F"] & \infcatname{D} \\
				\infcatname{C}^{\textrm{$k$-truncated}} \ar[u, hook] \ar[r, "F"] & \infcatname{D}^{\textrm{$k$-truncated}} \ar[u, hook].
			\end{diag}
		\end{prop}
		\begin{remark}[{\cite[Warning 1.2.1.9]{Lurie:HA}}]
			\label{remark:warning_truncation_stable}
			If $\infcatname{C}$ is a stable $\infty$-category then there are no truncated objects besides $0$. In fact, for $0 \neq X \in \infcatname{C}$ we can consider 
			\[
				0 \neq \id_{ X } \in \pi_{ n+1 }\left( \Map_{ \infcatname{C} }(X[-n-1], X) \right) \simeq \pi_{ 0 }(\Map_{ \infcatname{C} }(X, X)),
			\]
			giving us non-trivial homotopy classes (in different mapping spaces) for each $n$. 
		\end{remark}

		Let us now introduce the notion of almost compact objects, which is a weaker condition than being compact.

		\begin{defn}[{\cite[Definition 7.2.4.8]{Lurie:HA}}]
			\label{defn:almost_compact_object}
			Let $\infcatname{C}$ be a compactly-generated $\infty$-category and $X \in \infcatname{C}$.
			We say that $X$ is \emph{almost compact} if $\mathfrak{t}_{ \leq k }(X)$ is compact in $\infcatname{C}^{ \textrm{$k$-truncated} }$, for each $k \geq 0$.
			Denote the full subcategory of almost compact objects by $\infcatname{C}^{ \textrm{almost-cpt} }$.
		\end{defn}

		Let us now choose a stable $\infty$-category $\infcatname{C}$ endowed with a t-structure (our main example will be a module category over a connective ring spectrum with the standard t-structure). There is an interesting interplay between truncated objects in the connective part and general t-truncated objects (that is, objects that are truncated for the t-structure). It is mentioned in \cite[Warning 1.2.1.9]{Lurie:HA}, but we will state it as a proposition.

		\begin{prop}
			\label{prop:truncated_object_connective_part}

			Let $\infcatname{C}$ be a stable $\infty$-category endowed with a t-structure and $k\geq 0$. Then 
			\[
				(\infcatname{C}^{ \leq 0 })^{ \textrm{$k$-truncated} } \simeq \infcatname{C}^{ \leq 0} \cap \infcatname{C}^{ \geq -k }.
			\]
		\end{prop}
		\begin{proof}
			The idea is to observe that $X \in \infcatname{C}^{ \leq 0, \geq -k }$ if and only if
			\[
				\pi_{ 0 }\left( \Map_{ \infcatname{C}^{ \leq 0 }}(Y, X) \right) \simeq 0
			\]
			for each $Y \in \infcatname{C}^{ \leq -k-1 }$. Now we observe that
			\[
				\pi_{ 0 }\left( \Map_{ \infcatname{C}^{ \leq 0 }}(Y, X) \right) \simeq \pi_{ k+1 }\left( \Map_{ \infcatname{C}^{ \leq 0 } }(Y[-k-1], X) \right)
			\]
			and this last space being $0$ is basically $X$ being $k$-truncated in $\infcatname{C}^{ \leq 0 }$.
		\end{proof}

		An easy corollary is the following (be it for spectra or for $\C$-modules).
		\begin{coroll}
			\label{coroll:truncated_spectra}

			A connective $\C$-module $X$ is $k$-truncated if and only if $H^{ i }(X) \simeq 0$ for $i < -k$. 
		\end{coroll}

		We also observe, en passant, that there are no nontrivial truncated ``algebras'' in a stable $\infty$-category.

		\begin{prop}
			\label{prop:no_truncated_algebras}

			Let $\infcatname{C}$ be a presentable stable symmetric monoidal $\infty$-categories and $\mathcal{O}^{ \otimes }$ an $\infty$-operad. Then there are no non-trivial $k$-truncated objects in $\Alg_{ \mathcal{O} }(\infcatname{C})$.
		\end{prop}
		\begin{proof}
			The idea is that we can choose a color $X \in \mathcal{O}$, consider the ``evaluation functor at $X$''
			\[
				\ev_{ X }\colon \Alg_{ \mathcal{O} }(\infcatname{C}) \longrightarrow \Alg_{ \mathcal{T}riv }(\infcatname{C}) \simeq \infcatname{C}
			\]
			and apply \cite[Proposition 3.1.3.13]{Lurie:HA} to obtain that it admits a left adjoint, and therefore it preserves limits.
			By \cref{prop:left_exact_functor_preserve_truncatedness}, this means that all $k$-truncated $\mathcal{O}$-algebras in $\infcatname{C}$, when evaluated in $X$, are $0$, since $\tau_{ \leq k }(\infcatname{C}) \simeq 0$ thanks to \cref{remark:warning_truncation_stable}.
			Varying $X$, we obtain that each $k$-truncated $\mathcal{O}$-algebra is zero when seen as a functor $\mathcal{O} \to \mathcal{C}$. Therefore, we conclude.
		\end{proof}

		Some $\infty$-operads of interests are $\E_{ n }$, $0 \leq n \leq \infty$.
		Sometimes we will consider, though, $k$-truncated objects in \emph{connective} $\E_{ n }$-algebras in $\infcatname{C}$. Here we describe how they can be recognized at the level of their underlying object.
				
		\begin{prop}
			\label{prop:truncated_connective_algebras}
			
			Let $\infcatname{C}$ be a stable presentable symmetric monoidal $\infty$-category endowed with a t-structure compatible with the monoidal structure, so that $\infcatname{C}^{ \cn }$ inherits a symmetric monoidal structure.
			Then for each $1 \leq n \leq \infty$ we have
			\[
				\left( \Alg_{ \E_{ n } }(\infcatname{C}^{ \cn }) \right)^{ \textrm{$k$-truncated} } \simeq \Alg_{ \E_{ n } }(\infcatname{C}^{ \cn }) \times_{ \infcatname{C}^{ \cn } } (\infcatname{C}^{ \cn })^{ \textrm{$k$-truncated} };
			\]
			that is, $k$-truncated connective $\E_{ n }$-algebras are connective algebras whose underlying object in $\infcatname{C}$ is $k$-truncated, for each $k \geq 0$.

		\end{prop}

		\begin{proof}
			The forgetful functor $\oblv\colon \Alg_{ \E_{ n } }(\infcatname{C}^{ \cn }) \to \infcatname{C}^{ \cn }$ is right adjoint and therefore, by \cref{prop:left_exact_functor_preserve_truncatedness}, preserves $k$-truncated objects.
			We have then an induced functor
			\begin{diag}
				\Alg_{ \E_{ n } }(\infcatname{C}^{ \cn })^{\textrm{$k$-truncated}} \ar[drr, "\oblv", bend left=20] \ar[ddr, hook, bend right=20] \ar[dr, dashed, red] & & \\
																																									  & \Alg_{ \E_{ n } }(\infcatname{C}^{ \cn }) \times_{ \infcatname{C}^{ \cn } } (\infcatname{C}^{ \cn })^{ \textrm{$k$-truncated} } \ar[r] \ar[d, hook] \ar[dr, phantom, very near start, "\lrcorner"] & (\infcatname{C}^{\cn})^{\textrm{$k$-truncated}} \ar[d, hook] \\
																										              & \Alg_{\E_n}(\infcatname{C}^{\cn}) \ar[r, "\oblv"] & \infcatname{C}^{\cn}
			\end{diag}
			It is fully faithful by $2$-out-of-$3$, so to prove that it's an equivalence it suffices to show that it is essentially surjective. Equivalently, we can prove that each connective $\E_{ n }$-algebra with underlying object $k$-truncated is already $k$-truncated as an algebra. Let $B$ be such an algebra. We need to prove that for each $A \in \Alg_{ \E_{ n } }(\infcatname{C}^{ \cn })$ the mapping space 
			\[
				\Map_{ \Alg_{ \E_{ n } }(\infcatname{C}^{ \cn }) }(A, B)
			\]
			is $k$-truncated. Since $k$-truncated spaces are stable under limits and each algebra is a sifted colimits of free ones (where we denote $\Free\colon \infcatname{C}^{ \cn } \to \Alg_{ \E_{ n } }(\infcatname{C}^{ \cn })$ the left adjoint of $\oblv$), we can reduce to the case where $A = \Free(M)$ for $M \in \infcatname{C}^{ \cn }$. Then we have
			\[
				\Map_{ \Alg_{ \E_{ n } }(\infcatname{C}^{ \cn }) }(\Free(M), B) \simeq \Map_{ \infcatname{C}^{ \cn } }(M, \oblv(B))
			\]
			and we conclude since, by assumption, $\oblv(B)$ is a $k$-truncated object of $\infcatname{C}^{ \cn }$.
		\end{proof}
	
	\subsection{Some finiteness conditions on algebras}
		\label{subsubsection:finiteness_conditions_algebras}
		
		In this subsection we will recall some finiteness conditions on algebras, which will be useful later on. We will focus either on $\E_{ 1 }$-algebras or $\E_{ \infty }$-algebras. Our main interest is the connective case.

		\begin{defn}[{\cite[Definition 7.2.4.16, 7.2.4.30]{Lurie:HA}}]
			\label{defn:coherent_noetherian_rings}
			Let $R \in \Alg_{ \C }$; we say that $R$ is \emph{left coherent} if 
			\begin{enumerate}[label=(\roman*)]
				\item $R$ is connective;
				\item $H^{ 0 }(R)$ is a classical left coherent ring (i.e.\ all its finitely generated left ideals are finitely presented left $H^{ 0 }(R)$-modules);
				\item each $H^{ i }(R)$ is a finitely presented left $H^{ 0 }(R)$-module.
			\end{enumerate}
			Suppose now that $A \in \CAlg^{ \cn }_{ \C }$. We say that $A$ is \emph{Noetherian} if $A$ is left coherent (when $A$ is regarded as an $\E_{ 1 }$-ring) and $H^{ 0 }(A)$ is a classical commutative Noetherian ring (i.e.\ all its ideals are finitely generated $H^{ 0 }(A)$-modules).
			We denote by $\CAlg^{ \cn, \coh }_{ \C }$ and $\CAlg^{ \cn, \mathrm{Noeth} }_{ \C }$ the full subcategories of $\CAlg_{ \C }^{ \cn }$ spanned by coherent and Noetherian connective commutative $\C$-algebras, respectively.
		\end{defn}
		
		The condition of left coherence on a ring makes the $\infty$-category of almost perfect left modules (see \cref{defn:almost_perfect_modules}) well-behaved (see for example the statements of \cref{prop:almost_perfect_modules_over_coherent}, \cref{prop:almost_perfect_symmetric_monoidal}).\\
		Let us report \cite[Definition 7.2.4.26]{Lurie:HA} in the commutative case. It can be read exactly in the same way for $\E_{ n }$-algebras for $n \geq 1$.

		\begin{defn}[Commutative algebras locally of finite presentation]
			\label{defn:locally_finite_presentation}

			A morphism $A \to B$ in $\CAlg^{ \cn }_{ \C }$ is said to be \emph{locally of finite presentation}  if
			\[
				B \in \left( \CAlg_{ A/- }^{ \cn } \right)^{ \omega }.
			\]
			The full subcategory of locally finitely presented commutative algebras under $A$ is denoted by $\CAlg_{ A/- }^{ \lfp }$.
		\end{defn}

		\begin{defn}[Commutative algebras almost of finite presentations]
			\label{defn:almost_finite_presentation}

			A morphism $A \to B$ in $\CAlg^{\cn  }_{ \C }$ is said to be \emph{almost of finite presentation} if $B$ is an almost compact object of $\CAlg^{ \cn }_{ A/- }$.
			The full subcategory of almost finitely presented commutative algebras under $A$ is denoted by $\CAlg_{ A/- }^{ \afp }$.
		\end{defn}

		\begin{remark}
			\label{remark:connective_almost_finite_presentation}
			Observe how in \cref{defn:almost_finite_presentation} it is necessary to consider \emph{connective} algebras: in fact there are no non-trivial truncated objects in $\CAlg_{ A/- }$ thanks to \cref{prop:no_truncated_algebras}. We believe that, not restricting to connective algebras, is a typo in \cite[Definition 7.2.4.26]{Lurie:HA}. Moreover, we also believe that all statements in \cite[Proposition 7.2.4.27]{Lurie:HA} must be taken for connective objects: for example, the embedding $\CAlg^{ \cn } \hookrightarrow \CAlg$ commutes with filtered colimits and therefore $\Ind(\CAlg^{ \lfp })$ can only give us connective objects.
		\end{remark}

		The next proposition is again written for commutative algebras but it is valid, more generally, for $\E_{ n }$-algebras with $n \geq 1$.
		
		\begin{prop}[{\cite[Proposition 7.4.2.27]{Lurie:HA}}]
			\label{prop:connective_algebras_compactly_generated}
			Let $A \in \CAlg^{ \cn }_{ \C }$. The $\infty$-category of connective commutative algebras under $A$ is compactly generated by the locally finitely presented $A$-algebras. That is, we have
			\[
				\CAlg^{ \cn }_{ A/- } \simeq \Ind(\CAlg^{ \lfp }_{ A/- }).
			\]
		\end{prop}

		Let us remark how, a priori, it can be non trivial to verify if an algebra is locally of finite presentation or almost finitely presented. We will see later that, under certain circumstances, we can compute a relative cotangent complex (and hence ``linearize'' the problem) to be able to answer.
		Here's a criterion for identifying \emph{commutative} algebras almost of finite presentation (under $\C$ or, more generally, under any Noetherian $\E_{ \infty }$-algebra). It can be thought as an $\infty$-categorical version of the Hilbert basis theorem.

		\begin{prop}[{\cite[Proposition 7.2.4.31]{Lurie:HA}}]
			\label{prop:hilbert_basis_thm}
			Let $A \to B$ a morphism in $\CAlg_{ \C }^{ \cn }$ and suppose that $A$ is Noetherian. Then $B$ is almost of finite presentation over $A$ if and only if $B$ is Noetherian and $H^{ 0 }(B)$ is a finitely generated $H^{ 0 }(A)$-algebra.
		\end{prop}
	
		\begin{remark}
			The statement of \cref{prop:hilbert_basis_thm} explains why, when working over $\C$, we can identify commutative $\C$-algebras almost of finite presentation with those $A \in \CAlg^{ \cn }_{ \C }$ such that $H^{ 0 }(A)$ if a finitely presented classical $\C$-algebra and each $H^{ j }(A)$ is a finitely presented $H^{ 0 }(A)$-module.
			This characterization is often taken as a definition (and almost compact objects need not to be mentioned): see \cite[Chapter 2, 1.7.1]{GR:DAG1}.
		\end{remark}

		\begin{coroll}
			\label{coroll:chain_inclusion_algebras_C}
			We have the following chain of inclusion of full subcategories
			\[
				\CAlg^{ \lfp }_{ \C } \hookrightarrow \CAlg^{ \afp }_{ \C } \hookrightarrow \CAlg_{ \C }^{ \Noeth } \hookrightarrow \CAlg_{ \C }^{ \cn, \coh } \hookrightarrow \CAlg_{ \C }^{ \cn }.
			\]
		\end{coroll}

	\subsection{Some finiteness conditions on modules}
		\label{subsubsection:finiteness_conditions_modules}

		In this subsection we will recall some finiteness conditions on modules, which will be useful later on. We will introduce some notions for left modules over a connective $\E_{ 1 }$-ring spectrum (which can be read for right modules as well) and some others for modules over an $\E_{ \infty }$-ring spectrum. We always assume that the base is $\C$, but all constructions can be carried on over $\mathbb{S}$.
		For a more thorough study, which considers also $\E_{ n }$ for $2 \leq n < \infty$, we refer the reader to \cite[Chapter 4, 7]{Lurie:HA} and \cite[Chapter 2.7, Chapter 4]{Lurie:SAG}. 
		First of all, let us record the basic properties of module $\infty$-categories that we always use.

		\begin{prop}
			\label{prop:basic_properties_left_modules}
			Let $R \in \Alg_{ \C }$. The following assertions are true:
			\begin{enumerate}[label=(\roman*)]
				\item the $\infty$-category $\LMod_{ R }$ is stable and presentable;
				\item for any morphism $f\colon R \to T$ of $\E_{ 1 }$-algebras there exists a monadic pullback-pushforward adjunction 
					\[
						\adjunction{f^{ * }}{\LMod_{ R }}{\LMod_{ T }}{f_{ * }},
					\]
					where $f_{ * }$ preserves all (small) limits and colimits, and $f_{ * }f^{ * } \simeq {}_{ R }T_{ R } \otimes_{ R } -$;
				\item if $R$ is connective, then there exists a canonical t-structure $(\LMod_{ R }^{ \leq 0 }, \LMod_{ R }^{ \geq 0 })$ on $\LMod_{ R }$, defined by the homotopy groups of the underlying spectra, which is accessible, left and right complete and whose heart $(\LMod_{ R })^{ \heartsuit }$ is equivalent to the abelian $1$-category of left $H^{ 0 }(R)$-modules;
				\item if $R$ is connective then $\LMod_{ R }^{ \leq 0 }$ is the smallest full subcategory of $\LMod_{ R }$ that contains $R$ and is stable under colimits;
				\item if $f\colon R \to T$ is a morphism of connective $\E_{ 1 }$-algebras, then $f^{ * }$ is right t-exact and $f_{ * }$ is t-exact.
			\end{enumerate}
		\end{prop}
		\begin{proof}
			The stability of $\LMod_{ R }$ is proved in \cite[Corollary 7.1.1.5]{Lurie:HA}.	The second point is proved in \cite[Corollary 4.2.4.8]{Lurie:HA} and the formula comes from \cite[Proposition 4.6.2.17]{Lurie:HA}. The fact that the adjunction is monadic is an application of (the simple version of) Barr-Beck-Lurie, see \cite[Theorem 4.7.0.3]{Lurie:HA}.
			The existence and the properties of the t-structure on $\LMod_{ R }$, for a connective $R$, are proved in \cite[Proposition 7.1.1.13]{Lurie:HA}.
			For the final point, it is clear, by definition, that $f_{ * }$ is t-exact since it doesn't change the underlying object. This implies that $f^{ * }$ is right t-exact by an easy adjoint computation.
		\end{proof}

		\begin{prop}
			\label{prop:basic_properties_commutative_modules}
			Let $A \in \CAlg_{ \C }$. The following assertions are true:
			\begin{enumerate}[label=(\roman*)]
				\item there exists a symmetric monoidal structure on $\Mod_{ A }$, given by relative tensor product over $A$;
				\item if $f\colon A \to B$ is a morphism of $\E_{ \infty }$-algebras then $f^{ * }\colon \Mod_{ A } \to \Mod_{ B }$ is symmetric monoidal;
				\item if $A$ is connective, then the t-structure $(\Mod_{ A }^{ \leq 0 }, \Mod_{ A }^{ \geq 0 })$ is compatible with the symmetric monoidal structure.
			\end{enumerate}
		\end{prop}
		\begin{proof}
			The first point is \cite[Theorem 4.5.2.1]{Lurie:HA}. The second point is \cite[Theorem 4.5.3.1]{Lurie:HA}. The third point is proved in \cite[Lemma 7.1.3.10]{Lurie:HA}.
		\end{proof}

		Let us now introduce a common construction.

		\begin{defn}
			\label{defn:thick_subcategories}
			Let $\infcatname{C}$ be a stable $\infty$-category closed under retracts and $K \in \infcatname{C}$. The \emph{thick subcategory generated by $K$} is defined as the smallest stable full subcategory of $\infcatname{C}$ that contains $K$ and is closed under retracts. It is denoted by $\Thick(K)$.
		\end{defn}
		
		The most important example of thick subcategory will be perfect modules. 

		\begin{defn}[{\cite[Definition 7.2.4.1]{Lurie:HA}}]
			\label{defn:perfect_modules}
			Let $R \in \Alg_{ \C }$ and $M \in \LMod_{ R }$. We say that $M$ is \emph{perfect} if it belongs to $\Thick(R)$, as defined in \cref{defn:thick_subcategories}. We will denote this full subcategory with $\LMod_{ R }^{ \mathrm{perf} }$ and, if $R \in \CAlg_{ \C }$, with $\Perf(R)$ or $\Perf_{ R }$.
		\end{defn}

		Let us record some of its main properties.
		\begin{prop}
			\label{prop:stuff_about_perfect_left_modules}
			Let $R \in \Alg_{ \C }$. The following assertions are true:
			\begin{enumerate}[label=(\roman*)]
				\item the $\infty$-category $\LMod_{ R }$ is compactly generated and $(\LMod_{ R })^{ \omega } \simeq \LMod_{ R }^{ \mathrm{perf} }$;
				\item if $R$ is connective then any perfect left $R$-module $M$ lies in $\LMod_{ R }^{ - }$.
			\end{enumerate}
		\end{prop}
		\begin{proof}
			The first point is proved in \cite[Proposition 7.2.4.2]{Lurie:HA} and the second point in \cite[Corollary 7.2.4.5]{Lurie:HA}.
		\end{proof}

		Perfect modules for an $\E_{ \infty }$-ring admit a ``monoidal'' characterization.

		\begin{prop}
			\label{prop:perfect_modules_dualizable}
			Let $A \in \CAlg_{ \C }$. Then the dualizable objects of $(\Mod_{ A }, \otimes_{ A })$ coincide with $\Perf(A)$.
		\end{prop}
		\begin{proof}
			See the proof of \cite[Proposition 6.2.6.2]{Lurie:SAG} or \cite[Lemma 3.5]{BFN:IntegralTransformsDrinfeldCentersDAG}.
		\end{proof}

		If $A$ is a classical Noetherian commutative algebra and $M$ is a discrete $A$-module, then $M$ is perfect if and only if it admits a finite resolution by projective $A$-modules of finite rank; that is, $M$ is perfect if and only if it has finite projective dimension.
		What happens more often in practice is that $M$ is finitely presented over $A$, so that it admits a resolution by projective $A$-modules of finite type which might also be infinite.
		This condition corresponds to $M$ being an \emph{almost perfect} $A$-module, as we will now define.
		
		\begin{defn}[{\cite[Definition 7.2.4.10]{Lurie:HA}}]
			\label{defn:almost_perfect_modules}
			Let $R \in \Alg_{ \C }^{ \cn }$ and $M \in \LMod_{ R }$. We say that $M$ is \emph{almost perfect} if there exists $k \in \Z$ such that $M$ is an almost compact object, as in \cref{defn:almost_compact_object}, of $\LMod_{ R }^{ \leq k }$. 
			The full subcategory spanned by almost perfect modules is denoted by $\LMod_{ R }^{ \mathrm{aperf} }$; if $R$ is commutative then it is denoted by $\APerf(R)$ or $\APerf_{ R }$.
		\end{defn}	

		\begin{remark}
			\label{remark:pseudocoherent_modules}
			Almost perfect modules are sometimes called \emph{pseudocoherent}.
		\end{remark}

		Let us record some properties of almost perfect modules.

		\begin{prop}[{\cite[Proposition 7.2.4.11]{Lurie:HA}}]
			\label{prop:basic_properties_almost_perfect_modules}
			Let $R \in \Alg_{ \C }^{ \cn }$. The following assertions are true:
			\begin{enumerate}[label=(\roman*)]
				\item the $\infty$-category $\LMod_{ R }^{ \mathrm{aperf} }$ is a full subcategory of $\LMod_{ R }$, closed under retracts, which contains $\LMod_{ R }^{ \mathrm{perf} }$;
				\item the full subcategory $(\LMod_{ R }^{ \mathrm{aperf} })^{ \leq 0 } \hookrightarrow \LMod_{ R }$ is closed under geometric realizations of simplicial objects;
				\item each $M \in (\LMod_{ R }^{ \mathrm{aperf} })^{ \leq 0 }$ can be realized as a geometric realization of a simplicial left $R$-module $P_{ \bullet }$ where each $P_{ n }$ is a finite coproduct of copies of $R$.
			\end{enumerate}
		\end{prop}

		Let us observe how the third point of \cref{prop:basic_properties_almost_perfect_modules} generalizes the procedure of taking ``infinite free resolutions'' of finitely presented classical left modules.
		Checking if $M \in \LMod_{ R }^{ \leq k }$ is almost perfect amounts, by \cref{defn:almost_perfect_modules}, to check that each $\tau^{ \geq n }(M)$ is compact in $\LMod_{ R }^{ \geq n, \leq k }$ (where we used \cref{prop:truncated_object_connective_part}). This is non trivial in general.
		If $R$ respects some mild condition (for example $R$ is left Noetherian) then this task becomes an easier question that can be checked at the level of cohomologies.

		\begin{prop}[{\cite[Proposition 7.2.4.17, 7.2.4.18]{Lurie:HA}}]
			\label{prop:almost_perfect_modules_over_coherent}
			Let $R \in \Alg^{ \cn }$ be left coherent, as in \cref{defn:coherent_noetherian_rings}. Then $M \in \LMod_{ R }$ is almost perfect if and only if:
			\begin{enumerate}
				\item there exists $m \in \Z$ such that $M \in \LMod_{ R }^{ \leq m }$, i.e.\ $H^{ i }(M) \simeq 0$ for each $i > m$;
				\item each $H^{ j }(M)$ is a finitely presented left $H^{ 0 }(R)$-module.
			\end{enumerate}
			Moreover, the left coherency of $R$ implies (and is equivalent to) that the $\infty$-category $\LMod_{ R }^{ \mathrm{aperf} }$ inherits a t-structure $(\LMod_{ R }^{ \mathrm{aperf} } \cap \LMod_{ R }^{ \leq 0 }, \LMod_{ R }^{ \mathrm{aperf} } \cap \LMod_{ R }^{ \geq 0 })$, whose heart is the abelian $1$-category of finitely presented left $H^{ 0 }(R)$-modules.
		\end{prop}

		We now want to consider a connective commutative $\E_{ \infty }$-algebra and see how almost perfect modules relate with the symmetric monoidal structure.

		\begin{prop}
			\label{prop:almost_perfect_symmetric_monoidal}
			Let $A$ be a connective commutative algebra which is \emph{coherent}. Then $\APerf_{ A }$ is a full symmetric monoidal subcategory of $\Mod_{ A }$ and the t-structure is compatible with the symmetric monoidal structure.
		\end{prop}
		\begin{proof}
			It suffices to prove that $\APerf_{ A }$ is closed under relative tensor product over $A$. It clearly contains the unit $A$. By \cref{prop:almost_perfect_modules_over_coherent} we know that an almost perfect module $M$ over $A$ is such that $H^{ i }(M) = 0$ for $i \gg 0$ and $H^{j}(M)$ is a finitely presented module over $H^{ 0 }(A)$ for all $j \in \Z$. We want to prove that, for $M, N \in \APerf_{ A }$, $M \otimes_{ A } N$ is almost perfect.
			Up to suspension, we can assume that $M$ and $N$ are connective. This implies that $M \otimes_{ A } N$ is connective as well by \cref{prop:basic_properties_commutative_modules} (another reference is \cite[Corollary 7.2.1.23]{Lurie:HA}).
			To prove that each cohomology is finitely presented over $H^{ 0 }(A)$ we can use the converging (first quadrant) spectral sequence
			\[
				E_{ 2 }^{ p, q } \coloneqq \Tor_{ p }^{ H^{ -* }(A) }(H^{ -* }(M), H^{ -* }(N))(q) \Rightarrow H^{ -(p+q) }(M \otimes_{ A } N),
			\]
			which comes from \cite[Proposition 7.2.1.19]{Lurie:HA}.
			It suffices to prove that all terms in the $E_{ 2 }$-page are finitely presented $H^{ 0 }(A)$-modules: by coherence of $H^{ 0 }(A)$ we can conclude then that all the terms in the following pages and in the $E_{ \infty }$-page are finitely presented as well, and hence we conclude (an $H^{ 0 }(A)$-module endowed with a finite filtration whose associated graded terms are finitely presented, is finitely presented itself).
		\end{proof}

		Finally, with some boundedness condition on $A \in \CAlg^{\cn }_{ \C }$, we can also relate almost perfect $A$-modules with the duality. They won't be fully dualizable objects (they lack a co-evaluation map) but they will be reflexive (i.e.\ equivalent to their double $A$-linear dual).

		\begin{defn}
			\label{defn:eventually_coconnective}
			A module $M \in \Mod_{ \C }$ is said to be \emph{$n$-coconnective} if it belongs to $\Mod_{ \C }^{ \geq -n }$, i.e.\ if $H^{ j }(M) = 0$ for $j < -n$. A module $M \in \Mod_{ \C }$ is said to be \emph{eventually coconnective} if it is $n$-coconnective for some $n$. 
			We say that $M$ is \emph{bounded} if $M \in \Mod_{ \C }^{ \leq k }$ for some $k \in \Z$ and $M$ is eventually coconnective. We denote by $\Mod_{ \C }^{ \bdd }$ the full subcategory spanned by bounded modules.
			An algebra $A \in \Alg_{ \C }$ is eventually coconnective (or bounded) if the underlying module is so. We denote by $\Alg_{ \C }^{ \cn, \bdd }$ the full subcategory spanned by connective and bounded algebras.
		\end{defn}

		The following proposition is already known and used, for example in \cite[Corollary 4.2]{Nuiten:KoszulLieAlgbd}. Notice how we crucially use the fact that we are over $\C$ in the proof, which is made easy using a dg-model presentation. 

		\begin{prop}
			\label{prop:almost_perfect_duality_bounded}
			Let $A \in \CAlg^{ \cn }_{ \C }$ be coherent and assume that $A \in \Mod_{ \C }^{ \geq n, \leq 0}$ for some $n < 0$ (that is, $A$ is eventually coconnective). Then, the $A$-linear duality functor restricted to almost perfect $A$-modules 
			\[
				(-)^{ \vee } \coloneqq \IntHom_{ A }(-, A)\colon \APerf_{ A} ^{ \op } \longrightarrow \Mod_{ A }
			\]
			is symmetric monoidal.
		\end{prop}
		\begin{proof}
			The coherency of $A$ makes $\APerf_{ A }^{ \leq 0 }$ and $\APerf_{ A }$ full symmetric monoidal subcategories of $\Mod_{ A }$ thanks to \cref{prop:almost_perfect_symmetric_monoidal}.
			The functor $(-)^{ \vee }$ is lax symmetric monoidal (since $A$ is an algebra), i.e.\ we have natural morphisms
			\[
				\IntHom_{ A }(V, A) \otimes_{ A } \IntHom_{ A }(W, A) \to \IntHom_{ A }(V \otimes_{ A } W, A).
			\]
			We can assume, without loss of generality, that $V$ and $W$ are connective.
			Let us work in the (projective) model category of $\Mod_{ A }^{ \dg }$, where we chose a cofibrant cdga representative for $A$. Now, the $1$-functor of forgetting differentials is conservative, strong monoidal and closed (see \cref{remark:forget_differential}), so we can prove that the morphism above, once appropriate representatives are chosen everywhere, is an $1$-isomorphism of graded modules.
			We can choose explicit cofibrant representatives of $V$ and $W$ (without differentials) as
			\[
				V \simeq \bigoplus_{ k \geq 0 } A^{ n_{ k } }[k], \qquad W \simeq \bigoplus_{ k \geq 0 } A^{ m_{ k } }[k],
			\]
			where $n_{ k }, m_{ k } \in \N$ for each $k$. This comes from ``resolving'' $V$ and $W$ by free shifted copies of $A$, as in \cref{prop:basic_properties_almost_perfect_modules}.
			Observe now that, thanks to the boundedness of $A$, we have
			\[
				\IntHom_{ A }(V, A) \simeq \prod_{ k \geq 0 } A^{ n_{ k } }[-k] \simeq \bigoplus_{ k \geq 0 } A^{ n_{ k } }[-k], \qquad \IntHom_{ A }(W, A) \simeq \prod_{ k \geq 0 } A^{ m_{ k } }[-k] \simeq \bigoplus_{ k \geq 0 } A^{ m_{ k } }[-k].
			\]
			We conclude since tensor products commutes with direct sum.
		\end{proof}

		\begin{conj}
			\label{conj:almost_perfect_reflexive_spectra}
			We claim that the result of \cref{prop:almost_perfect_duality_bounded} also holds over $\mathbb{S}$.
		\end{conj}

		Finally, let us investigate how almost perfect modules interact with base change.
		\begin{prop}
			\label{prop:pullback_preserves_almost_perfect}
			Let $f\colon A \to B$ be a morphism in $\Alg^{ \cn }_{ \C }$. 
			Then the base change morphism preserves almost perfect (left) modules. That is, we have
			\[
				f^{ * } \simeq B \otimes_{ A } -\colon \LMod_{ A }^{ \mathrm{aperf} } \longrightarrow \LMod_{ B }^{ \mathrm{aperf} }.
			\]
		\end{prop}
		\begin{proof}
			Let $M \in \LMod_{ A }^{ \mathrm{aperf} }$ and let us assume, without loss of generality, that $M$ is connective.
			By \cref{prop:basic_properties_almost_perfect_modules}, we can approximate $M$ ``at each level'' with perfect left $A$-modules. That is, for each $k \in \Z$ there exists $E_{ k } \in \LMod^{ \mathrm{perf} }_{ A }$ equipped with a morphism $f_{ k }\colon E_{ k } \to M$ with $\fib(f_{ k }) \in \LMod_{ A }^{ < k }$.
			Tensoring with $B$ we obtain a fiber sequence of left $B$-modules
			\[
				\fib(B \otimes_{ A } f_{ k }) \simeq B \otimes_{ A } \fib(f_{ k }) \longrightarrow B \otimes_{ A } E_{ k } \longrightarrow B \otimes_{ A } M.
			\]
			We know that $B \otimes_{ A } E_{ k } \in \LMod^{\mathrm{perf} }_{ B }$ (because the right adjoint $f_{ * }$ preserves filtered colimits) and that $B \otimes_{ A } M \in \LMod_{ A }^{ \leq 0 }$ since $f^{ * }$ is right t-exact by \cref{prop:basic_properties_left_modules}. For the same reason, $\fib(B \otimes_{ A } f_{ k }) \in \LMod_{ B }^{ \leq k-1 }$.
			This is sufficient to conclude: we proved that $\tau^{ \geq k }(B \otimes_{ A } M) \simeq \tau^{ \geq k }(B \otimes_{ A } E_{ k })$ and this gives ``perfect approximations'' for each level.
		\end{proof}	
		Let us introduce the definition of ``Tor amplitude''. 
		\begin{defn}[{\cite[Definition 7.2.4.21]{Lurie:HA}}]
			\label{defn:tor_amplitude_left_module}
			Let $R \in \Alg_{ \C }^{ \cn }$ and $M \in \LMod_{ R }$. We say that $M$ has \emph{$\Tor$-amplitude $\leq n$}, for $n \geq 0$, if for any discrete right $R$-module $N$, we have 
			\[
				j < -n \Rightarrow H^{ j }(N \otimes_{ R } M) = 0.
			\]
		\end{defn}
		The main properties of $\Tor$-amplitude are explained in \cite[Proposition 7.2.4.23]{Lurie:HA}.\\
		Let us finally introduce the notion of coherent modules.
		\begin{defn}
			\label{defn:coherent_modules}
			Let $R \in \Alg_{ \C }^{ \cn }$ and $M \in \LMod_{ R }$. We say that $M$ is \emph{coherent} if it is almost perfect and eventually coconnective. We denote the (essentially small) full subcategory of coherent left $R$-modules by $\LMod_{ R }^{ \mathrm{coh} }$. If $R$ is commutative then we denote this subcategory by $\Coh(R)$ or $\Coh_{ R }$.
		\end{defn}

		\begin{prop}
			\label{prop:condition_for_coherence}
			Let $R \in \Alg_{ \C }^{ \cn }$ be a coherent algebra. A module $M \in \LMod_{ R }$ is coherent if and only if $M \in \LMod_{ R }^{ \geq n, \leq k }$ for some $n, k \in \Z$ and each $H^{ j }(M)$ is a finitely presented left $H^{ 0 }(R)$-module.
		\end{prop}
		\begin{proof}
			This is an immediate corollary of \cref{prop:almost_perfect_modules_over_coherent}.
		\end{proof}

		Perfect and coherent modules are two important small subcategories of the large category of all modules. 

		\begin{prop}
			\label{prop:perfect_modules_are_coherent_bounded_base}
			Let $R \in \Alg_{ \C }^{ \cn }$. All perfect left $R$-modules are coherent if and only if $R$ is bounded.
		\end{prop}
		\begin{proof}
			It is immediate by definition of perfect left $R$-modules as being the thick subcategory generated by $R$, see \cref{defn:perfect_modules}.
		\end{proof}

		\begin{example}
			\label{example:coherent_which_is_not_perfect}
			Any connective algebra $R$ which is not eventually coconnective is an example of perfect left module which is not coherent.
			An example of a $\C$-module which is coherent but not perfect is $\C[t]/(t^{ 2 })$: it is finitely presented but it has infinite projective dimension.
			An example of an almost perfect $\C$-module which is not coherent is $\oplus_{ n \in \N } \C[n]$.
		\end{example}

		\begin{remark}
			Observe that neither coherent nor eventually coconnective modules (over a commutative connective coherent algebra) are symmetric monoidal subcategories. They are not even non-unital symmetric monoidal. For example, $\C$ is a coherent (and not perfect) $\Sym_{ \C }(\C[1])$-module but 
			\[
				\C \otimes_{ \Sym_{ \C }(\C[1]) } \C \simeq \Sym_{ \C }(\C[2])
			\]
			is unbounded.
		\end{remark}

	\subsection{The cotangent complex}
		\label{subsubsection:algebraic_cotangent_complex}

		In this subsection we will briefly recall the construction and main properties of the cotangent complex for commutative algebras. We work over $\C$, so there will be no difference between the \emph{topological} cotangent complex and the \emph{algebraic} cotangent complex. For a comparison between the two, outside of characteristic $0$, see \cite[Chapter 25.3]{Lurie:SAG}.\\
		For a brief motivational note: the raison d'être of the cotangent complex is to give an ``homotopically correct'' version of the module of K{\"a}hler differentials, which is only well behaved on classical smooth commutative algebras. As with the K{\"a}hler differentials, one of the main interests of the cotangent complex is how it controls infinitesimal (derived) deformations of (derived) affine schemes.
		In fact, the cotangent complex can be interpreted as the non-abelian left derived functor of the functor sending a classical algebra $A$ to $\Omega^{ 1 }_{ A } \in \Mod_{ A }$. For a more extensive motivation and introduction to the cotangent complex, see \cite[§1]{Toen:DAG}.

		\begin{remark}
			\label{remark:history_cotangent_complex}
			Our exposition aims only to quickly recall the construction and main properties of the cotangent complex, and we will follow \cite[Chapter 7.3]{Lurie:HA}. Let us immediately point out that the constructions there are carried out in a greater generality, whereas we will specialize them to the case of commutative algebras. 
			The construction of the cotangent complex of commutative rings goes way back in history, before Lurie. Let us mention the works of André and Quillen (for the so called André-Quillen homology, see \cite{Andre:HomologieAlgebresCommutatives} and \cite{Quillen:CohomologyCommutativeRings} ) and of Grothendieck and Illusie (for the globalization to schemes and ringed topoi, see \cite{Grothendieck:CategoriesCofibreesComplexeCotangent} and \cite{Illusie:ComplexeCotangent1}).
		\end{remark}

		We will introduce first a ``point-wise'' definition of the cotangent complex, and then comment on how to recover the functoriality.
		Let us first recall the definition of split square zero extension of a commutative algebra.

		\begin{defn}
			\label{defn:split_square_zero_extension}
			Let $A \in \CAlg_{ \C }$ and $M \in \Mod_{ A }$. The composition
			\[
				A \oplus -\colon \Mod_{ A } \stackrel{\sim}{\longrightarrow} \Sp(\CAlg_{ -/A })\stackrel{\Omega^{ \infty }_{ \CAlg_{ -/A } }}{\longrightarrow} \CAlg_{ -/A }
			\]
			is called the \emph{split square-zero extension functor}.
		\end{defn}

		\begin{remark}
			\label{remark:split_square_zero_extension}
			The first equivalence in \cref{defn:split_square_zero_extension} comes from \cite[Corollary 7.3.4.14]{Lurie:HA}. The notation $A \oplus -$ is used because, once we forget the commutative algebra structure, we obtain $A \coprod -$ in spectra. That is, the composite acts as
			\[
				\Mod_{ A } \stackrel{A \oplus -}{\longrightarrow} \CAlg_{ -/A } \stackrel{\oblv_{ \CAlg} }{\longrightarrow} \Sp_{ -/A }, \qquad M \mapsto (A \oplus M \stackrel{pr_{ A }}{\longrightarrow} A).
			\]
			It it called \emph{split} square-zero extension because of the existence of the canonical section $A \to A \oplus M$ in $\CAlg_{ -/A }$. More abstractly, the canonical section arises from writing the functor $A \oplus -$ as the composite
			\[
				\Mod_{ A } \stackrel{\sim}{\longrightarrow} \Sp(\CAlg_{ -/A }) \simeq \Sp(\CAlg_{ A/-/A }) \stackrel{\Omega^{ \infty }}{\longrightarrow} \CAlg_{ A/-/A } \stackrel{\oblv_{ A/-} }{\longrightarrow} \CAlg_{ -/A },
			\]
			since the stabilisation of an $\infty$-category with a terminal object is equivalent to the stabilisation of its pointed version, by \cite[Remark 1.4.2.18]{Lurie:HA}. 
			Furthermore, it really deserves to be called ``split square zero extension'' by \cite[Remark 7.3.4.16]{Lurie:HA}; observe that, after passing to homotopy groups, we have an equivalence of graded commutative algebras
			\[
				\pi_{ * }(A \oplus M) \simeq \pi_{ * }(A) \oplus \pi_{ * }(M),
			\]
			where the right-hand side is the classical split square zero extension of the graded algebra $\pi_{ * }(A)$ by its graded module $\pi_{ * }(M)$.
		\end{remark}

		\begin{lemma}
			\label{lemma:split_square_zero_extension_preserves_limits}
			Let $A \in \CAlg_{ \C }$. The split square-zero extension functor
			\[
				A \oplus -\colon \Mod_{ A } \to \CAlg_{ -/A }
			\]
			is a right adjoint.
		\end{lemma}
		\begin{proof}
			This is immediate from the definition, since for any presentable $\infty$-category $\infcatname{C}$ the functor 
			\[
				\Omega^{ \infty }\colon \Sp(\infcatname{C}) \to \infcatname{C}
			\]
			is a right adjoint by \cite[Proposition 1.4.4.4]{Lurie:HA}.
		\end{proof}

		\begin{remark}
			\label{remark:split_square_zero_extension_in_families}
			The notion of split square zero extension from \cref{defn:split_square_zero_extension} can be made functorial in changing $A$. More precisely, we can consider the functor  
			\[
				(\CAlg_{ \C })^{ \Delta^{ 1 } } \stackrel{\ev_{ 1 }}{\longrightarrow} \CAlg_{ \C }, \qquad (A \to B) \mapsto B
			\]
			which is a presentable fibration (as in \cite[Definition 5.5.3.2]{Lurie:HTT}), whose fiber over $B$ is the presentable $\infty$-category $\CAlg_{ -/B }$.
			Stabilising all those fibers we obtain the \emph{tangent bundle} to $\CAlg_{ \C }$
			\begin{diag}
				T_{\CAlg} \ar[rr, "G"] \ar[dr, "p"] & & (\CAlg_{\C})^{\Delta^1} \ar[dl, "\ev_1"] \\
										  & \CAlg_{\C} &
			\end{diag}
			where $p$ is a presentable fibration satisfying the universal property described in \cite[Definition 7.3.1.9]{Lurie:HA}. The fiber of the tangent bundle over $B \in \CAlg_{ \C }$ is $\Sp(\CAlg_{ -/B })$. Recall the notation $\Mod(\Mod_{ \C })$ from \cite[Construction 4.5.1.3]{Lurie:HA}, which denotes the $\infty$-category whose objects are pairs $(A, M)$ where $A \in \CAlg_{ \C }$ and $M \in \Mod_{ A }$.
			The equivalence
			\[
				\Mod(\Mod_{ \C }) \simeq T_{ \CAlg_{ \C } }
			\]
			(of presentable fibrations over $\CAlg_{ \C }$) from \cite[Theorem 7.3.4.18]{Lurie:HA} allows us to conclude that the composite functor (over $\CAlg_{ \C }$)
			\[
				\Mod(\Mod_{\C}) \stackrel{\sim}{\longrightarrow} T_{\CAlg_{\C}} \to (\CAlg_{ \C })^{ \Delta^{ 1 } } 
			\]
			is the parametrized version of the split square zero extension functor. That is, the functor induced on the fibers at $A \in \CAlg_{ \C }$ is exactly the one introduced in \cref{defn:split_square_zero_extension}.
		\end{remark}
		
		\begin{defn}
			\label{defn:cotangent_complex_absolute_pointwise}
			Let $A \in \CAlg_{ \C }$. Consider the $\infty$-functor
			\[
				\CAlg_{ -/A } \stackrel{\Sigma^{ \infty }_{ \CAlg_{ -/A } }}{\longrightarrow} \Sp(\CAlg_{ -/A }) \stackrel{\sim}{\longrightarrow} \Mod_{ A }
			\]
			obtained as a left adjoint of the split square zero extension functor.
			Its value on $A \in \CAlg_{ -/A }$ is denoted by $\L_{ A/\C } \in \Mod_{ A }$ and is called the \emph{absolute cotangent complex of $A$} (here we use absolute to mean relative to $\C$).
		\end{defn}

		On more concrete terms, for any $M \in \Mod_{ A }$, we have the equivalence of spaces
		\[
			\Map_{ A }(\L_{ A/\C }, M) \simeq \Map_{ \CAlg_{ -/A } }(A, A \oplus M).
		\]
		As it happens in the classical case, the right-hand side can be interpreted as the space of (derived) $\C$-linear derivations from $A$ to $M$. We see then that $\L_{ A/\C }$ has the same (homotopical) universal property of the K{\"a}hler module of differential; for a more thorough discussion see the introduction of \cite[§7.4.1]{Lurie:HA}.

		\begin{remark}[Functoriality of the absolute cotangent complex]
			\label{remark:functoriality_absolute_cotangent_complex}
			As explained in \cref{remark:split_square_zero_extension_in_families} we can consider the diagram
			\begin{diag}
				\Mod(\Mod_{\C}) \ar[rr, "G"] \ar[dr, "p"] & & (\CAlg_{\C})^{\Delta^1} \ar[dl, "\ev_1"] \\
										  & \CAlg_{\C} &
			\end{diag}
			where the downward arrows are presentable fibrations. For the case of modules, this boils down to the fact that Cartesian lifts are obtained by restriction of scalars and coCartesian lifts by base change. As explained in \cite[Proposition 7.3.2.6, 7.3.2.14]{Lurie:HA}, $G$ admits a \emph{relative left adjoint} $F$. This simply means that $G$ admits a left adjoint on each fibers (proved in \cref{lemma:split_square_zero_extension_preserves_limits}) and that $G$ preserves locally Cartesian morphisms. The latter condition boils down to saying that, given a morphism of commutative algebras $f\colon A \to B$ and a $B$-module $M$, we have 
			\[
				(A \oplus f_{ * }(M)) \simeq A \times_{ B } (B \oplus M)
			\]
			as commutative algebras.
			The functoriality in $A$ of the absolute cotangent complex can then be obtained as a left adjoint of the composition
			\[
				\Mod(\Mod_{ \C }) \simeq T_{ \CAlg_{ \C } } \stackrel{G}{\longrightarrow} (\CAlg_{ \C })^{ \Delta^{ 1 } } \stackrel{\ev_{ 0 }}{\longrightarrow} \CAlg_{ \C },
			\]
			where we observed that $\CAlg_{ \C } \stackrel{\id}{\longrightarrow} (\CAlg_{ \C })^{ \Delta^{ 1 } }$ is the left adjoint of $\ev_{ 0 }$. That is, we have
			\[
				\L_{ (-) }\colon \CAlg_{ \C } \stackrel{\id}{\longrightarrow} (\CAlg_{ \C })^{ \Delta^{ 1 } } \stackrel{F}{\longrightarrow} \Mod(\Mod_{ \C }).
			\]
			This is recorded in \cite[Definition 7.3.2.14, 7.3.2.15]{Lurie:HA}.
		\end{remark}

		Having the fully fledged functoriality of the absolute cotangent complex, we can also define the relative one.

		\begin{defn}[{\cite[Definition 7.3.3.1]{Lurie:HA}}]
			\label{defn:relative_cotangent_complex}
			The \emph{relative cotangent complex} functor is defined as the composite
			\[
				(\CAlg_{ \C })^{ \Delta^{ 1 } } \stackrel{(\L_{ (-) })^{ \Delta^{ 1 } }}{\longrightarrow} (\Mod(\Mod_{ \C }))^{ \Delta^{ 1 } } \longrightarrow \Mod(\Mod_{ \C })
			\]
			where the last arrow takes the ``relative cofiber''.
			It acts by
			\[
				(A \to B) \mapsto \cofib(\L_{ A } \otimes_{ A } B \to \L_{ B }) \in \Mod_{ B }.
			\]			
		\end{defn}

		\begin{remark}
			\label{remark:exact_sequence_kahler}
			Observe how, by construction (\cref{defn:relative_cotangent_complex}), we have a generalization of the following well known fact in algebraic geometry that, given two morphisms $A \to B \to C$ of classical commutative algebras, there is always an exact sequence of $C$-modules
			\[
				\Omega^{ 1 }_{ B/A } \otimes_{ B } C \to \Omega^{ 1 }_{ C/A } \to \Omega^{ 1 }_{ C/B } \to 0.
			\]
			This is not always exact to the left (i.e.\ the first map is not always injective), unless some smoothness conditions are imposed.
			Let us appreciate how much cleaner the (homotopical) algebraic counterpart of this is: the above one becomes a fiber sequence of cotangent complexes, and no smoothness assumptions are needed.
			Even better: the way to ``continue to the left'' the sequence above is by considering the long exact sequence in cohomology arising from the fiber sequence of cotangent complexes, see \cref{prop:cotangent_complex_kahler}. This is related to the fact that the cotangent complex can be seen as a (non abelian) left derived functor of the K{\"a}hler differentials (on finite polynomial algebras).
		\end{remark}

		Let us record two basic properties of (relative) cotangent complexes.
		\begin{prop}[{\cite[Corollary 7.3.3.6]{Lurie:HA}}]
			\label{prop:fiber_sequence_cotangent_complex}
			Let $A \to B$ and $B \to C$ be two morphisms in $\CAlg_{ \C }$. There exists a fiber sequence in $\Mod_{ A }$ given by
			\[
				\L_{ B/C } \otimes_{ B } A \longrightarrow \L_{ A/C } \longrightarrow \L_{ A/B }.
			\]
		\end{prop}

		\begin{prop}[{\cite[Proposition 7.3.3.7]{Lurie:HA}}]
			Given a pushout diagram in $\CAlg_{ \C }$
			\begin{diag}
				A \ar[r] \ar[d] & B \ar[d] \\
				C \ar[r] & C \otimes_A B \ar[ul, very near start, phantom, "\ulcorner"]
			\end{diag}
			there exists a canonical equivalence $\L_{ C \otimes_{ A } B/B } \stackrel{\sim}{\longrightarrow} \L_{ C/A } \otimes_{ C } (C \otimes_{ A } B)$ in $\Mod_{ C \otimes_{ A } B }$.
		\end{prop}

		Let us now state more explicitly how the cotangent complex is related to the K{\"a}hler module of differentials.

		\begin{prop}[{\cite[Proposition 7.4.3.9]{Lurie:HA}}]
			\label{prop:cotangent_complex_kahler}
			Let $A \to B$ a morphism in $\CAlg^{ \cn }_{ \C }$. Then $\L_{ B/A }$ is a connective $B$-module and $H^{ 0 }(\L_{ B/A }) \simeq \Omega^{ 1 }_{ H^{ 0 }(B)/H^{ 0 }(A) }$ as classical $H^{ 0 }(B)$-modules.
		\end{prop}

		We will finish this subsection by mentioning how properties of the relative cotangent complex can be used to deduce global properties of the commutative algebras (for example their local finite presentation).

		\begin{thm}[{\cite[Theorem 7.4.3.18]{Lurie:HA}}]
			\label{thm:locally_finite_presentation_perfect_cotangent_complex}
			Let $f\colon A \to B$ be a map in $\CAlg_{ \C }^{ \cn }$ which is locally of finite presentation.
			Then the relative cotangent complex is perfect, i.e.\ 
			\[
				\L_{ B/A } \in \Perf(B).
			\]
			If instead $A \to B$ is almost of finite presentation, then 
			\[
				\L_{ B/A } \in \APerf(B).
			\]
			The converse implications hold whenever $H^{ 0 }(B)$ is a finitely presented $H^{ 0 }(A)$-algebra.
		\end{thm}
	
	\subsection{Stacks and prestacks}
		\label{subsubsection:stacks_prestacks}

		In this subsection we will recall the basic definitions of derived algebraic geometry: derived affine schemes, (derived) stacks and prestacks, formal completions and quasi-coherent sheaves. Our main reference is \cite[Chapter 2,3]{GR:DAG1} and \cite{TV:HAG2}. 
		Let us start by defining derived affine schemes and several useful subcategories.

		\begin{defn}[Derived affine schemes]
			\label{defn:derived_affine_schemes}
			\begin{enumerate}[label=(\arabic*)]
				\item The $\infty$-category of \emph{derived affine schemes} (over $\C$) is defined as
					\[
						\dAff_{ \C } \coloneqq (\CAlg^{ \cn }_{ \C })^{ \op }
					\]
					and the symbol $\Spec$ is used as expected.

				\item The full subcategory of \emph{locally finitely presented} derived affine schemes is defined as
					\[
						\dAff_{ \C }^{ \lfp } \coloneqq (\CAlg^{ \lfp }_{ \C })^{ \op }.
					\]
				\item The full subcategory of \emph{almost finitely presented} derived affine schemes is defined as 
					\[
						\dAff_{ \C }^{ \afp } \coloneqq (\CAlg^{ \afp }_{ \C })^{ \op }.
					\]
				\item The full subcategory of \emph{eventually coconnective} (or bounded) derived affine schemes is defined as
					\[
						\dAff_{ \C }^{ \bdd } \coloneqq (\CAlg^{ \cn, \bdd }_{ \C })^{ \op }.
					\]
			\end{enumerate}
			Combinations of the above superscripts are interpreted in the obvious ways.
		\end{defn}

		The most general objects we will consider are prestacks (the prefix ``derived'' is omitted, since we never consider non derived ones).

		\begin{defn}[Prestacks]
			\label{defn:prestacks}
			The $\infty$-category of \emph{prestacks} (over $\C$) is defined as the full subcategory
			\[
				\Fun_{ \mathrm{acc} }(\CAlg_{ \C }^{ \cn }, \S) \hookrightarrow \Fun(\CAlg_{ \C }^{ \cn }, \S)
			\]
			of accessible functors. It is denoted by $\PreStk_{ \C }$.
		\end{defn}

		To be able to talk about geometric objects we will need, as a bare minimum, to impose some étale descent condition. Let us recall how étale topology, in this homotopical context, can be defined.

		\begin{defn}[Étale topology]
			\label{defn:etale_topology}
			A morphism of derived affine schemes $f\colon \Spec A \to \Spec B$ is \emph{étale} if
			\begin{itemize}
				\item the morphism of classical commutative algebras $H^{ 0 }(B) \to H^{ 0 }(A)$ is étale (so, in particular, flat);
				\item the morphism $f$ is \emph{strong}, i.e.\ for each $j \leq 0$ it induces an equivalence
					\[
						H^{ j }(B) \otimes_{ H^{ 0 }(B) } H^{ 0 }(A) \stackrel{\sim}{\longrightarrow} H^{ j }(A).
					\]
			\end{itemize}
			An morphism $f\colon \Spec A \to \Spec B$ is an \emph{étale covering} if it is étale and $H^{ 0 }(f)$ is a surjective morphism of classical schemes.
		\end{defn}

		One can prove that étale coverings give a Grothendieck topology to $\dAff_{ \C }$, see \cite[Lemma 2.2.2.13]{TV:HAG2}.

		\begin{remark}
			\cref{defn:etale_topology} works also when replacing ``étale'' by ``flat'' or ``smooth'' or ``Zariski open immersion''.
		\end{remark}

		\begin{defn}[Stacks]
			The $\infty$-category of \emph{stacks} (over $\C$) is defined by 
			\[
				\Stk_{ \C } \coloneqq \Shv^{ \mathrm{hyp} }_{ \mathrm{et} }(\dAff_{ \C }).
			\]
		\end{defn}

		That is, stacks are prestacks that satisfy étale hyper-descent. The reason why hyper-descent is chosen over descent is to have an hypercomplete $\infty$-topos, i.e.\ where equivalences can be checked at the level of homotopy ``groups''; for more information see \cite{TV:HAG1} and \cite[Chapter 6]{Lurie:HTT}.
		Let us mention the crucial basic fact that all derived affine schemes are stacks.
		
		\begin{prop}
			\label{prop:affine_schemes_are_stacks}
			The Yoneda embedding
			\[
				\dAff_{ \C }^{ \op } \hookrightarrow \PreStk_{ \C }
			\]
			factors through $\Stk_{ \C }$.
		\end{prop}
		\begin{proof}
			See \cite[Chapter 2, Proposition 2.4.2]{GR:DAG1} and the references therein.
		\end{proof}

		We won't need more specific notions about stacks in this work; let us just mention that one can single out Deligne-Mumford stacks (equipped with an étale atlas from derived affine schemes) and $n$-Artin stacks (equipped with a smooth atlas from $(n-1)$-Artin stacks, defined in a recursive way). For more information about them, see \cite[Chapter 2.4]{GR:DAG1} and \cite[Chapter 2.2]{TV:HAG2} (where they are called \emph{geometric} stacks).
		Let us now introduce some finiteness conditions on prestacks, which will be of interest later one. For a more comprehensive explanation see \cite[Chapter 2.1]{GR:DAG1}.

		\begin{defn}[Prestacks locally almost of finite type]
			\label{defn:laft_prestacks}
			Consider the composite embedding
			\[
				\Fun(\CAlg_{ \C }^{ \cn, \bdd, \afp }, \S) \stackrel{\mathrm{LKE}}{\hookrightarrow} \Fun(\CAlg_{ \C }^{ \cn, \bdd }, \S) \stackrel{\mathrm{RKE}}{\hookrightarrow} \Fun(\CAlg_{ \C }^{ \cn }, \S).
			\]
			The prestacks (i.e.\ accessible functors) in the essential image are called \emph{locally almost of finite type}, often shortened to laft. They span a full subcategory of prestacks which is denoted by $\PreStk_{ \laft }$. 
		\end{defn}

		\begin{remark}
			\label{remark:other_definition_of_laft}
			A prestack $X$ is laft if and only if 
			\begin{enumerate}[label=(\arabic*)]
				\item $X$ is convergent, i.e.\ it respects Postnikov towers;
				\item the restriction of $X$ to $n$-coconnective derived affine schemes commutes with filtered colimits, for each $n$.
			\end{enumerate}
			See \cite[Chapter 2, 1.7]{GR:DAG1}.
		\end{remark}

		We can now define quasi-coherent sheaves (and their pullbacks) on prestacks.

		\begin{defn}[Quasi-coherent sheaves]
			\label{defn:quasi_coherent_sheaves_prestacks}
			The $\infty$-functor of quasi-coherent sheaves (and $*$-pullback) on prestacks is defined as right Kan extension
			\begin{diag}
				\CAlg_{\C}^{\cn} \ar[r, "\Mod"] \ar[d, hook] & \CAlg(\Pr^{L, st}_{\C}) \\
				\PreStk_{\C} \ar[ur, dashed, "\Qcoh^*"']
			\end{diag}
			where $\Mod$ sends $A$ to $\Mod_{ A }$ and $A \to B$ to the base change symmetric monoidal functor $B \otimes_{ A } -$ (see \cite[Corollary 4.8.5.21]{Lurie:HA}).
		\end{defn}

		\begin{remark}
			\label{remark:negative_qcoh_prestacks}
			
			Observe that, by right t-exactness of the functors $f^{ * }$ on modules (see \cref{prop:basic_properties_left_modules}), it makes sense to talk about $\Qcoh(X)^{ \leq n }$, for any $n$, on any prestack $X$. 
			Moreover, the same definition by right Kan extension, yields the $\infty$-category $\Perf(X)$. It corresponds to dualizable objects in $\Qcoh(X)$ by \cite[Chapter 3, Lemma 3.6.3]{GR:DAG1}.
		\end{remark}
		
		We will be content with this definition; for a more thorough study (for example of direct images, base change properties etc.) we direct the reader to \cite[Chapter 3]{GR:DAG1}.\\
		Let us define when a prestack admits a cotangent complex, which can be thought as a globalization of \cref{defn:cotangent_complex_absolute_pointwise} (which corresponds to the affine case).

		\begin{defn}[Cotangent complex of a prestack]
			\label{defn:cotangent_complex_prestack}
			Let $X \in \PreStk_{ \C }$. An object $\L_{ X } \in \Qcoh(X)^{ - }$ is a \emph{cotangent complex of $X$} if for every $x\colon \Spec A \to X$ the pullback $x^{ * }\L_{ X } \in \Mod_{ A }^{ - }$ satisfies the following universal property:
			\[
				\Map_{ A }(x^{ * }\L_{ X }, M) \simeq \Map_{ \PreStk_{ \Spec A/- } }(\Spec(A \oplus M), X)
			\]
			for every $M \in \Mod_{ A }^{ \leq 0 }$.
			Similarly, for a morphism of prestacks $f\colon X \to Y$, an object $\L_{ X/Y } \in \Qcoh(X)^{ - }$ is a \emph{relative cotangent complex} if, with the same notations as above, we have
			\[
				\Map_{ A }(x^{ * }\L_{ X/Y }, M) \simeq \Map_{ \PreStk_{ \Spec A/ - / Y } }(\Spec(A \oplus M), X)
			\]
			where we consider $\Spec(A \oplus M) \stackrel{d_{ 0 }}{\to} \Spec(A) \stackrel{f \circ x}{\to} Y$.
		\end{defn}

		\begin{remark}
			\label{remark:cotangent_complexes_unique}
			It is clear that a (relative) cotangent complex is unique (i.e.\ the space of them is contractible), if it exists.
		\end{remark}

		\begin{remark}[Pro-cotangent complex]
			\label{remark:pro_cotangent_complex}
			A prestack admits a cotangent complex as in \cref{defn:cotangent_complex_prestack} if, in particular, for each point $x\colon \Spec A \to X$, the functor
			\[
				\Mod_{ A }^{ \leq 0 } \stackrel{A \oplus -}{\longrightarrow} \CAlg^{ \cn }_{ -/A } \to \S, \qquad M \mapsto \Map_{ \PreStk_{ \Spec A/- } }(\Spec(A \oplus M), X)
			\]
			is ``co-represented'' by an object in $\Mod_{ A }^{ - }$.
			It can happen that the above functor is only pro-corepresented; that is, it's not co-represented by anything in $\Mod_{ A }^{ - }$, but preserves finite limits and therefore can be seen as an object of $\Pro(\Mod_{ A }^{ - })$. In this case, we say that $X$ admits a \emph{pro-cotangent space at $x$} (see \cite[Chapter 1, 2.2.3]{GR:DAG2}).
			If, moreover, all such pro-cotangent spaces are compatible with each other under $\Pro$-extension of $*$-pullbacks, we then say that $X$ admits a \emph{pro-cotangent complex}, that lives in a formal limit of $\infty$-categories of pro-objects (see \cite[Chapter 1, 4.1.4]{GR:DAG2}). 
		\end{remark}

		\begin{defn}[Infinitesimally cohesive]
			\label{defn:inf_cohesive}
			A prestack $X \in \PreStk_{ \C }$ is \emph{inf-cohesive} if it preserves pullback diagrams in $\CAlg_{ \C }^{ \cn }$ of the form 
			\begin{diag}
					A^{\eta} \ar[r, "f"] \ar[d] \ar[dr, phantom, very near start, "\lrcorner"] & A \ar[d, "d_{\eta}"] \\
					A \ar[r, "d_0"] & A \oplus M[1],
			\end{diag}
			where $A \in \CAlg_{ \C }^{ \cn }$, $M \in \Mod_{ A }^{ \leq 0 }$, $A \oplus M[1]$ is the split square-zero extension, $d_{ 0 }$ is the trivial derivation and $d_{ \eta } \in \pi_{ 0 }\Map_{ A }(\L_{ A }, M[1])$ (see \cref{defn:square_zero_extension}).
		\end{defn}	

		Let us now recall \cite[Chapter 1, 7.1.2]{GR:DAG2}.

		\begin{defn}[Prestacks with deformation theory]
			\label{defn:prestack_deformation_theory}
			A prestack $X \in \PreStk_{ \C }$ is said to \emph{admit deformation theory} if:
			\begin{enumerate}[label=(\alph*)]
				\item it is convergent, as defined in \cref{remark:other_definition_of_laft};
				\item it is inf-cohesive, as described in \cref{defn:inf_cohesive};
				\item it admits a pro-cotangent complex, as described in \cref{remark:pro_cotangent_complex}.
			\end{enumerate}
			We denote the full subcategory they span by $\PreStk_{ \mathrm{def} }$.
		\end{defn}

		\begin{remark}[How is this related to deformation theory?]
			\label{remark:explanation_deformation_theory}
			Let us consider $X \in \PreStk_{ \mathrm{def} }$ equipped with a point $x\colon \Spec A \to X$, and suppose we want to study (infinitesimal and derived) deformation of $X$ around $x$.
			That is, we want to study spaces of liftings like
			\begin{diag}
				\Spec A  \ar[d, hook, "\mathrm{nil-iso}"', "i"] \ar[r, "x"] & X \\
				\Spec B \ar[ur, dashed]
			\end{diag}
			where $i\colon \Spec A \hookrightarrow \Spec B$ is a nilpotent thickening. By convergence of $X$ we can assume that $A$ and $B$ are bounded.
			By \cite[Chapter 1, Proposition 5.5.3]{GR:DAG2}, we can write $\Spec A \hookrightarrow \Spec B$ as a finite composition of square-zero extensions. We assume then, without loss of generality, that $B$ is a square-zero extension of $A$ by $M \in \Mod_{ A }^{ \leq -1 }$ for a certain derivation $\eta \in \pi_{ 0 }\Map_{ A }(\L_{ A }, M[1])$.
			By inf-cohesiveness of $X$, we have a pullback diagram of spaces (over $X(A)$)
			\begin{diag}
				X(B) \ar[d, "X(i)"] \ar[r] \ar[dr, very near start, phantom, "\lrcorner"] & X(A) \ar[d, "X(d_0)"] \\
				X(A) \ar[r, "X(d_{\eta})"] & X(A \oplus M[1])
			\end{diag}
			and we are interested in the homotopy fiber at $x \in X(A)$. This can be seen as a pullback of spaces given by
			\begin{diag}
				\Map_{\PreStk_{\Spec A/-}}(\Spec B, X) \ar[r] \ar[d] \ar[dr, phantom, very near start, "\lrcorner"] & * \ar[d, "x \circ d_0"] \\
				* \ar[r, "x \circ d_{\eta}"] & \Map_{\PreStk_{\Spec A/-}}(\Spec(A \oplus M[1]), X).
			\end{diag}
			where we see the derivations as morphisms $\Spec(A \oplus M[1]) \to \Spec A$.
			We conclude by observing that the bottom right space is exactly the value of the pro-cotangent space at $x \in X$, denoted by $T_{ x }^{ * }(X) \in \Pro(\Mod_{ A }^{ - })$ evaluated at $M$.
			This explains exactly how, for prestacks that admit deformation theory, the pro-cotangent complex controls infinitesimal deformation theory.
		\end{remark}
		
		The final topics we want to review are formal completions of prestacks. Let us start defining the de Rham prestack, originally due to Simpson in \cite{Simpson:deRhamStack}.

		\begin{defn}[{\cite[Definition 18.2.1.1]{Lurie:SAG}}]
			\label{defn:de_rham_prestack}
			Let $X \in \PreStk_{ \C }$. Its \emph{de Rham prestack} $X_{ \dR } \in \PreStk_{ \C }$ is defined by
			\[
				\CAlg_{ \C }^{ \cn } \ni R \mapsto \colim_{ I \subseteq H^{ 0}(R) } X(H^{ 0 }(R)/I)
			\]
			where the colimit is indexed by the filtered poset of nilpotent ideals of $H^{ 0 }(R)$.
		\end{defn}

		The de Rham stack is particularly well behaved, and has a simpler formula, when restricted to laft prestacks. Since the prestacks of interest in \cite{GR:DAG2} are (almost always) laft, this explains the apparent discrepancy between \cref{defn:de_rham_prestack} and \cite[Chapter 4, 1.1.1]{GR:DAG2}.

		\begin{prop}
			\label{prop:de_rham_laft_prestack}
			Let $X \in \PreStk_{ \laft }$. Its de Rham prestack $X_{ \dR }$ is still locally almost of finite type and is defined by
			\[
				\CAlg_{ \C }^{ \cn, \bdd, \afp } \ni R \mapsto X(H^{ 0 }(R)^{ \red })
			\]
			where $H^{ 0 }(R)^{ \red } \simeq H^{ 0 }(R)/\sqrt{0}$.
			This yields a functor
			\[
				(-)_{ \dR }\colon \PreStk_{ \laft } \to \PreStk_{ \laft }
			\]
			which commutes with all (small) limits and colimits and which is, by construction, right adjoint to
			\[
				(-)^{ \red }\colon \PreStk_{ \laft } \to \PreStk_{ \laft }
			\]
			sending $Y$ to the left Kan extension along $\CAlg^{ \red }_{ \C } \hookrightarrow \CAlg^{ \cn, \afp, \bdd }_{ \C }$ of its restriction to $\CAlg^{ \red }_{ \C }$ (which is the $1$-category of classical reduced commutative $\C$-algebras).
		\end{prop}
		\begin{proof}
			Everything can be found in \cite[Chapter 4]{GR:DAG2}. For the first formula, it suffices to remark that if $R \in \CAlg^{ \cn, \bdd, \afp }$ then $H^{ 0 }(R)$ is finitely presented over $\C$ (hence Noetherian) and therefore has a maximal nilpotent ideal given by $\sqrt{0}$, which simplifies the colimit computation in \cref{defn:de_rham_prestack}.
		\end{proof}
		
		Let us record the observation that the derived part of $X$ (that is, its derived points) does not matter when taking its de Rham prestack. More precisely, not even the nilpotent elements of the classical part of $X$ matter as well. 

		\begin{prop}[{\cite[Lemma 1.1.7]{GR:CrystalsDMod}}]
			\label{prop:de_rham_prestack_is_classical}
			The functor $(-)_{ \dR }\colon \PreStk_{ \laft } \to \PreStk_{ \laft }$ is left Kan extended from its restriction to $(\CAlg_{ \C }^{ \red })^{ \op }$.
			Equivalently,
			\[
				(X)_{ \dR } \simeq (X^{ \red })_{ \dR }.
			\]
		\end{prop}

		Observe that, by \cref{prop:de_rham_laft_prestack}, we always have a canonical morphism
		\[
			p_{ Y }\colon Y \to Y_{ \dR }
		\]
		for $Y \in \PreStk_{ \laft }$. It is a nil-isomorphism by \cref{prop:de_rham_prestack_is_classical}. It is a regular epimorphism (i.e.\ $Y_{ \dR }$ is the geometric realization of the \v{C}ech nerve of $p$) whenever $Y$ is classically formally smooth, see \cite[Lemma 1.2.4]{GR:CrystalsDMod}.

		Let us finally define formal completions.

		\begin{defn}[{\cite[Definition 6.1.2]{GR:DGIndSch}}]
			\label{defn:formal_completion}
			Let $X \to Y$ a morphism of laft prestacks. The \emph{formal completion of $Y$ along $X$} is the prestack obtained by the following pullback
			\begin{diag}
				(Y)_X^{\wedge} \ar[r] \ar[d] \ar[dr, phantom, very near start, "\lrcorner"] & X_{\dR} \ar[d, "f_{\dR}"] \\
				Y \ar[r, "p_Y"] & Y_{\dR}.
			\end{diag}
			Note how it depends on $f$, although it is suppressed in the notation.
		\end{defn}

		Let us record its basic properties.

		\begin{prop}
			\label{prop:basic_properties_formal_completion}
			Let $f\colon X \to Y$ in $\PreStk_{ \laft }$ and consider $(Y)^{ \wedge }_{ X }$ as in \cref{defn:formal_completion}. The following assertions are true:
			\begin{enumerate}[label=(\roman*)]
				\item the morphism $f$ can be factored as 
					\[
						X \to (Y)^{ \wedge }_{ X } \to Y
					\]
					where the first morphism is a nil-isomorphism;
				\item $(Y)^{ \wedge }_{ X }$ is only sensible to $X^{ \red }$;
				\item if $f\colon X^{ \red } \to Y^{ \red }$ is a monomorphism then the prestack $(Y)^{ \wedge }_{ X }$ is laft;
				\item if $Y$ is a stack then $(Y)^{ \wedge }_{ X }$ is a stack too.
			\end{enumerate}
		\end{prop}
		\begin{proof}
			They can all be found in \cite[Chapter 6]{GR:DGIndSch}.
		\end{proof}

%% file: IndCoh.tex
\section{Ind-coherent sheaves}
	\label{subsection:ind_coherent_sheaves}
		Ind-coherent sheaves are a certain renormalisation (``at $-\infty$'') of quasi-coherent sheaves that is better suited for deformation theory and geometric representation theory. The main references are \cite{GR:DAG1} and \cite{GR:DAG2}, where all functorialities are studied in details. The seminal paper \cite{G:IndCoh} is also relevant, especially for Serre duality.
		Let us recall some of the main steps in its construction.

		\begin{defn}
			\label{defn:ind_coherent_modules}
			Consider a Noetherian derived scheme $S$, so that $\Coh(S) \hookrightarrow \Qcoh(S)$ is a well defined small stable subcategory, given by bounded quasi-coherent sheaves with coherent cohomologies (which are seen as quasicoherent sheaves on $t_{ 0 }(S)$, the classical truncation of $S$).
			The $\infty$-category of \emph{ind-coherent sheaves on $S$} is defined as $\Ind(\Coh(S))$.
		\end{defn}
		If $S = \Spec A$, we clearly have $\Coh(S) = \Coh_{ A } \hookrightarrow \Mod_{ A }$.
		What gets the theory started is the ``realisation'' morphism
		\[
			\psi_{ S }\colon \Ind(\Coh(S)) \longrightarrow \Qcoh(S)
		\]
		which is the left Kan extension of $\Coh(S) \hookrightarrow \Qcoh(S)$. Endowing $\Ind(\Coh(S))$ with the accessible t-structure whose connective part is $\Ind(\Coh(S)^{ \leq 0 })$ (see \cite[§1.2.1]{G:IndCoh}) one proves that $\psi_{ S }$ identifies $\Qcoh(S)$ with the left t-completion of $\Ind(\Coh(S))$.
		Of uttermost importance is \cite[Proposition 1.2.4]{G:IndCoh} that proves that 
		\[
			\psi_{ S }^{ \geq n }\colon \Ind(\Coh(S))^{ \geq n } \stackrel{\sim}{\longrightarrow} \Qcoh(S)^{ \geq n }
		\]
		is an equivalence for each $n$. This explains how $\Ind(\Coh(S))$ and $\Qcoh(S)$ only differs at $-\infty$, i.e.\ at $\Ind(\Coh(S))^{ \leq -\infty }$.
		This allows us to transfer functorialities from $\Qcoh$, the first of which is the direct image, see \cite[§3.1]{G:IndCoh}.

		\begin{defn}
			\label{defn:*_pushforward_ind_coh}
			There exists an $\infty$-functor
			\[
				\Ind(\Coh( - ))\colon \dSch_{ \C }^{ \Noeth } \longrightarrow \Pr^{ L }
			\]
			sending \[
				(f\colon X \to Y)  \mapsto f_{ * }^{ \IC }\colon \Ind(\Coh(X)) \longrightarrow \Ind(\Coh(Y)) 
			\]
			where $f_{ * }^{ \mathrm{IC} }$ is the unique colimit-preserving functor making the diagram
			\begin{diag}
				\Ind(\Coh(X)) \ar[r, "\psi_X"] \ar[d, "f_*^{\mathrm{IC}}"] & \Qcoh(X) \ar[d, "f_*"] \\
				\Ind(\Coh(Y)) \ar[r, "\psi_Y"] & \Qcoh(Y)
			\end{diag}
			commute.
		\end{defn}

		One of the main features of ind-coherent sheaves is that they are arguably the right framework for Grothendieck duality. That is, the $!$-pullback functor can be defined in full generality on ind-coherent sheaves, whereas some additional finiteness assumptions must be taken in order to be able to define it at the quasi-coherent level.
		The construction of $!$-pullbacks and their ($\infty$-categorical) functoriality is pretty hard to obtain: one book and a half (DAG 1, for schemes, and the first part of DAG 2 for ind-inf-schemes) has been devoted to this technical challenge. 
		Let us quickly review the main practical steps.

		\begin{defn}
			\label{defn:!_pullback_ind_coh_proper}
			Let $f\colon X \to Y$ be a proper map of Noetherian derived schemes (recall that $f$ is proper if $t_{ 0 }(f)$ is a proper map of classical schemes). Define
			\[
				f^{ ! }\colon \Ind(\Coh(Y)) \longrightarrow \Ind(\Coh(X))
			\]
			as the \emph{right} adjoint of $f_{ * }^{ \mathrm{IC} }$.
		\end{defn}
 		
		Observe that $f^{ ! }$ preserves colimits since it is exact and preserves filtered colimits; this holds because  $f_{ * }^{ \mathrm{IC} }$ preserves compact objects (i.e.\ coherent sheaves) by properness, see \cite[Corollary 3.3.6]{G:IndCoh}.
		Proper base change (interchanging $!$-pullbacks with ind-coherent $*$-pushforwards) is proved in \cite[Proposition 3.4.2]{G:IndCoh}.

		\begin{defn}
			\label{defn:!_pullback_ind_coh_open}
			Let $j\colon U \to X$ be an open immersion of Noetherian derived schemes. Define
			\[
				j^{ ! }\colon \Ind(\Coh(X)) \longrightarrow \Ind(\Coh(U))
			\]
			as the \emph{left} adjoint of $j_{ * }^{ \mathrm{IC} }$.
		\end{defn}

		\begin{remark}
			\label{remark:*_indcoh_pullback}
			Observe that \cref{defn:!_pullback_ind_coh_open} is well posed: in general, if $f\colon X \to Y$ (assumed to be almost of finite type) is of finite Tor amplitude, i.e.\ $f^{ * }$ is left t-exact up to a finite shift, then there exists a colimit (and compact objects-) preserving functor $f^{\mathrm{IC}, * }$ making the following diagram
			\begin{diag}
				\Ind(\Coh(Y)) \ar[d, "f^{\mathrm{IC}{,}*}"] \ar[r, "\psi_Y"] & \Qcoh(Y) \ar[d, "f^*"] \\
				\Ind(\Coh(X)) \ar[r, "\psi_X"] & \Qcoh(X)
			\end{diag}
			commute.
			It is left adjoint to $f_{ * }^{ \mathrm{IC} }$, as proved in \cite[Lemma 3.5.8]{G:IndCoh}. 
		\end{remark}

		The construction of the $!$-pullback for separated\footnote{For sake of readability, we will describe the process only for \emph{separated} morphisms. In general, one needs to right Kan extend to reach all almost finite type morphisms.} almost of finite type morphisms proceeds as follows:
		\begin{enumerate}
			\item let $f\colon X \to Y$ be such a morphism;
			\item it can be factored, in an essentially unique way, as
				  \begin{diag}
					& \overline{X} \ar[dr, "\overline{f}"] & \\
					  X \ar[ur, hook, "j"] \ar[rr, "f"] & & Y
				  \end{diag}
				  where $j$ is an open immersion and $\overline{f}$ is proper, by a derived version of Nagata compactification theorem (see \cite[§6.3]{G:IndCoh} or \cite[Theorem 8.15]{Scholze:SixFunctor});
		    \item then one defines
				\[
					f^{ ! } \coloneqq j^{ ! } \circ \overline{f}^{ ! },
				\]
				since we already know how to define $!$-pullbacks on open immersions and proper maps.
		\end{enumerate}

		Making this process $\infty$-functorial, encoding base change isomorphisms with respect to $*$-pushforwards (which is a feature that we want in Grothendieck duality) is not easy. 
		The difficulty is that we want to encode additional \emph{structure}, not just properties. This is because, since in general the $!$-pullback won't be adjoint to the ind-coherent $*$-pushforward, there are no Beck-Chevalley transformations a priori, so they need to be given (and in a homotopy-coherent way). 
		The idea to encode their datum, and the fact that they must be equivalences, is to build $\IndCoh$ as an $\infty$-functor from a certain $\infty$-category of \emph{correspondences} in (almost finite type) derived schemes, denoted by $\infcatname{Corr}(\dSch^{ \afp })$.
		The rough idea is that morphisms from $X$ to $Y$ are encoded by ``hats''
		\begin{diag}
			& \ar[dl] Z \ar[rd] & \\
			X & & Y
		\end{diag}
		and composition is done by pullbacks
		\begin{diag}
			& & \ar[dl] Z_3 \coloneqq Z_1 \times_{X_1} Z_2 \ar[dd, phantom, very near start, "\lrcorner"] \ar[dr] & & \\
			& \ar[dl] Z_1 \ar[dr] & & \ar[dl] Z_2 \ar[dr] & \\
			X_0 &  & X_1 &  & X_2
		\end{diag}
		We endow it with a symmetric monoidal structure given by the Cartesian product of derived schemes (beware that this is not the Cartesian product in the correspondence category).
		Notice that we have two canonically defined functors
		\begin{gather*}
			\dSch^{ \afp }_{ \C } \longrightarrow \Corr(\dSch_{ \C }^{ \afp }), \qquad (f\colon X \to Y) \mapsto (X \stackrel{\id_{ X }}{\longleftarrow} X \stackrel{f}{\longrightarrow} Y),\\
			(\dSch^{ \afp }_{ \C })^{ \op } \longrightarrow \Corr(\dSch_{ \C }^{ \afp }), \qquad (f\colon X \to Y) \mapsto  (Y \stackrel{f}{\longleftarrow} X \stackrel{\id_{ X }}{\longrightarrow} X). 
		\end{gather*}
		We observe that both inclusion are symmetric monoidal (where we consider the Cartesian monoidal structure on $\dSch_{ \C }^{ \afp }$ and the coCartesian one on its opposite category).
		For a precise account and construction, see \cite[Chapter 7]{GR:DAG1} and \cite[§2]{HeyerMann:6FF_smooth_representations}. 
		As explained in the first reference mentioned above, the correspondence category admits an $(\infty, 2)$-categorical enhancement (compatible with the tensor product), denoted by $\infcatname{Corr}(\dSch^{ \mathrm{aft} })^{ \mathrm{proper} }$. Roughly, we add non-invertible $2$-morphisms given by commutative diagrams
		\begin{diag}
			& \ar[dl] Z_0 \ar[dd, "h"] \ar[dr] & \\
			X & & Y \\
			  & \ar[ul] Z_1 \ar[ur] &
		\end{diag}
		where $h$ is a proper map.
		The main result is the following (recall that the $\Pr^{ L, st }_{ \C } \coloneqq \Mod_{ \Mod_{ \C } }(\Pr^{ L, st })$ is naturally seen as an $(\infty, 2)$-category adding non invertible natural transformations).
		
		\begin{thm}
			\label{thm:ind_coh_correspondence_schemes}
			There exists an $(\infty, 2)$-functor
			\[
				\IndCoh\colon \infcatname{Corr}(\dSch_{ \C }^{ \mathrm{aft} })^{\mathrm{proper}  } \longrightarrow (\Pr^{ L, st }_{ \C })^{ \mathrm{2-cat} }.
			\]
			Moreover, it has a canonical symmetric monoidal structure.
		\end{thm}
		\begin{proof}
			The construction of the $(\infty, 2)$-functor is done in \cite[Theorem 2.1.4]{GR:DAG1}.
			The symmetric monoidal structure comes from \cite[Theorem 4.1.2]{GR:DAG1}.
			Morally, the action on $1$-morphism is given by
			\[
				(X \stackrel{f}{\leftarrow} Z \stackrel{g}{\rightarrow} Y) \mapsto \left( \Ind(\Coh(X)) \stackrel{f^{ ! }}{\longrightarrow} \Ind(\Coh(Z)) \stackrel{g_{ * }^{ \IC }}{\longrightarrow} \Ind(\Coh(Y)) \right).
			\]
			The homotopy-coherent compatibilities with compositions are exactly the homotopy-coherent base change isomorphisms.
			The symmetric monoidal structure tells us, in particular, that there are ``box product'' morphisms
			\[
				- \boxtimes -\colon \Ind(\Coh(X)) \times \Ind(\Coh(Y)) \longrightarrow \Ind(\Coh(X \times Y))
			\]
			inducing equivalences when they are seen from $\Ind(\Coh(X)) \otimes \Ind(\Coh(Y))$, where $\otimes$ is the symmetric monoidal structure of $\Pr^{ L, st }_{ \C }$ (notice that here it matters to work on a perfect base field, see \cite[Proposition 6.3.4]{GR:DAG1}). 
		\end{proof}
		Several remarks are in order.
		\begin{remark}
			\label{remark:2_categories_ind_coh}
			As explained in \cite[Chapter 5, 0.1.5]{GR:DAG1}, $(\infty, 2)$-categories are necessary for the following reason: consider a proper map $f\colon X \to Y$ and the composition of correspondences
			\begin{diag}
				& & \ar[dl, "g"]  X \times_Y X \ar[dr, "g"] & & \\
				& \ar[dl, "\id_X"] X \ar[dr, "f"] & & \ar[dl, "f"'] X \ar[dr, "\id_X"] \\
				X & & Y & & X
			\end{diag}
			which gives, by $(\infty, 1)$-functoriality of $\IndCoh$, the following natural equivalence
			\[
				g^{ \IC }_{ * } \circ g^{ ! } \stackrel{\sim}{\longrightarrow} f^{ ! }\circ f^{ \IC }_{ * }.
			\]
			The point is that, since $f$ is proper and therefore $f^{ ! }$ is right adjoint to $f^{ \IC }_{ * }$, there exists another natural transformation (Beck-Chevalley) $g^{ \IC }_{ *}g^{ ! } \to f^{ ! }f^{ \IC }_{ * }$, coming from the adjunction (see \cite[Definition 4.7.4.13]{Lurie:HA}). A priori, there is no reason for the first equivalence to be given by the Beck-Chevalley transformation, although we really would like it to be true.
			The introduction of non-invertible $2$-morphisms (given by proper maps) ensures that this is the case, as explained in \cite[§5.4]{G:IndCoh}.
			It is then natural to wonder what happens, instead, with open morphisms, where the $!$-pullback is left adjoint to the ind-coherent $*$-pushforward. There, the adjunction is already encoded by the $(\infty, 1)$-functor since the derived self-intersection of an open immersion is itself, see \cite[§5.3]{G:IndCoh}.
		\end{remark}

		\begin{remark}
			\label{remark:mannscholze_vs_gr}
			In \cite[Remark 8.25]{Scholze:SixFunctor}, the $(\infty, 1)$-functor $\IndCoh$ is built using alternative techniques that do not rely on $(\infty, 2)$-categorical tools\footnote{For a good reason: the theory of $(\infty, 2)$-categories is not yet fully developed.}, see \cite{HeyerMann:6FF_smooth_representations}.	What is not clear in this $(\infty,1)$-categorical version, as opposed to the $(\infty,2)$-one, is what is explained in \cref{remark:2_categories_ind_coh}; that is, it is not clear (to the author) that the natural equivalence of ``proper base change'' for proper maps are obtained as Beck-Chevalley transformation. In other words, it is not clear how to recover the adjunction $(f^{ \IC }_{ * }, f^{ ! })$ for proper maps $f$ which are \emph{not} truncated (and of which there's plenty).
		\end{remark}

		\begin{remark}
			\label{remark:extension_ind_inf_schemes}
			There exists a more general version of \cref{thm:ind_coh_correspondence_schemes}, where not only schemes are considered by \emph{ind-inf-schemes}, see \cite[Chapter 3]{GR:DAG2}. Ind-inf-schemes are prestacks that are ``laft'', admit deformation theory and whose reduced classical truncation is a quasi-compact ind-scheme.
			The more general statement is the existence of an $(\infty, 2)$-functor
			\[
				\IndCoh\colon \Corr(\PreStk_{ \C }^{ \mathrm{laft} })_{ \mathrm{ind-inf-schematic}, \mathrm{all}  }^{ \mathrm{ind-inf-schematic}, \mathrm{ind-proper} } \longrightarrow (\Pr^{ L, st }_{ \C })^{ \mathrm{2-cat} },
			\]
			satisfying the expected properties (i.e.\ encoding proper base change isomorphisms etc).
			This means that we can always consider $!$-pullbacks, ind-coherent $*$-pushforward can be defined for morphisms of laft prestacks which are represented by ind-inf-schemes and the adjunction $(f^{\IC }_{ * }, f^{ ! })$ holds for all $f$ which are ind-inf-schematic and ind-proper. 
			Also symmetric monoidal structures are respected (when restricted to $\Corr(\infcatname{indinfSch}_{ \mathrm{laft} })$), see \cite[Chapter 3, Corollary 6.1.2]{GR:DAG2}.
		\end{remark}

		\begin{example}
			\label{example:indcoh_de_rham}
			We introduced the more general construction in \cref{remark:extension_ind_inf_schemes} because we will talk about ind-coherent sheaves on the de Rham stack $X_{ \dR }$, for $X$ an almost finitely presented derived scheme. Then, $t_{ 0 }(X_{ \dR })^{ \red } \simeq t_{ 0 }(X)^{ \red }$ so it is an inf-scheme and we can use the fully fledged functoriality of $\IndCoh$ (which would not have been available if we had restricted ourselves only to derived schemes), see \cite[Theorem 2.4.5]{Ber:DerivedDMod}.
		\end{example}

		We can now, finally, define the $!$-pullback on ind-coherent sheaves as an $\infty$-functor.

		\begin{defn}
			\label{defn:!_pullback_ind_coh}
			Consider the symmetric monoidal functor 
			\[
				\IndCoh^{ ! }\colon (\dSch^{ \mathrm{aft} }_{ \C })^{ \op } \to \infcatname{Corr}(\dSch_{ \C }^{ \mathrm{aft} }) \to \Pr^{ L, st }_{ \C }.
			\]
			It sends $f\colon X \to Y$ to a functor 
			\[
				f^{ ! }\colon \Ind(\Coh(Y)) \to \Ind(\Coh(X))
			\]
			which we define to be the \emph{$!$-pullback along $f$}.
		\end{defn}

		Analogously, observe that the symmetric monoidal restriction
		\[
			\IndCoh_{ * }\colon \dSch^{ \afp }_{ \C } \to \Corr(\dSch^{ \afp }_{ \C }) \to \Pr^{ L, st }_{ \C }
		\]
		recovers the functor of \cref{defn:*_pushforward_ind_coh}.

		\begin{coroll}
			\label{coroll:symm_monoidal_structure_ind_coh}
			Given $X$ an (almost finite type) derived $\C$-scheme, the $\infty$-category $\Ind(\Coh(X))$ is endowed with a symmetric monoidal structure, denoted by $- \otimes^{ \mathrm{IC} } -$. The monoidal unit is $\omega_{ X } \simeq f^{ ! }(\C)$ for $f\colon X \to \Spec(\C)$.
		\end{coroll}
		\begin{proof}
			This is a corollary of the fact that 
			\[
				\IndCoh^{ ! }\colon (\dSch_{ \C }^{ \afp })^{ \op } \longrightarrow \Pr^{ L, st}_{ \C }
			\]
			is symmetric monoidal and thus, in particular, it preserves commutative algebras.
			Since $(\dSch_{ \C }^{ \afp })^{ \op }$ is endowed with the coCartesian symmetric monoidal structure, the functor above can be upgraded to
			\[
				\IndCoh^{ ! }\colon (\dSch_{ \C }^{ \afp })^{ \op } \simeq \CAlg(\dSch_{ \C }^{ \afp, \op }) \longrightarrow \CAlg(\Pr^{ L, st }_{ \C }).
			\]
			The first equivalence is thanks to \cite[Corollary 2.4.3.10]{Lurie:HA}.
			This means that the tensor product on $\Ind(\Coh(X))$ is computed as
			\[
				\Ind(\Coh(X)) \times \Ind(\Coh(X)) \stackrel{- \boxtimes -}{\longrightarrow} \Ind(\Coh(X \times X)) \stackrel{\Delta^{ ! }}{\longrightarrow} \Ind(\Coh(X)), \qquad (\mathcal{F}, \mathcal{G}) \mapsto \Delta^{ ! }(\mathcal{F} \boxtimes \mathcal{G}).
			\]
			Observe finally that we have
			\[
				\omega_{ X } \simeq f^{ ! }(\C),
			\]
			for $f\colon X \to \Spec(\C)$ simply because of the symmetric monoidal equivalence $\Ind(\Coh_{ \C }) \simeq \Mod_{ \C }$ coming from \cref{prop:upsilon_symm_monoidal}.
		\end{proof}

		A final piece of the puzzle is the (symmetric monoidal) interaction of $\Qcoh$ and $\IndCoh$.

		\begin{prop}
			\label{prop:upsilon_symm_monoidal}
			There exists a symmetric monoidal natural transformation
			\[
				\Upsilon\colon \Qcoh^{ * } \longrightarrow \IndCoh^{ ! }
			\]
			in $\Fun_{ \CAlg(\Pr^{ L, st }_{ \C }) }(\Qcoh^{ * }, \IndCoh^{ ! })$. Given $X \in \dSch_{ \C }^{ \afp }$, it sends $\mathcal{E} \in \Qcoh(X)$ to $\mathcal{E} \odot \omega_{ X }$, where
			\[
				- \odot -\colon \Qcoh(X) \times \Ind(\Coh(X)) \to \Ind(\Coh(X))
			\]
			is the action of $\Qcoh(X)$ on $\Ind(\Coh(X))$ in $\Pr^{ L, st }_{ \C }$ coming from the standard action of $\Perf(X)$ on $\Coh(X)$, see \cite[§1.4]{G:IndCoh}.
			If $A$ is a classical regular $\C$-algebra, then $\Upsilon_{ A }$ is a symmetric monoidal equivalence.
		\end{prop}
		\begin{proof}
			The existence of $\Upsilon$ comes from \cite[Chapter 6, 3.2]{GR:DAG1}.
			For the last claim, observe that
			\[
				\psi_{ A }\colon \Ind(\Coh_{ A }) \longrightarrow \Mod_{ A }
			\]
			is an equivalence by \cite[Lemma 1.1.6]{G:IndCoh}. Since both $\infty$-categories are self-dual (using, respectively, Serre duality and the ``naive'' duality) we obtain that also
			\[
				\psi_{ A }^{ \vee }\colon \Mod_{ A } \simeq \Mod_{ A }^{ \vee } \longrightarrow (\Ind(\Coh_{ A }))^{ \vee } \simeq \Ind(\Coh_{ A })
			\]
			is an equivalence. By \cite[Chapter 6, Theorem 4.4.3]{GR:DAG1}, this is identified with $\Upsilon_{ A }$, and we obtain therefore the sought for symmetric monoidal equivalence.
		\end{proof}

%% file: ProCoherentSheaves.tex
\section{Pro-coherent sheaves}
 	\label{subsection:pro-coherent_sheaves}

	In this section, we will introduce pro-coherent sheaves, first on derived affine schemes and then on prestacks, and recall some of their main properties. The main references are \cite[Appendix A.2]{BrantnerMagidsonNuiten:FormalIntegrationDerivedFoliations} and \cite[§2]{BrantnerCamposNuiten:PDOperadsExplicitPartitionLieAlgebras}. Beware of the different grading conventions (cohomological for us and homological for them).

	\subsection{Definition and main properties}
		\label{subsubsection:pro_coherent_modules_definitions_main_prop}

		Let us start with the affine case.

		\begin{defn}
			\label{defn:pro_coherent_modules}
			Let $A \in \CAlg^{ \cn }_{ \C }$. The $\infty$-category of \emph{pro-coherent modules} over $A$ is the full subcategory of $\Fun_{ \mathrm{ex} }(\Mod_{ A }, \Sp)$ spanned by those $F\colon \Mod_{ A } \to \Sp$ which are
			\begin{enumerate}[label=(\arabic*)]
				\item \emph{convergent}, meaning that for each $M \in \Mod_{ A }$ we have natural equivalences
					\[
						F(M) \stackrel{\sim}{\longrightarrow} \lim_{ n } F(\tau^{ \geq n }M);
					\]
				\item \emph{almost finitely presented}, meaning that each restriction
					\[
						F\colon \Mod_{ A }^{ \geq n } \to \Sp
					\]
					commutes with filtered colimits.
			\end{enumerate}
			It is denoted by $\ProCoh_{ A } \coloneqq \Fun_{ \ex, \conv, \afp }(\Mod_{ A }, \Sp)$.
		\end{defn}

		\begin{remark}
			We chose to follow the notation of \cite[§7.1]{Nuiten:KoszulLieAlgbd}; notice that $\ProCoh_{ A }$ is denoted by $\infcatname{QC}^{ \vee }_{ A }$ in \cite{BrantnerMagidsonNuiten:FormalIntegrationDerivedFoliations}.
			We believe no confusion will arise: we will never consider the $!$-pullback functor on quasi-coherent sheaves.
		\end{remark}

		This $\infty$-category can be presented in a plethora of equivalent ways, see \cite[Lemma A.26]{BrantnerMagidsonNuiten:FormalIntegrationDerivedFoliations}. For general $A$, it can be badly behaved; we will focus only on \emph{coherent} $A$. The name, in the latter case, is explained by the following lemma.

		\begin{lemma}
			\label{lemma:pro_coherent_modules_presentations}
			Let $A \in \CAlg^{ \cn }$ be a coherent commutative algebra. Then the following $\infty$-categories are equivalent:
			\begin{enumerate}[label=(\arabic*)]
				\item $\ProCoh_{ A }$;
				\item $\Fun_{ \ex, \afp }(\Mod_{ A }^{ + }, \Sp)$;
				\item $\Fun_{ \ex }(\Coh_{ A }, \Sp)$;
				\item $\Fun_{ \lex }(\Coh_{ A }, \S)$;
				\item $\Fun_{ \ex, \conv }(\APerf_{ A }, \Sp)$.
			\end{enumerate}
		\end{lemma}
		\begin{proof}
			This is a mix of \cite[Lemma A.42]{BrantnerMagidsonNuiten:FormalIntegrationDerivedFoliations} and \cite[Lemma A.26]{BrantnerMagidsonNuiten:FormalIntegrationDerivedFoliations}. Observe that the compact objects of $\Mod_{ A }^{ + }$ are exactly coherent $A$-modules (i.e.\ bounded almost perfect $A$-modules). See \cite[Corollary 2.20]{BrantnerCamposNuiten:PDOperadsExplicitPartitionLieAlgebras} as well.
			Each time we want to impose convergence property we right Kan extend (for example along $\Coh_{ A } \hookrightarrow \APerf_{ A }$) and each time we want to impose finiteness conditions we left Kan extend.
			The equivalence between $(2)$ and $(3)$ stems from the fact that $\Coh_{ A }$ is stable, see \cite[Corollary 1.4.2.23]{Lurie:HA}. 
		\end{proof}

		This lemma allows us to use, interchangeably, the notations $\ProCoh_{ A }$ and $\Ind(\Coh_{ A }^{ \op })$. Notice that the latter is, technically, the opposite of what one would call pro-coherent sheaves, but we will stick with this terminology to be coherent\footnote{Pun intended.} with the literature.

		\begin{coroll}
			\label{corollary:pro_coherent_stable_presentable}
			For a coherent $A \in \CAlg^{ \cn }$, the $\infty$-category $\ProCoh_{ A }$ is stable and compactly generated.
		\end{coroll}

		\begin{prop}
			\label{prop:modules_into_pro_coherent}
			There exists a left adjoint functor
			\[
				i\colon \Mod_{ A } \longrightarrow \ProCoh_{ A }
			\]
			sending $M$ to $M \otimes_{ A } -\colon \Coh_{ A } \to \Sp$.
		\end{prop}
		\begin{proof}
			The tensor product on $\Mod_{ A }$ gives a functor in $\Pr^{ L, st }$
			\[
				- \otimes_{ A } -\colon \Mod_{ A } \otimes \Mod_{ A } \to \Sp.
			\]
			We can then consider the composition
			\[
				\Mod_{ A } \longrightarrow \Fun^{ L }_{ \ex }(\Mod_{ A }, \Sp) \longrightarrow \Fun_{ \ex }(\Coh_{ A }, \Sp),
			\]
			where the last arrow is precomposition along $\Coh_{ A } \hookrightarrow \Mod_{ A }$. It is formal to see that $i$ is a left adjoint and has the prescribed form. Observe that we can equivalently see $i(M)$ as \[
				\APerf_{ A } \ni N \mapsto \lim_{ n } M \otimes_{ A } \tau^{ \geq n }(N).
			\]
			One easy end computation gives us the formula for the right adjoint
			\[
				\nu\colon \ProCoh_{ A } \simeq \Fun_{ \ex, \conv }(\APerf_{ A }, \Sp) \longrightarrow \Mod_{ A }, \qquad F \mapsto F(A) \simeq \lim_{ n } F(\tau^{ \geq n }(A) ). \qedhere
			\]
		\end{proof}

		Let us now record some properties of $i$: see \cite[Proposition 2.26, Example 2.28]{BrantnerCamposNuiten:PDOperadsExplicitPartitionLieAlgebras}.

		\begin{coroll}
			\label{coroll:A_bounded_i_left_adjoint}
			Then the functor
			\[
				i\colon \Mod_{ A } \longrightarrow \ProCoh_{ A }
			\]
			is fully faithful when restricted to $\Mod_{ A }^{ - }$. It is fully faithful on the whole $\Mod_{ A }$ if and only if $A$ is bounded. 
			It is an equivalence if and only if $A$ is a classical smooth commutative algebra.
		\end{coroll}
		\begin{proof}
			The counit morphism of the adjunction in \cref{prop:modules_into_pro_coherent} is given by
			\[
				M \longrightarrow \nu(i(M)) \simeq \lim_{ n } M \otimes_{ A } \tau^{ \geq n }(A).
			\]
			Let $M \in \Mod_{ A }^{ \leq m }$ and consider the cofiber sequence of $A$-modules
			\[
				M \simeq A \otimes_{ A } M \longrightarrow \tau^{ \geq n }(A) \otimes_{ A } M \longrightarrow \cofib(A \to \tau^{ \geq n }(A)) \otimes_{ A } M.
			\]
			The rightmost term lies in $\Mod_{ A }^{ \leq (n+m-1) }$ so passing to the limit we obtain $0$, since the t-structure on $\Mod_{ A }$ is left complete. This allows us to conclude that 
			\[
				i\colon \Mod_{ A }^{ - } \hookrightarrow \ProCoh_{ A }.
			\]
			Assume now that $A$ is bounded, so that $A \simeq \tau^{ \geq n_{ 0 } }(A)$ for some $n_{ 0 } \in \Z$; then the counit morphism is an equivalence for any $M$ since the limit on the right-hand side is constant.
			To prove the converse implication, consider $M = \oplus_{ k \geq 0 } A[-k]$ and observe that the counit is identified with the map
			\[
				\bigoplus_{ k \geq 0  } A[-k] \longrightarrow \prod_{ k \geq 0 } A[-k]
			\]
			which is not an equivalence unless $A$ is bounded.

			For the final claim, we refer to \cite[Proposition 1.6.4]{G:IndCoh} (recalling that ind-coherent and pro-coherent modules are equivalent, in our context).
		\end{proof}
		
		Pro-coherent $A$-modules carry a left complete t-structure such that $i$ identifies $\Mod_{ A }$ (with its standard left and right complete t-structure) as a right completion of $\ProCoh(A)$, see \cite[Observation A.30]{BrantnerMagidsonNuiten:FormalIntegrationDerivedFoliations}.
		There is another functor from (the opposite of) $A$-modules to pro-coherent ones, given by the restricted Yoneda embedding 
		\[
			j\colon \Mod_{ A }^{ \op } \longrightarrow \ProCoh(A), \qquad M \mapsto \IntHom_{ A }(M, -) \in \Fun_{ \ex, \conv }(\APerf_{ A }, \Sp).
		\]
		The convergence condition is automatic but the almost finite presentation is not: this is why the mapping spectrum formula is valid only on almost perfect modules.

		\begin{defn}[{\cite[Definition A.33]{BrantnerMagidsonNuiten:FormalIntegrationDerivedFoliations}}]
			\label{defn:dually_almost_perfect}
			Let $A \in \CAlg^{ \cn }$ be coherent. The functor $j$ introduced above restricts to a fully faithful embedding
			\[
				\APerf_{ A }^{ \op } \hookrightarrow \ProCoh(A).
			\]
			The full subcategory of \emph{dually almost perfect $A$-modules} is defined as the essential image of the above embedding. 
		\end{defn}

		See \cref{prop:dually_almost_perfect} for an explanation of the terminology, which will be more clear once we endow pro-coherent modules with the appropriate symmetric monoidal structure.

	\subsection{Symmetric monoidal structure}
		\label{subsubsection:tensor_product_pro_coherent_modules}

		We will now introduce the symmetric monoidal structure on pro-coherent $A$-modules (where $A$ will tacitly be assumed to be coherent, unless specified otherwise). 
		First, a technical lemma.

		\begin{lemma}
			\label{lemma:left_Day_convolution_stable}
			Let $\infcatname{C}$ be a small stable $\infty$-category, endowed with a symmetric monoidal structure which commutes with finite colimits in each variable. We then have equivalences of symmetric monoidal $\infty$-categories
			\[
				\Ind(\infcatname{C}) \simeq \Fun_{ \lex }(\infcatname{C}^{ \op }, \S) 
			\]
			where $\Ind(\infcatname{C})$ is endowed with the ``Ind-extension'' of the tensor product on $\infcatname{C}$ and the right-hand side with (the restriction of) left Day convolution\footnote{In \cite[§2.2.6]{Lurie:HA}, the term ``Day convolution'' is used instead. We add ``left'' to recall that we use the left Kan extension.}.
		\end{lemma}
		\begin{proof}
			The statement is well known for the underlying $\infty$-categories (that is, forgetting the tensor products), thanks to the stability of $\infcatname{C}^{ \op }$ (see \cite[Corollary 1.4.2.23]{Lurie:HA}). The symmetric monoidal structure on $\Ind(\infcatname{C})$ comes from \cite[Corollary 4.8.1.14]{Lurie:HA}. The left Day convolution symmetric monoidal structure on the (whole) functor categories $\P(\infcatname{C}) \simeq \Fun(\infcatname{C}^{ \op }, \S)$ comes from \cite[Example 2.2.6.17]{Lurie:HA}.
			What we need to prove first is that $\Fun_{ \lex }(\infcatname{C}^{ \op }, \S)$ is a symmetric monoidal subcategory of $\P(\infcatname{C})$.
			Let us now recall a fundamental result: the left Day convolution on $\P(\infcatname{C})$ makes the Yoneda embedding $\infcatname{C} \hookrightarrow \P(\infcatname{C})$ symmetric monoidal and commutes with colimits in each variable (see \cite[Proposition 2.14]{Glasman:DayConvolution}). This two features are sufficient to uniquely determine the symmetric monoidal structure on $\P(\infcatname{C})$, see \cite[Corollary 4.8.1.12]{Lurie:HA}. 
			This means that, to check that left exact functors are closed under the Day convolution product, it suffices to check that 
			\[
				\Map_{ \infcatname{C} }(-, X) \circledast \Map_{ \infcatname{C} }(-, Y) 
			\]
			is left exact (since the left Day convolution commutes with colimits in each variable). This is true because, via the symmetric monoidal Yoneda embedding, it corresponds to $\Map_{ \infcatname{C} }(-, X \otimes Y)$.
		\end{proof}

		\begin{coroll}
			\label{coroll:Day_convolution_exact_functors}
			Let $\infcatname{C}$ be as above. The $\infty$-category $\Fun_{ \ex }(\infcatname{C}^{ \op }, \Sp)$ inherits a symmetric monoidal structure which
			\begin{enumerate}[label=(\arabic*)]
				\item commutes with colimits in each variable;
				\item makes the stable Yoneda embedding \[
						X \mapsto \IntHom_{ \infcatname{C} }(-, X)
				\]
				symmetric monoidal.
			\end{enumerate}
		\end{coroll}
		\begin{proof}
			This is an immediate consequence of the fact that the stability of $\infcatname{C}$ is equivalent to an $\Sp$-enrichment such that
			\[
				\Omega^{ \infty }(\IntHom_{ \infcatname{C} }(-, X)) \simeq \Map_{ \infcatname{C} }(-, X),
			\]
			see \cite[Notation 1.1.2.17]{Lurie:HA}.
			By stability of $\infcatname{C}$, \cite[Corollary 1.4.2.23]{Lurie:HA} gives a commuting diagram
			\begin{diag}
				\infcatname{C} \ar[r, hook] \ar[dr, hook] & \Fun_{\lex}(\infcatname{C}^{\op}, \S) \\
														  & \Fun_{\ex}(\infcatname{C}^{\op}, \Sp) \ar[u, "\Omega^{\infty} \circ -"', "\sim"]
			\end{diag}
			where the two arrows are stable and unstable Yoneda.
			We transfer the restricted Day convolution on $\Fun_{ \lex }(\infcatname{C}^{ \op }, \S)$ to $\Fun_{ \ex }(\infcatname{C}^{ \op }, \Sp)$ and conclude by \cref{lemma:left_Day_convolution_stable}.
		\end{proof}

		\begin{example}
			\label{example:tensor_product_modules_is_day_convolution}
			The symmetric monoidal structure on $\Mod_{ A } \simeq \Ind(\Perf_{ A })$ is obtained by ``Ind-extending'' the tensor product on $\Ind(\Perf_{ A })$. By \cref{lemma:left_Day_convolution_stable}, we deduce that it coincides with the restriction of Day convolution on $\Fun_{ \lex }(\Perf_{ A }^{ \op }, \Sp)$.
		\end{example}

		\begin{remark}
			We don't know if the symmetric monoidal structure on $\Fun_{ \ex }(\infcatname{C}^{ \op }, \Sp)$ is the restriction of left Day convolution on $\Fun(\infcatname{C}^{ \op }, \Sp)$. In fact, one key point in the proof of the symmetric monoidal Yoneda embedding in \cite[§3]{Glasman:DayConvolution} is to use that each space $Y$ is the colimit of the constant $*$-valued functor $Y \to \S$, to be able to compute the left Kan extension required in Day convolution.
			We cannot see how to generalize his proof to spectra, due to this specific point in particular. Nevertheless, \cref{coroll:Day_convolution_exact_functors} suffices for our purposes.
		\end{remark}
	
		\begin{prop}
			\label{prop:convergent_functor_symm_monoidal_localization_day_convolution}
			Consider an $\infty$-category $\infcatname{C}$ such that
			\begin{enumerate}[label=(\roman*)]
				\item $\infcatname{C} \in \CAlg(\Pr^{ L, st })$;
				\item $\infcatname{C}$ is endowed with a left complete accessible t-structure which is compatible with the tensor products (i.e.\ $\infcatname{C}^{ \leq 0 }$ is a symmetric monoidal full subcategory);
				\item $\infcatname{C} \simeq \infcatname{C}^{ - }$, i.e.\ each object is right-bounded;
				\item each object of $\infcatname{C}$ can be written as a filtered colimit of objects in $\Thick(1_{ \infcatname{C} })$.
			\end{enumerate}
			Let $\infcatname{D} \in \CAlg(\Pr^{ L })$. The inclusion
			\[
				\Fun_{ \conv }(\infcatname{C}, \infcatname{D}) \hookrightarrow \Fun(\infcatname{C}, \infcatname{D})
			\]
			is a \emph{symmetric monoidal localisation}, where the right-hand side is endowed with left Day convolution product.
			If $\infcatname{D} = \S$, then the induced symmetric monoidal structure on $\Fun_{ \conv }(\infcatname{C}, \S)$ commutes with colimits in each variable and makes the Yoneda embedding
			\[
				\infcatname{C}^{ \op } \hookrightarrow \Fun_{ \conv }(\infcatname{C}, \S)
			\]
			symmetric monoidal. It is, therefore, the unique symmetric monoidal structure having such properties.
		\end{prop}
		\begin{proof}
			The definition of convergent functor is \cite[Definition 2.17]{BrantnerCamposNuiten:PDOperadsExplicitPartitionLieAlgebras}. 
			The equivalence 
			\[
				\Fun_{ \conv }(\infcatname{C}, \infcatname{D}) \simeq \Fun(\infcatname{C}^{ + }, \infcatname{D})
			\]
			given by \cite[Lemma 2.19]{BrantnerCamposNuiten:PDOperadsExplicitPartitionLieAlgebras} suffices to deduce the existence of the left adjoint (the ``localization'' functor). It is given by restriction along $\infcatname{C}^{ + } \hookrightarrow \infcatname{C}$ followed by right Kan extension.

			We only need to prove that this localisation is compatible with the symmetric monoidal structure, as in \cite[Example 2.2.1.7]{Lurie:HA}. That is, given a natural transformation $\eta\colon F \to G$ which is an equivalence when restricted to $\infcatname{C}^{ + }$, we need to prove that
			\[
				\eta \circledast H\colon F \circledast H \longrightarrow G \circledast H
			\]
			satisfies the same property for any $H\colon \infcatname{C} \to \infcatname{D}$ (where $\circledast$ is the left Day convolution product). For $M \in \infcatname{C}^{ + }$, the left Kan extension formula gives us
			\[
				(F \circledast H)(M) \simeq \colim_{ (X, Y, X \otimes Y \to M) } F(X) \otimes H(Y)
			\]
			where the index category is 
			\[
				\mathcal{I}_{ M } \coloneqq (\infcatname{C} \times \infcatname{C}) \times_{ \infcatname{C} } \infcatname{C}_{ -/M }.
			\]
			It suffices to prove that 
			\[
				\mathcal{J}_{ M } \coloneqq (\infcatname{C}^{ + } \times \infcatname{C}) \times_{ \infcatname{C} } \infcatname{C}_{ -/M } \hookrightarrow \mathcal{I}_{ M }
			\] 
			is cofinal (notice that the only difference is that $X$ needs to be eventually coconnective in $\mathcal{J}_{ M }$). By (the $\infty$-categorical version of) Quillen's theorem A, it suffices to prove that for each $(X, Y, X \otimes Y \to M) \in \mathcal{I}_{ M }$, the slice category
			\[
				\mathcal{E} \coloneqq \mathcal{J}_{ M } \times_{ \mathcal{I}_{ M } } (\mathcal{I}_{ M })_{ (X, Y, X \otimes Y \to M) / - }
			\]
			is weakly contractible, see \cite[Theorem 4.1.3.1]{Lurie:HTT}. This means that the $\infty$-groupoid obtained by inverting all morphisms, denoted by $|\mathcal{E}|$, is a contractible space. By Whitehead's theorem, it suffices to prove that all its homotopy groups are zero. One way of proving it is to show that all possible diagrams like
			\begin{diag}
				K \ar[r] \ar[d, hook] & \mathcal{E} \\
				K^{\vartriangleright} \ar[ur, dashed] & 
			\end{diag}
			admit a lift, for $K$ any finite simplicial set. In fact, lifting for $K = \partial\Delta^{ n+1 }$ proves that $\pi_{ n }(|\mathcal{E}|) = 0$ for $n \geq 0$, and lifting for $K = \emptyset$ proves that $|\mathcal{E}|$ is not the empty space.
			Let us assume without loss of generality that $M \in \infcatname{C}^{ \geq n }$. Since $\tau^{ \geq n}$ is the left adjoint of the (notationally silent) inclusion $\infcatname{C}^{ \geq n } \hookrightarrow \infcatname{C}$, we can then factor
			\[
				X \otimes Y \longrightarrow \tau^{\geq n}( X \otimes Y ) \to M.
			\]
			We claim that $X \otimes Y \to \tau^{ \geq n }(X \otimes Y)$ can be factored as
			\[
				X \otimes Y \longrightarrow \tau^{ \geq n-p }(X) \otimes Y \to \tau^{ \geq n }(X \otimes Y)	
			\]
			for some $p$ such that $Y \in \infcatname{C}^{ \leq p }$. The key idea is that, thanks to the accessibility of the t-structure and the assumption (iv) on $\infcatname{C}$, each object  in $\infcatname{C}^{ \leq p }$ can be written as a filtered colimit of $1_{ \infcatname{C} }[-p]$, up to retracts (which won't be relevant since they are preserved by each functor). Observe now that  
			\[
				\tau^{ \geq n-p }(X) \otimes 1_{ \infcatname{C} }[-p]  \simeq \tau^{ \geq n }(X \otimes 1_{ \infcatname{C} }[-p]).
			\]
			The existence of the factorization claimed above is then due to the simple fact that tensor product commutes with colimits in each variable and to the universal property of colimits (that is, the left-hand side is a colimit in $\infcatname{C}$ and each term has a natural arrow towards the right-hand side). This proves the claim, which allows us to deduce that $|\mathcal{E}|$ is not the empty space.

			Consider now $K$ a non-empty finite simplicial set and denote by $(C_{ k }, D_{ k }, C_{ k } \otimes D_{ k } \to M)$ a given $K$-shaped diagram in $\mathcal{E}$. We want to add a final object $(C, D, C \otimes D \to M)$ (and compatible homotopies) in the diagram.
			Let us consider $D \coloneqq \colim_{ k } D_{ k }$ computed in $\infcatname{C}_{ Y/- }$, and observe that, without loss of generality, the initial $K$-diagram maps to the new $K$-diagram in $\mathcal{E}$ given by $(C_{ k }, D, C_{ k } \otimes D \to M)$. Moreover, this new $K$-diagram lives under $(X, D, X \otimes D \to M)$. This proves that we can assume the second object to be constant in each entry of the diagram. That is, renaming $Y$ to $D$, we reduced ourselves to consider a new diagram
			\[
				K \longrightarrow \mathcal{F}_{ X } \coloneqq \left( (\infcatname{C}^{ + } \times \{Y\}) \times_{ \infcatname{C} } \infcatname{C}_{ -/M } \right)_{ (X, Y, X \otimes Y \to M ) / -}.
			\]
			Notice the slight abuse of notation: it is not really a comma category but it is obtained as a pullback of one, since $X$ is not necessarily in $\infcatname{C}^{ + }$.
			Focus now on the first variable, and let $N \coloneqq \min(n-p, \min_{ k }(n_{ k }))$ where each $C_{ k }$ lives in $\infcatname{C}^{ \geq n_{ k } }$ and $Y \in \infcatname{C}^{ \leq p }$. Consider now the truncation map (tensored by $Y$)
			\[
				X \otimes Y \longrightarrow \tau^{ \geq N }(X) \otimes Y.
			\]
			It is clear, by choice of $N$ and by adjunction, that our $K$-shaped diagram factors through the induced map $\mathcal{F}_{ \tau^{ \geq N }(X) } \to \mathcal{F}_{ X }$. This means we can assume that $X$ lives in $\infcatname{C}^{ + }$.
		    Consider
			\[
				(C \coloneqq \colim_{ k } C_{ k }, Y, C \otimes Y \to M),
			\]
			where the first colimit is computed in $\infcatname{C}_{ X/- }$. Since $K$ is finite, all categories in play are stable and $X$ is eventually coconnective, such colimit is actually computed in $\infcatname{C}^{ + }_{ X/- }$. This means that the above object lives in $\mathcal{F}_{  X }$. It is also clear that it lives under $K$; that is, it defines one extension along $K \hookrightarrow K^{ \vartriangleright }$.

			This concludes the proof of the cofinality of $\mathcal{J}_{ M } \hookrightarrow \mathcal{I}_{ M }$ and shows that $\eta \circledast H$ is a $\conv$-local equivalence.
			We can conclude that $\Fun_{ \conv }(\infcatname{C}, \infcatname{D})$ inherits the convergent renormalization of the Day convolution product.

			Let us now assume $\infcatname{D} = \S$. By left completion of the t-structure on $\infcatname{C}$, it is clear that the Yoneda embedding factors as
			\[
				\infcatname{C}^{ \op } \hookrightarrow \Fun_{ \conv }(\infcatname{C}, \S).
			\]
			Notice that, by properties of localisations, this map is also the composition
			\[
				\infcatname{C}^{ \op } \hookrightarrow \Fun(\infcatname{C}, \S) \longrightarrow \Fun_{ \conv }(\infcatname{C}, \S).
			\]
			This proves that it is symmetric monoidal.
			Finally, since the convergent tensor product is obtained by renormalising the Day convolution, by properties of localisations and by the fact that left Day convolution on $\Fun(\infcatname{C}, \S)$ commutes with colimits in each variable, we deduce that also convergent Day convolution does so (beware that $\Fun_{ \conv }(\infcatname{C}, \S) \hookrightarrow \Fun(\infcatname{C}, \S)$ does not commute with general colimits).
		\end{proof}

		\begin{example}
			The main examples we have in mind are $\infcatname{C} = \Mod_{ B }^{ - }$ or $\infcatname{C} = \APerf_{ B }$ for any $B$ (coherent) connective commutative algebra.
			Concretely, the symmetric monoidal structure on convergent functors from \cref{prop:convergent_functor_symm_monoidal_localization_day_convolution} can be described object-wise as 
			\[
				F \circledast^{ \conv } G \simeq (F \circledast G)^{ \conv }
			\]
			where $(-)^{\conv}$ is the left adjoint to the inclusion of convergent functors into all functors. This means that, to compute $F \circledast^{ \conv } G$ on $\infcatname{C}^{ + }$, we can simply compute the value of the standard Day convolution.
		\end{example}

		\begin{remark}
			\label{remark:convergent_day_convolution_yoneda}
			The cofinality argument in \cref{prop:convergent_functor_symm_monoidal_localization_day_convolution} also applies for (left) exact functors, when $\infcatname{D} = \S$. That is,
			\[
				\Fun_{ \lex, \conv }(\infcatname{C}, \S) \hookrightarrow \Fun_{ \lex }(\infcatname{C}, \S)
			\]
			is a symmetric monoidal localisation for the left Day convolution. The same holds for exact functors to $\Sp$, using \cref{coroll:Day_convolution_exact_functors}. The tensor product commutes with colimits in each variable and makes the Yoneda embedding symmetric monoidal. As such, it is uniquely determined (i.e.\ any symmetric monoidal structure with those two properties is canonically equivalent to the convergent Day convolution).
		\end{remark}

		\begin{prop}
			\label{prop:pro_coherent_tensor_product}
			There exists a unique symmetric monoidal structure on $\ProCoh_{ A }$, denoted by $- \otimes^{ ! }_{ A } -$, such that:
			\begin{enumerate}[label=(\roman*)]
				\item it commutes with colimits in each variable;
				\item the Yoneda embedding 
					\[
						j\colon \APerf^{ \op }_{ A } \hookrightarrow \ProCoh_{ A }
					\]
					is symmetric monoidal.
			\end{enumerate}
			The functor $i\colon \Mod_{ A } \to \ProCoh_{ A }$ inherits a canonical symmetric monoidal structure.
		\end{prop}
		\begin{proof}
			Let us observe that $ \APerf_{ A } \hookrightarrow \Mod_{ A }$ is a symmetric monoidal full subcategory by \cref{prop:almost_perfect_symmetric_monoidal} and, by coherence of $A$, inherits an (accessible, left complete and bounded to the right) t-structure.
			Now consider the inclusion
			\[
				\ProCoh_{ A } \simeq \Fun_{ \ex, \conv }(\APerf_{ A }, \Sp) \hookrightarrow \Fun_{ \ex }(\APerf_{ A }, \Sp).
			\]
			Applying \cref{prop:convergent_functor_symm_monoidal_localization_day_convolution} and \cref{remark:convergent_day_convolution_yoneda}, the convergent Day convolution defines a symmetric monoidal structure on $\ProCoh_{ A }$ satisfying the two assumptions.
			It remains to verify that $i\colon \Mod_{ A } \to \ProCoh_{ A }$ is symmetric monoidal. Using \cref{example:tensor_product_modules_is_day_convolution}, it suffices to prove that 
			\[
				i\colon \Perf_{ A } \hookrightarrow \Mod_{ A } \to \ProCoh_{ A }
			\]
			is symmetric monoidal. 
			Let $P \in \Perf_{ A }$ and observe that
			\[
				i(P) \simeq P \otimes_{ A } - \simeq \IntHom_{ A }(P^{ \vee }, -) \simeq j(P^{ \vee })
			\]
			in $\Fun_{ \ex }(\Coh_{ A }, \Sp)$. This means that we can factor the restriction of $i$ to perfect $A$-modules as
			\[
				\Perf_{ A } \stackrel{(-)^{ \vee }}{\longrightarrow} \Perf_{ A }^{ \op } \hookrightarrow \APerf_{ A }^{ \op } \stackrel{j}{\longrightarrow} \Fun_{ \ex, \conv }(\APerf_{ A }, \Sp) \simeq \ProCoh_{ A },
			\]
			where $(-)^{ \vee }$ is $A$-linear duality and $j$ is the Yoneda embedding. All functors above are symmetric monoidal and hence we conclude.
		\end{proof}

		\begin{remark}
			\label{remark:link_with_super_abstract_pro_coherent_tensor_product}
			We decided to explicitly prove that $\ProCoh_{ A }$ is a symmetric monoidal localisation of the Day convolution product. This is claimed in \cite[Example 2.54]{BrantnerCamposNuiten:PDOperadsExplicitPartitionLieAlgebras} and \cite[Corollary A.45]{BrantnerMagidsonNuiten:FormalIntegrationDerivedFoliations}.
		\end{remark}

		\begin{remark}
			The pro-coherent tensor product introduced here agrees with the one used in \cite[§8.6]{Scholze:SixFunctor}, which commutes with colimits in each variable and is explicitly described on coherent modules. The explanation given there is quite intuitive: the $!$-tensor product of two coherent modules is computed first in $A$-modules and then it is considered as (formal) limit of its (coherent) truncations. It is easy to verify that such a description holds also for \cref{prop:pro_coherent_tensor_product}.
		\end{remark}

		We are now able to better explain the name ``dually almost perfect'' from \cref{defn:dually_almost_perfect}.
		\begin{prop}[{\cite[Proposition 2.55]{BrantnerCamposNuiten:PDOperadsExplicitPartitionLieAlgebras}}]
			\label{prop:dually_almost_perfect}
			Consider $\ProCoh_{ A }$ endowed with its (closed) symmetric monoidal structure. 
			The following diagram commutes
			\begin{diag}
				\APerf_A \ar[d, hook, "i"] \ar[dr, hook, "j^{\op}"] & \\
				\ProCoh_A \ar[r, "\IntHom(-{,}i(A))"] & (\ProCoh_{A})^{\op}.
			\end{diag}
			More precisely, there is an induced equivalence
			\[
				\IntHom(-,i(A))\colon i(\APerf_{ A }) \stackrel{\sim}{\longrightarrow} j^{ \op }(\APerf_{ A }).
			\]
			Concretely, this means that all functors $K \otimes_{ A }-$ with $K$ almost perfect are dualizable in $\ProCoh_{ A }$ and their dual is given by $\Map_{ A }(K, -)\colon \Coh_{ A } \to \S$.
		\end{prop}

		Let us give now an explicit example of computation of the pro-coherent tensor product; in particular, the symmetric monoidal functor $i\colon \Mod_{ A } \to \ProCoh_{ A }$ gives a $\Mod_{ A }$-module structure to $\ProCoh_{ A }$.
		\begin{coroll}
			\label{coroll:qcoh_module_structure_on_procoh}
			Let $M \in \APerf_{ A }$ and $F \in \ProCoh_{ A }$. Then
			\[
				F \otimes^{ ! }_{ A } (M \otimes_{ A } -) \simeq F(M \otimes_{ A } -)\colon \Coh_{ A } \to \S,
			\]
			where we recall that $- \otimes^{ ! }_{ A } -$ is the symmetric monoidal structure on $\ProCoh_{ A }$ introduced in \cref{prop:pro_coherent_tensor_product}.
		\end{coroll}
		\begin{proof}
			Let $K \in \Coh_{ A }$ and let us use Yoneda lemma to evaluate the left-hand side on $K$:
			\begin{gather*}
				(F \otimes^!_{ A } i(M))(K) \simeq \Map_{ \ProCoh_{ A } }(\Map_{ A }(K, -), F \otimes^{ ! }_{ A } i(M)) \simeq \\
				\simeq \Map_{ \ProCoh_{ A } }(\Map_{ A }(K, -) \otimes^{ ! }_{ A } (i(M))^{ \vee }, F), 
			\end{gather*}
			where $(-)^{ \vee }$ denotes pro-coherent duality and we used that $i(M)$ is dualizable thanks to \cref{prop:dually_almost_perfect}. Since $i(M)^{ \vee } \simeq j(M)$ and $j\colon \APerf_{ A }^{ \op } \hookrightarrow \ProCoh_{ A }$ is symmetric monoidal, we obtain
			\begin{gather*}
				(F \otimes^{ ! }_{ A } i(M))(K) \simeq \Map_{ \ProCoh_{ A } }(\Map_{ A }(K, -) \otimes^{ ! }_{ A } \Map_{ A }(M, -), F) \simeq \\
				\simeq \Map_{ \ProCoh_{ A } }(\Map_{ A }(K \otimes_{ A } M, -), F)  \simeq F(M \otimes_{ A } -)(K). \qedhere
			\end{gather*}
		\end{proof}

	\subsection{Functoriality}
		\label{subsubsection:pro_coherent_functoriality}
		
		Let us study the functoriality of $\ProCoh_{ A }$ changing $A$. The main functor will be the $!$-pullback. 

		\begin{defn}
			\label{defn:!-pullback_pro_coherent_modules}
			Let $f\colon A \to B$ be a morphism in $\CAlg^{ \cn }_{ \C }$. Define
			\[
				f^{ ! } \coloneqq (- \circ f_{ * })\colon \ProCoh_{ A } \longrightarrow \ProCoh_{ B }
			\]
			as precomposition with $f_{ * }\colon \Mod_{ B }^{ + } \to \Mod_{ A }^{ + }$. Its right adjoint is denoted by $f_{ * }^{ \mathrm{PC} }$.
		\end{defn}

		\begin{remark}
			What we denote by $f^{ ! }$ is denoted by $f^{ * }$ in \cite[Definition A.37]{BrantnerMagidsonNuiten:FormalIntegrationDerivedFoliations}. We chose to denote it as $!$-pullback because it corresponds to the $!$-pullback of ind-coherent sheaves under Serre duality\footnote{We reckon that this can be seen as a reason to instead stick with the $*$-notation.} and because there is no risk of confusion with $*$-pullback of quasi-coherent sheaves. Moreover, we will not need to talk about $!$-pullback of quasi-coherent sheaves, so there should be no confusion.
		\end{remark}

		\begin{prop}
			\label{prop:!_pullback_makes_sense}
			Let $f\colon A \to B$ as above. Then $f^{ ! }$ from \cref{defn:!-pullback_pro_coherent_modules} is well defined, exact, $\C$-linear and admits a right adjoint.
		\end{prop}
		\begin{proof}
			Exactness and $\C$-linearity are automatic from the definition, since $f_{ * }$ is exact and $\C$-linear.
			Let us start observing that the connectivity assumptions on $A$ and $B$ imply that 
			\[
				f^{ * } \simeq B \otimes_{ A } - \colon \Mod_{ A } \to \Mod_{ B }
			\] 
			is right t-exact. Therefore, its right adjoint $f_{ * }$ is left t-exact and induces a functor
			\[
				f_{ * }^{ + }\colon \Mod_{ B }^{ + } \to \Mod_{ A }^{ + }.
			\]
			Remark that $f_{ * }$ is obtained by the above by right Kan extending along $\Mod_{ B }^{ + } \hookrightarrow \Mod_{ B }$, so that precomposing along $f_{ * }$ or $f_{ * }^{ + }$ is the same for our purposes. Since $\ProCoh_{ A } \simeq \Fun_{ \ex, \mathrm{afp} }(\Mod_{ A }^{ + }, \Sp)$ by \cref{lemma:pro_coherent_modules_presentations}, it suffices to prove that precomposition by $f_{ * }^{ + }$ preserves almost finitely presented functors (exactness is immediate).
			Let $F$ be almost finitely presented: it suffices to prove that for any $n \in \N$, the restriction
			\[
				F(f_{ * }^{ + }(-))\colon \Mod_{ B }^{ \geq n } \to \Sp
			\]
			commutes with filtered colimits. Again, this is immediate because $f_{ * }$ commutes with (all) colimits, because the conservative forgetful functor $\Mod_{ B } \to \Sp$ preserves and creates colimits.
			We conclude that precomposition with $f^{ + }_{ * }$ induces a well defined functor
			\[
				f^{ ! }\colon \ProCoh_{ A } \longrightarrow \ProCoh_{ B }.
			\]
			It preserves colimits because $\ProCoh_{ A } \hookrightarrow \Fun_{ \ex }(\Mod_{ A }^{ + }, \Sp)$ does so, and colimits of (exact) functors are computed term-wise. By the adjoint functor theorem, there exists a right adjoint
			\[
				f_{ * }^{ \mathrm{PC} }\colon \ProCoh_{ B } \longrightarrow \ProCoh_{ A }. \qedhere
			\]
		\end{proof}
		
		We will need a technical lemma, about construction of certain functors out of modules and pro-coherent modules. It is a modified and easier version of \cite[Proposition 2.40]{BrantnerCamposNuiten:PDOperadsExplicitPartitionLieAlgebras}, whose proof is immediately adapted from theirs (observing that a functor preserves finite colimits if and only if it preserves finite coproducts and finite geometric realisations). In fact, we only need to identify the subcategories on the right-hand side which correspond to colimit-preserving functors; that is, we just need to add the requirement that finite coproducts are preserved.

		\begin{defn}[{\cite[Definition 2.38]{BrantnerCamposNuiten:PDOperadsExplicitPartitionLieAlgebras}}]
			\label{defn:regular_functors}
			Let $A \in \CAlg_{ \C }^{ \cn, \coh }$ and let $F \in \Fun_{ \ex }(\APerf_{ A }^{ \op }, \infcatname{V})$. We say that $F$ is \emph{regular} if
			\[
				F^{ \op }\colon \APerf_{ A } \to \infcatname{V}^{ \op }
			\]
			is convergent, i.e.\ preserves limits of Postnikov towers. They span a full subcategory denoted by $\Fun_{ \ex, \mathrm{reg} }(\APerf_{ A }^{ \op }, \infcatname{V})$.
		\end{defn}

		\begin{lemma}
			\label{lemma:extending_functor_mod_procoh}
			Let $A \in \CAlg^{ \cn, \coh}_{ \C }$ and $\infcatname{V} \in \Pr^{ L, st }$.
			The following diagram commutes
			\begin{diag}
				\Fun^L_{\ex}(\ProCoh_A, \infcatname{V}) \ar[r, "\sim"] \ar[d, "- \circ i_A"] & \Fun_{\ex, \mathrm{reg}}(\APerf_A^{\op}, \infcatname{V}) \simeq \Fun_{\ex}(\Coh_A^{\op}, \infcatname{V}) \ar[d] \\
				\Fun^L_{\ex}(\Mod_A, \infcatname{V}) \ar[r, "\sim"] & \Fun_{\ex}(\Perf_A, \infcatname{V}).
			\end{diag}
			The horizontal arrows are equivalences given by restriction (and whose inverse is left Kan extension). The vertical rightmost arrow is given by precomposing along the exact functor
			\[
				\Perf_{ A } \stackrel{(-)^{ \vee }}{\simeq} \Perf_{ A }^{ \op } \hookrightarrow \APerf_{ A }^{ \op }.
			\]
		\end{lemma}

		Let us record here a symmetric monoidal version of the above lemma, which will be important later on. We will denote by $\Fun^{ L, \otimes }(\infcatname{C}, \infcatname{D})$ the $\infty$-category of symmetric monoidal colimit-preserving functors from $\infcatname{C}$ and $\infcatname{D}$ (which are elements of $\CAlg(\Pr^{ L })$).

		\begin{lemma}
			\label{lemma:symm_monoidal_functor_mod_procoh}
			Let $A \in \CAlg^{ \cn, \coh}_{ \C }$ and $\infcatname{V} \in \CAlg(\Pr^{ L, st })$.
			The following diagram commutes
			\begin{diag}
				\Fun^{L, \otimes}_{\ex}(\ProCoh_A, \infcatname{V}) \ar[r, "\sim"] \ar[d, "- \circ i_A"] & \Fun_{\ex, \mathrm{reg}}^{\otimes}(\APerf_A^{\op}, \infcatname{V}) \ar[d] \\
				\Fun^{L, \otimes}_{\ex}(\Mod_A, \infcatname{V}) \ar[r, "\sim"] & \Fun_{\ex}^{\otimes}(\Perf_A, \infcatname{V}).
			\end{diag}
		\end{lemma}
		\begin{proof}
			Let us first observe that the diagram is well defined, since all arrows are given by restriction along symmetric monoidal functors.
			The bottom horizontal arrow, given by restriction along the symmetric monoidal inclusion $\Perf_{ A } \hookrightarrow \Mod_{ A }$, is an equivalence because $\Mod_{ A } \simeq \Ind(\Perf_{ A })$ is equipped with the symmetric monoidal structure $\Ind$-extended from $\Perf_{ A }$ (and hence satisfying the corresponding universal property explained in \cite[Corollary 4.8.1.14]{Lurie:HA}).
			We can use that 
			\[
				\adjunction{(-)^{ \mathrm{conv} }}{\Ind(\APerf_{ A }^{ \op })}{\ProCoh_{ A }}{i}
			\]
			is a symmetric monoidal localization (i.e.\ $(-)^{\conv }$ is symmetric monoidal, and therefore $i$ is lax symmetric monoidal) by \cref{prop:convergent_functor_symm_monoidal_localization_day_convolution}. In particular, this implies that
			\[
				\adjunction{( (-)^{ \conv } )^{ \otimes }}{\Ind(\APerf_{ A }^{ \op })^{ \otimes }}{(\ProCoh_{ A })^{ \otimes }}{i^{ \otimes }}
			\]
			is also a Bousfield localization (where we considered the total categories over $\E_{ \infty }^{ \otimes }$ giving the symmetric monoidal structures to the ones before).
			This gives us
			\[
				\Fun^{ L, \otimes }_{ \ex }(\ProCoh_{ A }, \infcatname{V}) \stackrel{- \circ (-)^{ \conv }}{\hookrightarrow} \Fun^{ L, \otimes }_{ \ex }(\Ind(\APerf_{ A }^{ \op }), \infcatname{V}).
			\]
			using \cite[Proposition 5.2.7.12]{Lurie:HTT}. The same proposition allows us to understand the essential image: it is spanned by those functors which send to equivalences all the morphisms in $\Ind(\APerf_{ A }^{ \op })^{ \otimes }$ that are sent to equivalences by $( (-)^{ \conv } )^{ \otimes }$. 
			Like before, we have
			\[
				\Fun^{ L, \otimes }_{ \ex }(\Ind(\APerf_{ A }^{ \op }), \infcatname{V}) \stackrel{\sim}{\longrightarrow} \Fun^{ \otimes }_{ \ex }(\APerf_{ A }^{ \op }, \infcatname{V}).
			\]
			Unravelling the definitions and using \cref{lemma:extending_functor_mod_procoh}, we see that this essential image is exactly the full subcategory of $\Fun_{ \ex }^{ \otimes }(\APerf_{ A }^{ \op }, \infcatname{V})$ spanned by the regular functors (observe that regularity is a property that only depends on the underlying functor). To be more explicit, we observe that for any $K \in \APerf_{ A }$ the following morphism in $\Ind(\APerf_{ A }^{ \op })$
			\[
				\colim_{ n \in \N } \Map_{ A }(\tau^{ \geq -n }(K), -) \longrightarrow \Map_{ A }(K, -)
			\]
			is not an equivalence but becomes one when $(-)^{ \conv }$ is applied. This is how the regularity condition arises.	This allows us to conclude that the top arrow of the diagram is an equivalence.
			Finally, the fact that the diagram in the statement commutes follows from \cref{prop:dually_almost_perfect} which tells us that the diagram
			\begin{diag}
				\Perf_A \ar[r, hook] \ar[d, "(-)^{\vee}"] & \Mod_A \ar[d, "i_A"] \\
				\APerf_A^{\op} \ar[r, "j_A"] & \ProCoh_A
			\end{diag}
			commutes in $\CAlg(\Pr^{ L, st })$.\qedhere
		\end{proof}

		\begin{prop}
			\label{prop:pro_coh_!_pullback_compatibilities}
			Let $f\colon A \to B$ be a morphism of coherent connective commutative $\C$-algebras.
			Then the following diagram commutes
			\begin{diag}
				\Mod_A \ar[d, "i_A"] \ar[r, "f^*"] & \Mod_B \ar[d, "i_B"] \\
				\ProCoh_A \ar[r, "f^!"] & \ProCoh_B \\
				\APerf_A^{\op} \ar[u, hook, "j_A"] \ar[r, "(f^*)^{\op}"] & \APerf_B^{\op} \ar[u, hook, "j_B"]
			\end{diag}
			It can, furthermore, be seen as a diagram in $\CAlg(\Pr^{ L, st })$.
		\end{prop}
		\begin{proof}
			Let us start with the bottom square and let $K \in \APerf_{ A }$; going up left gives
			\[
				\Map_{ A }(K, f_{ * }(-)) \in \Fun_{ \lex, \afp }(\Mod_{ A }^{ + }, \S)
			\]
			By the $(f^{ * }, f_{ * })$ adjunction, it is equivalent to 
			\[
				\Map_{ B }(f^{ * }(K), -) \in \Fun_{ \lex, \afp }(\Mod_{ B }^{ + }, \S),
			\]
			which is exactly what we obtain by going right-up in the diagram.
			Since this equivalence is natural in $K$ (coming from the adjunction), we can conclude that the bottom diagram commutes in $\Pr^{ L, st }$. Observe that all the $\infty$-categories are symmetric monoidal, $f^{ * }$ is symmetric monoidal and both $j_{ A }$ and $j_{ B }$ are symmetric monoidal by \cref{prop:pro_coherent_tensor_product}. We just need to prove that $f^{ ! }$ is symmetric monoidal to be able to say that it is a commutative square in $\CAlg(\Pr^{ L, st })$. This holds by \cref{lemma:symm_monoidal_functor_mod_procoh}.
			Let us now prove the statement for the top square, observing that by \cref{lemma:symm_monoidal_functor_mod_procoh} we can restrict the top left corner to $\Perf_{ A }$. Let $E \in \Perf_{ A }$ and consider
			\[
				f^{ ! }(E \otimes_{ A } -) \simeq E \otimes_{ A } f_{ * }(-) \in \Fun_{ \lex, \afp }(\Mod_{ B }^{ + }, \S).
			\]
			By the projection formula (see \cite[Chapter 3, Remark 3.4.2.6]{Lurie:SAG} or the classical \cite[\href{https://stacks.math.columbia.edu/tag/0944}{Tag 0944}]{stacks-project}) we obtain
			\[
				f^{ ! }(E \otimes_{ A } -) \simeq f^{ * }(E) \otimes_{ B } - \in \Fun_{ \lex, \afp }(\Mod_{ B }^{ + }, \S).
			\]
			The naturality in $E$ of the projection formula allows us to conclude.
		\end{proof}

		\begin{prop}
			\label{prop:functoriality_!_pro_coherent_modules}
			The object-wise assignation 
			\[
				(f\colon A \to B) \rightsquigarrow (f^{ ! }\colon \ProCoh_{ A } \to \ProCoh_{ B })
			\]
			is the shadow of an $\infty$-functor 
			\[
				\ProCoh \colon \CAlg^{ \cn, \mathrm{coh} }_{ \C } \to \Pr^{ L, st }.
			\]	
			Moreover there exists a natural transformation
			\[
				i\colon \Qcoh^{ * } \longrightarrow \ProCoh
			\]
			in $\Fun(\CAlg_{ \C }^{ \cn }, \Pr^{ L, st })$ which gives back, on objects, the $i$ introduced in \cref{prop:modules_into_pro_coherent}.
			Finally, the two functors and the natural transformation actually take value in $\CAlg(\Pr^{ L, st })$; that is, everything is compatible with the tensor products.
		\end{prop}
		\begin{proof}
			Recall that we denote by $\Qcoh^{ * }$ the symmetric monoidal $\infty$-functor
			\[
				\CAlg^{ \cn }_{ \C } \to \Pr^{ L, st }
			\]
			sending $A$ to $\Mod_{ A }$ and $A \to B$ to the base change functor $B \otimes_{ A } -\colon \Mod_{ A } \to \Mod_{ B }$.\\
			The proof comes from a slight modification of \cite[Theorem 2.52]{BrantnerCamposNuiten:PDOperadsExplicitPartitionLieAlgebras}. We borrow their notations.
			Let us replace $\infcatname{Add}^{ \mathrm{coh}, \mathrm{poly} }$ by the (non full) subcategory $\infcatname{Add}^{ \mathrm{coh}, \mathrm{lin} }$ where we only consider \emph{linear} $\infty$-functors (i.e.\ that preserves finite products). It is equipped with a symmetric monoidal structure making the inclusion into $\infcatname{Add}^{ \mathrm{coh}, \mathrm{poly}}$ lax symmetric monoidal.
			It is easy to see (by modifying \cite[Proposition 2.40]{BrantnerCamposNuiten:PDOperadsExplicitPartitionLieAlgebras}), that then the target can be restricted from $\Pr^{ st }_{ \Sigma }$ to $\Pr^{ L, st }$, which is again equipped with a symmetric monoidal structure making $\Pr^{ L, st } \hookrightarrow \Pr^{ st }_{ \Sigma }$ lax symmetric monoidal.
			The restriction of the main diagram of the mentioned theorem to these subcategories continues to respect the symmetric monoidal structures. Imposing $\C$-linearity simply amounts to consider additive categories endowed with an action of the commutative algebra $\Vect^{ \omega }_{ \C }$. That is, we have a symmetric monoidal natural transformation 
			\begin{diag}
				\infcatname{Add}^{\mathrm{coh}, \mathrm{lin}}_{\C} \ar[r, bend left=50, ""{name=U, below}, "\Mod"] \ar[r, bend right=50, ""{name=D}, "\infcatname{QC}^{\vee}"'] & \Pr^{L, st}_{\C} \arrow[Rightarrow, from=U, to=D, "k"]
			\end{diag}
			which is obtained using \cref{lemma:extending_functor_mod_procoh} (where instead of $A$ one uses $\infcatname{A} \in \infcatname{Add}^{ \coh, \mathrm{lin} }$).
			Let us precompose the above diagram with the symmetric monoidal $\infty$-functor
			\[
				\Vect^{ \omega }\colon \CAlg^{ \cn, \mathrm{coh}  }_{ \C } \longrightarrow \infcatname{Add}^{ \mathrm{coh}, \mathrm{lin} }_{\C }, \qquad A \mapsto \Vect^{ \omega }_{ A }
			\]
			sending $A \to B$ to the restricted base change.
			Observe that $\Qcoh^{ * }$ is identified with the top composition by \cite[Example 2.2]{BrantnerCamposNuiten:PDOperadsExplicitPartitionLieAlgebras} and by abuse of notation denote again by $\infcatname{QC}^{ \vee }$ the bottom composition. We then have a symmetric monoidal natural transformation
			\begin{diag}
				\CAlg^{\cn, \mathrm{coh}}_{\C} \ar[r, bend left=50, ""{name=U, below}, "\Qcoh^{*}"] \ar[r, bend right=50, ""{name=D}, "\infcatname{QC}^{\vee}"'] & \Pr^{L, st}_{\C} \arrow[Rightarrow, from=U, to=D, "k"].
			\end{diag}
			We claim that $\infcatname{QC}^{ \vee }$ is the symmetric monoidal $\infty$-functor enhancing the object-wise assignation $\ProCoh$ and that $k$ is the (symmetric monoidal) natural transformation enhancing the object-wise $i$.
			It is clear by construction that $\infcatname{QC}^{ \vee }_{ A } \simeq \ProCoh_{ A }$ as $\infty$-categories and the fact that $k_{ A } \simeq i_{ A }$ is clear by construction. 
			Let us check what happens on $1$-morphisms, say $f\colon A \to B$. Precisely, we need to prove that 
			\[
				f^{ ! }\colon \ProCoh_{ A } \to \ProCoh_{ B }
			\]
			is obtained by right-left extension of $B \otimes_{ A } -\colon \Vect_{ A }^{ \omega } \to \Vect_{ B }^{ \omega }$, according to \cite[Corollary 2.49]{BrantnerCamposNuiten:PDOperadsExplicitPartitionLieAlgebras}. It is, again, immediate thanks to \cref{lemma:extending_functor_mod_procoh} (notice also \cref{prop:pro_coh_!_pullback_compatibilities}).\\
			Let us finally verify that the symmetric monoidal structure defined by renormalized Day convolution on $\ProCoh_{ A }$ in \cref{prop:pro_coherent_tensor_product} coincides with the one coming from the (exterior) symmetric monoidal structure on the functor $\infcatname{QC}^{ \vee }$.
			This boils down to verify that the former tensor product is obtained by $!$-pullback along the diagonal morphism. That is, given $A \in \CAlg^{ \cn, \coh }_{ \C }$ with multiplication $m\colon A \otimes_{ \C } A \to A$, we need to verify that the composition
			\[
				\ProCoh_{ A } \times \ProCoh_{ A } \to \ProCoh_{ A } \otimes \ProCoh_{ A } \simeq \ProCoh_{ A \otimes_{ \C } A } \stackrel{m^{ ! }}{\longrightarrow} \ProCoh_{ A }
			\] 
			coincides with the previously defined $\otimes^{ ! }$ (and all coherent homotopies). It suffices to check it for $M, N \in \Coh_{ A }^{ \op }$ (seen through Yoneda embedding in $\ProCoh_{ A }$). By compatibility of $!$-pullback and $*$-pullback on almost perfect modules through Yoneda embedding (see \cref{prop:pro_coh_!_pullback_compatibilities}) we deduce that the image of $(M, N)$ is
			\[
				\Map_{ A }(m^{ * }(M \otimes_{ \C } N), -) \simeq \Map_{ A }(M \otimes_{ A } N, -) \in \Fun_{ \ex, \afp }(\Mod_{ A }^{ + }, \S).
			\]
			The right-hand side coincides exactly with $j(M) \otimes^{ ! } j(N)$ thanks to \cref{prop:pro_coherent_tensor_product}.
		\end{proof}

		We are finally ready to define pro-coherent sheaves on general prestacks (although they will be badly behaved). 

		\begin{defn}[{\cite[Definition B.15]{BrantnerMagidsonNuiten:FormalIntegrationDerivedFoliations}}]
			\label{defn:pro_coherent_sheaves_prestacks}

			For each prestack $X\colon \CAlg^{ \cn }_{ \C } \to \S$, the (symmetric monoidal) presentable stable $\infty$-category of pro-coherent sheaves is defined as
			\[
				\ProCoh_{ X } \coloneqq \lim_{ S \in \dAff_{ -/X } } \ProCoh_{ S } \in \CAlg(\Pr^{ L, st }_{ \C })
			\]
			where $!$-pullbacks are considered. More generally, we define an $\infty$-functor
			\[
				\ProCoh\colon \PreStk_{ \C }^{ \op } \longrightarrow \CAlg(\Pr^{ L, st }_{ \C })
			\]
			as right Kan extension along the Yoneda embedding $\dAff_{ \C }^{ \op } \hookrightarrow \PreStk_{ \C }^{ \op }$ of the functor $\ProCoh$ from \cref{prop:functoriality_!_pro_coherent_modules}.
			The right Kan extension of $i\colon \Qcoh^{ * } \to \ProCoh$ along the Yoneda embedding gives a symmetric monoidal transformation
			\[
				i_{ X }\colon \Qcoh(X) \longrightarrow \ProCoh_{ X }.
			\]
		\end{defn}

		This allows us to talk about pro-coherent sheaves, their monoidal structure and their (monoidal) $!$-pullbacks on any prestack. As before, this is very general and quite often ill-behaved. 
		For example it is not always compactly generated; it is the case for \emph{locally coherent} derived schemes, see \cite[Proposition B.16]{BrantnerMagidsonNuiten:FormalIntegrationDerivedFoliations}.

	\subsection{Relation with ind-coherent sheaves}
		\label{subsubsection:pro_coh_vs_ind_coh}

		In this subsection, we will introduce Serre duality on ind-coherent sheaves and relate them with pro-coherent sheaves. Moreover, we will compare how they interact with quasi-coherent sheaves and how this comparison is functorial with respect to pullbacks.\\
		Let us start with some observations for derived schemes\footnote{They also hold for ind-inf-schemes, but we focus on derived schemes for readability.}. By construction,
		\[
			\IndCoh(X) \simeq \Ind(\Coh(X))
		\]
		is compactly generated, and hence dualizable in $\Pr^{ L, st }_{ \C }$ and with dual given by $\Ind(\Coh(X)^{ \op })$. This holds thanks to \cite[Chapter 1, Proposition 7.3.2]{GR:DAG1}.
		Moreover, observe that $\Corr(\dSch_{ \C }^{ \afp })$ is rigid, i.e.\ each object $X$ is dualizable and self-dual thanks to the two diagonal correspondences
		\[
			\begin{tikzcd}
				& \ar[dl, "\Delta"'] X \ar[dr]   &  & & \ar[dl] X \ar[dr, "\Delta"] & \\
				X \times X & & \Spec(\C) &   \Spec(\C) & & X \times X
			\end{tikzcd}
		\]
		which are, respectively, evaluation and co-evaluation. As all small rigid $\infty$-categories, its categorical dual is its opposite and it is self dual, i.e.\ passing to the dual object defines an equivalence
		\[
			\varpi \colon \Corr(\dSch_{ \C }^{ \afp })^{ \op } \stackrel{\sim}{\longrightarrow} \Corr(\dSch_{ \C }^{ \afp }).
		\]
		By \cite[Chapter 7, Proposition 2.3.4]{GR:DAG1}, such an equivalence is just ``swapping the arrows'', i.e.\ on $1$-morphisms is given by 
		\[
			(X \leftarrow Z \rightarrow Y) \mapsto (Y \leftarrow Z \rightarrow X).	
		\]
		Combining these observations with the symmetric monoidal structure on $\IndCoh$ (as a functor from correspondences) coming from \cref{thm:ind_coh_correspondence_schemes}, one formally obtains the following theorem.

		\begin{thm}[{Serre duality, \cite[Chapter 3, Theorem 6.2.2]{GR:DAG2}}]
			\label{thm:serre_duality}
			There is a commutative diagram of functors
			\begin{diag}
				(\Corr(\indinfSch_{\laft}))^{\op} \ar[r, "(\IndCoh)^{\op}"] \ar[d, "\varpi"] & (\Pr^{L, st}_{\C})^{\mathrm{dualizable}, \op} \ar[d, "(-)^{\vee}"] \\
				\Corr(\indinfSch_{\laft}) \ar[r, "\IndCoh"] & (\Pr^{L, st}_{\C})^{\mathrm{dualizable}}.
			\end{diag}
			This implies, in particular, that for each $X$ there exists an ``involutive'' equivalence
			\[
				\mathbb{D}^{ \mathrm{Serre} }_{ X }\colon \Ind(\Coh(X)^{ \op }) \stackrel{\sim}{\longrightarrow} \Ind(\Coh(X))
			\]
			and for each $f\colon X \to Y$ we have a natural ``involutive'' identifications
			\[
				(f^{ ! })^{ \vee } \simeq f_{ * }^{ \IC }, \qquad (f_{ * }^{ \IC })^{ \vee } \simeq f^{ ! }.
			\]
		\end{thm}

		To be precise, in the above equivalences we identify $(f^{ ! })^{ \vee }$ with the composition
		\[
			\Ind(\Coh(X)) \stackrel{\mathbb{D}^{ \mathrm{Serre} }_{ X }}{\longrightarrow} \Ind(\Coh(X)^{ \op }) \simeq (\Ind(\Coh(X)))^{ \vee } \stackrel{(f^{ ! })^{ \vee }}{\longrightarrow} (\Ind(\Coh(X)))^{ \vee } \stackrel{\mathbb{D}^{ \mathrm{Serre} }_{ Y }}{\longrightarrow} \Ind(\Coh(Y)).
		\]

		It is called ``Serre duality'' thanks to the following proposition.

		\begin{prop}[{\cite[Proposition 9.5.7]{G:IndCoh}}]
			\label{prop:link_classical_Serre_duality}
			For $X$ a derived scheme, we can identify the equivalence
			\[
				\mathbb{D}^{ \mathrm{Serre} }_{ X }\colon \Ind(\Coh(X)^{ \op }) \longrightarrow \Ind(\Coh(X))
			\]
			with the $\Ind$-extension of the ``classical'' Serre duality functor
			\[
				\mathbf{D}^{ \mathrm{Serre} }_{ X } \coloneqq \IntHom(-, \omega_{ X })\colon \Coh(X)^{ \op } \to \Coh(X),
			\]
			given by the $\Qcoh(X)$-enriched hom objects of $\Ind(\Coh(X))$ (described in \cref{prop:upsilon_symm_monoidal}).
		\end{prop}

		Let us now compare ind-coherent sheaves with pro-coherent ones.
		Let us start on derived affine schemes.

		\begin{prop}
			\label{prop:functorial_equivalence_procoh_indcoh}
			Let $A \in \CAlg^{ \cn }_{ \C }$ be coherent. The Serre duality equivalence can be thought as an equivalence between pro-coherent and ind-coherent modules over $A$, i.e.
			\[
				\Serre_{ A }\colon \ProCoh_{ A } \stackrel{\sim}{\longrightarrow} \Ind(\Coh_{ A }).
			\]
			Moreover, the following diagram is commutative
			\begin{diag}
				\Mod_A \ar[d, "i_A"] \ar[dr, "\Upsilon_A"] & \\
				\ProCoh_A \ar[r, "\sim"', "\Serre_A"] & \Ind(\Coh_A).
			\end{diag}
			and natural in $A$, i.e.\ compatible with $*$-pullbacks on modules and $!$-pullbacks on pro-coherent and ind-coherent modules. It can be enhanced to a commutative diagram of \emph{symmetric monoidal} functors $\CAlg^{ \cn, \mathrm{coh} } \to \Pr^{ L, st }_{ \C }$; in particular, it is compatible with the previously defined tensor products.
		\end{prop}
		\begin{proof}
			The naturality in $A$ of $i\colon \Mod_{ - } \to \ProCoh$, together with its symmetric monoidal structure, comes from \cref{prop:functoriality_!_pro_coherent_modules}. The naturality of $\Upsilon\colon \Mod_{ - } \to \IndCoh$ comes from \cref{prop:upsilon_symm_monoidal}.
			To study the naturality of Serre duality, observe that \cref{thm:serre_duality} implies, in particular, that for a map $f\colon A \to B$ in $\CAlg^{ \cn, \coh }_{ \C }$ the following diagram commutes
			\begin{diag}
				\ProCoh_A \ar[d, "(f_*^{\IC})^{\vee}"] \ar[r, "\sim"', "\Serre_A"] & \Ind(\Coh_A) \ar[d, "f^!"] \\
		\ProCoh_B \ar[r, "\sim", "\Serre_B"'] & \Ind(\Coh_B).
			\end{diag}

			Observe that the categorical dual of $f_{ * }^{ \IC }$ is identified with the composition
			\[
				\ProCoh_{ A } \simeq \Fun^{ L }(\Ind(\Coh_{ A }), \S) \stackrel{- \circ f_{ * }^{ \IC }}{\longrightarrow} \Fun^{ L }(\Ind(\Coh_{ B }), \S) \simeq \ProCoh_{ B }
			\]
			It is easy to see that, by first restricting to coherent modules and then left Kan extending to eventually coconnective modules (recall that $\Ind(\Coh_{ A })^{ + } \simeq \Mod_{ A }^{ + }$ and the restriction of $f_{ * }^{ \IC }$ to here is simply $f_{ * }$), that it is also identified with
			\[
				\ProCoh_{ A } \simeq \Fun_{ \lex, \afp }(\Mod_{ A }^{ + }, \S) \stackrel{- \circ f_{ * }}{\longrightarrow} \Fun_{ \lex, \afp }(\Mod_{ B }^{ + }, \S).
			\]
			This corresponds to the pro-coherent $!$-pullback defined in \cref{defn:!-pullback_pro_coherent_modules}.

			The fact that Serre duality is symmetric monoidal is now a corollary of the fact that, on both sides, the tensor product is computed as $!$-pullback along the diagonal and we just proved naturality with respect to $!$-pullbacks on both sides (one final necessary technical observation is that $\Serre_{ A \otimes_{ \C } B } \simeq \Serre_{ A } \otimes \Serre_{ B }$ and it stems from the fact that the commuting diagram in \cref{thm:serre_duality} is actually a diagram of symmetric monoidal $\infty$-categories).
		\end{proof}

		This allows us to finally identify the symmetric monoidal $\infty$-categories of ind-coherent and pro-coherent sheaves on any prestack.

		\begin{coroll}
			\label{coroll:serre_duality_prestack}
			For any laft prestack $X$, (the right Kan extension of) Serre duality gives a natural equivalence in $\CAlg(\Pr^{ L, st })_{ \Qcoh(X)/- }$
			\[
				\Serre_{ X }\colon \ProCoh_{ X } \stackrel{\sim}{\longrightarrow} \IndCoh(X).
			\]
			It is natural in $X$, under $!$-pullbacks on both sides.
		\end{coroll}

		Let us do a final comment on all this story before moving on.

		\begin{remark}
			\label{remark:pro_coh_are_better}
			
			Pro-coherent sheaves offer two major advantages compared to ind-coherent ones. Firstly, they can be defined also for commutative algebras which are not coherent.
			Last, and definitely not least, the $!$-pullback and $!$-tensor product is way easier to compute here (via a precomposition and a renormalized Day convolution, respectively).
			Moreover, colimit-preserving symmetric monoidal functors out of pro-coherent modules are equivalent to regular symmetric monoidal functors from (the opposite of) almost perfect modules (see \cref{prop:pro_coherent_tensor_product} and \cref{lemma:extending_functor_mod_procoh}). This will be crucial later on.
		\end{remark}

%% file: filtered.tex
\section{Filtered objects}
			\label{subsection:filtered_objects}
			In this section, we will review the definitions and some important properties of filtered and graded objects in a stable presentable $\infty$-category $\infcatname{C}$. Our main references are \cite{GwilliamPavlov:FilteredDerivedCat} and \cite{Moulinos:GeometryFiltrations}. Let $\infcatname{C}$ be a stable (presentable) $\infty$-category throughout this part. We will denote by $\Z^{ \leq }$ the category associated with the poset $\Z$ with the standard order, and by $\Z^{ \mathrm{disc} }$ the category associated with the discrete set $\Z$. As usual, we will think of them automatically as $\infty$-categories.

			\begin{defn}
				\label{defn:filtered_objects}
				The $\infty$-category of \emph{filtered objects} in $\infcatname{C}$ is 
				\[
					\infcatname{C}^{\Fil} \coloneqq \Fun( \Z^{ \leq }, \infcatname{C})
				\]
				and the $\infty$-category of \emph{graded objects} in $\infcatname{C}$ is 
				\[
					\infcatname{C}^{ \gr } \coloneqq \Fun(\Z^{ \mathrm{disc} }, \infcatname{C}).
				\]
				If $\infcatname{C}$ admits $\Z$-indexed limits, we can define the (stable) $\infty$-category of \emph{(left) filtered complete objects} in $\infcatname{C}$ as the full subcategory $\infcatname{C}^{ \cpl } \subset \infcatname{C}^{ \Fil }$ spanned by the objects $(\dots \to F^{ n }X \to F^{ n+1 }X \to \dots)$ such that $\varprojlim_{ n } F^{ n }X = 0$.
			\end{defn}
			Observe that both $\infcatname{C}^{ \Fil }$ and $\infcatname{C}^{ \gr }$ are presentable and stable as well, being functor categories (see \cite[Prop 1.1.3.1]{Lurie:HA} and \cite[Prop 5.5.3.6]{Lurie:HTT}).
			Moreover, $\infcatname{C}^{ \cpl }$ is presentable and stable too, by \cite[Thm 2.5]{GwilliamPavlov:FilteredDerivedCat}.
			Finally, filtered and graded objects are obviously functorial, meaning they define $\infty$-functors
			\[
				(-)^{ \Fil } \coloneqq \Fun(\Z^{ \leq }, -)\colon \Pr^{ L, st } \to \Pr^{ L, st }, \qquad (-)^{ \gr }\coloneqq \Fun(\Z^{ \mathrm{disc} }, -)\colon \Pr^{ L, st} \to \Pr^{ L, st }.
			\]
			
			\begin{remark}
				Our terminology differs from \cite{GwilliamPavlov:FilteredDerivedCat}: their ``filtered objects'' are what we call \emph{complete} filtered objects, and their ``sequences'' are what we call \emph{filtered} objects.
			\end{remark}

			\begin{example}
				\label{example:filtered_vector_space}
				One classical example of a filtered object is simply a filtered vector space $V = \bigcup_{ i \in \Z } V_{ i }$, corresponding to the increasing filtration
				\[
					\dots \hookrightarrow  F^{ n }V = \bigcup_{ i \geq -n } V_{ i } \hookrightarrow F^{ n+1 }V = \bigcup_{ i \geq -n-1 } V_{ i } \hookrightarrow \dots
				\]
				The completeness condition, in the classical world, is exactly the separatedness condition, that is, $\bigcap_{ i \in \Z } V_{ i } = 0$.
				We need to be careful, though: the limit in \cref{defn:filtered_objects} is a homotopical limit, and therefore it can differ from a classical limit (in a model category, for example).
				See \cite[Remark 2.13]{GwilliamPavlov:FilteredDerivedCat}.
			\end{example}

			\begin{remark}
				\label{remark:right_filtered_complete_objects}
				A (right) filtered complete object is simply a left filtered complete object in $\infcatname{C}^{ \op }$. That is, a diagram in $\infcatname{C}$ 
				\[
					\dots \to F_{ n }X \to F_{ n-1 }X \to F_{ n-2 }X \to \dots
				\]
				such that $\colim_{ n } F_{ n }X = 0$. The object of interest, in this case, will be then $\lim_{ m } F_{ m }X$. This might be closer to the intuition from classical $I$-adic completion, where one example of such tower for $\Z$ with ideal $(p) \subset \Z$ is
				\[
					\dots \twoheadrightarrow \Z/p^{ n }\Z \twoheadrightarrow \Z/p^{ n-1 }\Z \twoheadrightarrow \dots \twoheadrightarrow \Z/p^{ 2 }\Z \twoheadrightarrow \Z/p\Z \twoheadrightarrow 0 \to \dots
				\]
				More about the relation with the classical point of view can be found at \cite[Lemma 2.9]{GwilliamPavlov:FilteredDerivedCat} and \cite[Remark 3.11]{Antieau:SpectralSequencesDecalage}.
				In this text, we will only work with left complete objects.
			\end{remark}

			\begin{remark}
				\label{remark:pass_left_right_complete_filtered_object}
				Each left filtered complete object $FX$ (which we think as an increasing filtration on $X \coloneqq \colim_{ n } F^{ n }X$) gives rise to a right filtered complete object $FY$ defined as
				\[
					F^{ n }Y \coloneqq \cofib(F^{ -n }X \to X).
				\]
				It is immediate to verify that $\colim_{ m } F^{ m }Y \simeq 0$ and $\lim_{ m } F^{ m }Y \simeq X$. In some sense, this procedure endows $X$ with the decreasing filtration corresponding to the ``associated system of quotients'', and exchanges positive and negative filtered objects. It is clearly reversible.
			\end{remark}
	
			Let us now assume $(\infcatname{C}, \otimes)$ is a presentably monoidal stable $\infty$-category, that is, an object of $\Alg(\Pr^{L, st })$. In this case, by \cite[Example 2.2.6.10]{Lurie:HA}, the category $\infcatname{C}^{ \Fil }$ inherits a monoidal structure by left Day convolution (where the symmetric monoidal structure on $N(\Z^{ \leq })$ is given by the sum $+$), denoted by $\otimes^{ \Fil }$. 
			Object-wise, given $FX, FY \in \infcatname{C}^{ \Fil }$ their tensor product is given by the following left Kan extension
			\begin{diag}
				N(\Z^{\leq}) \times N(\Z^{\leq}) \ar[d, "(-) + (-)"] \ar[r, "FX \times FY"] & \infcatname{C} \times \infcatname{C} \ar[r, "(-) \otimes (-)"] & \infcatname{C} \\
				N(\Z^{\leq}) \ar[urr, "FX \otimes^{\Fil} FY" below, dashed]
			\end{diag}
			which, by the left Kan extension formula, means 
			\[
				F^{ n }(FX \otimes^{ \Fil } FY) \simeq \colim_{ p+q \leq n } F^{ p }X \otimes F^{ q }Y.
			\]
			We immediately see that the monoidal unit is given by 
			\[
				\langle 0, 1_{ \infcatname{C} }\rangle \coloneqq \left( \dots \to 0 \to 1_{ \infcatname{C} } \to 1_{ \infcatname{C} } \to \dots \right)
			\] 
			which is $0$ in negative degrees and $1_{ \infcatname{C} }$ (monoidal unit of $\infcatname{C}$) in non-negative ones, with identity maps.\\	
			Similarly, equipping $N(\Z^{ \mathrm{disc} })$ with the sum $+$ as the tensor product, we can equip $\infcatname{C}^{ \gr }$ with the left Day convolution as well (similar diagram as above), denoting the resulting monoidal structure by $\otimes^{ \gr }$. For two graded objects $X, Y \in \infcatname{C}^{ \gr }$ this simply means
			\[
				(X \otimes^{ \gr } Y)(n) \simeq \bigoplus_{ p + q = n } X(p) \otimes Y(q).
			\]
			The monoidal unit is $(1_{ \infcatname{C} })_{ 0 }$, that is, the graded object having $1_{ \infcatname{C} }$ in weight $0$ and $0$ elsewhere.\\
			As remarked in \cite[Remark 2.3]{Fu:DualityLieAlgbdFoliations}, by functoriality of left Day convolution (and by the fact that $\Pr^{L, st}\hookrightarrow \Pr^{ L }$ is a symmetric monoidal localization), we can upgrade the filtered and graded functors at the level of algebras
			\[
				(-)^{ \Fil }\colon \Alg(\Pr^{ L, st }) \to \Alg(\Pr^{ L, st }), \quad (-)^{ \gr }\colon \Alg(\Pr^{ L, st }) \to \Alg(\Pr^{ L, st }).
			\]
			We will soon describe a (monoidal) natural transformation between them $(-)^{ \Fil } \to (-)^{ \gr }$. \\
			Here are some useful functors and properties.
			\begin{enumerate}[label=(\emph{\alph*}), ref=(\emph{\alph*})]
				\item \label{filtered:constant_functor} The \emph{constant functor} $\const_{ \infcatname{C} }\colon \infcatname{C} \simeq \Fun(*, \infcatname{C}) \to \infcatname{C}^{ \Fil }$, i.e. restriction along the unique functor $\Z^{ \leq } \to *$, is right adjoint to the colimit functor 
					\[
						(-)_{ \infty }\colon \infcatname{C}^{ \Fil }\to \infcatname{C}, \qquad (\dots \to F^{ n }X \to F^{ n+1 }X \to \dots) \mapsto \colim_{ n }F^{ n }X
					\]
					often referred to as the ``underlying object'' (or as ``realization''). Since $\infcatname{C}$ is presentably monoidal, this functor is strongly monoidal. Moreover, everything is natural in $\infcatname{C}$.

				\item \label{filtered:associated_graded} The \emph{associated graded functor}, denoted by $\gr_{ \infcatname{C} }\colon \infcatname{C}^{ \Fil } \to \infcatname{C}^{ \gr }$, is the composition 
					\begin{gather*}
						  \infcatname{C}^{ \Fil }=\Fun(\Z^{ \leq }, \infcatname{C}) \stackrel{\Delta}{\longrightarrow} \prod_{ n \in \Z } \Fun(\Z^{ \leq }, \infcatname{C}) \stackrel{\prod_{ n \in \Z } \textrm{restr} }{\longrightarrow} \prod_{ n \in \Z } \Fun([n-1] \to [n], \infcatname{C}) \\ 
						  	\simeq  \prod_{ n \in \Z } \Fun([n-1] \to [n], \infcatname{C})  \stackrel{\prod_{ n } \textrm{cofib}}{\longrightarrow} \prod_{ n \in \Z } \infcatname{C} = \infcatname{C}^{ \gr }.
					  \end{gather*}
					  That is, for a filtered object $FX$, the associated graded is given, in weight $m$, by 
					  \[
					  	\gr(FX)(m) \simeq \cofib(F^{ m-1 }X \to F^{ m }X).
					  \]
					  The associated graded functor preserves limits and colimits (which are computed term-wise in $\infcatname{C}^{ \Fil }$ and are preserved by all functors above, since $\infcatname{C}$ is stable). The morphisms in $\infcatname{C}^{ \Fil }$ which are sent to equivalences in $\infcatname{C}^{ \gr }$ are called ``graded equivalences''. 
					  Moreover, it is natural in $\infcatname{C}$; that is, it describes a natural transformation $\gr\colon (-)^{ \Fil } \to (-)^{ \gr }$.

				  \item \label{filtered:right_triv} The right adjoint to the associated graded functor is the \emph{(right) trivial} functor, denoted by $\triv_{ \infcatname{C} }\colon \infcatname{C}^{ \gr } \to \infcatname{C}^{ \Fil }$, defined by the composition
					  \begin{gather*}
						  \infcatname{C}^{ \gr } = \prod_{ n \in \Z } \infcatname{C} \stackrel{\prod_{ n \in \Z } (0 \to \id)}{\longrightarrow} \prod_{ n \in \Z } \Fun([n-1] \to [n], \infcatname{C}) \stackrel{\prod_{ n \in \Z } \mathrm{Ran}_{ n }}{\longrightarrow} \prod_{ n \in \Z } \Fun(\Z^{ \leq }, \infcatname{C}) \\
					\simeq \prod_{ n \in \Z } \Fun(\Z^{ \leq }, \infcatname{C})  \stackrel{\lim_{ n \in \Z }}{\longrightarrow} \Fun(\Z^{ \leq }, \infcatname{C}) = \infcatname{C}^{ \Fil }.
					  \end{gather*}
					  An easy computation, where by $\mathrm{Ran}_{ n }$ we mean the right Kan extension along $\{[n] \to [n+1]\} \subseteq \Z^{ \leq }$, gives us the explicit formula
					  \[
						  \triv(X) = (\dots \stackrel{0}{\longrightarrow} X(n) \stackrel{0}{\longrightarrow} X(n+1) \stackrel{0}{\longrightarrow} \dots),
					  \]
					  that is, it equips the graded object $X$ with the trivial filtration, which is $X(n)$ at filtered-degree $n$ and has zero maps between them.
					  Let us observe that, since limits and colimits of graded (and filtered) objects are computed degree-wise, $\triv\colon \infcatname{C}^{ \gr } \to \infcatname{C}^{ \Fil }$ commutes both with limits and colimits.

				  \item \label{filtered:left_triv} The left adjoint to the associated graded functor is the \emph{(left) trivial} functor, denoted by $\triv_{ \infcatname{C} }^{ L }\colon \infcatname{C}^{ \gr } \to \infcatname{C}^{ \Fil }$, defined by the composition
					  \begin{gather*}
						  \infcatname{C}^{ \gr } = \prod_{ n \in \Z } \infcatname{C} \stackrel{\prod_{ n \in \Z } (\id[-1] \to 0)}{\longrightarrow} \prod_{ n \in \Z } \Fun([n-1] \to [n], \infcatname{C}) \stackrel{\prod_{ n \in \Z } \mathrm{Lan}_{ n }}{\longrightarrow} \prod_{ n \in \Z } \Fun(\Z^{ \leq }, \infcatname{C}) \\
						  \simeq \prod_{ n \in \Z } \Fun(\Z^{ \leq }, \infcatname{C}) \stackrel{\colim_{ n \in \Z }}{\longrightarrow} \Fun(\Z^{ \leq }, \infcatname{C}) = \infcatname{C}^{ \Fil }.
					  \end{gather*}
				 	  An easy computation, where by $\mathrm{Lan}_{ n }$ we mean the left Kan extension along $\{[n] \to [n+1]\} \subseteq \Z^{ \leq }$, gives us the explicit formula
				      \[
						  \triv^{ L }(X) = \left( \dots \stackrel{0}{\longrightarrow} X(n+1)[-1] \stackrel{0}{\longrightarrow} X(n+2)[-1] \stackrel{0}{\longrightarrow} \dots \right), 
				      \]
					  that is, it equips the graded object $X$ with the shifted (both in weights and in cohomological degrees) trivial filtration, which is $X(n+1)[-1]$ at filtered-degree $n$ and has zero maps between them.
					  As before, $\triv^{ L }\colon \infcatname{C}^{ \gr } \to \infcatname{C}^{\Fil }$ commutes both with limits and colimits.

				  \item \label{filtered:completion_functor} The \emph{completion functor} $\widehat{(-)}\colon \infcatname{C}^{ \Fil } \to \infcatname{C}^{ \cpl }$ is the left adjoint of the inclusion. That is, complete filtered objects form a stable localization for filtered objects (along the so-called graded equivalences, see \cite[Lemma 2.15]{GwilliamPavlov:FilteredDerivedCat}). Concretely, on objects it sends $X \in \infcatname{C}^{ \Fil }$ to the filtered object $\widehat{X}$ such that
					\[
						F^{ n }\widehat{X} \simeq \cofib(X(-\infty) \to F^{ n }X),
					\]
					for $X(-\infty) \coloneqq \lim_{ m }F^{ m }X$ (see \cite[Lemma 2.9]{GwilliamPavlov:FilteredDerivedCat}). That is, we just need to mod out the part at ``$-\infty$''.

				\item \label{filtered:warning_complete} \textbf{Warning:} Since we are doing everything directly in the $\infty$-categorical world, the limit defined above for $X(-\infty)$ is to be interpreted as an homotopy limit (when working, for example, with model categories). This means, also, that our definition of complete filtration (an ``homotopical'' version) doesn't always agree with the classical definition of separatedness of a filtration: see \cite[Remark 2.13]{GwilliamPavlov:FilteredDerivedCat} and \cref{example:filtered_vector_space}. 

				\item \label{filtered:triv_gr_complete_adjunctions} Observe that trivial (and shifted trivial) filtered objects are always complete. This implies that, restricting $\gr$ to complete filtered objects, we still have the adjunctions
					\[
						\adjunction{\triv^{ L }}{\infcatname{C}^{ \gr }}{\infcatname{C}^{\cpl}}{\gr}, \qquad \adjunction{\gr}{\infcatname{C}^{ \cpl }}{\infcatname{C}^{ \gr }}{\triv}.
					\]

				\item \label{filtered:complete_tensor_product} There is a (symmetric) monoidal structure on $\infcatname{C}^{ \cpl }$ given, object-wise, by the formula 
					\[
						X \otimes^{ \cpl } Y \coloneqq \widehat{X \otimes^{ \Fil } Y}, \qquad \forall X, Y \in \infcatname{C}^{ \cpl }.
					\] 
					This follows from general facts about localizations in (symmetric) monoidal $\infty$-categories: see \cite[Theorem 2.25]{GwilliamPavlov:FilteredDerivedCat}. This implies that the completion functor is strong monoidal, by definition of completed tensor product, and its right adjoint (inclusion) is lax (symmetric) monoidal. 
					\begin{remark}
						We wrote ``(symmetric)'' because the same result holds for any $\O^{ \otimes }$-monoidal structure.
					\end{remark}

				\item \label{filtered:filtered_internal_hom} If $\infcatname{C}$ is also closed (i.e.\ admits internal mapping objects), then $\infcatname{C}^{ \Fil }$ is closed. Given $X, Y \in \infcatname{C}^{ \Fil }$, their filtered mapping object is defined, at weight $n$, as the object of ``natural transformations of degree $n$'', given by 
					\[
						F^{ n }\IntHom(X, Y) \coloneqq \int_{ m \in \Z } \IntHom_{ \infcatname{C} }\left( F^{ m }X, F^{ m+n }Y \right),
					\]
					where we used the classical integral notation for ends in $\infcatname{C}$. The fact that this is indeed an internal mapping object, i.e.\ a right adjoint for the filtered tensor product, comes from \cite[Lemma 2.27]{GwilliamPavlov:FilteredDerivedCat}: they state it only in the case of a complete filtered object, but we can observe that the proof works equally well for filtered objects (remark the heavy usage of \cite[Remark 2.24]{GwilliamPavlov:FilteredDerivedCat}).

				\item \label{filtered:internal_hom_complete} As already mentioned, the above filtered hom object is also the internal hom for $\infcatname{C}^{ \cpl }$. Indeed, let us observe that if $Y$ is complete, then also the filtered mapping object above is so: in fact
					\begin{gather*}
						\lim_{ n } F^{ n }\IntHom(X, Y) \simeq \lim_{ n } \int_{ m \in \Z } \IntHom_{ \infcatname{C} }(F^{ m }X, F^{ m+n }Y) \simeq \int_{ m \in \Z }\IntHom_{ \infcatname{C} }(F^{ m }X, \lim_{ n }F^{ m+n }Y) \\
						\simeq \int_{ m \in \Z } \IntHom_{ \infcatname{C} }(F^{ m }X, 0) \simeq 0,
					\end{gather*}
					since limits commute between them.

				\item \label{filtered:graded_internal_hom} The $\infty$-category of graded objects $\infcatname{C}^{ \gr }$ is closed as well. It is easy to see that the internal mapping object is given by the explicit formula 
					\[
						\IntHom^{ \gr }(X, Y)(n) \coloneqq \prod_{ m \in \Z} \IntHom(X(m), Y(m+n)).
					\]
					Observe that if we consider the full (symmetric) monoidal subcategory of positively graded objects $\infcatname{C}^{ \gr, \geq 0 }$ (see \cref{remark:possible_confusion_positively_filtered}), the internal hom there is the restriction of the above one to positive $m$. That is, 
					\[
						\IntHom^{ \gr, \geq 0 }(X, Y)(n) \coloneqq \prod_{ m \in \N } \IntHom(X(m), Y(m+n))
					\]
					for $n \in \N$.
				
			\item \label{filtered:properties_gr} Let us now consider the restriction of the associated graded functor to complete filtered object $\gr\colon \infcatname{C}^{ \cpl } \to \infcatname{C}^{ \gr }$. It is conservative (simply because the first $\infty$-category is the localization of filtered objects along graded equivalences)\footnote{This is not true for filtered objects, where instead the family $(\gr, (-)_{ -\infty})$ is jointly conservative.}. By the explicit cofiber formula for the completion functor, we obtain that \[
						\gr \circ \widehat{(-)} \simeq \gr(-).
				\]
				This implies that it preserves limits (because $\infcatname{C}^{ \cpl } \hookrightarrow \infcatname{C}^{ \Fil }$ does so) and colimits (because they are computed as completion of standard filtered colimits). Moreover, $\gr$ is a strong monoidal functor by \cite[Proposition 2.26]{GwilliamPavlov:FilteredDerivedCat} and the following formula holds 
				\[
					\gr\left( X \otimes^{ \Fil } Y \right)(n) \simeq \bigoplus_{ p + q = n } \gr(X)(p) \otimes \gr(Y)(q), \qquad \gr(X \otimes^{ \Fil } Y) \simeq \gr(X) \otimes^{ \gr } \gr(Y),
				\]
				that is, it interchanges $\otimes^{ \Fil }$ (and $\otimes^{ \cpl }$) with $\otimes^{ \gr }$. Finally, $\gr$ is also strong closed thanks to \cite[Proposition 2.28]{GwilliamPavlov:FilteredDerivedCat}; that is, we have 
				\[
					\gr(\IntHom^{ \Fil }(FX, FY)) \simeq \IntHom^{ \gr }(\gr(FX), \gr(FY)).
				\]

			\item \label{filtered:triv_lax_monoidal} By formal reasons, the right adjoint $\triv\colon \infcatname{C}^{ \gr } \to \infcatname{C}^{ \Fil }$ is lax symmetric monoidal. Similarly, $\triv$ is also lax symmetric monoidal when seen as a functor towards complete filtered objects.
				  For example, we have 
				  \[
					  \triv( (1_{ \infcatname{C} })_{ 0 } ) \otimes^{ \Fil } \triv( (1_{ \infcatname{C} })_{ 0 } ) \simeq \left(  \dots \stackrel{0}{\longrightarrow} 0 \stackrel{0}{\longrightarrow} 1_{ \infcatname{C} } \stackrel{0}{\longrightarrow} 1_{ \infcatname{C} }[1] \stackrel{0}{\longrightarrow} 0 \stackrel{0}{\longrightarrow} \dots \right)
				  \]
				  where $1_{ \infcatname{C} }$ is in filtered-degree $0$, see \cite[Example 3.2.9]{Lurie:RotationInvarianceKTheory} where $\infcatname{C} = \Sp$ and $\mathbb{A} = \triv(\mathbb{S}_{ 0 })$.

			\item \label{filtered:triv_L_oplax_monoidal} By formal reasons, the left adjoint $\triv^{ L }\colon \infcatname{C}^{ \gr } \to \infcatname{C}^{ \Fil }$ is op-lax symmetric monoidal. Similarly, $\triv^{ L }$ is also op-lax symmetric monoidal when seen as a functor towards complete filtered objects.

			\item \label{filtered:evaluation_degree} The \emph{evaluation in degree $n$} is the functor $(n) = \ev_{ n }\colon \infcatname{C}^{ \gr } \to \infcatname{C}$ selecting the object in degree $n$. It preserves limits and colimits and its left and right adjoint is simply given by (left = right) Kan extension along $* \stackrel{n}{\to} \Z^{ \textrm{disc} }$, denoted by $(-)_{ n }\colon \infcatname{C} \to \infcatname{C}^{ \gr }$. Let us observe that, for $n=0$, the functor $\ev_{ 0 }$ is lax symmetric monoidal. This is due to formal properties of left Day convolution (see \cite[Theorem 2.2.6.2]{Lurie:HA}), since we are restricting along the symmetric monoidal functor $\{0\} \hookrightarrow \Z^{ \mathrm{disc} }$.

				\item \label{filtered:step_functor} The \emph{$n$-th step functor} $\langle n, -\rangle \colon \infcatname{C} \to \infcatname{C}^{ \Fil }$ is the left Kan extension along $* \stackrel{n}{\to} \Z^{ \leq }$ (so the right adjoint is simply selecting the $n$-th filtered component). It sends $A$ to the filtration which is $0$ in degrees smaller than $n$ and $A$ in degrees bigger or equal to $n$, with identity maps. Let us immediately observe, by the explicit formula, that 
					\[
						\langle m, A \rangle \otimes^{ \Fil } \langle n, B \rangle \simeq \langle m+n, A \otimes B \rangle.
					\]
					This implies that $\langle 0, -\rangle$ is strong monoidal.
					We will see in \cref{prop:left_adjoint_step_0} that it admits a left adjoint.

				\item \label{filtered:0_step_functor_monoidal} Let us now consider the symmetric monoidal functor $\langle 0, - \rangle \colon \infcatname{C} \to \infcatname{C}^{ \Fil }$ (which is actually a natural transformation $\id \to (-)^{ \Fil }$ between endomorphisms of $\CAlg(\Pr^{ L, st })$). This gives $\infcatname{C}^{ \Fil }$ a $\infcatname{C}$-module structure: that is, $\infcatname{C}^{ \Fil }$ is also $\infcatname{C}$-enriched. By an easy computation, the $\infcatname{C}$-enriched mapping objects are given by 
				\[
					F^{ 0 }\IntHom^{ \Fil }(FX, FY) = \int_{ m \in \Z } \IntHom(F^{ m }X, F^{ m }Y),
				\]
				which is simply the object of ``internal natural transformations''.

			\end{enumerate}

			Let us now introduce a useful subcategory of filtered objects in a stable $\infty$-category.
			\begin{defn}
				\label{defn:positively_filtered}
				Let $\infcatname{C} \in \Pr^{ L, st}$; the $\infty$-category of \emph{positively filtered objects} of $\infcatname{C}$ is the full subcategory of $\infcatname{C}^{ \Fil }$ spanned by those $FX$ such that $F^{ n }X \simeq 0$ for $n < 0$. It is denoted by $\infcatname{C}^{ \Fil, \geq 0 }$.
			\end{defn}

			\begin{defn}
				\label{defn:positively_graded}
				Let $\infcatname{C} \in \Pr^{ L, st }$; the $\infty$-category of \emph{positively graded objects} of $\infcatname{C}$ is the full subcategory of $\infcatname{C}^{ \gr }$ spanned by those $(X_{ n })_{ n \in \Z}$ such that $X_{ n } \simeq 0$ for $n < 0$. It is denoted by $\infcatname{C}^{ \gr, \geq 0 }$.
				Similarly, the $\infty$-category of \emph{negatively graded objects} of $\infcatname{C}$ is the full subcategory of $\infcatname{C}^{ \gr }$ spanned by those $(X_{ n })_{ n \in \Z }$ such that $X_{ n } \simeq 0$ for $n > 0$. It is denoted by $\infcatname{C}^{ \gr, \leq 0 }$.
			\end{defn}

			\begin{remark}
				\label{remark:possible_confusion_positively_filtered}
				Consider $\infcatname{C}^{ \Fil } \times_{ \infcatname{C}^{ \gr } } \infcatname{C}^{ \gr, \geq 0 }$, i.e.\ the full subcategory of $\infcatname{C}^{ \Fil }$ spanned by filtered objects whose associated graded is $0$ in negative weights.
				This is a different $\infty$-category than $\infcatname{C}^{ \Fil, \geq 0 }$: for example any non-zero constant filtered object belongs to it but not to $\infcatname{C}^{ \Fil, \geq 0 }$. 
			\end{remark}

			\begin{prop}
				\label{prop:positively_graded_objects}
				Let $\infcatname{C} \in \Pr^{ L, st }$. There are natural embeddings 
				\[
					\infcatname{C}^{ \gr, \leq 0 } \hookrightarrow \infcatname{C}^{ \gr }, \qquad \infcatname{C}^{ \gr, \geq 0 } \hookrightarrow \infcatname{C}^{ \gr }.
				\]
				If $\infcatname{C}$ is presentably (symmetric) monoidal then $\infcatname{C}^{ \gr, \leq 0 }$ and $\infcatname{C}^{ \gr, \geq 0 }$ are (symmetric) monoidal full subcategories of $\infcatname{C}^{ \gr }$. If $\infcatname{C}$ is closed, then $\infcatname{C}^{ \gr, \leq 0 }$ and $\infcatname{C}^{ \gr, \geq 0 }$ are closed as well with internal hom given by
				\[
					\IntHom^{ \gr, \geq 0 }(X, Y)(n) \simeq \delta_{ n, \N } \cdot \IntHom^{ \gr }(X, Y)(n), \qquad \IntHom^{ \gr, \leq 0 }(X, Y)(n) \simeq \delta_{ n, \Z_{ \leq 0 } } \cdot \IntHom^{ \gr }(X, Y)(n)
				\]
				for $n \in \Z$ (where $\delta_{ x, A }$ is $1$ if $x \in A$ and $0$ otherwise).
			\end{prop}
			\begin{proof}
				Same reasoning as the proof of \cref{prop:positively_filtered_objects_properties}.
			\end{proof}

			\begin{prop}
				\label{prop:positively_filtered_objects_properties}
				Let $\infcatname{C} \in \Pr^{ L, st }$. There is a natural factorization
				\[
					\infcatname{C}^{ \Fil, \geq 0 } \hookrightarrow \infcatname{C}^{ \cpl } \hookrightarrow \infcatname{C}^{ \Fil }.
				\]
				Moreover, if $\infcatname{C}$ is presentably (symmetric) monoidal then $\infcatname{C}^{ \Fil, \geq 0 }$ is a (symmetric) monoidal subcategory of $\infcatname{C}^{ \Fil }$ and of $\infcatname{C}^{ \cpl }$. If $\infcatname{C}$ is closed, then $\infcatname{C}^{ \Fil, \geq 0 }$ is closed as well with internal hom given by 
				\[
					F^{ n }\IntHom^{ \Fil, \geq 0 }(FX, FY) \simeq \int_{ m \in \N } \IntHom(F^{ m }X, F^{ n+m }Y),
				\]
				for $n \in \N$.
			\end{prop}
			\begin{proof}
				Let us first observe that the inclusion $\infcatname{C}^{ \Fil, \geq 0 } \hookrightarrow \infcatname{C}^{ \Fil }$ can be identified with the left Kan extension 
				\[
					\Fun(\N^{ \leq }, \infcatname{C}) \hookrightarrow \Fun(\Z^{ \leq }, \infcatname{C})
				\]
				along the inclusion of posets $\N^{ \leq } \hookrightarrow \Z^{ \leq }$.
				In particular, this implies that $\infcatname{C}^{ \Fil, \geq 0 }$ is presentable and stable, being a functor category.
				Since the inclusion of posets $\Z_{ < 0 } \hookrightarrow \Z$ is final, one easily sees that all positively filtered objects are complete. That is, we have the sought-for factorization.
				\[
					\infcatname{C}^{ \Fil, \geq 0 } \hookrightarrow \infcatname{C}^{ \cpl } \hookrightarrow \infcatname{C}^{ \Fil }.
				\]
				The explicit formula for the tensor product in $\infcatname{C}^{ \Fil }$ implies that the $n$-th filtration term of the tensor product depends only on previous terms. In particular, $\infcatname{C}^{ \Fil, \geq 0 }$ is stable under such monoidal structure.
				Finally, the internal hom formula is verified by using the formula from \cref{filtered:filtered_internal_hom} and the fact that the right adjoint of $\infcatname{C}^{ \Fil, \geq 0 } \hookrightarrow \infcatname{C}^{ \Fil }$ is simply restriction along $\N^{ \leq } \hookrightarrow \Z^{ \leq }$.
			\end{proof}

			Let us now record two useful adjunctions.

			\begin{prop}
				\label{prop:adjoint_positively_filtered}
				Let $\infcatname{C} \in \Pr^{ L,st }$. There exists a left adjoint of the inclusion $\infcatname{C}^{ \Fil, \geq 0 } \hookrightarrow \infcatname{C}^{\Fil}$, denoted by $\tau^{ \geq 0 }$, and it is given on objects by 
				\[
					(\dots \to F^{ n }X \to F^{ n+1 }X \to \dots) \mapsto (F^{ 0 }X/F^{ -1 }X \to F^{ 1 }X/F^{ -1 }X \to \dots)
				\]
				i.e.\ quotienting\footnote{i.e.\ taking cofibers} by $F^{ -1 }X$ the positive terms of the filtration. The right adjoint is simply restriction along $\N^{ \leq } \hookrightarrow \Z^{ \leq }$.
			\end{prop}
			\begin{proof}
				The last claim is immediate from \cref{prop:positively_filtered_objects_properties}.
				Observe first that $\infcatname{C}^{ \Fil, \geq 0 } \hookrightarrow \infcatname{C}^{ \Fil }$ preserves limits and colimits, since they are computed component-wise in functor categories by \cite[Lemma 5.1.2.1, Corollary 5.1.2.3]{Lurie:HTT}.
				Since both $\infty$-categories are presentable, the adjoint functor theorem tells us that there exists a left adjoint $\tau^{ \geq 0 }\colon \infcatname{C}^{ \Fil } \to \infcatname{C}^{ \Fil, \geq 0 }$.
				The explicit formula follows from the following end computation (see \cref{filtered:filtered_internal_hom}). Let $FX \in \infcatname{C}^{ \Fil }$ and $FY \in \infcatname{C}^{ \Fil, \geq 0 }$:
				\begin{gather*}
					\Map_{ \infcatname{C}^{ \Fil, \geq 0 } }\left( FX/F^{ -1 }X, FY \right) \simeq \int_{ n \in \N } \Map_{ \infcatname{C} }\left( F^{ n }X/F^{ -1 }X, F^{ n }Y \right) \simeq \\
					\simeq \int_{ n \in \N} \left( \Map_{ \infcatname{C} }(F^{ n }X, F^{ n }Y) \times_{ \Map_{ \infcatname{C} }(F^{ -1 }X, F^{ n }Y ) } \{0\} \right) \simeq \\
					\simeq \left( \int_{ n \in \N } \Map_{ \infcatname{C} }(F^{ n }X, F^{ n }Y) \right) \times_{ \left( \int_{ n  \in \N} \Map_{ \infcatname{C} }(F^{ -1 }X, F^{ n }Y) \right) } \{0\} \simeq \\
					\simeq \left( \int_{ n  \in \N} \Map_{ \infcatname{C} }(F^{ n }X, F^{ n }Y) \right) \times_{ \Map_{ \infcatname{C} }(F^{ -1 }X, F^{ 0 }Y) } \{0\} \simeq \\
					\simeq \Map_{ \infcatname{C}^{ \Fil } }(FX, FY)
				\end{gather*}
				where the last equivalence comes from the fact that $FY$ is positively filtered, i.e.\ $F^{ n }Y \simeq 0$ for $n < 0$.
			\end{proof}
			
			\begin{prop}
				\label{prop:left_adjoint_step_0}
				Let $\infcatname{C} \in \Pr^{ L, st }$ and consider the ``$0$-step'' functor 
				\[
					\langle 0, - \rangle\colon \infcatname{C} \to \infcatname{C}^{ \Fil }.
				\]
				It admits a left adjoint whose action on objects is given by 
				\[
					FX \mapsto \cofib(F^{ -1 }X \to (FX)_{ \infty } ).
				\]
			\end{prop}
			\begin{proof}
				Observe that the step functor can be written as 
				\[
					\langle 0, - \rangle \colon \infcatname{C} \stackrel{\const}{\longrightarrow} \infcatname{C}^{ \Fil } \stackrel{\textrm{restr}}{\longrightarrow} \infcatname{C}^{ \Fil, \geq 0 } \hookrightarrow \infcatname{C}^{ \Fil }
				\] 
				so its left adjoint is given by first truncating by $\tau^{ \geq 0 }$ (see \cref{prop:adjoint_positively_filtered}) and then computing colimit (see \cref{filtered:constant_functor}).
				The formula on objects is verified, since
				\[
					\cofib(F^{ -1 }X \to (FX)_{ \infty }) \simeq \left( F^{ 0 }X/F^{ -1 }X \to F^{ 1 }X/F^{ -1 }X \to \dots \to F^{ n }X/F^{ -1 }X \to  \dots \right)_{ \infty }. \qedhere
				\]
			\end{proof}
					
			\begin{remark}
				There might be a notational confusion when talking about connective/coconnective objects for the Beilinson t-structure and negatively/positively filtered objects as defined above. We'll use the notations introduced before for the latter.
			\end{remark}

			Here are some useful lemmas, to characterize graded modules over a certain graded algebra which are concentrated in weight $0$.

			\begin{lemma}
				\label{lemma:graded_modules_weight_0}
				Let $\infcatname{C}$ be as above and $B \in \Alg(\infcatname{C}^{ \gr, \geq 0 })$. We have an adjunction
				\[
					\adjunction{i^{ L }}{\LMod_{ B }(\infcatname{C}^{ \gr })}{\LMod_{ B(0) }(\infcatname{C})}{i},
				\]
				where the right adjoint is fully faithful.
				Moreover, we have the following pullback diagram in $\CAT$
				\begin{diag}
					\LMod_{B(0)}(\infcatname{C}) \ar[d, "\oblv_{B(0)}"] \ar[r, "i", hook] \ar[dr, phantom, very near start, "\lrcorner"] & \LMod_B(\infcatname{C}^{\gr}) \ar[d, "\oblv_B"] \\
					\infcatname{C} \ar[r, "(-)_0", hook] & \infcatname{C}^{\gr}.
				\end{diag}
				The same statement holds if $B \in \Alg(\infcatname{C}^{ \gr, \leq 0 })$.
			\end{lemma}
			\begin{proof}
				Let us first observe that we have an adjunction
				\[
					\adjunction{\ev_{ 0 }}{\infcatname{C}^{ \gr, \geq 0 }}{\infcatname{C}}{(-)_{ 0 }}
				\]
				where the unit is easily seen to be an isomorphism, implying that $(-)_{ 0 }\colon \infcatname{C} \to \infcatname{C}^{ \gr, \geq 0 }$ is fully faithful.
				Both functors are symmetric monoidal, since we are considering \emph{positively} graded objects. By \cite[Corollary 7.3.2.12]{Lurie:HA}, the adjunction lifts at the algebra level
				\[
					\adjunction{\ev_{ 0 }}{\Alg(\infcatname{C}^{ \gr, \geq 0 })}{\Alg(\infcatname{C})}{(-)_{ 0 }}.
				\]
				We obtain a canonical augmentation morphism
				\[
					\phi\colon B \to B(0)_{ 0 }
				\]
				of positively graded algebras (which has a section $B(0)_{ 0 } \to B$).
				Composing with the symmetric monoidal embedding $\infcatname{C}^{ \gr, \geq 0 } \hookrightarrow \infcatname{C}^{ \gr }$, we can then consider the functor
				\[
					i\colon \LMod_{ B(0) }(\infcatname{C}) \stackrel{(-)_{ 0 }}{\longrightarrow} \LMod_{ B(0)_{ 0 } }(\infcatname{C}^{ \gr }) \stackrel{\oblv_{ \phi }}{\longrightarrow} \LMod_{ B }(\infcatname{C}^{ \gr })
				\]
				whose left adjoint is given by the composition
				\[
					i^{ L }\colon \LMod_{ B }(\infcatname{C}^{ \gr }) \stackrel{B(0)_{ 0 } \otimes^{ \gr }_{ B } -}{\longrightarrow} \LMod_{ B(0)_{ 0 } }(\infcatname{C}^{ \gr }) \stackrel{\ev_{ 0 }}{\longrightarrow} \LMod_{ B(0) }(\infcatname{C}).
				\]
				Here we used that $\ev_{ 0 }\colon \infcatname{C}^{ \gr } \to \infcatname{C}$ is lax symmetric monoidal, see \cref{filtered:evaluation_degree}.
				Since evaluation at $0$ preserves limits and colimits (which are computed weight-wise), thinking about the relative tensor product as the geometric realization of the simplicial bar construction, we see that the counit morphism
				\[
					i^{ L }(i(M)) \simeq \ev_{ 0 }\left( B(0)_{ 0 } \otimes^{ \gr }_{ B } \oblv_{ \phi }(M_{ 0 }) \right) \to M
				\]
				is an equivalence for each $M \in \LMod_{ B(0) }(\infcatname{C})$, since $\oblv_{ \phi }$ doesn't change the underlying $\infcatname{C}^{ \gr }$-object.
				This implies that $i$ is fully faithful.
				This proves that the diagram of the statement is well-defined (and evidently commutative). Let us now prove it is a pullback diagram.
				Let $\infcatname{M}$ be the pullback category, so that we have 
				\begin{diag}
					\LMod_{B(0)}(\infcatname{C}) \ar[dr, dashed, red, "H"] \ar[drr, hook, "i"] \ar[ddr, "\oblv_{B(0)}"'] & & \\
																													   & \infcatname{M} \ar[r, hook, "G"] \ar[dr, phantom, very near start, "\lrcorner"] \ar[d, "F"] & \LMod_B(\infcatname{C}^{\gr}) \ar[d, "\oblv_B"] \\
					& \infcatname{C} \ar[r, "(-)_0", hook] & \infcatname{C}^{\gr} 
				\end{diag}
				where the functor $H$ exists simply because the diagram in the statement commutes (easily seen in weight $0$, where both functors are forgetting to $\infcatname{C}$, and in other weights where both compositions are $0$).
				The upper commutative triangle gives us the fully faithfulness of $H$, by $2$-out-of-$3$.
				We need now to prove that $H$ is essentially surjective, to be able to conclude that it is an equivalence. Fix $m \in \infcatname{M}$ and consider $G(m) \in \LMod_{ B }(\infcatname{C}^{ \gr })$; we want to prove that $G(m)$ is contained in the essential image of $i \simeq G \circ H$, so that by fully faithfulness of $G$ we can conclude that $m$ is in the essential image of $H$.
				We know that 
				\[
					\oblv_{ B }(G(m)) \simeq F(m)_{ 0 }
				\]
				and hence the datum of $G(m) \in \LMod_{ B }(\infcatname{C}^{ \gr })$ can be equivalently thought of as an algebra morphism
				\[
					B \to \IntHom^{ \gr }(F(m)_{ 0 }, F(m)_{ 0 })
				\]
				in $\infcatname{C}^{ \gr }$, using the self-enrichment of $\infcatname{C}^{ \gr }$ explained in \cref{filtered:graded_internal_hom}.
				An easy computation of graded hom (see \cref{filtered:graded_internal_hom}) gives
				\[
					\IntHom^{ \gr }(F(m)_{ 0 }, F(m)_{ 0 }) \simeq \IntHom^{\gr, \geq 0 }(F(m)_{ 0 }, F(m)_{ 0 }) \simeq \IntHom(F(m), F(m))_{ 0 }
				\]
				as objects in $\infcatname{C}^{ \gr, \geq 0 }$. To identify the algebra structure on them, we can apply \cref{coroll:adjunction_enriched_hom_easy_version} to the $(\ev_{ 0 }, (-)_{ 0 })$ adjucntion on positively graded objects.
				Hence, the map $B \to \IntHom(F(m), F(m))_{ 0 }$ corresponds by the above adjunction to				 
				\[
					B(0) \to \IntHom(F(m), F(m)) \in \Alg(\infcatname{C}).
				\]
				This is a left $B(0)$-module structure on $F(m)$, which we will denote by $x \in \LMod_{ B(0) }(\infcatname{C})$.
				The situation is depicted by the following commutative diagram in $\Alg(\infcatname{C}^{ \gr })$
				\begin{diag}
					B \ar[r, "G(m)"] \ar[d, "\phi"] & \IntHom^{\gr}(F(m)_0, F(m)_0) \ar[d, "\sim"] \\
					B(0)_0 \ar[r, "x"] & \IntHom(F(m), F(m))_0
				\end{diag}
				which proves that $i(x) \simeq G(m)$ (that is, the left $B$-module structure from $G(m)$ comes from restricting scalars from the left $B(0)$-module $x$, inserted in weight $0$). This implies that $m \simeq H(x)$, proving that $H$ is essentially surjective and therefore an equivalence.\\
				It is clear that the argument goes through also when $B$ is in $\infcatname{C}^{ \gr, \leq 0 }$, replacing positively with negatively filtered objects everywhere.
			\end{proof}
	
			Finally, let us record a similar fact that will be useful later.
			\begin{lemma}
				\label{lemma:filtered_modules_degree_0}
				Let $A \in \Alg(\infcatname{C}^{ \Fil, \geq 0 })$ and let $A_{ \infty } \in \Alg(\infcatname{C})$ be its realization. We have an adjunction 
				\[
					\adjunction{j^{ L }}{\LMod_{ A }(\infcatname{C}^{ \Fil, \geq 0 })}{\LMod_{ A_{ \infty } }(\infcatname{C})}{j}
				\]
				where the right adjoint $j$ is fully faithful.
				Moreover, we have a pullback square in $\CAT$
				\begin{diag}
				\LMod_{A_{\infty}}(\infcatname{C}) \ar[d, "\oblv_{A_{\infty}}"] \ar[r, "j", hook] \ar[dr, phantom, "\lrcorner", very near start] & \LMod_A(\infcatname{C}^{\Fil, \geq 0}) \ar[d, "\oblv_A"] \\
					\infcatname{C} \ar[r, "\langle 0{,} -\rangle", hook] & \infcatname{C}^{\Fil, \geq 0}.
				\end{diag}
			\end{lemma}
			\begin{proof}
				As observed in \cref{prop:left_adjoint_step_0}, we have an adjunction
				\[
					\adjunction{(-)_{ \infty }}{\infcatname{C}^{ \Fil, \geq 0 }}{\infcatname{C}}{\langle 0, - \rangle}.
				\]
				The counit 
				\[
					\left( \langle 0, X \rangle \right)_{ \infty } \to X
				\]
				is easily seen to be an equivalence, implying that the adjunction above is a localization, i.e.\ $\langle 0, -\rangle$ is fully faithful.
				Let us observe that $\langle 0, - \rangle$ is also symmetric monoidal by \cref{filtered:step_functor} (since the filtered tensor product of positively-filtered objects is positively filtered by \cref{prop:positively_filtered_objects_properties}) and also that the realization functor $(-)_{ \infty }$ is symmetric monoidal as well, being a colimit (thanks to the fact that $\infcatname{C}$ is presentably monoidal).
				This implies that, for $A \in \Alg(\infcatname{C}^{ \Fil, \geq 0 })$ we can upgrade those functors to 
				\[
					(-)_{ \infty }\colon \LMod_{ A }(\infcatname{C}^{ \Fil, \geq 0 }) \to \LMod_{ A_{ \infty } }(\infcatname{C}), \qquad \langle 0, - \rangle\colon \LMod_{ A_{ \infty } }(\infcatname{C}) \to \LMod_{ \langle 0, A_{ \infty } \rangle }(\infcatname{C}^{ \Fil, \geq 0 }),
				\]
				which are still adjoint. Using the filtered algebra map $A \to \langle 0, A_{ \infty } \rangle$ (coming from the unit of the first symmetric monoidal adjunction upgraded to algebra-level) we can compose with the forgetful to obtain
				\[
					j \colon \LMod_{ A_{ \infty } }(\infcatname{C}) \stackrel{\langle 0, - \rangle}{\longrightarrow} \LMod_{ \langle 0, A_{ \infty } \rangle }(\infcatname{C}^{ \Fil, \geq 0 }) \stackrel{\oblv}{\longrightarrow} \LMod_{ A }(\infcatname{C}^{ \Fil, \geq 0 })
				\]
				whose left adjoint is then 
				\[
					j^{ L }\colon \LMod_{ A }(\infcatname{C}^{ \Fil, \geq 0 }) \stackrel{\langle 0, A_{ \infty }\rangle \otimes^{ \Fil }_{ A } -}{\longrightarrow} \LMod_{ \langle 0, A_{ \infty } \rangle }(\infcatname{C}^{ \Fil, \geq 0 }) \stackrel{(-)_{ \infty }}{\longrightarrow} \LMod_{ A_{ \infty } }(\infcatname{C})
				\]
				which, once we forget to $\infcatname{C}$, is again the realization $(-)_{ \infty }$, since $A_{ \infty }$ is the realization of both $A$ and $\langle 0, A_{ \infty } \rangle$.
				The fully faithfulness of $j$ is verified by observing that the counit is given by
				\[
					j^{ L }(j(M)) \simeq j^{ L }(\oblv \langle 0, M \rangle) \simeq \left( \langle 0, A_{ \infty } \rangle \otimes^{ \Fil }_{ A } \oblv(\langle 0, M \rangle) \right)_{ \infty } \simeq (\langle 0, A_{ \infty } \rangle)_{ \infty } \otimes_{ A_{ \infty } } (\langle 0, M \rangle)_{ \infty } \simeq M,
				\]
				since the realization is symmetric monoidal and commutes with colimits (hence with relative tensor products) and both $A$ and $\langle 0, A_{ \infty } \rangle$ have $A_{ \infty }$ as underlying object.
				We verified that the diagram of the beginning is well-defined and commutative, so let us now prove it is indeed a pullback. Calling $\infcatname{M}$ the pullback category, observe that we have 
				\begin{diag}
					\LMod_{A_{\infty}}(\infcatname{C}) \ar[dr, dashed, red, "H"] \ar[ddr, "\oblv_{A_{\infty}}"'] \ar[drr, "j", hook] & & \\
					& \infcatname{M} \ar[d, "F"] \ar[r, "G", hook] \ar[dr, phantom, "\lrcorner", very near start] & \LMod_A(\infcatname{C}^{\Fil, \geq 0}) \ar[d, "\oblv_A"] \\
					& \infcatname{C} \ar[r, "\langle 0{,} - \rangle", hook] & \infcatname{C}^{\Fil, \geq 0}.
				\end{diag}
				By $2$-out-of-$3$, $H$ is fully faithful (since $G$ and $\langle 0, - \rangle$ are fully faithful). We just need to prove that it is essentially surjective as well.
				Let $m \in \infcatname{M}$; we know that the underlying filtered object of $G(m)$ is given by $\langle 0, F(m) \rangle$ and its left $A$-module structure is encoded by a morphism
				\[
					A \stackrel{G(m)}{\longrightarrow} \IntHom^{ \Fil, \geq 0 }(\langle 0, F(m) \rangle, \langle 0, F(m) \rangle)
				\]
				in $\Alg(\infcatname{C}^{ \Fil, \geq 0 })$ (see \cref{prop:positively_filtered_objects_properties} for the internal hom formula). 
				Applying \cref{coroll:adjunction_enriched_hom_easy_version} to the adjunction $\adjunction{(-)_{ \infty }}{\infcatname{C}^{ \Fil, \geq 0 }}{\infcatname{C}}{\langle 0, - \rangle}$ (since the left adjoint is symmetric monoidal and the right one is fully faithful) we obtain
				\[
					\IntHom^{ \Fil, \geq 0 }(\langle 0, F(m) \rangle, \langle 0, F(m) \rangle) \simeq \langle 0, \IntHom(F(m), F(m)) \rangle 
				\]
				as algebras in $\infcatname{C}^{ \Fil, \geq 0 }$.
			As the end of the proof of \cref{lemma:graded_modules_weight_0}, we have an induced left $A_{ \infty }$-module structure on $F(m)$, which we call $x \in \LMod_{ A_{ \infty } }(\infcatname{C})$, and a commutative diagram of filtered algebras
				\begin{diag}
					A \ar[r, "G(m)"] \ar[d, "\mathrm{unit}"] & \IntHom^{\Fil, \geq 0}(\langle 0, F(m) \rangle, \langle 0, F(m) \rangle) \ar[d, "\sim"] \\
					\langle 0, A_{\infty} \rangle \ar[r, "j(x)"] & \langle 0, \IntHom(F(m), F(m)) \rangle.
				\end{diag}
				This implies that $G(m) \simeq j(x) \simeq GH(x)$ and, by fully faithfulness of $G$, we may conclude $m \simeq H(x)$, which proves the essential surjectivity of $A$.
			\end{proof}
			
			\begin{remark}
				In \cref{lemma:filtered_modules_degree_0}, the positivity of the filtration is essential: without this assumption, the realization functor $(-)_{ \infty }\colon \infcatname{C}^{ \Fil } \to \infcatname{C}$ is not the left adjoint of $\langle 0, -\rangle$.
			\end{remark}

			\begin{coroll}
				\label{coroll:filtered_complete_modules_gr_0}
				Let $A \in \Alg(\infcatname{C}^{ \Fil, \geq 0 })$ and let $A_{ \infty } \in \Alg(\infcatname{C})$ be its realization. We have a pullback diagram
				\begin{diag}
					\LMod_{A_{\infty}}(\infcatname{C}) \ar[r, hook] \ar[d, "\oblv_{A_{\infty}}"] \ar[dr, phantom, very near start, "\lrcorner"] & \LMod_A(\infcatname{C}^{\cpl}) \ar[d, "\gr"] \\
					\infcatname{C} \ar[r, "(-)_0"] & \infcatname{C}^{\gr}.
				\end{diag}
			\end{coroll}
			\begin{proof}
				Let us start by observing that we have a pullback diagram
				\begin{diag}
					\infcatname{C}^{\Fil, \geq 0} \ar[r, hook] \ar[d, "\gr"] \ar[dr, phantom, very near start, "\lrcorner"] & \infcatname{C}^{\cpl} \ar[d, "\gr"] \\
					\infcatname{C}^{\gr, \geq 0} \ar[r, hook] & \infcatname{C}^{\gr}
				\end{diag}
				where the horizontal arrows are symmetric monoidal inclusions (see \cref{defn:positively_filtered} and \cref{defn:positively_graded}) and the vertical arrows are symmetric monoidal and conservative (see \cref{filtered:properties_gr}). The commutativity of the diagram is evident, let us now prove that it is a pullback diagram.
				For this, it suffices to observe that given any $M \in \infcatname{C}^{ \cpl }$ whose associated graded lies in $\infcatname{C}^{ \gr }$ we can consider $\tau^{ \geq 0 }(M) \in \infcatname{C}^{ \Fil, \geq 0 }$ (defined in \cref{prop:adjoint_positively_filtered}) which is equipped with a natural morphism
				\[
					M \longrightarrow \tau^{ \geq 0 }(M).
				\]
				An easy cofiber computation yields $\gr(M) \simeq \gr(\tau^{\geq }(M))$ (since both live in $\infcatname{C}^{ \gr, \geq 0 }$) and hence we conclude that $M \simeq \tau^{ \geq 0 }(M)$ since $\gr$ is conservative on complete filtered objects (see \cref{filtered:properties_gr}).
				Let us remark that completion is fundamental here: see \cref{remark:possible_confusion_positively_filtered}.
				Observe, similarly, that we also have pullback diagrams
				\begin{diag}
					\infcatname{C} \ar[r, "\langle 0{,}- \rangle", hook] \ar[d, "\sim"] \ar[dr, phantom, very near start, "\lrcorner"] & \infcatname{C}^{\Fil, \geq 0} \ar[d, "\gr"] & & \LMod_A(\infcatname{C}^{\Fil, \geq 0}) \ar[r, hook] \ar[d, "\oblv_A"] \ar[dr, phantom, very near start, "\lrcorner"] & \LMod_A(\infcatname{C}^{\cpl}) \ar[d, "\oblv_A"] \\
					\infcatname{C} \ar[r, "(-)_0", hook] & \infcatname{C}^{\gr, \geq 0} & & \infcatname{C}^{\Fil, \geq 0} \ar[r, hook] & \infcatname{C}^{\cpl}.
				\end{diag}
				Pasting them together with the pullback diagram coming from \cref{lemma:filtered_modules_degree_0}, we conclude.
			\end{proof}

			We already mentioned that the associated graded functor $\gr\colon \infcatname{C}^{ \cpl } \to \infcatname{C}^{ \gr }$ is symmetric monoidal and strong closed. We will need a slightly stronger version of the latter property: consider a complete algebra $A \in \Alg(\infcatname{C}^{ \cpl })$, and the $\infty$-categories $\LMod_{ A }(\infcatname{C}^{ \cpl })$ and $\LMod_{ \gr(A) }(\infcatname{C}^{ \gr })$. The former is (canonically) a module over $\infcatname{C}^{ \cpl }$ and the latter over $\infcatname{C}^{ \gr }$, and, by \cref{lemma:internal_hom_categories} they admit enriched homs.
			The content of the following proposition is that $\gr$ interchanges them.

			\begin{prop}
				\label{prop:gr_closed}
				Let $A \in \Alg(\infcatname{C}^{ \cpl })$ and consider the induced functor 
				\[
					\gr\colon \LMod_{ A }(\infcatname{C}^{ \cpl }) \to \LMod_{ \gr(A) }(\infcatname{C}^{ \gr });
				\]
				it ``respects enriched homs'', meaning we have natural isomorphisms
				\[
					\gr\left(\IntHom^{ \cpl }_{ A }(M, N)\right) \simeq \IntHom^{ \gr }_{ \gr(A) }(\gr(M), \gr(N))
				\]
				in $\infcatname{C}^{ \gr }$, for each $M, N \in \LMod_{ A }(\infcatname{C}^{ \cpl })$.\footnote{Observe that we abused notation here, by denoting $\gr$ both the induced functor on left modules and the original one.}
			\end{prop}
			\begin{proof}
				Let us give a natural map in $\infcatname{C}^{ \gr }$
				\[
						\gr\left(\IntHom^{ \cpl }_{ A }(M, N)\right) \dashrightarrow \IntHom^{ \gr }_{ \gr(A) }(\gr(M), \gr(N)).
				\]
				It is equivalent to giving a map of graded left $\gr(A)$-modules
				\[ 
					\gr\left( \IntHom^{ \cpl }_{ A }(M, N) \right) \odot \gr(M) \dashrightarrow \gr(N)
				\]
				where we denoted by $\odot$ the canonical action of $\infcatname{C}^{ \gr }$ on $\LMod_{ \gr(A) }(\infcatname{C}^{ \gr })$.
				By strong monoidality of $\gr$, which intertwines $\odot$ with the action of $\infcatname{C}^{ \cpl }$ on $\LMod_{ A }(\infcatname{C}^{ \cpl })$, we have
				\[
					\gr(\IntHom_{ A }^{ \cpl }(M, N)) \odot \gr(M) \simeq \gr\left( \IntHom_{ A }^{ \cpl }(M, N) \odot M \right) \to \gr(N)
				\]
				where we used evaluation maps in the last passage.\\
				Let us fix $N \in \LMod_{ A }(\infcatname{C}^{ \cpl }) $ and let us prove that for all $M \in \LMod_{ A }(\infcatname{C}^{ \cpl })$ the map built above is an equivalence of graded objects.
				Since the free-forgetful adjunction
				\[
					\adjunction{A \otimes^{ \cpl } -}{\infcatname{C}^{ \cpl }}{\LMod_{ A }(\infcatname{C}^{ \cpl })}{\oblv_{ A }}
				\]
				is monadic, every $M$ is a geometric realization of free left $A$-modules.
				Since $\gr$ commutes with colimits and  
				\begin{gather*}
					\IntHom_{ A }^{ \cpl }(-, N)\colon \LMod_{ A }(\infcatname{C}^{ \cpl })^{ \op } \to \infcatname{C}^{ \cpl },\quad \IntHom^{ \gr }_{ \gr(A) }(-, \gr(N))\colon \LMod_{ \gr(A) }(\infcatname{C}^{ \gr })^{ \op } \to \infcatname{C}^{ \gr }
				\end{gather*}
				commute with limits, we can assume that $M \simeq A \otimes^{ \cpl } V$ is a free left $A$-module without loss of generality in the first map of the proof.
				Using now that $A \otimes^{ \cpl } -$ is a $\infcatname{C}^{ \cpl }$-linear left adjoint of the forgetful functor, by \cref{lemma:enriched_adjunctions} we can write
				\[
					\gr\left( \IntHom^{ \cpl }_{ A }(A \otimes^{ \cpl } V, N) \right) \simeq \gr\left( \IntHom^{ \cpl }(V, \oblv_{ A }(N)) \right) \simeq \IntHom^{ \gr }(\gr(V), \gr(\oblv_{ A }(N))),
				\]
				where the last passage holds since $\gr$ is strong closed, see \cref{filtered:properties_gr}.
				The following diagram is commutative
				\begin{diag}
					\LMod_A(\infcatname{C}^{\cpl}) \ar[d, "\oblv_A"] \ar[r, "\gr"] & \LMod_{\gr(A)}(\infcatname{C}^{\gr}) \ar[d, "\oblv_{\gr(A)}"] \\
					\infcatname{C}^{\cpl} \ar[r, "\gr"] & \infcatname{C}^{\gr},
				\end{diag}
				and therefore we can continue the above equivalences by
				\begin{gather*}
					\IntHom^{ \gr }(\gr(V), \oblv_{ \gr(A) }(\gr(N))) \simeq \IntHom^{ \gr }_{ \gr(A) }(\gr(A) \otimes^{ \gr } \gr(V), \gr(N)) \\
					\simeq \IntHom^{ \gr }_{ \gr(A) }(\gr(A \otimes^{ \cpl } V), \gr(N)),
				\end{gather*}
				where we used, similarly, that $\gr(A) \otimes^{ \gr } -$ is a $\infcatname{C}^{ \gr }$-linear left adjoint functor.
				This proves the statement for free $M$, and therefore for all $M$.
				Varying $N$, we conclude the proof.
			\end{proof}
			
			\begin{remark}
			    The statement of \cref{prop:gr_closed} also holds if $A \in \CAlg^{ \cpl }$, and the equivalence is an equivalence of graded $\gr(A)$-modules.
			\end{remark}

%% file: filtered_first_definition.tex
\chapter{Derived crystals}
		\label{section:first_definition_derived_DModules}

		In this chapter, we will introduce the first definition of derived D-modules (which are an extension of D-modules which remembers derived information). It is a special case of crystals over the terminal derived foliation, see \cite{TV:AlgFoliations2}.       
		Classical left D-modules on a smooth scheme $X$ can be thought of as quasicoherent sheaves $M$ equipped with a flat connection $\nabla\colon M \to M \otimes_{ \O_{ X } } \Omega^{ 1 }_{ X }$, see \cite[Lemma 1.2.1]{Hotta:GreenBookDMod}.
		As mentioned in the introduction, the interpretation of left D-modules as quasicoherent sheaves equipped with a flat connection allows us to define their \emph{de Rham complex} as
		\[
			M \stackrel{\nabla}{\longrightarrow} M \otimes_{ \O_{ X } } \Omega^{ 1 }_{ X } \stackrel{\nabla}{\longrightarrow} M \otimes_{ \O_{ X } } \Omega^{ 2 }_{ X } \longrightarrow \dots
		\]
		as a dg-module over the de Rham commutative differential graded algebra $(\Omega^{ * }_{ X }, d_{ \dR }, \wedge)$.
		If we want to generalize this interpretation to one which keeps into account the derived structure on $X$, the first step is clearly to replace the cotangent sheaf $\Omega^{ 1 }_{ X }$ with the cotangent complex $\L_{ X }$ (see \cref{subsubsection:algebraic_cotangent_complex}).
		It is then subtle how to generalize the exterior algebra on $\Omega^{ 1 }_{ X }$ in a homotopy-coherent way. 
		Such an object, in the guise of a complete filtered algebra (see \cref{remark:filtered_instead_mixed_graded} for the relation with homotopy coherent cdga's), will be introduced in the next section, after which one definition of derived D-modules will be easy to formulate. 

        \section{The de Rham filtered commutative algebra}
			\label{subsection:deRham_mixed_commutative_algebra}

			The object of interest of this section is the ``commutative de Rham algebra'' of a derived scheme $X$ (although we will mainly focus on derived affine schemes in this text, for simplicity).
			It is more or less clear what to expect: as a graded (dg-)algebra we want $\Sym^{ \gr }_{ \O_{ X } }(\L_{ X }[-1])$, which is the exterior powers of $\L_{ X }$ (where we consider $\L_{ X }[-1]$ in weight $-1$). It is hard, and in fact unfeasible,  to define ``by hand'' the analogue of the de Rham differential\footnote{This is, instead, the standard approach in classical algebraic geometry.}. Another deeper problem is that the standard condition $d^{ 2 } = 0$ will need to be relaxed in a homotopy-coherent way; that is, we will want to have an homotopy between the maps $d^{ 2 }$ and $0$, with additional homotopies making it compatible with all meaningful compositions.
			
			Luckily, there are several equivalent $\infty$-categories (admitting some model categorical presentations) in which such ``homotopy coherent cochain complexes'' live. 
			Our choice will be to use the language of complete filtered objects, see \cref{subsection:filtered_objects}. In the relevant literature, like \cite{TV:AlgFoliations2} and \cite{CPTVV:ShiftedPoisson}, it is common to find mixed graded objects. Let us remind the reader of their definition (for all other properties, see the aforementioned references).

			\begin{defn}[{\cite[§1.1]{TV:AlgFoliations1}}]
				\label{defn:mixed_graded_object}
				Let $\infcatname{C}$ be a presentable stable $\C$-linear $\infty$-category, presented by a (nice enough)\footnote{See \cite[§1.1]{CPTVV:ShiftedPoisson}.} $\Mod_{ \C }^{ \dg }$-model category $\mathcal{C}$.
				A \emph{graded mixed} object of $\mathcal{C}$ is the datum of:
				\begin{enumerate}
					\item a $\Z$-indexed family of objects of $\mathcal{C}$, $(E(n))_{ n \in \Z }$;
					\item for each $n \in \Z$, morphisms 
						\[
							\epsilon_{ n }\colon E(n) \to E(n+1)[-1]
					    \]
						such that $\epsilon_{ n+1 }[-1] \circ \epsilon_{ n } = 0$.
				\end{enumerate}
				A morphism of graded mixed objects $E \to F$ is the datum for each weight $n \in \Z$ of a map $f_{ n }\colon E(n) \to F(n)$ such that $f_{ n+1 }[-1] \circ \epsilon_{ n } = \epsilon_{ n+1 }\circ f_{ n+1 }$. This forms a $1$-category denoted by $\mathcal{C}^{ \egr }$.
				A morphism of graded object is a \emph{weak equivalence} if each $f_{ n }$ is a weak equivalence in $\mathcal{C}$. 
				The $\infty$-category of mixed graded objects in $\infcatname{C}$ is defined as the $\infty$-categorical localization
				\[
					\infcatname{C}^{ \egr } \coloneqq \mathcal{C}^{ \egr }[ \textrm{weak equiv}^{ -1 }].
				\]
			\end{defn}

			\begin{example}
				\label{example:mixed_graded_complex}
				Let us consider $\mathcal{C} = \Mod_{ \C }^{ \dg }$, the $1$-category of chain complexes of $\C$-vector spaces (with quasi-isomorphisms as weak equivalences).
				Then one can picture a mixed graded chain complex $E$ as follows
				\[
				\begin{tikzcd}
					\dots \ar[r] & E(2)^{-2} \ar[r] & E(2)^{-1} \ar[r] & E(2)^0 \ar[r] & E(2)^1 \ar[r] & \dots  \\
					\dots \ar[r] & E(1)^{-2} \ar[r] & E(1)^{-1} \ar[r] \ar[ul, "\epsilon", blue] & E(1)^0 \ar[r] \ar[ul, "\epsilon", blue] & E(1)^1 \ar[r] \ar[ul, "\epsilon", blue] & \dots  \\
					\dots \ar[r] & E(0)^{-2} \ar[r] & E(0)^{-1} \ar[r] \ar[ul, "\epsilon", blue] & E(0)^0 \ar[r] \ar[ul, "\epsilon", blue] & E(0)^1 \ar[r] \ar[ul, "\epsilon", blue] & \dots \\
					\dots \ar[r] & E(-1)^{-2} \ar[r] & E(-1)^{-1} \ar[r] \ar[ul, "\epsilon", blue] & E(-1)^0 \ar[r] \ar[ul, "\epsilon", blue] & E(-1)^1 \ar[r] \ar[ul, "\epsilon", blue] & \dots
				\end{tikzcd}
				\]
				The horizontal black arrows are the ``cohomological'' differentials for each $E(n)$ while the blue diagonal arrows are the mixed differentials. In a way, we can think of $E$ as a certain shifted bicomplex, although weak equivalences of mixed graded objects work differently.
				In fact, aside from shifts, the main difference with bicomplexes and their quasi-isomorphisms (both in the vertical and horizontal direction) is that we consider quasi-isomorphisms only in the horizontal direction; that is, we never consider cohomology in the $\epsilon$-direction. 
			\end{example}

			See the following remark and example to have an idea about the relation between mixed graded objects and complete filtered ones.
				\begin{remark}
					\label{remark:filtered_instead_mixed_graded}

					Although the mentioned reference uses the language of \emph{mixed graded objects}, we will only talk about complete filtered objects.
					There is a symmetric monoidal equivalence between the two: see \cite[Theorem 2.3.1]{Pavia:TStructure} or \cite[Prop 1.3.1]{TV:AlgFoliations2}. Moreover, the same $\infty$-category can also be thought of as \emph{coherent chain complexes}, see \cite{Ariotta:CoherentChainComplexes}.
				\end{remark}

				\begin{example}
					\label{example:derivation_map_mixed_graded}

					Let $A \in \CAlg^{ \cpl }_{ \C } \coloneqq \CAlg(\Mod_{ \C }^{ \cpl })$ be a complete filtered commutative algebra. The mixed graded algebra (or, equivalently, commutative algebra in coherent cochain complexes) we can associate with it consists, on its most basic level, in keeping track of the various (shifted) graded pieces of $A$ together with some additional morphisms called ``mixed differentials'' (and other homotopy coherent data we won't mention here).
					One way to extract such a morphism is to consider the following pushout diagram (compare with \cite[§1.2.2]{Lurie:HA}):
					\begin{diag}
						A^{-2} \ar[r] \ar[d] & A^{-1} \ar[r] \ar[d] & A^0 \ar[d] & \\
						0 \ar[r] & \gr(A)(-1) \ar[d] \ar[r] \ar[ul, very near start, "\ulcorner", phantom]  & ? \ar[ul, phantom, very near start, "\ulcorner"]  \ar[r] \ar[d] & 0 \ar[d] \\
								 & 0 \ar[r] & \gr(A)(0) \ar[r, "\epsilon"] \ar[ul, very near start, phantom, "\ulcorner"] & \gr(A)(-1)[1] \ar[ul, phantom, very near start, "\ulcorner"] .
					\end{diag}
					The morphism $\epsilon\colon \gr(A)(0) \to \gr(A)(-1)[1] = \gr(A)(-1)[2][-1]$ is then a $\C$-linear derivation, i.e.\ it respects the Leibniz rule (to be interpreted in a homotopy-coherent way)\footnote{More precisely: the construction of all the $\epsilon$ maps given here can be promoted to an $\infty$-functor from (complete) filtered objects to graded mixed ones. It can be proved that it gives a symmetric monoidal equivalence and, therefore, also an equivalence between commutative algebras. Then one proves that the algebra structure on $A$ makes those $\epsilon$ respect Leibniz rule.}. 
					Moreover, we can also start to see why those are ``differentials'', meaning why their square is nullhomotopic (such homotopy is part of the data that needs to be given).
					Consider the bigger diagram
					\[
						\begin{tikzcd}[row sep=large, column sep=small]
						A^{-3} \ar[d] \ar[r] & A^{-2} \ar[r] \ar[d] & A^{-1} \ar[r] \ar[d] & A^0 \ar[d] & & \\
						0 \ar[r] & \gr(A)(-2) \ar[ul, phantom, very near start, "\ulcorner"] \ar[r] \ar[d] & ? \ar[ul, phantom, very near start, "\ulcorner"] \ar[r] \ar[d] & ? \ar[ul, phantom, very near start, "\ulcorner"] \ar[r] \ar[d] & 0 \ar[d] & \\
						 & 0 \ar[r] & \gr(A)(-1) \ar[ul, phantom, very near start, "\ulcorner"] \ar[r] \ar[d] & ? \ar[ul, phantom, very near start, "\ulcorner"] \ar[r] \ar[d] & \gr(A)(-2)[1] \ar[ul, phantom, very near start, "\ulcorner"] \ar[r] \ar[d] & 0 \ar[d] \\
						 & & 0 \ar[r] & \gr(A)(0) \ar[ul, phantom, very near start, "\ulcorner"] \ar[r, "\epsilon"] & ? \ar[ul, phantom, very near start, "\ulcorner"] \ar[d] \ar[r, "\epsilon"] & \gr(A)(-1)[1] \ar[ul, phantom, very near start, "\ulcorner"] \ar[d, "\varepsilon"] \\
						 & & & & 0 \ar[r] & \gr(A)(-2)[2] \ar[ul, phantom, very near start, "\ulcorner"]
						\end{tikzcd}
					\]
					where $?$ denotes the corresponding pushout. We then see that the composite 
					\[
						\gr(A)(0) \longrightarrow \gr(A)(-1)[1] \longrightarrow \gr(A)(-2)[2]
					\]
					is canonically nullhomotopic to zero.
					All of this can be made precise by studying a model presentation of mixed graded objects, see \cite[§1]{CPTVV:ShiftedPoisson} (being careful about their different conventions and shifts).
				\end{example}

				\begin{remark}
					\label{remark:conventions_signs_mixed_graded}
					
					There can be different sign and weight conventions on mixed graded objects: we will stick with the ones adopted in \cite{TV:AlgFoliations2}, where everything is cohomological\footnote{By this we mean that all differentials raise by one the degree/weight.} and mixed differentials are like $\epsilon\colon X(n) \to X(n+1)[-1]$. All formulas for (Tate) realizations are then taken from there. In characteristic $0$, anyway, all these sign and grading conventions are equivalent thanks to the symmetric monoidal \emph{shearing} functor.\\
					Some care must also be given when comparing such objects with complete filtered ones. Starting from a complete filtered object $X$ we hinted at how we can get a mixed graded object $X^{ \egr }$. Confusingly enough, in weight $1$ of $X^{ \egr }$, we will find $\gr(X)(-1)[2]$ (this is coherent with the diagram above, recalling how mixed differentials already bring a cohomological shift), and, more generally,
					\[
						X^{ \egr }(n) \simeq \gr(X)(-n)[2n].
					\]
					This means, for example, that mixed graded objects concentrated in positive weights (like the de Rham mixed algebra we will see shortly) will correspond to (complete) filtered objects whose associated graded is concentrated in non-positive weights, so they are constant in $[0, +\infty]$.
					An explicit proof can be found in \cite[Proposition 2.3.9]{Pavia:TStructure}, although he uses a different weight grading convention.
				\end{remark}

				Let us now introduce the main object of this section: the complete filtered de Rham algebra of a derived affine scheme. We will follow \cite[Proposition 1.3.8, 1.3.12]{CPTVV:ShiftedPoisson}, where the language of mixed graded objects is used\footnote{The article just mentioned proves everything in a great generality; for derived affine schemes, the construction is already known, see \cref{remark:DR_illusie_bhatt}.}.
				For the sake of readability we will reproduce their proof translated in the filtered setting.
				\begin{prop}
					\label{prop:deRham_graded_commutative_algebra_structure}
                	The functor $\gr(-)(0)\colon \CAlg^{ \cpl }_{ \C } \to \CAlg_{ \C }$ admits a left adjoint $\DR\colon \CAlg_{ \C } \to \CAlg^{ \cpl }_{ \C }$. Moreover, for each $A \in \CAlg_{ \C }$, there exists a canonical equivalence of graded commutative algebras 
                	\[
                    	\phi_A\colon \Sym^{\gr}_A(\L_A[-1]) \simeq \gr(\DR(A)),
                	\]
					where in the left-hand side $\L_{ A }[-1]$ is of weight $-1$. Moreover, this is natural in $A$. \footnote{Observe how this proposition doesn't require $A$ to be connective, i.e.\ to be the functions of a derived affine scheme. Indeed, one can define $\L_{ A }$ for a general commutative algebra $A$; it is denoted by $\L_{ A }^{ \mathrm{int} }$ in the cited article.}
            	\end{prop}
				\begin{proof}
					\textbf{Existence of the Hodge-completed de Rham algebra}\\
					Observe first that both $\gr$ and $\ev_{ 0 }$ are lax symmetric monoidal so that the functor mentioned above is well-defined. Consider the commuting square (with vertical arrows being the forgetful functors)
					\begin{diag}
						\CAlg^{\cpl}_{\C} \ar[r, "\gr(-)(0)"] \ar[d] & \CAlg_{\C} \ar[d] \\
						\Mod^{\cpl}_{\C} \ar[r, "\gr(-)(0)"] & \Mod_{\C}.
					\end{diag}
					All $\infty$-categories are presentable, and the vertical arrows are simply forgetting the $\E_{ \infty }$-algebra structures (and hence they preserve limits and sifted colimits). By the adjoint functor theorem, it suffices to prove that the top arrow is accessible and preserves limits. This is true if and only if its composition with the vertical right arrow has the same property, but then, thanks to the commutative diagram above, we can first go down left and then use the bottom arrow. Since also the vertical left arrow is accessible and preserves limits, we are reduced to studying the bottom arrow. Since it can be written as the composition
					\[
						\Mod_{ \C }^{ \cpl } \stackrel{\gr}{\longrightarrow} \Mod_{ \C }^{ \gr } \stackrel{\ev_{ 0 }}{\longrightarrow} \Mod_{ \C },
					\]
					we have an explicit left adjoint by \cref{filtered:triv_gr_complete_adjunctions} and \cref{filtered:evaluation_degree}. It is given by 
					\[
						M \mapsto \triv^{ L }(M_{ 0 }) = \left( \dots \longrightarrow 0 \longrightarrow M[-1] \longrightarrow 0 \longrightarrow 0 \longrightarrow \dots \right),
					\]
					where $M[-1]$ is in filtered degree $-1$.
					Therefore we conclude that the left adjoint $\DR$ exists.\\	
					For the second part, let us fix $A \in \CAlg$. By \cref{example:derivation_map_mixed_graded}, the map $\epsilon\colon \gr(\DR(A))(0) \to \gr(\DR(A))(-1)[1]$ respects the Leibniz rule (i.e.\ it is a derivation). Then, by definition of the cotangent complex (relative to $\C$), this corresponds to a morphism
					\[
						\L_{ \gr(\DR(A))(0) } \to \gr(\DR(A))(-1)[1] 
					\]
					in $\Mod_{ \gr(\DR(A))(0) }$.

					Considering the morphism of commutative algebras $A \to \gr(\DR(A))(0)$ (coming from the unit of the above adjunction) we obtain, by general properties of the cotangent complex, a morphism 
					\[ 
						\L_{ A } \otimes_{ A } \gr(\DR(A))(0) \to \L_{ \gr(\DR(A))(0) }
					\]
					of $\gr(\DR(A))(0)$-modules. Precomposing it to the above one and using the tensor-forgetful adjunction, we get a morphism of $A$-modules
					\[
						\L_{ A }[-1] \to \gr(\DR(A))(-1).
					\]

					Using first the ambidextrous adjunction $\adjunction{(-)(-1)}{\Mod_{ A }^{ \gr }}{\Mod_{ A }}{(-)_{ -1 }}$ and then the free-forgetful adjunction for graded commutative algebras, this corresponds to a morphism in graded commutative $A$-algebras 
					\[
						\phi_{ A }\colon \Sym_{ A }^{ \gr }(\L_{ A }[-1]) \to \gr(\DR(A))
					\]
					where $\L_{ A }[-1]$ is seen as a graded $A$-module concentrated in weight $-1$. One verifies that all this construction is functorial in $A$, in the appropriate sense. Now our goal is to prove that $\phi_{ A }$ is an isomorphism in $\CAlg^{ \gr }$ (since the forgetful functor $\CAlg_{ A }^{ \gr } \to \CAlg^{ \gr }$ is conservative).

					Our strategy of proof will be the following: first, prove that $\phi_{ A }$ is an equivalence for free commutative algebras $A$ and then conclude, using that $\phi_{ (-) }$ commutes with sifted colimits and every $A$ can be written as a sifted colimit of free commutative algebras (using the bar construction for the (monadic) free-forgetful adjunction for commutative algebras).

					\textbf{Assertion for $A$ a free commutative algebra}

					Let $A = \Sym_{ \C }(X)$ for $X \in \Mod_{ \C }$, where $\Sym_{ \C }$ is the left adjoint of the forgetful functor $\CAlg_{ \C } \to \Mod_{ \C }$. Then, passing to left adjoints in the initial square, we find that 
					\[
						\DR(\Sym_{ \C }(X)) \simeq \Sym^{ \cpl }_{ \C }(\Free^{ \cpl }(X))
					\]
					where we named $\Free^{ \cpl }$ the explicit left adjoint of $\gr(-)(0)\colon \Mod^{ \cpl } \to \Mod$. The symmetric monoidal structure of $\gr$ gives us the following commutative square 
					\begin{diag}
						\CAlg^{\cpl}_{\C} \ar[d, "\gr"] \ar[r] & \Mod^{\cpl}_{\C} \ar[d, "\gr"] \\
						\CAlg^{\gr}_{\C} \ar[r] & \Mod^{\gr}_{\C}.
					\end{diag}

					We claim that it is left adjointable (i.e.\ it still commutes when passing to left adjoints of the horizontal forgetful functors). This is used also in \cite[§1.2.13]{Pavia:MixedCE}.
					Let us prove this claim: consider the left Beck-Chevalley transformation $\Sym^{ \gr } \circ \gr \to \gr \circ \Sym^{ \cpl }$ in $\CAlg^{ \gr }_{ \C }$; we can check it is an equivalence after forgetting the commutative algebra structure.
					Let us observe that for a complete filtered module $X$ we have 
					\[
						\Sym^{ \gr }_{ \C }(\gr(X) ) \simeq \bigoplus_{ n \in \N } \left(\gr(X)^{ \otimes^{ \gr, n }} \right)_{ \Sigma_{ n } }
					\]
					while the right-hand side (before taking $\gr$) is
					\[
						\Sym^{ \cpl }_{ \C }(X) \simeq \bigoplus_{ n \in \N } \left(X^{ \otimes^{ \cpl, n }}\right)_{ \Sigma_{ n } }.
					\]
					Now we can conclude because $\gr$ commutes both with (completed) tensor products and with taking co-invariants (which in characteristic $0$ need not to be derived), since $G$-co-invariants are the colimit of the corresponding functor from $BG$ and $\gr$ preserves colimits.
						
					This implies that 
					\[
						\gr(\DR(\Sym_{ \C }(X))) \simeq \Sym^{ \gr }_{ \C }(\gr( \Free^{ \cpl }(X) ) )
					\]
					By explicit construction, $\gr(\Free^{ \cpl }(X)) = X_{ 0 } \oplus X[-1]_{ -1 }$, where subscripts denote the weights, and since $\Sym^{ \gr }$ commutes with colimits (being a left adjoint) we obtain that 
					\[
						\Sym_{ \C }^{ \gr }(X_{ 0} \oplus X[-1]_{ -1 }) \simeq \Sym_{ \C }^{ \gr }(X_{ 0 }) \otimes^{ \gr }_{ \C } \Sym_{ \C }^{ \gr }(X[-1]_{ -1 })
					\]
					where we used that the coproduct in $\CAlg^{ \gr }_{ \C }$ is given by the standard tensor product, by \cite[Prop 3.2.4.7]{Lurie:HA}. Since $\Sym_{ \C }^{ \gr }(X_{ 0 }) \simeq \Sym_{ \C }(X)_{ 0 }$ (explicit computations using the formulas above), by formal properties of Sym\footnote{Again this comes from writing a commutative diagram using the symmetric monoidal structure of base change $- \otimes_{ \C }^{ \gr } \Sym_{ \C }(X)_{ 0 }$ and then proving that it is left adjointable, using that $\Sym^{ \gr }_{ \C }(-)$ is left adjoint to the forgetful functor $\CAlg_{ \C }^{ \gr } \to \Mod_{ \C }^{ \gr }$. Being in characteristic $0$ seems fundamental to such computations.}  we obtain 
					\[
						\Sym_{ \C }^{ \gr }(X_{ 0 }) \otimes^{ \gr }_{ \C } \Sym_{ \C }^{ \gr }(X[-1]_{ -1 }) \simeq \Sym^{ \gr }_{ \Sym(X) }(\Sym_{ \C }(X)_{ 0 } \otimes^{ \gr }_{ \C } X[-1]_{ -1 }).
					\]
					It finally remains to observe that the cotangent complex of a free commutative algebra is computed by $\L_{ \Sym_{ \C }(X) } \simeq \Sym_{ \C }(X) \otimes_{ \C } X \in \Mod_{ \Sym(X) }$ (to verify it, it suffices to prove that the right-hand side satisfies the correct universal property) and hence we conclude that \[
						\gr(\DR(\Sym(X))) \simeq \Sym_{ \Sym(X) }^{ \gr }(\L_{ \Sym(X) }[-1]),
					\]
					where $\L_{ \Sym(X) }[-1]$ is in weight $-1$.

					We proved our assumption for all free commutative algebras.

					\textbf{General case}

					To finally conclude the proof, it suffices to show that all functors involved commute with sifted colimits, since each commutative algebra is a sifted colimit of free ones (by the bar resolution). The functors $\DR$ and $\gr$, being left adjoints, commute with all colimits, so the right-hand side is okay. For the left-hand side we need to prove that the functor 
					\[
						\CAlg_{ \C } \ni A \mapsto \Sym_{ A }^{ \gr }(\L_{ A }[-1]) \in \CAlg^{ \gr }_{ \C }
					\]
					commutes with (sifted) colimits. One way to see this is to notice that for any $B \in \CAlg^{ \gr }_{ \C }$, by precomposing with the canonical algebra map $A_{ 0 } \to \Sym_{ A }^{ \gr }(\L_{ A }[-1])$, we have a natural map 
					\[
						\Map_{ \CAlg^{ \gr } }(\Sym_{ A }^{ \gr }(\L_{ A }[-1]), B) \to \Map_{ \CAlg }(A, B(0))
					\]	
					where the fiber at any map $A \to B(0)$ is equivalent to $\Map_{ \Mod_{ A } }(\L_{ A }, B(-1)[1])$. This latter space, by definition of cotangent complex, is equivalent to $\Map_{ \CAlg_{ /B(0) } }(A, B(0) \oplus B(-1)[1])$ (where we implicitly used that $A \oplus B(-1)[1] \simeq A 	\times_{ B(0) } (B(0) \oplus B(-1)[1])$, as split square-zero extensions). Now this space (corresponding to each fiber) sends colimits in $A$ to limits, as well as the target space $\Map_{ \CAlg }(A, B(0))$, which implies that also $\Map_{ \CAlg^{ \gr } }(\Sym_{ A }^{ \gr }(\L_{ A }[-1]), B)$ does so. By Yoneda, we finally conclude that our initial functor commutes with colimits in $A$. The construction of the arrow $\phi$ does the same as well, and hence we conclude.
				\end{proof}

				This proposition justifies the name ``de Rham complex of $A$'', especially once we represent it as a mixed graded object or coherent chain complex: we got the object we expected, that is, $\Sym^{ \gr }_{ A }(\L_{ A }[-1]_{ -1 })$, and the mixed differential (whose existence is only ``weak'', see \cite[Remark 1.3.15]{CPTVV:ShiftedPoisson}) represents the analogue of the de Rham differential. This approach is different from the classical one, where we can build the algebraic de Rham complex by hand on affine spaces. The commutative algebra structure corresponds to the classical wedge product of algebraic differential forms.
				More precisely we should think of $\DR(A)$ as the completed Hodge-filtered de Rham algebra of $A$. The Hodge filtration is taken into account by the $\G_{ m }$-action, i.e.\ the weight grading.
            	We can then consider modules over $\DR(A)$ in the symmetric monoidal $\infty$-category of $\Mod^{\cpl}$, and this is how we're going to generalize D-modules in the following sections.
				
				\begin{remark}
					\label{remark:DR_arbirtrary base}
					The same proof of \cref{prop:deRham_graded_commutative_algebra_structure} works on any base $A \in \CAlg_{ \C }$, giving us a \emph{relative de Rham} filtered algebra functor, denoted by
					\[
						\DR(-/A)\colon \CAlg_{ A } \longrightarrow \CAlg_{ A }^{ \cpl }
					\]
					and with the property that
					\[
						\gr(\DR(B/A)) \simeq \Sym_{ B }^{ \gr }(\L_{ B/A }[-1]),
					\]
					where $\L_{ B/A }[-1]$ is in weight $-1$.
				\end{remark}

				\begin{remark}
					\label{remark:DR_illusie_bhatt}
					The construction of the ``derived de Rham algebra'' of a map $A \to B$ goes back to Illusie, who proved that derived de Rham cohomology coincides with Hartshorne's algebraic de Rham cohomology for lci morphisms (see \cite{Illusie:ComplexeCotangent1}).	Bhatt then extended this result to more general (finite type and Noetherian) morphisms of $\Q$-schemes, see \cite{Bhatt:deRhamDerivedCohomologyCompletions}. 
					There, the (Hodge-completed) derived de Rham algebra of a morphism $A \to B$ is denoted by $\widehat{\mathrm{dR}}_{ B/A }$ and is exactly $\DR(B/A)$ endowed with the ``quotient filtration'' (see \cref{remark:pass_left_right_complete_filtered_object}).
				\end{remark}
        
		\section{Constant modules over a complete filtered algebra}
			\label{subsection:crystals_derived_foliation}

			We continue our recall of \cite{TV:AlgFoliations2} and introduce the definition and main properties of constant modules over a complete filtered commutative algebra\footnote{They only consider the case of a derived foliation, where they call ``quasicoherent crystals'' what we call constant modules, but the definition can be given in full generality.}.

			\begin{defn}
				\label{defn:constant_graded_mixed_modules}
				Let $B \in \CAlg^{ \cpl }_{ \C }$, a complete filtered commutative $\C$-algebra. The $\infty$-category of \emph{constant modules over $B$} is defined as the following pullback
				\begin{center}
					\begin{tikzcd}[row sep=huge, column sep = huge]
						\Mod^{\const}_{B} \ar[r, hook] \ar[d] \ar[dr, phantom, "\lrcorner", very near start] & \Mod_{B}(\Mod^{\cpl}_{\C}) \ar[d, "\gr"] \\
						\Mod_{\gr(B)(0)} \ar[r, "\gr(B) \otimes^{\gr}_{\gr(B)(0)} (-)_0", hook] & \Mod_{\gr(B)}(\Mod^{\gr}_{\C})
					\end{tikzcd}
				\end{center}
				in $\CAT$. 
			\end{defn}

			\begin{remark}
				\label{remark:constant_modules_fully_faithful}
				Let us give a proof that the bottom arrow in the diagram of \cref{defn:constant_graded_mixed_modules} is fully faithful. Let $A \coloneqq \gr(B)(0) \in \CAlg_{ \C }$. For $M, N \in \Mod_{ A }$ we have 
				\begin{gather*}
					\Map_{ \Mod^{ \gr }_{ \gr(B) } }\left( M_{ 0 } \otimes^{ \gr }_{ A_{ 0 } }  \gr(B), N_{ 0 } \otimes^{ \gr }_{ A_{ 0 } }  \gr(B) \right) \\
					\simeq \Map_{ \Mod_{ A }^{ \gr } }\left(M_{ 0 }, N_{ 0 } \otimes^{ \gr }_{ A_{ 0 } } \gr(B) \right) 	\simeq \Map_{ \Mod_{ A }}\left(M, N\right),
				\end{gather*}
				where we used the tensor-forgetful adjunction, followed by the fact that evaluation in weight $0$ is right adjoint to insertion in weight $0$ and finally that $(N_{ 0 } \otimes_{ A }^{ \gr } \gr(B))(0) \simeq N$.
			\end{remark}
			
			The next proposition shows that, if $\gr(B)$ respects certain perfection assumptions, then constant modules over $B$ are both a localization and a colocalization of all (complete filtered) $B$-modules.
			\begin{prop}
				\label{prop:constant_modules_Pr_L_R}
				Let $B \in \CAlg_{ \C }^{ \cpl }$ and consider the cartesian diagram in \cref{defn:constant_graded_mixed_modules}. It is a pullback diagram in $\Pr^{ L }$.
				Moreover, if each $\gr(B)(p)$ is a perfect $\gr(B)(0)$-module, it is also a pullback diagram in $\Pr^{ R }$.
			\end{prop}
			\begin{proof}
				Let us first recall that the inclusions $\Pr^{ L } \to \CAT$ and $\Pr^{ R } \to \CAT$ both preserve and create limits by \cite[Proposition 5.5.3.13, Theorem 5.5.3.18]{Lurie:HTT}.
				This means that it suffices to prove that the functors involved are left (and right, if $\gr(B) \in \Perf_{ \gr(B)(0) }^{ \gr }$) adjoints. \\
				We know that $\gr\colon \Mod_{ \C }^{ \cpl } \to \Mod^{ \gr }_{ \C }$ has both a left and a right adjoint (see \cref{filtered:left_triv} and \cref{filtered:right_triv}); moreover, since $\Mod_{ A }(\infcatname{C}) \to \infcatname{C}$ preserves all limits and colimits (for example for $\infcatname{C} \in \Pr^{ L }$, by \cite[§4.2.3]{Lurie:HA}), we deduce that 
				\[
					\gr\colon \Mod_{ B }(\Mod^{ \cpl }_{ \C }) \longrightarrow \Mod_{ \gr(B) }(\Mod_{ \C }^{ \gr })
				\]
				commutes with both limits and colimits and therefore, by the adjoint functor theorem \cite[Corollary 5.5.2.9]{Lurie:HTT}, admits both a left and a right adjoint. On the other hand, the bottom functor, being a tensor product, commutes with colimits and so is a left adjoint.
				We conclude that the diagram is also a pullback in $\Pr^{ L }$. This implies that $\Mod^{ \const }_{ B } \hookrightarrow \Mod^{ \cpl }_{ B }$ is a \emph{colocalization}.\\
				Suppose now that for each $n \in \Z$, $\gr(B)(n) \in \Perf_{ \gr(B)(0) }$; then $\gr(B) \otimes^{ \gr }_{ \gr(B)(0) } (-)_{ 0 }$ preserves limits: it suffices to observe that limits are computed weight-wise (since graded objects are a functor category, see \cite[Lemma 5.1.2.1]{Lurie:HTT}) and in weight $p$ it corresponds to tensoring with $\gr(B)(p)$, which is perfect by assumption, and hence dualizable. The analogue of the reasoning goes through to obtain that this same diagram is a pullback in $\Pr^{ R }$.
				This implies that $\Mod^{ \const }_{ B } \hookrightarrow \Mod^{ \cpl }_{ B }$ is a \emph{localization} as well.
			\end{proof}
		
			\begin{remark}
				\label{remark:constant_modules_monoidal}
				Let us observe that all $\infty$-categories in the pullback diagram of \cref{defn:constant_graded_mixed_modules} are (presentably) symmetric monoidal and stable, and all functors are symmetric monoidal themselves (and also exact and left adjoints). That is, it is a diagram in $\CAlg(\Pr^{ L, st })$ and a pullback there as well (since forgetting to $\CAT$ commutes with limits by \cite[Corollary 3.2.2.4]{Lurie:HA}). Therefore we deduce that $\Mod^{ \const }_{ B }$ is a symmetric monoidal subcategory of $\Mod^{ \cpl }_{ B }$ (equipped with the obvious monoidal structure).
			\end{remark}

			We can now define derived D-modules on $\Spec(A)$ as constant modules over $\DR(A)$, see \cite[§2.2]{TV:AlgFoliations2} and \cref{prop:deRham_graded_commutative_algebra_structure}. This generalizes the ``de Rham complex'' of a (left) D-module. We will refer to them as ``derived left crystals''.

			\begin{defn}
				\label{defn:toen_vezzosi_derived_d_modules}
				Let $X = \Spec A$ be a derived affine scheme. The symmetric monoidal $\infty$-category of \emph{derived left crystals over $X$} is $\Mod^{ \const }_{ \DR(A) }$, with the standard monoidal structure (tensor product over $\DR(A)$).
			\end{defn}
			The following proposition will be useful later on.

			\begin{prop}
				\label{prop:constant_modules_forgetful_conservative}
				Let $B \in \CAlg^{ \cpl }_{ \C }$ be such that for each $n \in \Z$, $\gr(B)(n)$ is a perfect $\gr(B)(0)$-module.
				Then the associated graded in weight $0$ functor
				\[
					\gr(-)(0)\colon \Mod^{ \const }_{ B } \to \Mod_{ \gr(B)(0) }, \qquad M \mapsto \gr(M)(0)
				\]
				is a monadic right adjoint.
			\end{prop}
			\begin{proof}
				It is a right adjoint thanks to \cref{prop:constant_modules_Pr_L_R}, so we just need to prove it is monadic. It suffices to prove that it preserves colimits and is conservative, so then we conclude by applying Barr-Beck-Lurie, \cite[Theorem 4.7.3.5]{Lurie:HA}. Preservation of colimits comes again from \cref{prop:constant_modules_Pr_L_R}, so we just need to prove it is conservative.
				Let us now consider the commutative diagram
				\begin{diag}
					\Mod_B(\Mod^{\cpl}_{\C}) \ar[d, "\gr"] \ar[r, "\oblv_B"] & \Mod^{\cpl}_{\C} \ar[d, "\gr"] \\
					\Mod_{\gr(B)}(\Mod^{\gr}_{\C}) \ar[r, "\oblv_{\gr(B)}"] & \Mod^{\gr}_{\C}
				\end{diag}
				and observe that the rightmost $\gr$ is conservative by \cref{filtered:properties_gr}, since we are working with \emph{complete} filtered modules. Forgetting the module structure (both for $B$ and $\gr(B)$) is conservative as well, see \cite[Corollary 4.2.3.2]{Lurie:HA}. This implies that the leftmost $\gr$ is conservative too. 
				We conclude that our initial functor $\gr(-)(0)$ is conservative since it arises as the pullback of a conservative functor along a fully faithful one, $\gr(B) \otimes_{ \gr(B)(0) }^{ \gr } (-)_{ 0 }$ (see \cref{remark:constant_modules_fully_faithful}).
			\end{proof}

			\begin{remark}
				\label{remark:question_left_adjoint_gr_0}
				It is then a legit question to wonder what is the left adjoint of $\gr(-)(0)\colon \Mod^{ \const }_{ B } \to \Mod_{ \gr(B)(0) }$. We will give an answer at the very end of this manuscript, see \cref{coroll:left_adjoint_gr_0}, for those $B$'s arising as filtered Chevalley-Eilenberg cohomology complexes of perfect Lie algebroids.
			\end{remark}

			\begin{coroll}
				\label{coroll:constant_modules_de_rham_conservative}
				Let $A \in \CAlg_{ \C }^{ \cn }$ be locally finitely presented, as in \cref{defn:locally_finite_presentation}. Then the functor
				\[
					\gr(-)(0)\colon \Mod_{ \DR(A) }^{ \const } \to \Mod_{ A }
				\]
				is a monadic right adjoint.
			\end{coroll}
			\begin{proof}
				Let us consider the complete filtered commutative de Rham algebra $\DR(A)$.
				By \cref{prop:deRham_graded_commutative_algebra_structure}, we know that 
				\[
					\gr(\DR(A)) \simeq \Sym^{ \gr }_{ A }(\L_{ A }[-1]),
				\]
				with $\L_{ A }[-1]$ in weight $-1$, so that $\gr(\DR(A))(0) \simeq A$. Moreover, since $A$ is locally finitely presented, $\L_{ A }$ is a perfect $A$-module by \cref{thm:locally_finite_presentation_perfect_cotangent_complex}. We deduce that in each weight $p$, $\Sym^{ p }_{ A }(\L_{ A }[-1])$ is a perfect $A$-module too, being the retract of the $p$-th tensor power of $\L_{ A }$ via the symmetrization map
				\[
					\Sym^{ p }_{ A }(\L_{ A }[-1]) \to (\L_{ A }[-1])^{ \otimes_{ A }^{ p } },\qquad x_{ 1 } \otimes^{ s } x_{ 2 } \otimes^{ s } \dots \otimes^{ s } x_{ p } \mapsto \frac{1}{p!} \sum_{ \sigma } x_{ \sigma(1) } \otimes x_{ \sigma(2) } \otimes \dots \otimes x_{ \sigma(p) },
				\]
				where $\sigma \in \Sigma_{ n }$, described, for example, in the proof of \cite[Theorem 2.2]{Hennion:TangentLieAlg}. Notice how characteristic $0$ of the base field is essential, since we are relying on cdga's as a model for commutative algebras. We conclude applying \cref{prop:constant_modules_forgetful_conservative}.
			\end{proof}

%% file: lie_algebroids.tex
\section{Dg-Lie Algebroids}
			\label{subsection:dg_lie_algebroids}
			The goal of this section is to review the main definitions and properties of dg-Lie algebroids from \cite{Nuiten:KoszulLieAlgbd}.
			We emphasize that the adjective ``dg'' signals that we will be working with $1$-categories, unless otherwise specified.

			\subsection{Definition}
				\label{subsubsection:definition_dg_lie_algebroids}

				As the name suggests, (dg)-Lie algebroids are a generalization of (dg)-Lie algebras, where there are two different bases that interact with each other. In fact, to be precise, one talks about dg-Lie algebroids over a base $A \in \CAlg$.
				A lot of results in dg-Lie algebras, for example \cite[Chapter 2]{Lurie:DAGX}, also hold in this context, up to slight modifications, thanks to \cite{Nuiten:KoszulLieAlgbd}.\\
				Let us fix a base $A \in \CAlg^{ \cn }_{ \C }$, by which we will also denote a cofibrant cdga representing it (recall the rectification result $\infcatname{cdga}_{ \C }^{ \cn }[\textrm{qis}^{ -1 }] \simeq \CAlg_{ \C }^{ \cn }$, see \cite[Theorem 4.5.4.7]{Lurie:HA}). As already mentioned in \cref{section:second_definition_derived_DModules}, its tangent complex $ \T_{ A } $ can be interpreted as $\C$-linear (derived) derivations of $A$ (hence it acts on $A$) and it has the standard commutator bracket making it a dg-Lie algebra over $\C$. 

				\begin{defn}
					A dg-Lie algebroid over $A$ is the datum of $E \in \Mod_{ A }^{ \dg }$ equipped with a $\C$-linear dg-Lie algebra structure and an \emph{anchor map} $\rho\colon E \to \T_{ A }$ which is both a map of dg-$A$-modules and of dg-Lie algebras over $\C$. Moreover, the following must hold:
					\[
						[e, a \cdot e'] = (-1)^{ |e||a| } a \cdot [e, e'] + \rho(e)(a) \cdot e'
					\]
					for all homogeneous $e, e' \in E$ and $a \in A$.
					A morphism of dg-Lie algebroids is a map $E \to E'$ preserving all this structure. They form a $1$-category denoted by $ \LieAlgbd^{ \dg }_{ A }$.
				\end{defn}

				Two easy examples of dg-Lie algebroids are the initial one $0 \to \T_{ A }$ and the final one, i.e.\ the tangent complex itself $\id\colon \T_{ A } \to \T_{ A }$. Moreover all dg-Lie algebras over $A$ can be seen as dg-Lie algebroids with zero anchor map.

				As homotopically-minded people, whenever we work with differential graded objects, what we really want to understand is the $\infty$-category obtained by localizing along quasi-isomorphisms (by which we mean morphisms of dg-Lie algebroids inducing quasi-isomorphisms on the underlying dg-modules), i.e.\ \[
					\LieAlgbd_{ A } \coloneqq \LieAlgbd^{ \dg }_{ A }\left[ \textrm{qis}^{ -1 } \right].
				\]

				\begin{remark}
					For a non cofibrant cdga $A'$ we can define $\LieAlgbd_{ A' } \coloneqq \LieAlgbd_{ Q(A') }$, for a cofibrant replacement $Q(A')$. In fact, derivations of $Q(A)$ are the correct object to consider for the tangent complex, see \cite[Remark 2.2]{Vezzosi:LieAlgebroidsModel}. 
				\end{remark}
			
				\begin{remark}
					\label{remark:forget_differential}
					We gave the definition in the \emph{differential} graded setting because this allows us to do homotopy theory by localizing along quasi-isomorphisms. We could have written an equal definition simply in the \emph{graded} setting, ignoring differentials. 
					In fact, all these constructions can simply be performed at the graded level first and then adding (or checking/requiring compatibility) with the differentials.
					In other words, the $1$-functor that forgets differentials
					\[
						\oblv_{ \partial }\colon \Mod_{\C}^{ \dg } \to \left( \Mod_{ \C }^{ \heartsuit}\right)^{  \gr }
					\]
					is conservative, symmetric monoidal, closed and preserves both limits and colimits (where everything is meant $1$-categorically).
				\end{remark}			
				The main result of \cite{Nuiten:HomotopicalLieAlgbd} is the existence of a semi-model structure\footnote{ It is not a model structure because there is no functorial fibrant replacement functor.} on $\LieAlgbd^{ \dg }$ with quasi-isomorphisms as weak equivalences and degree-wise surjections as fibrations. Cofibrant dg-Lie algebroids and their properties are studied there.\\
				Let us introduce some terminology for later.
				\begin{defn}
					\label{defn:free_lie_algebroid_functor}
					The \emph{free Lie algebroid over $A$} is the $\infty$-functor
					\[
						\FreeLieAlgbd\colon (\Mod_{ A })_{ -/\T_{ A } } \longrightarrow \LieAlgbd_{ A }
					\]
					defined as the left adjoint of the forgetful $\oblv\colon \LieAlgbd_{ A } \to (\Mod_{ A })_{ -/\T_{ A } }$.
				\end{defn}
				The existence of the free Lie algebroid functor comes from \cite[Theorem 3.3]{Nuiten:HomotopicalLieAlgbd}. The same result tells us that the adjunction
				\[
					\adjunction{\FreeLieAlgbd}{(\Mod_{ A })_{ -/\T_{ A } }}{\LieAlgbd_{ A }}{\oblv}
				\]
				is monadic.\\
				Another monadic adjunction relates Lie algebroids over $A$ with Lie algebras in $\Mod_{ A }$.

				\begin{prop}[{\cite[Proposition 3.5]{Nuiten:HomotopicalLieAlgbd}}]
					\label{prop:lie_algebras_and_lie_algebroids}
					Let $A \in \CAlg^{ \cn }_{ \C }$. There is a monadic $\infty$-categorical adjunction
					\[
						\adjunction{\mathrm{diag}}{\infcatname{LieAlg}_{ A }}{\LieAlgbd_{ A }}{\ker}
					\]
					where the composition
					\[
						\LieAlgbd_{ A } \stackrel{\ker}{\longrightarrow} \infcatname{LieAlg}_{ A } \stackrel{\oblv_{ \Lie }}{\longrightarrow} \Mod_{ A }
					\]
					sends $\g \to \T_{ A }$ to $\fib(\g \to \T_{ A }) \in \Mod_{ A }$.
				\end{prop}

				Let us remark that, although the $1$-categorical version of $\mathrm{diag}$ is fully faithful, the $\infty$-functor $\mathrm{diag}$ is not. That is, Lie algebroids over $A$ can be thought as $A$-linear Lie algebras with additional \emph{structure}, see \cite[Remark 3.7]{Nuiten:HomotopicalLieAlgbd}

			\subsection{Universal enveloping algebra}
				\label{subsubsection:universal_enveloping_algebra}
								
				Similarly to the case of dg-Lie algebras, we can define an associative dg-algebra $U(\g)$ such that being a $\g$-representation is equivalent to being a \emph{left} dg-$U(\g)$-module. The first construction of such an object, for classical Lie algebroids, comes from \cite[§2]{Rinehart:LieAlgbd} (where the term $(\C, A)$-Lie algebras is used) but it is easily adapted to the dg context. Let us rapidly recall its ($1$-categorical) construction:
				\begin{enumerate}
					\item Consider $A \oplus \g \in \Mod_{ \C }^{ \dg }$ equipped with the Lie bracket given on homogeneous elements by \[
						[a + g, a' + g'] \coloneqq \left(\rho(g)(a') - (-1)^{ |g'||a|} \rho(g')(a)\right) + [g, g']
					\]
					which makes it into a dg-Lie algebra over $\C$.
					\item Form the universal enveloping dg-algebra (defined as the usual quotient of the tensor dg-algebra) $U(A \oplus \g)$ of the aforementioned dg-Lie algebra over $\C$. 
					\item Consider the non-unital dg-algebra $U(A \oplus \g)^{ + } \coloneqq U(A\oplus \g)/\C$ and its two-sided ideal $P$ generated by \[
						\overline{a} \cdot \overline{(a' + g')} - \overline{(a\cdot a' + a\cdot g')}
					\]
					where $\overline{(-)}\colon A \oplus \g \to U(A \oplus \g)^{ + }$.
					\item Finally, define \[
						U(\g) \coloneqq U(A \oplus \g)^{ + }/P \in \Alg^{ \dg }_{ \C }.
					\]
				\end{enumerate}

				Let us immediately make some observations.
				\begin{itemize}
					\item This construction defines a $1$-functor
						\[
							U(-)\colon \LieAlgbd_{ A }^{ \dg } \to \Alg_{ \C }^{ \dg }.
						\]
					\item There is an inclusion of dg-algebras $A \hookrightarrow U(\g)$, obtained by composing 
						\[
							A \hookrightarrow A \oplus \g \hookrightarrow U(A \oplus \g)^{ + } \to U(\g).
								\]
						It is not central, though. This means $U(\g)$ is an algebra in the monoidal category of dg-$(A, A)$-bimodules (i.e.\ $A$ can act from the left and from the right, and those two actions are different). We can therefore upgrade the $1$-functor to 
						\[
							U(-)\colon \LieAlgbd_{ A }^{ \dg } \to \Alg^{ \dg }( (A, A)-\BiMod^{ \dg} ),
						\]
						where the monoidal structure on dg-$(A,A)$-bimodules is the obvious $- \otimes_{ A } -$ which uses the right $A$-action on the left term and the left $A$-action on the right (and the inherited $(A,A)$-bimodule structure is using the left action on the left term and the right action on the right term).
					\item There is an inclusion of left dg-$A$-modules\footnote{Remark that forgetting the right action $(A,A)-\BiMod \to \LMod_{ A }$ is not monoidal, so $U(\g)$ is not an algebra in left $A$-modules.} $\g \hookrightarrow U(\g)$, obtained in a similar way as above. It is not a morphism of dg-$(A,A)$-bimodules, though, since \[
					 	ag \neq ga
					 \]
					 in $U(\g)$ (whereas it is true in $\g$, seen as diagonal dg-$(A,A)$-bimodule).
	
				 \item There exists an increasing filtration, called the \emph{PBW filtration}, on the dg-$(A, A)$-bimodule $U(\g)$. The $n$-th step $U(\g)^{ \leq n }$ is the submodule spanned by ``words'' with at most $n$ elements of $\g$. The filtration is exhaustive because $\Sym_{ \C }(A \oplus \g) \to U(\g)$ is surjective. That is, we actually have a positively filtered dg-$(A,A)$-bimodule
					 \[
					 	U(\g)^{ \PBW } \coloneqq \left( 0 \hookrightarrow U(\g)^{ \leq 0 } \simeq A \hookrightarrow U(\g)^{ \leq 1 } \hookrightarrow \dots \hookrightarrow U(\g)^{ \leq n } \hookrightarrow \dots \right)
					 \]
					 whose underlying object (i.e.\ colimit) is $U(\g)$.
			     
			     \item It is easily seen that the filtration is functorial in $\g$ and $U(\g)^{ \PBW }$ is actually a filtered associative algebra (that is, the product sends $U(\g)^{ \leq n } \otimes_{ A } U(\g)^{ \leq m }$ to $U(\g)^{ \leq n+m }$). The $1$-functor becomes
					 \[
						 U(-)^{ \PBW }\colon \LieAlgbd_{ A }^{ \dg } \longrightarrow \Alg^{ \dg }\left( (A,A)-\BiMod^{ \dg, \Fil } \right).
					 \]
				 
			     \item We can consider the associated graded $1$-functor (defined with the same formula as \cref{filtered:associated_graded} but quotients are $1$-quotients). It is again monoidal, so that $\gr(U(\g)^{ \PBW })$ is an associative algebra in graded dg-$(A,A)$-bimodules. Observe that, like in the classical case, for $x \in U(\g)^{ \leq p }$ and $a \in A$ we have 
			    \[
			     	ax - xa \in U(\g)^{ \leq p-1 };
				\]
				this implies that $\gr(U(\g)^{ \PBW })$ is an actual (graded) $A$-algebra, i.e.\ $A$ is contained in its center. That is, $\gr(U(\g)^{ PBW }) \in \Alg^{ \dg }(\Mod_{ A }^{ \dg, \gr })$.
				\end{itemize}

				We will work with graded objects, so let us recall the notation introduced in \cref{filtered:evaluation_degree}. In particular we will often use the functor
				\[
					(-)_{ 1 }\colon \infcatname{C} \to \infcatname{C}^{ \gr }
				\]
				that sends $X \in \infcatname{C}$ to the graded object of $\infcatname{C}$ which is $X$ in weight $1$ and $0$ elsewhere.

				\begin{remark}
					This construction works verbatim in the \emph{graded} setting, that is, where $A$ and $\g$ do not have differentials. Moreover, calling $\oblv_{ \partial }$ the symmetric monoidal $1$-functor that forgets differentials, we have 
					\[
						\oblv_{ \partial }( U(\g) ) \simeq U( \oblv_{ \partial }(\g) )
					\]
					as filtered (graded) algebras. The symmetric algebra construction works similarly; that is, we have
					\[
						\oblv_{ \partial }( \Sym^{ \gr }_{ A }(\g_{ 1 }) ) \simeq \Sym^{ \gr }_{ \oblv_{ \partial }(A) }(\oblv_{ \partial }(\g))
					\]
					as graded (commutative graded) $\oblv_{ \partial }(A)$-algebras. See \cref{remark:forget_differential}.
				\end{remark}

				The Poincaré-Birkhoff-Witt theorem for dg-Lie algebras (over a characteristic $0$ field $k$) gives a canonical equivalence of dg-algebras between $\gr(U(\h)^{ PBW })$ and $\Sym_{ k }(\h)$, for $\h$ a dg-Lie algebra over $k$. A classical reference is \cite[Theorem 2.3]{Quillen:RationalHomotopyTheory}.
				The same theorem holds for a specific class of dg-Lie algebroids. 
				The language here is $1$-categorical, so ``surjections'' and ``isomorphisms'' between dg-modules are to interpret degree-wise.

				\begin{thm}
					\label{thm:pbw_free_graded}
					There exists a canonical surjection $\Sym^{ \gr }_{ A }(\g) \to \gr(U(\g)^{ \PBW })$ of graded dg-$A$-algebras (where $\g$ is in weight $1$). If $\g$ is free as a \emph{graded} $A$-module (that is, forgetting the differential), then this map is an isomorphism. Moreover, in such a case, $U(\g)$ has an underlying free \emph{graded left} $A$-module.
				\end{thm}
				\begin{proof}
					To build the surjection start with the canonical morphism of graded dg-$A$-modules
					\[
						\g_{ 1 } \to \gr(U(\g)^{ \PBW })
					\]
					and apply the (graded) symmetric-forgetful adjunction. It is surjective by construction.
					To prove it is an isomorphism, we can forget the differentials, which is conservative and all operations interact well with such an action (see \cref{remark:forget_differential}).
					We will now implicitly forget differentials everywhere.\\
					Let us then assume that $\g$ has an underlying free graded $A$-module, then one can show as in the first part of \cite[Theorem 3.1]{Rinehart:LieAlgbd} that $U(\g)$ is a free \emph{left} $A$-module as well: the idea is to build a ``faithful'' representation of $U(\g)$ on $\Sym_{ A }(\g)$ to prove the $A$-linear independence of the basis elements given by increasing sequences of the basis elements of $\g$. Once we know the explicit basis for each $U(\g)^{ \leq p }$, then we can see that a basis of $\gr(U(\g)^{ \PBW })(p)$ is given by increasing sequences of $p$ elements (taken from a chosen ordered basis of $\g$), and therefore we can conclude that it is isomorphic to $\Sym^{ p }_{ A }(\g)$ (using the symmetrization map).
				\end{proof}
				\begin{coroll}
					\label{coroll:pbw_cofibrant_lie_algebroids}
					Let $\g$ be a cofibrant dg-Lie algebroid over $A$. Then the canonical surjection
					\[
						\Sym^{ \gr }_{ A }(\g) \to \gr(U(\g)^{ \PBW })
					\]
					of graded dg-$A$-algebras (where $\g$ is in weight $1$) is a ($1$-categorical, i.e.\ degreewise) isomorphism. Moreover, in such a case, $U(\g)$ has an underlying graded projective left $A$-module.
				\end{coroll}
				\begin{proof}
					Let us first recall that by the proof of \cite[Lemma 4.20]{Nuiten:HomotopicalLieAlgbd}, all cofibrant dg-Lie algebroids are retracts of free Lie algebroids on graded vector space maps like $M \to \T_{ A }$, once differentials are forgotten everywhere (see \cite[Proposition 4.18]{Nuiten:HomotopicalLieAlgbd} and observe that condition $(\star)$ is stable by retract and by attaching cells via elementary cofibrations of dg-Lie algebroids). 
					That is, we have a commutative diagram of graded Lie algebroids over $A$
					\begin{diag}
						& \mathrm{FreeGradedLieAlgbd}(M \to \T_A) \ar[dr] & \\
						\g \ar[ur] \ar[rr, "\id"] & & \g.
					\end{diag}
					where $M$ is a graded $\C$-vector space equipped with a map to $\T_{ A }$.
					By \cite[Remark 2.17]{Nuiten:HomotopicalLieAlgbd}, the free (graded) Lie algebroid over $A$ of $M \to \T_{ A }$ is obtained by
					\[
						A \otimes_{ \C } \mathrm{FreeGradedLie}(M) \to \T_{ A }
					\]
					where $\mathrm{FreeLie}$ is the free (graded) Lie algebra $1$-functor (over $\C$), the Lie bracket is the ``action Lie algebroid''  from \cite[Example 2.10]{Nuiten:HomotopicalLieAlgbd} (so it is not $A$-linear) and the map above is obtained by successive adjunction and scalar extensions.
					Since $M$ is a (graded) free vector space, also $\mathrm{FreeLie}(M)$ has an underlying (graded) free vector space (although way bigger than $M$) and thus $A \otimes_{ \C } \mathrm{FreeLie}(M)$ has an underlying (graded) free $A$-module.
					Since the construction of the PBW morphism is functorial, the diagram above implies that the arrow $\Sym_{ A }^{ \gr }(\g) \to \gr(U(\g)^{ \PBW })$ is a retract of 
					\[
						\Sym_{ A }^{ \gr }(\mathrm{FreeLieAlgbd}(M \to \T_{ A })) \to \gr(U( \mathrm{FreeLieAlgbd}(M \to T_{ A }) )^{ \PBW } ).
					\]
					By \cref{thm:pbw_free_graded}, the latter map is an isomorphism, and since isomorphisms are stable under retract we can conclude.
					The last statement comes from the fact that $U(\g)$ is a retract of a free graded left $A$-module.
				\end{proof}

				\begin{remark}
					\label{remark:joost_pbw}
					A stronger version of \cref{thm:pbw_free_graded} is used in the proof of \cite[Lemma 2.10]{Nuiten:KoszulLieAlgbd}, where the result is used not only for free graded $A$-modules but for all $\g$ with an underlying cofibrant dg-$A$-module. 
					We gave a different proof, with stronger assumptions, that is sufficient of our purposes.
				\end{remark}

				\begin{remark}
					\label{remark:stronger_pbw}
					For dg-Lie algebras, there exists a different theorem that is also known as (strong) PBW theorem. It identifies the cocommutative coalgebra (and, in particular the underlying dg-modules) $U(\g)$ with $\Sym(\g)$, see \cite[Theorem 2.2]{Hennion:TangentLieAlg}.
					We don't believe that the same holds for dg-Lie algebroids.
				\end{remark}

				Let us report here a technical fact that will be useful later.
				\begin{lemma}
					\label{lemma:projective_cofibrations_dg_A_mod}
					Let $A$ be a cdga over $\C$ and consider the projective model structure on $\Mod_{ A }^{ \dg }$. A map $i\colon V \to W$ of dg-$A$-modules is a cofibration if and only if it is a monomorphism, it splits as a map of \emph{graded} $A$-modules (i.e.\ forgetting differentials) and $C = W/V$ is a cofibrant dg-$A$-module.
				\end{lemma}
				\begin{proof}
					Observe that, since all objects of $\Mod_{ \C }^{ \dg }$ are projectively cofibrant, $\Mod_{ A }^{ \dg }$ inherits the projective model structure by \cite[Proposition 4.3.3.15]{Lurie:HA}. It is cofibrantly generated by \cite[Theorem 4.1]{SchwedeShipley:AlgebrasModulesMonoidalModelCategories}, and generating (trivial) cofibrations are obtained simply by taking free $A$-module functor on those generating (trivial) cofibrations for $(\Mod_{ \C }^{ \dg })_{ \textrm{proj} }$.
					If $A$ is a discrete commutative algebra this is proved in \cite[Proposition 18.5.3]{MayPonto:MoreConciseAT}.
					Observe now that the proof of \cite[Proposition 18.5.2]{MayPonto:MoreConciseAT} can be read also in dg setting: it says that any cofibrant dg-$A$-module is a projective graded $A$-module (that is, once we forget differentials). 
					The proof of \cite[Proposition 18.5.3]{MayPonto:MoreConciseAT} then goes through, by using the above fact (so all the relevant splittings still happen) and by observing that the aciclicity of the dg-$A$-module of morphisms still holds (this is what allows to ``fix'' the graded lift $\tilde{f_{ 2 }}$).
				\end{proof}

				\begin{prop}
					\label{prop:universal_enveloping_algebra_cofibrant}
					Let $\g$ be a cofibrant dg-Lie algebroid over $A$ and consider $U(\g)^{ \PBW }$. The underlying \emph{left} dg-$A$-module $U(\g)$ is cofibrant for the projective model structure on $\Mod_{ A }^{ \dg }$.
				\end{prop}
				\begin{proof}
					Let us observe that $U(\g)^{ \leq 0 } \simeq A$ is a cofibrant left dg-$A$-module.
					Let us prove, by induction on $n \in \N$, that $U(\g)^{ \leq n }$ is a cofibrant left dg-$A$-module. Assume this is true for $n$ and consider the following short exact sequence of left dg-$A$-modules coming from \cref{coroll:pbw_cofibrant_lie_algebroids}:
					\[
						0 \longrightarrow U(\g)^{ \leq n } \longrightarrow U(\g)^{ \leq n+1 } \longrightarrow \Sym_{ A }^{ n+1 }(\g) \longrightarrow 0.
					\]
					By assumption $\g$ is a cofibrant dg-Lie algebroid, which implies that its underlying $A$-module is cofibrant and therefore also $\Sym_{ A }^{ n+1 }(\g)$ is so (this holds because, being in characteristic $0$, it is a retract of $\g^{ \otimes_{ A }^{ n } }$, which is cofibrant because the projective model structure on $\Mod_{ A }^{ \dg }$ is symmetric monoidal model).
					By \cite[Proposition 18.5.3]{MayPonto:MoreConciseAT}, $\Sym_{ A }^{ n+1 }(\g)$ is a projective graded $A$-module and therefore the monomorphism $U(\g)^{ \leq n } \hookrightarrow U(\g)^{ \leq n+1 }$ is graded-split.
					By inductive assumption, $U(\g)^{ \leq n }$ is a cofibrant left dg-$A$-module.
					This proves that the monomorphism \[
						U(\g)^{ \leq n } \hookrightarrow U(\g)^{ \leq n+1 }
					\]
					has a cofibrant cokernel, so it is a cofibration. We conclude that $U(\g)^{ \leq n+1 }$ is cofibrant as well by \cref{lemma:projective_cofibrations_dg_A_mod}.
					Since cofibrations are stable under filtered colimits, we conclude that $U(\g)$ has an underlying cofibrant left dg-$A$-module.
				\end{proof}

				\begin{remark}
					\label{remark:relation_enveloping_operad}
					The underlying dg-$(A,A)$-bimodule of $U(\g)$ is cofibrant as well, thanks to \cite[Theorem 4.22]{Nuiten:HomotopicalLieAlgbd} (plugging in $p=1$ and observing that $\mathrm{Env}_{ \g }(1) \simeq U(\g)$ as bimodules). This is stronger: all cofibrant dg-$(A,A)$-bimodules are cofibrant left dg-$A$-modules but not viceversa.
				\end{remark}

				We can now lift this construction at the level of $\infty$-categories. 

				\begin{prop}
					\label{prop:pbw_infinity_categories}
					There exists a \emph{filtered universal enveloping algebra} $\infty$-functor 
					\[
						U(-)^{ \PBW }\colon \LieAlgbd_{ A } \longrightarrow \Alg( (A,A)-\BiMod^{ \cpl }),
					\]
					equipped with a natural equivalence
					\[
						\Sym_{ A }^{ \gr }(\oblv(-)_{ 1 }) \stackrel{\sim}{\longrightarrow} \gr(U(-)^{ \PBW })
					\]
					in $\Fun(\LieAlgbd_{ A }, \CAlg^{ \gr }_{ A })$, where $\oblv\colon \LieAlgbd_{ A } \to \Mod_{ A }$ is the forgetful functor.
				\end{prop}
				\begin{proof}
					Let us first recall that, thanks to the (projective) semi-model structure on dg-Lie algebroids on $A$, we have the equivalence
					\[
						\LieAlgbd_{ A } \coloneqq \LieAlgbd_{ A }^{ \dg }[\textrm{qis}^{ -1 }] \simeq \LieAlgbd_{ A }^{ \dg, \mathrm{cofib} }[\textrm{qis}^{ -1 }]
					\]
					of $\infty$-categories; this simply means that each dg-Lie algebroid is weakly equivalent to a cofibrant one.
					Let us now observe that since $A$ is a cofibrant commutative dg-$\C$-algebra, we have the following rectification result
					\[
						\Alg\left( (A, A)-\BiMod^{ \cpl } \right) \simeq \Alg^{ \dg }\left( (A,A)-\BiMod^{ \dg, \cpl }  \right)[\mathrm{qis}^{ -1 }]
					\]
					as explained in \cite[Lemma 2.13]{Nuiten:KoszulLieAlgbd} (here we add the filtration on bimodules: see \cite[Theorem 3.9, Theorem 3.50, Proposition 3.57]{GwilliamPavlov:FilteredDerivedCat} for properties of the model structure on complete filtered dg-$(A,A)$-bimodules, where we observe that the proofs go through for non-symmetric monoidal structures).\\
					Restricting ourselves to \emph{cofibrant} dg-Lie algebroids, we have defined a $1$-functor of universal enveloping algebra
				\[
					U(-)^{ \PBW }\colon \LieAlgbd_{ A }^{ \dg, \mathrm{cofib} } \longrightarrow \Alg^{ \dg }( (A,A)-\BiMod^{ \dg, \Fil \geq 0 } ),
				\]
				where we observed that the PBW filtration is $0$ in negative degrees.
				Observe that, given a weak equivalence (i.e.\ quasi-isomorphism) $\g \to \h$ of cofibrant dg-Lie algebroids, this induces a morphism \[
					U(\g)^{ \PBW } \to U(\h)^{ \PBW }
				\]
				in non-negatively filtered algebras in $(A,A)$-dg-bimodules. Passing to associated graded, by \cref{coroll:pbw_cofibrant_lie_algebroids}, we obtain the quasi-isomorphism of graded commutative algebras over $A$ 
				\[
					\Sym_{ A }^{ \gr }(\g) \stackrel{\sim}{\longrightarrow} \Sym_{ A }^{ \gr }(\h).
				\]
				Since $\gr$ preserves and detects weak equivalences\footnote{Basically by definition, since weak equivalences in complete filtered objects are \emph{graded} weak equivalences. In our case, since the filtration is only positive, this means, by successive extensions, a standard weak equivalence as sequences.} for complete filtered objects, we deduce that $U(\g)^{ \PBW } \to U(\h)^{ \PBW }$ was a weak equivalence (meaning a quasi-isomorphism on the underlying complexes) to begin with.
				This means that the above $1$-functor induces an $\infty$-functor between the ($\infty$-categorical) localizations inverting quasi-isomorphisms on both sides. That is, we have an $\infty$-functor (that we will denote in the same exact way):
				\[
					U(-)^{ \PBW }\colon \LieAlgbd_{ A } \longrightarrow \Alg( (A,A)-\BiMod^{ \cpl } ).
				\]
				The natural equivalence 
				\[
					\Sym_{ A }^{ \gr }( \oblv(-)) \stackrel{\sim}{\longrightarrow} \gr(U(-)^{ \PBW })
				\] 
				goes through similarly, since it is natural for cofibrant dg-Lie algebroids and the left-hand side is easily seen to preserve quasi-isomorphisms.
				\end{proof}
						
		\subsection{Representations of dg-Lie algebroids}
				\label{subsubsection:representations_dg_lie_algebroids}
				
				Similarly to representations of dg-Lie algebras, there exists a notion of representation of a dg-Lie algebroids.

				\begin{defn}
					\label{defn:representation_lie_algebroid}

					A representation of a dg-Lie algebroid $\g \to \T_{ A }$ is a dg-$A$-module $E$ equipped with a $\C$-linear action map $\nabla\colon \g \otimes_{ \C } E \to E$ satisfying 
					\[
						\nabla_{ a X }(e) = a \cdot \nabla_{ X }(e), \qquad \nabla_{ X }(a\cdot e) = \rho(X)(a) \cdot e + a \cdot \nabla_{ X }(e)
					\]
					for all $a \in A$, $X \in \g$ and $e \in E$.
					A morphism of $\g$-representations is an $A$-linear map preserving the action. They form a $1$-category $\infcatname{Rep}_{ \g }^{ \dg }$.
				\end{defn}

				\begin{remark}
					\label{remark:representations_as_left_modules}

					As for dg-Lie algebras, the ``universal enveloping algebra'' of a dg-Lie algebroid $\g$, denoted by $U(\g)$ (see \cref{subsubsection:universal_enveloping_algebra}), is built such that a (dg) $\g$-representation is the same thing as a left dg-$U(\g)$-module; that is, we have an equivalence of categories 
					\[
						\Rep_{ \g }^{ \dg } \simeq \LMod_{ U(\g) }^{ \dg }.
					\]
					Notice that by the right-hand side we mean left dg-modules over $U(\g)$ in the symmetric monoidal category of dg-modules over $\C$.	From now on, we will only think about left dg-$U(\g)$-modules.
				\end{remark}

				Here follows some observations.
				\begin{enumerate}[label=(\roman*), ref=(\roman*)]
					\item \label{liealgbd:1_forgetful_mod_A_conservative} The dg-algebra inclusion $A \hookrightarrow U(\g)$ induces a forgetful $1$-functor 
						\[
							\oblv_{ \g }\colon \LMod_{ U(\g) }^{ \dg } \to \Mod_{ A }^{ \dg },
						\]
						which is conservative and right adjoint.
					\item \label{liealgbd:bialgebra_1_categories} As explained in \cite{MoerdijkMrcun:UniversalEnvelopingAlgLieAlgbd} (which goes through also in the $1$-categorical dg-setting), $U(\g)$ is a dg-$(A/\C)$-bialgebra. In particular, this says that $U(\g)$ is a cocommutative dg-coalgebra in the symmetric monoidal category $\Mod_{ A }^{ \dg }$ (where $A$ acts from the left). In particular, calling $\Delta$ the coproduct, we have
						\[
							\Delta(a) = a \otimes 1 = 1 \otimes a, \qquad \Delta(x) = x \otimes 1 + 1 \otimes x
						\]
						for $a \in A \subset U(\g)$ and $x \in \g \subset U(\g)$.
						Moreover, this coalgebra structure is ``compatible'' with the multiplication in $U(\g)$, in the sense of \cite[Definition 2.1]{MoerdijkMrcun:UniversalEnvelopingAlgLieAlgbd}.

						\begin{remark}
							We currently don't know whether an $\infty$-categorical generalization of $(A/\C)$-bialgebras (which are defined through explicit formulas) exists.
						\end{remark}

					\item \label{liealgbd:dg_representation_on_A} There is a canonical representation of $\T_{ A }$ on $A$, where the tangent complex acts by derivation. For every dg-Lie algebroid $\g \to \T_{ A }$, this gives $A$ a canonical structure of $\g$-representation, i.e.\ of left dg-$U(\g)$-module.
						\begin{remark}
							\label{remark:A_not_Ug_A_bimodule}
							Since $A$ is not central in $U(\g)$, the action of $U(\g)$ on $A$ is not compatible with the action of $A$ on itself. That is, $A$ is a left dg-$U(\g)$-module but not a $(U(\g), A)$-bimodule.
						\end{remark}

					\item \label{liealgbd:dg_linear_dual_representation} As for the case of dg-Lie algebras, whose universal enveloping algebra is a cocommutative Hopf dg-algebra, the category of representations inherits a symmetric monoidal structure $(\LMod_{ U(\g) }^{ \dg }, - \otimes_{ A } -)$. That is, if $V, W$ are two $\g$-representations then $V \otimes_{ A } W$ inherits an action of $\g$ by 
						\[
							X.(v \otimes w) \coloneqq (X.v) \otimes w + v \otimes (X.w),
						\]
						ignoring eventual Koszul signs. The monoidal unit is $A$, seen with the action of $\g$ by derivations (through the anchor map).
						In other words, the forgetful $1$-functor
						\[
							\oblv_{ \g }\colon \LMod_{ U(\g) }^{ \dg } \longrightarrow \Mod_{ A }^{ \dg }
						\]
						is symmetric monoidal.

					\item \label{liealgbd:dg_internal_hom_representation} Another explicit construction equips $\LMod_{ U(\g) }^{ \dg }$ with an internal hom: it simply consists in considering the internal hom in dg-$A$-modules between two $\g$-representations and making $\g$ act ``by conjugation''. That is, the forgetful $1$-functor $\oblv_{ \g }$ is also closed.

					\item \label{liealgbd:dg_representations} Consider the projective model structure on $\Mod_{ \C }^{ \dg }$. It is combinatorial, cofibrantly generated and symmetric monoidal model: see \cite[Proposition 7.1.2.8, 7.1.2.11]{Lurie:HA}. Moreover, any object is cofibrant thanks to \cite[Remark 7.1.2.10]{Lurie:HA}, so in particular $U(\g)$ is a cofibrant dg-$\C$-module. 
						Using \cite[Proposition 4.3.3.15]{Lurie:HA}, $\LMod_{ U(\g) }^{ \dg }$ can be equipped with the so called projective model structure (as in \cite{Nuiten:KoszulLieAlgbd}). Weak equivalences are maps inducing quasi-isomorphisms of the underlying complexes and fibrations are degree-wise surjections. We will then be interested in the $\infty$-categorical localization
						\[
							\Rep_{ \g } \coloneqq \LMod_{ U(\g) }^{ \dg }[\textrm{qis}^{ -1}].
						\]

					\item \label{liealgbd:rectification_representation} By \cite[Theorem 4.3.3.17]{Lurie:HA}, since $U(\g)$ has an underlying cofibrant dg-$\C$-module, we have the following rectification result
						\[
							\Rep_{ \g } \coloneqq \LMod_{ U(\g) }^{ \dg }[\textrm{qis}^{ -1 }] \simeq \LMod_{ U(\g) },
						\]
						where we recall that, since we work over $\C$, our ``base'' is $\Mod_{ \C } \simeq \Mod_{ \C }^{ \dg }[\textrm{qis}^{ -1 }]$.
						\begin{remark}
							\label{remark:filtered_Ug_cofibrant}
							The proof of \cref{prop:universal_enveloping_algebra_cofibrant} actually yields that $U(\g)^{ \PBW }$ is (projectively) cofibrant as a filtered left dg-$A$-module (and hence as a filtered dg-$\C$-module), see \cite[Lemma 3.1, Example 3.3]{GwilliamPavlov:FilteredDerivedCat}. This gives the following rectification result for filtered representations
							\[
								\LMod_{ U(\g)^{ \PBW } }^{ \dg }( \Mod_{ \C }^{ \dg, \Fil } )[\textrm{qis}^{ -1 }] \simeq \LMod_{ U(\g)^{ \PBW } }(\Mod_{ \C }^{ \Fil }).
							\]
						\end{remark}

					\item \label{liealgbd:tensor_product_representations} Let us assume now that $\g$ is a cofibrant dg-Lie algebroid. To upgrade the tensor product built ``by hand'' on the $1$-category $\LMod_{ U(\g) }^{ \dg }$ into a symmetric monoidal structure on the $\infty$-category $\LMod_{ U(\g) }$, one can use another model structure on $\LMod_{ U(\g) }^{ \dg }$, the \emph{relative injective} model structure (which is simply called ``injective'' in the original source). This is because the projective model structure introduced previously is not compatible with the tensor product, since $A$, the monoidal unit, is not (projectively)-cofibrant, see \cite[Remark 2.12]{Nuiten:KoszulLieAlgbd}.\\
						The relative injective weak equivalences are the same (so the underlying $\infty$-category doesn't change) and the relative injective cofibrations are morphisms whose underlying map of dg-$A$-modules is a cofibration\footnote{Observe how the only model structure considered on $\Mod_{ A }^{ \dg }$ is the projective model structure.}. The existence of such model structure is proved in \cite[Lemma 2.11]{Nuiten:KoszulLieAlgbd} (observe how, in the proof, one needs the cofibrancy of $U(\g)$ as left dg-$A$-module, proved in \cref{prop:universal_enveloping_algebra_cofibrant}, to say that all projective cofibrations of $\LMod_{ U(\g) }^{ \dg }$ are cofibrations of dg-$A$-modules).
						As mentioned above, the interest of the relative injective model structure is its compatibility with the tensor product, and allows us to endow the $\infty$-category $\LMod_{ U(\g) }$ with a symmetric monoidal structure.

					\item \label{liealgbd:injective_structure_filtered_representation} Similarly, one can endow the $1$-category $\LMod_{ U(\g)^{ \PBW } }^{ \dg }(\Mod_{ \C }^{ \dg, \Fil, \geq 0 })$ with a relative injective model structure which is compatible with the symmetric monoidal structure (which is given again as tensor products over $\langle 0, A \rangle$). Let us first recall the projective model structure on $\Mod_{ \C }^{ \dg, \Fil, \geq 0 }$ and on $\Mod_{ A }^{ \dg, \Fil, \geq 0 }$ where weak equivalences and fibrations are defined objectwise. Both are symmetric monoidal model and combinatorial.
						The cofibrancy of $U(\g)^{ \PBW }$ in $\Mod_{ A }^{ \dg, \Fil, \geq 0 }$, by \cref{coroll:pbw_cofibrant_lie_algebroids} and \cite[Lemma 3.1]{GwilliamPavlov:FilteredDerivedCat}, allows the argument to go through. This also holds for generic $\Z$-indexed filtration.

					\item \label{liealgbd:injective_structure_complete_representation} Again, one can endow the $1$-category $\LMod_{ U(\g)^{ \PBW } }^{ \dg }(\Mod_{ \C }^{ \dg, \cpl })$ with a relative injective symmetric monoidal model structure, with respect to the completed tensor product. Recall that weak equivalences in $\Mod_{ A }^{ \dg, \cpl }$ are the graded equivalences and the projective model structure on $\Mod_{ A }^{ \dg, \cpl }$ is again well behaved (see \cite[Theorem 3.9, 3.50]{GwilliamPavlov:FilteredDerivedCat}). The only missing ingredient is again the cofibrancy of $U(\g)^{ \PBW }$ in $\Mod_{ A }^{ \dg, \cpl }$, which is true for the same reasons as above.

					\item \label{liealgbd:forgetful_conservative_strong_monoidal} We have then a symmetric monoidal, closed, conservative and left adjoint\footnote{This forgetful functor is induced by tensoring with the $(A, U(\g))$-bimodule $U(\g)$, see \cref{subsection:morita_adjunctions}.} $\infty$-functor
						\[
							\oblv_{ \g }\colon \LMod_{ U(\g) } \to \Mod_{ A }.
						\]

					\item \label{liealgbd:rep_functor}	Thanks to \cite[Lemma 2.10]{Nuiten:KoszulLieAlgbd}, restricted to cofibrant dg-Lie algebroids over $A$, there exists an $\infty$-functor
						\[
							\Rep_{ - }\colon \LieAlgbd_{ A }^{ \op } \longrightarrow \CAlg(\Pr^{ L })_{ -/\Mod_{ A } }, \qquad \g \mapsto \oblv_{ \g }\colon \LMod_{ U(\g) } \to \Mod_{ A }
						\]
						where a morphism $\g \to \h$ is mapped to the symmetric monoidal ``forgetful'' functor $\LMod_{ U(\h) } \to \LMod_{ U(\g) }$. This can be made precise using the injective model structures and Quillen adjunctions, see \cite[Remark 2.12]{Nuiten:KoszulLieAlgbd}.

					\item \label{liealgbd:rep_commutes_limits} By \cite[Lemma 2.13]{Nuiten:KoszulLieAlgbd}, the functor $\Rep_{ - }$ introduced above commutes with limits. Concretely, this means that 
						\[
							\Rep_{ \colim_{ i } \g_{ i } } \simeq \lim_{ i } \Rep_{ \g_{ i } }
						\]
						where the limit is computed in $\CAlg(\Pr^{ L })_{ -/\Mod_{ A } }$ (so the underlying $\infty$-category is the limit in $\Pr^{ L }_{ -/\Mod_{ A } }$) and the morphisms are the forgetful functors.

					\item Consider the free-forgetful adjunction
						\[
							\adjunction{U(\g) \otimes_{ A } -}{\Mod_{ A }}{\LMod_{ U(\g) }}{\oblv_{ \g }},
						\]
						which gives a map of left $U(\g)$-modules 
						\[
							U(\g) \to A,
						\]
						which is obtained by acting on $1 \in A$. Beware that this is not an algebra map (in $\Mod_{ \C }$), for similar reasons as \cref{remark:A_not_Ug_A_bimodule}.

					\item Considering $U(\g)$ as an $(A, U(\g))$-bimodule, we observe that 
					\[
						\oblv_{ \g } \simeq U(\g) \otimes_{ U(\g) } -
					\]
					so that we have a Morita right adjoint given by $\IntHom_{ A }(U(\g), -)$, see \cref{subsection:morita_adjunctions}.
				\end{enumerate}

%% file: BeraldoVsNuiten.tex
	\section{Derived D-Modules as some subcategory of IndCoh}
			\label{subsection:derived_DModules_IndCoh}

			In this section, we will give a summary of the definitions of interest from \cite{Ber:DerivedDMod} for the $\infty$-category of derived D-modules. We will only consider the (bounded) derived affine case.

			\subsection{Bounded case}
				\label{subsubsection:derived_DModules_IndCoh_bounded_case}
				Here the category of derived D-modules is simple to define as a full subcategory of $\IndCoh(X_{ \dR }) = \Crys^{ r }(X)$. Let us recall \cite[Definition 3.1.1]{Ber:DerivedDMod} in the ``absolute case'', i.e.\ for $f\colon X \to \Spec \C$ where we can identify the formal completion as the de Rham prestack, i.e.\ 
				\[
					(\Spec \C)^{ \wedge }_{ X } \simeq X_{ \dR }.
				\]

				\begin{defn}
					\label{defn:beraldo_bounded_case}

					Let $X$ be a bounded derived scheme (or algebraic stack) over $\C$. The $\infty$-category of \emph{Beraldo's derived D-modules on $X$} is defined as the pullback
					\begin{diag}
						\infcatname{D}^{\der}(X) \ar[r, hook, "\Upsilon^{\mathfrak{D}}_X"] \ar[d] \ar[dr, phantom, very near start, "\lrcorner"] & \IndCoh(X_{\dR}) \ar[d, "\oblv^r"] \\
						\Qcoh(X) \ar[r, hook, "\Upsilon_X"] & \IndCoh(X)
					\end{diag}
					in $\CAT$.
				\end{defn}
				Let us record some of its properties.

				\begin{prop}
					\label{prop:properties_beraldo_d_modules}
					Let $X$ be as above. The following assertions are true:
					\begin{enumerate}[label=(\roman*)]
						\item the $\infty$-category $\infcatname{D}^{ \der }(X)$ lives in $\CAlg(\Pr^{ L, st })$;
						\item if $X$ has a perfect cotangent complex $\L_{ X }$ (i.e.\ it is locally finitely presented) then the left arrow
							\[
								\infcatname{D}^{ \der }(X) \to \Qcoh(X)
							\]
							is a monadic right adjoint;
						\item if $X$ has a perfect cotangent complex then the $\infty$-category $\infcatname{D}^{ \der }(X)$ is compactly generated and self-dual in $\Pr^{ L }$.
					\end{enumerate}
				\end{prop}
				\begin{proof}
					Consider the pullback diagram in \cref{defn:beraldo_bounded_case}: the rightmost arrow is the $\IndCoh$ $!$-pullback along the canonical morphism $X \to X_{ \dR }$, see \cref{defn:!_pullback_ind_coh}, and the bottom arrow is the action on $\Qcoh(X)$ on $\omega_{ X } \in \IndCoh(X)$ (which is the monoidal unit for the $\otimes^{ \IC }$ product), see \cref{prop:upsilon_symm_monoidal}. Both of these arrows are symmetric monoidal, so the pullback can be equally seen as a pullback diagram in $\CAlg(\Pr^{ L, st })$, since $\CAlg(\Pr^{ L,st }) \to \Pr^{ L,st } \to \CAT$ commutes with limits.
					Moreover the boundedness of $X$ (in the affine case this simply means that the connective commutative algebra $\Gamma(X, \O_{ X })$ is bounded) implies that $\Upsilon_{ X }$ is a fully faithful functor, see \cite[Corollary 9.6.3]{G:IndCoh}.
					The last two properties are proved in \cite[Corollary 3.1.7, Proposition 3.2.4]{Ber:DerivedDMod}.
				\end{proof}
	
				For more results and functorialities for such $\infty$-category: see \cite[§3.3]{Ber:DerivedDMod}.
				Let us finally record here the following fact, proving that derived D-modules on bounded derived affine schemes are, in the end, equivalent to classical D-modules on them. That is, unboundedness is fundamental to have a theory of derived D-modules that really differs from the classical one.
			
				\begin{prop}
					\label{conj:beraldo_bounded_negative_result}
					The symmetric monoidal functor 
					\[
						\Upsilon^{ \mathfrak{D} }_{ X }\colon \infcatname{D}^{ \der }(X) \hookrightarrow \IndCoh(X_{ \dR })
					\]
					is an equivalence whenever $X$ is a bounded derived affine scheme.
				\end{prop}
				\begin{proof}
					This is proved in \cite[Corollary 4.3.13]{Ber:DerivedDMod} on affine schemes of the form
					\[
						Y_{ n } = \Spec(\Sym_{\C}(\C[n]))
					\]
					for $n \in \N_{ >0 }$ odd. 
				\end{proof}

				\begin{remark}
					\label{remark:beraldo_unbounded}
					There exists a more complicated definition of derived D-modules on unbounded derived schemes (which are locally finitely presented): see \cite[Definition 4.2.9]{Ber:DerivedDMod}. They are again defined as algebras in $\Qcoh(X)$ for a certain monad, obtained as a ``filtered renormalization'' of the ind-coherent universal enveloping algebra of the ind-coherent tangent Lie algebroid of $X$.
					There is still a morphism 
					\[
						\Upsilon_{ X }^{ \mathfrak{D} }\colon \infcatname{D}^{ \der }(X) \to \IndCoh(X_{ \dR })
					\]
					which is not fully faithful anymore. 
					This fails to be an equivalence for the simplest unbounded derived affine scheme $Y_{ 2 } = \Spec(\Sym_{ \C }(\C[2]))$, see \cite[Corollary 4.3.13]{Ber:DerivedDMod}.
				\end{remark}
			
		\section{Relation to dg-Lie algebroids}
			\label{subsection:relation_lie_algebroids}

			There is an interpretation of \cref{defn:beraldo_bounded_case} in terms of dg-Lie algebroids and their representations. Let us immediately point out that, although Lie algebroids are mentioned already in \cite{Ber:DerivedDMod}, the author refers there to the $\infty$-category of Lie algebroids defined in \cite[Chapter 8]{GR:DAG2} (as ``ind-coherent'' formal moduli problems under the fixed base). 
			On the other hand, we will work with dg-Lie algebroids (see \cref{subsection:dg_lie_algebroids}).

			\subsection{Formal moduli problems and Koszul duality}
				\label{subsubsection:formal_moduli_problems_koszul}

				In this subsection we will state the main definitions and properties of formal moduli problems, with the goal of being able to state the Koszul duality result for dg-Lie algebroids in \cite{Nuiten:KoszulLieAlgbd}. The latter result is essential, since it allows us to understand dg-Lie algebroids over $A$ by restricting to the much more manageable \emph{good} dg-Lie algebroids.
				The definitions of small objects and formal moduli problems can be given in a formal setting: see \cite[Chapter 12]{Lurie:SAG} or \cite{CalaqueGrivaux:FormalModuli}. See also \cite{Toen:FMP}.

				The first notion we need to recall is that of (not necessarily split) square-zero extension of commutative algebras.
				\begin{defn}[{\cite[Definition 7.4.1.6]{Lurie:HA}}]
					\label{defn:square_zero_extension}
					Let $A \in \CAlg_{ \C }$ and $M \in \Mod_{ A }$. A \emph{square-zero extension} of $A$ by $M$ is the datum of a derivation $\eta \in \pi_{ 0 }\Map_{ A }(\L_{ A }, M[1])$ and of the morphism $f\colon A^{ \eta } \to A$ sitting in the pullback diagram
					\begin{diag}
						A^{\eta} \ar[r, "f"] \ar[d] \ar[dr, phantom, very near start, "\lrcorner"] & A \ar[d, "d_{\eta}"] \\
						A \ar[r, "d_0"] & A \oplus M[1].
					\end{diag}
					We denoted by $d_{ 0 }$ and $d_{ \eta }$ the morphisms in $\pi_{ 0 }\Map_{ \CAlg_{ -/A } }(A, A \oplus M[1])$ corresponding to $0$ and $\eta$ in $\pi_{ 0 }\Map_{ A }(\L_{ A }, M[1])$. The bottom right term is the split square-zero extension of $A$ by $M[1]$, as defined in \cref{defn:split_square_zero_extension}.
				\end{defn}	

				\begin{remark}
					\label{remark:square_zero_extension}
					As explained in \cite[Warning 7.4.1.10]{Lurie:HA}, it is not clear how to give a general \emph{intrinsic} characterization of square-zero extensions. That is, the functor
					\[
						\Phi\colon \infcatname{Der} \to \Fun(\Delta^{ 1 }, \CAlg_{ \C }), \qquad (A, M, \eta\colon \L_{ A } \to M[1]) \mapsto (A^{ \eta } \to A)
					\]
					is not always fully faithful. The situation is remedied when some boundedness conditions are imposed, see \cite[Theorem 7.4.1.26]{Lurie:HA} and \cite[Proposition 7.4.2.5]{Lurie:HA}.
				\end{remark}
			
				Square-zero extensions can be used to introduce two full subcategories of commutative algebras over $A$.

				\begin{defn}[{\cite[Definition 1.2]{Nuiten:KoszulLieAlgbd}}]
					\label{defn:small_commutative_algebras_over_base}
					Let $A \in \CAlg_{ \C }^{ \cn }$ be a connective commutative $\C$-algebra. The $\infty$-category $\CAlg^{ \smal }_{ -/A}$ is the smallest full subcategory of $(\CAlg_{ \C }^{ \cn })_{ -/A }$ which contains $A$ and is closed under square-zero extensions by $A[n]$ for $n \geq 0$.
					We refer to them as \emph{small commutative algebras over $A$}.
				\end{defn}

				\begin{defn}
					\label{defn:small_coherent_commutative_algebras_over_base}
					Let $A \in \CAlg_{ \C }^{ \cn, \coh }$ be a connective commutative $\C$-algebra which is coherent. The $\infty$-category $\CAlg^{ \smal, \coh }_{ -/A}$ is the smallest full subcategory of $(\CAlg_{ \C }^{ \cn })_{ -/A }$ which contains $A$ and is closed under square-zero extensions by \emph{coherent} $A$-modules lying in $\Mod_{ A }^{ \leq 0 }$.
					We refer to them as \emph{small-coherent commutative algebras over $A$}.
				\end{defn}
	
				\begin{remark}
					\label{remark:small_calg_deformation_context}
					Equivalently, the $\infty$-category $\CAlg_{ -/A }^{ \smal }$ is obtained as the small (or artinian) objects in the \emph{deformation context} $(\CAlg_{ \C/-/A }, \{A\})$ where $A \in \Mod_{ A } \simeq \Sp(\CAlg_{ -/A })$, see \cite[Definition 12.1.1.1, Example 12.1.2.2]{Lurie:SAG}.
					It is essentially small by \cite[Remark 12.1.2.8]{Lurie:SAG}.
				\end{remark}

				\begin{remark}
					\label{remark:small_coherent_calg_essentially_small}
					The $\infty$-category $\CAlg^{ \smal, \coh }_{ -/A }$ is essentially small as well.
				\end{remark}

				Let us record some easy observation about these two categories.

				\begin{lemma}
					\label{lemma:different_small_algebras}
					Let $A \in \CAlg^{ \cn }_{ \C }$ be a bounded coherent commutative algebra. There are fully faithful embeddings
					\[
						\CAlg^{ \smal }_{ -/A } \hookrightarrow \CAlg^{ \smal, \coh }_{ -/A } \hookrightarrow \CAlg^{ \cn, \bdd}_{ -/A }.
					\]
					Moreover, each $B \to A$ in $\CAlg^{ \smal, \coh }_{ -/A }$ induces an isomorphism $H^{ 0 }(B)^{ \red } \simeq H^{ 0 }(A)^{ \red }$.
				\end{lemma}
				\begin{proof}
					The first claim comes from the inclusion
					\[
						\Perf_{ A } \hookrightarrow \Coh_{ A }
					\]
					which holds by boundedness of $A$. The second claim is immediate from \cref{defn:small_coherent_commutative_algebras_over_base}.
				\end{proof}

				\begin{remark}
					\label{remark:why_small_coherent_algebras}
					The essentially small $\infty$-category $\CAlg^{ \smal, \coh }_{ -/A }$ naturally arises in derived algebraic geometry. In fact, consider $B \in \CAlg^{ \cn, \coh}$ and its convergent Postnikov tower
					\[
						\dots \longrightarrow \tau^{ \geq -2 }(B) \longrightarrow \tau^{ \geq -1 }(B) \longrightarrow \tau^{ \geq 0 }(B) \simeq H^{ 0 }(B).
					\]
					It is well known that $\tau^{ \geq n-1 }(B)$ is a square-zero extension of $\tau^{ \geq n }(B)$ by $H^{ n }(B)[n]$, see \cite[Corollary 7.4.1.28]{Lurie:HA}. Since $H^{ n }(B)[n]$ is a coherent $H^{ 0 }(B)$-module, we obtain that each $\tau^{ \geq n }(B)$ lives in $\CAlg^{ \smal, \coh }_{ -/H^{ 0 }(B)}$. For a more extensive explanation, see \cite[Remark 7.4.1.29]{Lurie:HA} and \cite[Chapter 1]{GR:DAG2}. Finally, let us observe that each nilpotent embedding of derived affine schemes is controlled by iterated square-zero extensions by connective coherent modules, see \cite[Chapter 1, 5.5]{GR:DAG2}.
				\end{remark}

				We are now ready to introduce ($\E_{ \infty }$) \emph{formal moduli problems under $\Spec A$}, specializing the general definition of formal moduli problems given in \cite[Definition 12.1.3.1]{Lurie:SAG}. This is denoted, in general, by $\infcatname{Moduli}^{ \CAlg_{ -/A } }$ by Lurie, but we decided to keep the notation from Nuiten.

				\begin{defn}[{\cite[Definition 1.2]{Nuiten:KoszulLieAlgbd}}]
					\label{defn:fmp_under_A}
					Let $A \in \CAlg^{ \cn, \bdd }_{ \C }$. The $\infty$-category of \emph{formal moduli problems} under $\Spec A$, denoted by $\FMP_{ \Spec A/- }$, is the full subcategory of $\Fun(\CAlg_{ -/A }^{ \smal }, \S)$ spanned by those $X$ such that
					\begin{enumerate}[label=(\roman*)]
						\item $X(A) \simeq *$ is a contractible space;
						\item $X$ preserves any pullback diagram of the form
							\begin{diag}
								B^{\eta} \ar[d] \ar[r] \ar[dr, phantom, very near start, "\lrcorner"] & A \ar[d, "d_0"] \\
								B \ar[r, "d_{\eta}"] & A \oplus A[n+1]
							\end{diag}
							which makes $B^{ \eta } \to B$ a square-zero extension of $B$ by $A[n]$ for $n \geq 0$, where $d_{ \eta } \in \pi_{ 0 }\Map_{ A }(\L_{ B } \otimes_{ B } A, A[n+1])$.
					\end{enumerate}
				\end{defn}
			
				\begin{remark}
					\label{remark:fmp_presentable_localisation}
					Observe that $\FMP_{ \Spec A/- }$ is presentable, by \cite[Remark 12.1.3.5]{Lurie:SAG}. Moreover, 
					\[
						\FMP_{ \Spec A/- } \hookrightarrow \Fun(\CAlg_{ -/A }^{ \smal }, \S)
					\]
					commutes with limits and filtered colimits, so that we obtain a localisation.
				\end{remark}
				
				There exists a second notion of formal moduli problems, which we call \emph{ind-coherent formal moduli problems}.

				\begin{defn}[{\cite[Chapter 5, 1.3.1]{GR:DAG2}}]
					\label{defn:ind_coherent_fmp}
					Let $A \in \CAlg^{ \cn, \coh }_{ \C }$. The $\infty$-category of \emph{ind-coherent formal moduli problems} under $\Spec A$, denoted by $\FMP^{ ! }_{ \Spec A/- }$, is the full subcategory of $(\PreStk_{ \laft })_{ \Spec A/- }$ spanned by those $\Spec A \to X$ where $X \in \PreStk_{ \laft, \defo }$ and $\Spec A \to X$ is a nil-isomorphism.
				\end{defn}

				\begin{remark}
					\label{remark:ind_coherent_fmp_as_thick}
					Remark that what we call $\FMP^{ ! }_{ \Spec A/- }$ is denoted by $\infcatname{ModuliStk}_{ A/\C }$ in \cite{BrantnerMagidsonNuiten:FormalIntegrationDerivedFoliations} and by $\infcatname{Thick}(\Spec A)$ in \cite[§4.1]{CalaqueGrivaux:FormalModuli}.
				\end{remark}

				Let us observe that our definition of ind-coherent formal moduli problems coincide with the ``pro-coherent formal moduli problems'' considered in \cite{Nuiten:KoszulLieAlgbd}. 

				\begin{lemma}
					\label{lemma:ind_coh_fmp_from_small_coherent}
					Let $A \in \CAlg^{ \cn, \coh }_{ \C }$. There exists a fully faithful embedding
					\[
						\FMP^{ ! }_{ \Spec A/- } \hookrightarrow \Fun(\CAlg^{ \smal, \coh }_{ -/A }, \S)
					\]
					whose essential image is spanned by those $X\colon \CAlg^{ \smal, \coh }_{ -/A }\to \S$ satisfying analogous conditions of \cref{defn:fmp_under_A}. 
				\end{lemma}
				\begin{proof}
					Let us first assume that $A$ is bounded. Then the arrow is obtained as composite of
					\[
						\FMP^{ ! }_{ \Spec A/-} \hookrightarrow (\PreStk_{ \laft })_{ \Spec A/- } \simeq \Fun_{ \mathrm{acc} }(\CAlg_{ -/A }^{ \bdd, \afp, \cn }, \S) \to \Fun(\CAlg^{ \smal, \coh }_{ -/A }, \S)
					\]
					where the last arrow is restriction along $\CAlg^{ \smal, \coh }_{ -/A } \hookrightarrow \CAlg_{ -/A }^{ \bdd, \afp, \cn }$.
					By \cref{remark:why_small_coherent_algebras}, we have
					\[
						(\dAff^{ \bdd, \afp }_{ \textrm{nil-isom from $\Spec A$ } })^{ \op } \simeq \CAlg^{ \smal, \coh }_{ -/A }
					\]
					and we can conclude by applying \cite[Chapter 5, Proposition 1.4.2]{GR:DAG2}. 
					The statement holds also in the general case, i.e.\ where $A$ is not bounded, by \cite[Theorem 5.3]{BrantnerMagidsonNuiten:FormalIntegrationDerivedFoliations}.
				\end{proof}

				\begin{coroll}
					\label{coroll:fmp_gr_presentable_localisation}
					The $\infty$-category $\FMP^{ ! }_{ \Spec A/- }$ is presentable and the fully faithful embedding
					\[
						\FMP^{ ! }_{ \Spec A/- } \hookrightarrow (\PreStk_{ \laft })_{ \Spec A/- }
					\]
					commutes with limits and filtered colimits, so that we obtain a localisation.
				\end{coroll}

				\begin{defn}
					\label{defn:formal_spectrum}
					The Yoneda embedding restricts to a functor
					\[
						\Spf\colon (\CAlg_{ -/A })^{ \op } \longrightarrow \FMP_{ \Spec A/- }, \qquad B \mapsto (\Map_{ \CAlg_{ -/A } }(B, -) \colon \CAlg_{ -/A}^{ \smal } \to \S)
					\]
					which is called the \emph{formal spectrum functor}.
				\end{defn}	
				We immediately observe that the restriction of $\Spf$ to $\CAlg^{ \smal }_{ -/A }$ is fully faithful, by the Yoneda lemma.	Let us explain the relation between these two version of formal moduli problems.

				\begin{prop}
					\label{prop:fmp_in_geometry}
					Let $A \in \CAlg^{ \cn, \lfp }_{ \C }$ be bounded. 
					There exists a commutative diagram
					\begin{diag}
						\ar[d, "\Spf"] (\CAlg^{\cn, \lfp}_{-/A})^{\op} \ar[r, "\Spec"] & (\dAff)_{\Spec A/-} \ar[d, "(-)^{\wedge}_{\Spec A}"] \\
						\FMP_{\Spec A/-} \ar[r, "\mathrm{LKE}", hook] & \FMP^!_{\Spec A/-},
					\end{diag}
					where the bottom arrow is the left Kan extension along the embedding $\CAlg_{ -/A }^{ \smal } \hookrightarrow \CAlg^{ \smal, \coh}_{ -/A }$ and the rightmost arrow is formal completion along $\Spec A$.
				\end{prop}
				\begin{proof}
					Observe that, for any $\Spec A \to \Spec B$, we have
					\[
						(\Spec B)^{ \wedge }_{ \Spec A } \in \FMP^{ ! }_{ \Spec A/- }.
					\]
					This is true because $\Spec A \to (\Spec B)^{ \wedge }_{ \Spec A }$ is a nil-isomorphism, and this formal completion is an inf-scheme by \cite[Chapter 2, 3.1.3]{GR:DAG2}.
					The bottom arrow given by left Kan extension gives a fully faithful embedding
					\[
						\FMP_{ \Spec A/- } \hookrightarrow \Fun(\CAlg_{ -/A }^{ \smal, \coh }, \S).
					\]
					It suffices now to verify that any element in the essential image verifies the pullback preserving conditions, see \cref{lemma:ind_coh_fmp_from_small_coherent}.
					The left Kan extension above is equivalent to first left Kan extend along $\CAlg^{ \smal }_{ -/A } \hookrightarrow \CAlg^{ \cn }_{ -/A }$ and then restrict along $\CAlg^{ \smal, \coh }_{ -/A } \hookrightarrow \CAlg^{ \cn }_{ -/A }$. 
					Observe that, if $X \in \FMP_{ \Spec A/- }$, then its left Kan extension will still send $A$ to $X(A) \simeq *$, because $A \in \CAlg^{ \smal }_{ -/A }$.
					It suffices to prove that the (big) left Kan extension of $X$ preserves pullbacks of the form 
					\begin{diag}
						C \ar[r] \ar[d] \ar[dr, phantom, very near start, "\lrcorner"] & A \ar[d, "d_0"] \\
						B \ar[r, "d_{\eta}"] & A \oplus V
					\end{diag}
					where $V \in \Coh_{ A }$ and $d_{ \eta }$ is a derivation from $B$ to $V$ (seen as a $B$-module thanks to the algebra map $B \to A$).
					We divide this proof into three steps.\\
					\textbf{Claim 1: Each $B \in \CAlg^{ \smal, \coh }_{ -/A }$ is a filtered colimit of small commutative algebras over $A$.}\\
					By assumption, we have
					\[
						B = B_{ n } \to B_{ n-1 } \to \dots \to B_{ 1 } \to B_{ 0 } \simeq A
					\]
					and each $B_{ i } \to B_{ i-1 }$ is a square-zero extension by a coherent $A$-module.
					Let us first prove the claim for $B_{ 1 }$. We have a pullback diagram
					\begin{diag}
						B_1 \ar[d] \ar[r] \ar[dr, phantom, very near start, "\lrcorner"] & A \ar[d, "d_0"] \\
						A \ar[r, "d_{\eta_1}"] & A \oplus W(1),
					\end{diag}
					in $\CAlg^{ \cn }_{ -/A }$ where $W(1) \in \Coh_{ A } \cap \Mod_{ A }^{ \leq -1 }$. By \cite[Proposition 7.2.4.11]{Lurie:HA}, we can write 
					\[
						W(1) \simeq \colim_{ n } W(1)_{ n }	
					\]
					with each $W(1)_{ n } \in \Perf_{ A }^{ \leq -1 }$. Let us now observe that the split-square zero extension functor $A \oplus -\colon \Mod_{ A }^{ \leq 0 } \to \CAlg_{ -/A }^{ \cn }$ commutes with filtered colimits: in fact, its left adjoint is given by
					\[
						\L_{ (-) } \otimes_{ (-) } A\colon \CAlg_{ -/A }^{ \cn } \to \Mod_{ A }^{ \leq 0 }
					\]
					which preserves compact objects by \cref{thm:locally_finite_presentation_perfect_cotangent_complex} (after having observed that $\CAlg_{ -/A }^{ \cn } \to \CAlg_{ \C }^{ \cn }$ preserves compact objects since its right adjoint $- \times A$ preserves filtered colimits). We then have
					\[
						A \oplus W(1) \simeq \colim_{ n } (A \oplus W(1)_{ n }) \in \CAlg^{ \cn }_{ -/A }.
					\]
					Since $A$ is a locally finitely presented commutative algebra over $\C$ (i.e.\ a compact object in $\CAlg^{ \cn }_{ \C }$), we can find $m$ such that $d_{ \eta_{ 1 } }$ factors as 
					\[
						A \stackrel{d_{ \eta_{ 1, m } }}{\longrightarrow} A \oplus W(1)_{ m } \longrightarrow A \oplus W(1).
					\]
					Since filtered colimits commutes with pullbacks in $\Mod_{ \C }$ and so does the forgetful functor $\oblv\colon \CAlg^{ \cn }_{ -/A } \to \Mod_{ \C }^{ \cn }$, we conclude that
					\[
						B_{ 1 } \simeq \colim_{ n \geq m } (A \times_{ A \oplus W(1)_{ n } } A) \in \CAlg^{ \cn }_{ -/A }.
					\]
					By perfection of each $W(1)_{ n }$, it is clear that the right-hand side is a filtered colimit of small $A$-algebras, since $\CAlg_{ -/A }^{ \smal } \hookrightarrow \CAlg_{ -/A }^{ \cn }$ commutes with pullbacks (see \cite[Remark 12.1.2.11]{Lurie:SAG}).
					By induction on $i$, we can now assume that $B_{ i } \simeq \colim_{ j } C_{ j }$ for a filtered poset $J$ (see \cite[Proposition 5.1.3.18]{Lurie:HTT}) and $C_{ j } \in \CAlg^{ \smal }_{ -/A }$. We have a pullback diagram
					\begin{diag}
						B_{i+1} \ar[d] \ar[r] \ar[dr, phantom, very near start, "\lrcorner"] & A \ar[d, "d_0"] \\
						B_i \ar[r, "d_{\eta}"] & A \oplus W(i+1),
					\end{diag}
					with $W(i+1)$ a coherent $A$-module. As before, it can be approximated by perfect $A$-modules $W(i+1)_{ n }$ so that
					\[
						A \oplus W(i+1) \simeq \colim_{ n \in \N } (A \oplus W(i+1)_{ n }) \in \CAlg^{ \cn }_{ -/A }.
					\]
					Observe now that we can define a morphism of posets $\phi\colon J \to \N$ that sends $j$ to the minimum $\phi(j) \in \N$ such that we have a factorization
					\[
						C_{ j } \to A \oplus W(i+1)_{ \phi(j) } \to A \oplus W(i+1) 
					\]
					in $\CAlg^{ \cn }_{ -/A }$.
					Using the ``rectification result'' of \cite[Proposition 5.3.5.15]{Lurie:HTT}, where $A$ is the finite poset given by the cospan shape, we deduce that the cospan $B_{ i } \rightarrow A \oplus W(i+1) \leftarrow A$ can be obtained as a filtered colimit of cospans
					\begin{diag}
						& A \ar[d, "d_0"] \\
						C_j	\ar[r, "d_{\eta_{j{,} n}}"] & A \oplus W(i+1)_n
					\end{diag}
					where the indexing poset is $\Delta = \{(n, j) | n \geq \phi(j) \} \subseteq J \times \N$. It is small and filtered, being a finite product of small filtered posets, and we claim that the two projections $\Delta \to J$ and $\Delta \to \N$ are cofinal. It is evident for both, as we can always choose a big enough $n \in \N$ once $j$ is fixed. A morphism over $A$ from $C_{ j }$ to $A \oplus W(i+1)_{ n }$ corresponds to a derivation $\eta_{ j, n }$ from $C_{ j }$ to $W(i+1)_{ n }$.
					As before, we deduce that
					\[
						B_{ i+1 } \simeq \colim_{ (n, j) \in \Delta } (C_{ j } \times_{ A \oplus W(i+1)_{ n } } A) \in \CAlg^{ \cn }_{ -/A }.
					\]
					Since $C_{ j } \in \CAlg^{ \smal }_{ -/A }$ and $W(i+1)_{ n } \in \Perf_{ A }^{ \leq -1 }$, we deduce that the right-hand side is a filtered colimits of small commutative $A$-algebras.
					This proves our first claim.\\
					\textbf{Claim 2: $X$ commutes with such filtered colimits.}\\
					It now suffices to consider the (big) left Kan extension of $X$. We claim that if $C \in \CAlg^{ \smal, \coh }_{ -/A }$ can be written as a filtered colimit of small algebras over $A$ as 
					\[
						C \simeq \colim_{ j \in J } C_{ j },
					\]
					then $X(C) \simeq \colim_{ j } X(C_{ j })$. This comes from the left Kan extension formula:
					\[
						X(C) \simeq \colim_{ K \to C } X(K)
					\]
					where the indexing $\infty$-category is $\CAlg^{ \smal }_{ -/A } \times_{ \CAlg^{ \smal, \coh }_{ -/A } } (\CAlg^{ \smal, \coh }_{ -/A })_{ -/C }$. Our claim holds simply because we have a commutative diagram
					\begin{diag}
						J \ar[r] \ar[dr] & \CAlg^{ \smal }_{ -/A } \times_{ \CAlg^{ \smal, \coh }_{ -/A } } (\CAlg^{ \smal, \coh }_{ -/A })_{ -/C } \ar[d, hook] \\
										 & (\CAlg^{\smal, \coh}_{-/A})_{-/C}
					\end{diag}
					where the diagonal arrow is cofinal by assumption, and the vertical arrow is fully faithful. It is then an easy application of Quillen's theorem A (see \cite[Theorem 4.1.3.1]{Lurie:HTT}) that gives that also the horizontal arrow is cofinal.\\
					\textbf{Conclusion of the first part of the proof.}\\
					To conclude, we can write, exactly like the proof of the first claim, the initial pullback diagram of small-coherent algebras over $A$ as 
					\[
						C = B \times_{ A \oplus V } A \simeq \colim_{ i } (B_{ i } \times_{ A \oplus V_{ i } } A) \in \CAlg^{ \cn }_{ -/A }
					\]
					for some filtered poset $I$, $B_{ i } \in \CAlg^{ \smal }_{ -/A }$ and $V_{ i } \in \Perf_{ A }^{ \leq -1 }$. We then have
					\[
						X(C) \simeq \colim_{ i } X(B_{ i } \times_{ A \oplus V_{ i } } A) \simeq \colim_{ i } (X(B_{ i }) \times_{ X(A \oplus V_{ i }) } X(A)) \in \S
					\]
					by the Schlessinger condition on $X \in \FMP_{ \Spec A/- }$.
					Since filtered colimits in spaces are left exact by \cite[Proposition 5.3.3.3]{Lurie:HTT}, we conclude that
					\[
						X(C) \simeq \colim_{ i } (X(B_{ i }) \times_{ X(A \oplus V_{ i }) } X(A)) \simeq (\colim_{ i } X(B_{ i })) \times_{ \colim_{ i } X(A \oplus V_{ i }) } X(A) \simeq X(B) \times_{ X(A \oplus V) } X(A),
					\]
					where the last equivalence is claim 2 proved above.\\
					\textbf{The diagram commutes.}\\
					Let us finally prove that the diagram commutes. Let us immediately notice that the forgetful functor
					\[
						\CAlg^{ \cn }_{ -/A } \to \CAlg^{ \cn }_{ \C }	
					\]
					preserves compact objects: in fact, it is left adjoint to the functor sending $B$ to $A \times B \to A$, which preserves filtered colimits (which are computed at the level of underlying spectra, where filtered colimits are left exact by \cite[Proposition 1.4.3.7]{Lurie:HA}). This removes any possible ambiguity from the notation $\CAlg^{ \cn, \lfp }_{ -/A }$. Let us fix an object $B \to A$ in this $\infty$-category. To conclude the proof it suffices now to prove that the formal completion $(\Spec B)^{ \wedge }_{ \Spec A }$ is left Kan extended from its restriction to \emph{small} commutative algebras over $A$.
					Notice that this is stronger than \cite[Chapter 5, Proposition 1.4.2]{GR:DAG2} since we are restricting to a smaller $\infty$-category. 
					Observe now that for each $C \to A \in \CAlg^{ \smal, \coh }_{ -/A }$ we obtain a nil-isomorphism $\Spec A \to \Spec C$ by \cref{lemma:different_small_algebras}; this means that 
					\[
						(\Spec B)_{ \Spec A }^{ \wedge }(C) \simeq \Map_{ \CAlg_{ -/A } }(B, C).
					\]
					The question boils down to prove that the natural morphism
					\[
						\mathrm{LKE}(\Map_{ \CAlg_{ -/A } }(B, -))(C) \longrightarrow \Map_{ \CAlg_{ -/A } }(B, C)
					\]
					is an equivalence of spaces, where the left-hand side is the left Kan extension of the functor $\Map_{ \CAlg_{ -/A } }(B, -)$ restricted to \emph{small} algebras.
					Using the compactness of $B$ and writing $C \in \CAlg^{ \smal, \coh }_{ -/A }$ as a filtered colimit of small $A$-algebras as in claim $1$, we conclude.
				\end{proof}

				\begin{remark}
					\label{remark:relation_with_thick_calaque_grivaux}
						The proposition above gives a rigorous proof of the fully faithfulness of the embedding mentioned in \cite[Proposition 4.1]{CalaqueGrivaux:FormalModuli}. We do not believe that it is an equivalence as stated there, though. That is, we have only a fully faithful embedding
						\[
							\FMP_{ \Spec A/- } \hookrightarrow \infcatname{Thick}(\Spec A).
						\]
						One object in the right-hand side which does not belong to the essential image is any formal completion $(\Spec B)_{ \Spec A }^{ \wedge }$ where $B$ is not a compact connective commutative algebra (that is, $\L_{ B }$ is not perfect).
						We claim that the essential image is given by full subcategory of $\FMP^{ ! }_{ \Spec A/- }$ spanned by those prestacks with a \emph{pro-perfect} cotangent complex.
				\end{remark}

				Let us now introduce good Lie algebroids over $A$, defined in \cite[§5.2]{Nuiten:KoszulLieAlgbd}.

				\begin{defn}
					\label{defn:good_Lie_algebroids}
					Let $A \in \CAlg^{ \cn}_{ \C }$. A dg-Lie algebroid $\g \in \LieAlgbd_{ A }$ is called \emph{good} if it admits a finite filtration
					\[
						0 \to \g^{ (1) } \to \g^{ (2) } \to \dots \to \g^{ (n-1) } \to \g^{ (n) } \simeq \g
					\]
					such that for each $1 \leq i \leq n$ there exists a pushout diagram of Lie algebroids on $A$
					\begin{diag}
						\FreeLieAlgbd(A[n_i] \stackrel{0}{\to} \T_A) \ar[r] \ar[d] & \g^{(i-1)} \ar[d] \\
						0 \ar[r] & \g^{(i)} \ar[ul, phantom, very near start, "\ulcorner"],
					\end{diag}
					for $n_{ i } \leq -2$ for each $i$.
					They span an essentially small full subcategory of $\LieAlgbd_{ A }$ denoted by $\LieAlgbd_{ A }^{ \good }$.
				\end{defn}

				Good Lie algebroids are particularly manageable. Recall the $1$-functor that forgets differential
				\[
					\oblv_{ \partial }\colon \Mod_{ \C }^{ \dg } \longrightarrow (\Mod_{ \C }^{ \heartsuit })^{ \gr }
				\]
				from \cref{remark:forget_differential}.

				\begin{prop}[{\cite[Lemma 5.9]{Nuiten:KoszulLieAlgbd}}]
					\label{defn:underlying_module_good_lie_algebroids}
					Let $\g \in \LieAlgbd_{ A }$ be good (choose a cofibrant representative). The following assertions hold:
					\begin{enumerate}[label=(\arabic*)]
						\item there exists a strictly positively graded $\C$-vector space $M$ equipped with a $\C$-linear morphism towards $\oblv_{ \partial }(\T_{ A })$ such that there exists an \emph{isomorphism}
							\[
								\oblv_{ \partial }(\g) \simeq \mathrm{FreeGradedLieAlgbd}(M \to \oblv_{ \partial }(\T_{ A }) );
							\]
						\item $\oblv_{ \partial }(\g)$ is a free graded $\oblv_{ \partial }(A)$-module isomorphic to $\bigoplus_{ n > 0 } A^{ k_{ n } }[-n]$ for $k_{ n } \in \N$ for each $n$.
					\end{enumerate}
				\end{prop}

				A key ingredient is the Chevalley-Eilenberg cohomology of dg-Lie algebroids, introduced in \cite[§3]{Nuiten:KoszulLieAlgbd}. It will be part of an adjunction that will give a \emph{deformation theory}, as defined in \cite[§12.3]{Lurie:SAG}. Let us recall the main properties.

				\begin{prop}
					\label{prop:classical_CE_deformation_theory}
					There exists an $\infty$-adjunction
					\[
						\adjunction{C^{ * }(-, A)}{\LieAlgbd_{ A }}{(\CAlg_{ -/A })^{ \op }}{\mathfrak{D}}
					\]
					such that
					\begin{enumerate}[label=(\roman*)]
						\item $C^{ * }(\g, A) \simeq \IntHom_{ U(\g) }(A, A)$ for each $\g \in \LieAlgbd_{ A }$;
						\item the Chevalley-Eilenberg cohomology complex of a free Lie algebroid over $A$ generated by $\rho\colon V \to \T_{ A }$ is the square-zero extension of $A$ by $V^{ \vee }[-1]$ given by the pullback diagram
							\begin{diag}
								C^*(\FreeLieAlgbd(V \to \T_A), A) \ar[r] \ar[d] \ar[dr, phantom, very near start, "\lrcorner"] & A \ar[d, "d_0"] \\
								A \ar[r, "d_{\rho^{\vee}}"] & A \oplus V^{\vee};
							\end{diag}			
						\item The composition of $\mathfrak{D}$ with the forgetful $\LieAlgbd_{ A } \to (\Mod_{ A })_{ -/\T_{ A } }$ acts by
							\[
								(B \to A) \mapsto (\T_{ A/B } \to \T_{ A });
							\]
							
						\item if $A$ is bounded, the adjunction above restricts to
							\[
								\adjunction{C^{ * }(-, A)}{\LieAlgbd_{ A }^{ \good }}{(\CAlg_{ -/A }^{ \smal })^{ \op }}{\mathfrak{D}}
							\]
							where the left adjoint is fully faithful;
						\item if $A$ is bounded, the adjunction $(C^{ * }(-, A), \mathfrak{D})$ is a deformation context.
					\end{enumerate}
				\end{prop}
				\begin{proof}
					The existence of the adjunction and the first three assertions are proved in \cite[§3]{Nuiten:KoszulLieAlgbd}. For the first point, see \cite[Remark 3.12]{Nuiten:KoszulLieAlgbd}.
					When $A$ is bounded, then \cite[Corollary 4.2]{Nuiten:KoszulLieAlgbd} proves that $C^{ * }(-, A)$ is fully faithful on good Lie algebroids over $A$. It is easily verified that the adjunction restricts to small and good objects using \cite[Proposition 12.3.2.2]{Lurie:SAG}. For example, using the second point together with \cref{defn:good_Lie_algebroids}, it is immediate to see that the Chevalley-Eilenberg cohomology complex of a good Lie algebroid over $A$ gives a small algebra as in \cref{defn:small_commutative_algebras_over_base} (see the beginning of the proof of \cite[Theorem 5.1]{Nuiten:KoszulLieAlgbd}).\\
					The last claim is readily verified, since $\LieAlgbd_{ A }^{ \good } \hookrightarrow \LieAlgbd_{ A }$ satisfies all desired conditions from \cite[Definition 12.3.1.1]{Lurie:SAG}. Moreover, by the second point, it is immediate that
					\[
						C^{ * }(\FreeLieAlgbd(A[n] \stackrel{0}{\to} \T_{ A }), A) \simeq A \oplus A[-n-1].
					\]
					It only remains to observe that it is a deformation theory as in \cite[Definition 12.3.3.2]{Lurie:SAG} because
					\[
						\LieAlgbd_{ A } \stackrel{\oblv}{\longrightarrow} (\Mod_{ A })_{ -/\T_{ A }} \stackrel{\fib}{\longrightarrow} \Mod_{ A } \stackrel{\oblv_{ A }}{\longrightarrow} \Sp
					\]
					is conservative and preserves sifted colimits by \cite[Proposition 2.4]{Nuiten:KoszulLieAlgbd}. Up to shift, this is the right functor to consider thanks to the computation
					\begin{gather*}
						\Map_{ \LieAlgbd_{ A } }(\mathfrak{D}(A \oplus A[n]), \h) \simeq \Map_{ \LieAlgbd_{ A } }(\mathfrak{D}(C^{ * }(\FreeLieAlgbd(A[-1-n] \stackrel{0}{\to} \T_{ A }), A)), \h)  \\
						\simeq \Map_{ \LieAlgbd_{ A } }(\FreeLieAlgbd(A[-1-n] \stackrel{0}{\to} \T_{ A }), \h) \simeq \Map_{ (\Mod_{ A })_{ -/\T_{ A } } }(A[-1-n] \stackrel{0}{\to} \T_{ A }, \oblv(\h))  \\
						\simeq \Map_{ A }(A[-1-n], \fib(\oblv(\h) \to \T_{ A })) \simeq \Map_{ \Sp }(\mathbb{S}[-1-n], \oblv_{ A }\fib(\oblv(\h) \to \T_{ A }) ) \\
						\simeq \Omega^{ \infty-n }(\oblv_{ A }\fib(\oblv(\h) \to \T_{ A })[1]),
					\end{gather*}
					which holds for each $n \geq 0$ and $\h \in \LieAlgbd_{ A }$.
				\end{proof}
				
				The main feature of (the essentially small $\infty$-category of) good Lie algebroids is that they can be used to understand the (large) whole $\infty$-category of Lie algebroids.

				\begin{prop}[{\cite[Proposition 5.4]{Nuiten:KoszulLieAlgbd}}]
					\label{prop:lie_algebroids_yoneda_into_good}
					The restricted Yoneda functor
					\[
						\phi\colon \LieAlgbd_{ A } \to \Fun(\LieAlgbd_{ A }^{ \good, \op }, \S), \qquad \g \mapsto \Map_{ \LieAlgbd_{ A } }(-, \g)
					\]
					is right adjoint and fully faithful. Its essential image is spanned by presheaves $F$ on good Lie algebroids such that
					\begin{enumerate}
						\item $F(0) \simeq *$;
						\item $F$ sends a pushout square as in \cref{defn:good_Lie_algebroids} to a pullback square of spaces.
					\end{enumerate}
				\end{prop}
				\begin{proof}
					The proof is given in the quoted reference. Let us just record that a left adjoint is given by the co-end
					\[
						\Fun(\LieAlgbd_{ A }^{ \good, \op }, \S) \longrightarrow \LieAlgbd_{ A }, \qquad F \mapsto \int^{ \h \in \LieAlgbd_{ A }^{ \good } } \h \odot F(\h),
					\]
					where $\odot$ is the tensoring of $\LieAlgbd_{ A }$ on $\S$ (by presentability).
				\end{proof}

				The main result of \cite{Nuiten:KoszulLieAlgbd} is that formal moduli problems under $\Spec A$ are equivalent to Lie algebroids over $A$, when $A$ is bounded.

				\begin{thm}[{\cite[Theorem 5.1]{Nuiten:KoszulLieAlgbd}}]
					\label{thm:koszul_duality_lie_algbd}
					Let $A \in \CAlg^{ \cn }_{ \C }$ be \emph{bounded}. There exists an adjoint equivalence
					\[
						\adjunction{\T_{A/-}}{\FMP_{ \Spec A/- }}{\LieAlgbd_{ A }}{\mathrm{MC}}.
					\]
				\end{thm}
				\begin{proof}
				    The statement is proved in the reference given above.
					Let us just record here the explicit formulas for the two functors.
					The right adjoint functor $\mathrm{MC}$ is obtained from the deformation theory $(C^{ * }(-, A), \mathfrak{D})$ following the general procedure described in \cite[§12.3.2]{Lurie:SAG}, where it is denoted by $\psi$.
					That is, it is obtained by the composition
					\[
						\LieAlgbd_{ A } \stackrel{\mathrm{Yoneda}}{\longrightarrow} \Fun(\LieAlgbd_{ A }^{ \op }, \S) \stackrel{- \circ \mathfrak{D}^{ \op }}{\longrightarrow} \Fun(\CAlg_{ -/A }, \S) \stackrel{\mathrm{restrict}}{\longrightarrow} \Fun(\CAlg_{ -/A }^{ \smal }, \S);
					\]
					it lands in $\FMP_{ \Spec A/- }$ by \cite[Proposition 12.3.2.1]{Lurie:SAG}.
					To completely spell out the definition, on objects it sends
					\[
						(\g \to \T_{ A }) \mapsto (\Map_{ \CAlg_{ -/A } }(\mathfrak{D}(-), \g)\colon \CAlg^{ \smal }_{ -/A } \to \S),
					\]
					see \cite[Remark 5.13]{Nuiten:KoszulLieAlgbd}. It is equivalent to the composition
					\[
						\LieAlgbd_{ A } \hookrightarrow \Fun(\LieAlgbd_{ A }^{ \good, \op }, \S) \stackrel{- \circ \mathfrak{D}^{ \op }}{\longrightarrow} \Fun(\CAlg^{ \smal }_{ -/A }, \S),
					\]
					where the first arrow is restricted Yoneda (see \cite[Proposition 5.4]{Nuiten:KoszulLieAlgbd}) and then we use the adjunction $(\mathfrak{D}^{ \op }, C^{ * }(-, A)^{ \op })$ on small/good objects from \cref{prop:classical_CE_deformation_theory}. Observe that the precomposition $- \circ \mathfrak{D}^{ \op }$ can be identified with right Kan extension along $C^{ * }(-,A)^{ \op }$, by adjunction.
					The left adjoint functor $\T_{ A/- }$ is then given by the composite 
					\[
						\FMP_{ \Spec A/- } \hookrightarrow \Fun(\CAlg^{ \smal }_{ -/A }, \S) \stackrel{- \circ C^{ * }(-, A)^{ \op }}{\longrightarrow} \Fun(\LieAlgbd_{ A }^{ \good, \op }, \S) \longrightarrow \LieAlgbd_{ A }, 
					\]
					where the last arrow is the left adjoint from \cref{prop:lie_algebroids_yoneda_into_good}. Such proposition implies that, starting from a formal moduli problem, we already land into the essential image of the restricted Yoneda (because $C^{ * }(-, A)$ preserves colimits). 
				\end{proof}

				The name $\T_{ A/-}$ for the right adjoint functor in \cref{thm:koszul_duality_lie_algbd} is not a coincidence, as explained in the following proposition. 

				\begin{prop}
					\label{prop:tangent_complex_small_fmp}
					Let $A \in \CAlg^{ \cn, \coh }$ be bounded. Interpreting each formal moduli problem under $A$ as a prestack, as in \cref{prop:fmp_in_geometry}, the composite functor
					\[
						\FMP_{ \Spec A/- } \stackrel{\T_{ A/- }}{\longrightarrow} \LieAlgbd_{ A } \stackrel{\oblv}{\longrightarrow} \Mod_{ A }
					\]
					gives a relative tangent complex, in the classical sense, to each formal moduli problem.
				\end{prop}
				\begin{proof}
					It suffices to verify that, for each formal moduli problem $\Spec A \to F$, the result of the composition above satisfies the universal property of the relative tangent complex (which is the $A$-linear dual of the relative cotangent complex).
					Observe that for each $n > 0$ we have
					\begin{gather*}
						\Map_{ \PreStk_{ \Spec A/- } }(\Spec(A \oplus A[n]), F) \simeq \\
						\simeq F(A \oplus A[n]) \simeq (F \circ C^{ * }(-, A))(\FreeLieAlgbd(A[-n-1] \stackrel{0}{\to} \T_{ A })).
					\end{gather*}
					We know by definition of the functor $\T_{ A/- }$ on formal moduli problems that
					\[
						F \circ C^{ * }(-, A) \simeq \Map_{ \LieAlgbd_{ A } }(-, \T_{ A/F }) \in \Fun(\LieAlgbd_{ A }^{ \good, \op }, \S);
					\]
					so that we have
					\begin{gather*}
						\Map_{ \PreStk_{ \Spec A/- } }(\Spec(A \oplus A[n]), F) \simeq \Map_{ \LieAlgbd_{ A } }(\FreeLieAlgbd(A[-n-1] \stackrel{0}{\to} \T_{ A }), \T_{ A/F }) \\
						\simeq \Map_{ (\Mod_{ A })_{ -/\T_{ A } } }(A[-n-1] \stackrel{0}{\to} \T_{ A }, \oblv(\T_{ A/F})) \simeq \Map_{ A }(A[-n], \fib(\T_{ A/F } \to \T_{ A })[1]) \simeq \\
						\simeq \Map_{ A }(A[-n], \cofib(\T_{ A/F } \to \T_{ A })).
					\end{gather*}
					To understand the classical relative tangent complex, we need to consider the mapping spaces
					\begin{gather*}
						\Map_{ \PreStk_{ \Spec A/ - / F } }(\Spec(A \oplus A[n]), \Spec A) \simeq \\
						\simeq \Map_{ \PreStk_{ \Spec A/- } }(\Spec(A \oplus A[n]), \Spec A) \times_{ \Map_{ \PreStk_{ \Spec A/- } }(\Spec(A \oplus A[n]), F) } \{d_{ 0 }\},
					\end{gather*}
					where $d_{ 0 }$ is the morphism $\Spec(A \oplus A[n]) \to \Spec A \to F$ induced by the $0$ derivation $0\colon \L_{ A } \to A[n]$.
					This space is equivalent to
					\[
						\Map_{ A }(A[-n], \T_{ A }) \times_{ \Map_{ A }(A[-n], \cofib(\T_{ A/F } \to \T_{ A })) } \{0\} \simeq \Map_{ A }(A[-n], \T_{ A } \times_{ \cofib(\T_{ A/F } \to \T_{ A }) } 0)
					\]
					so we conclude, by stability of $\Mod_{ A }$, that
					\[
						\Map_{ \PreStk_{ \Spec A/ - / F } }(\Spec(A \oplus A[n]), \Spec A) \simeq \Map_{ A }(A[-n], \T_{ A/F }). 
					\]
					This suffices to conclude that $\T_{ A/F } \in \Mod_{ A }$ satisfies the universal property of the relative tangent complex.
				\end{proof}

				\begin{remark}
					\label{remark:tangent_complex_fmp}
					We proved (the already known fact) that each formal moduli problem under $A$, seen as a particular nilpotent thickening $f\colon \Spec A \to F$, admits a relative tangent complex. We don't know if $F$ admits a global tangent complex $\T_{ F }$, although we think so. In fact, the $A$-module $\cofib(\T_{ A/F } \to \T_{ A })$ needs to be $f^{ * }(\T_{ F })$.
					We think that using an argument of nil-isomorphic proper descent (as is typical with ind-coherent sheaves) for the proper map $f\colon \Spec A \to F$ should yield a proof of the existence of $\T_{ F }$ as totalization of a certain cosimplicial module.
					We plan to pursue this direction in future work, as it seems a key step in proving the conjecture in \cref{remark:relation_with_thick_calaque_grivaux}.
				\end{remark}

				Let us finish this subsection about formal moduli problems with another comparison lemma.

				\begin{coroll}
					\label{lemma:tangent_complex_small_fmp_vs_ind_coherent_fmp}
					Let $A \in \CAlg^{ \cn, \coh }_{ \C }$ be bounded. 
					The following diagram commutes
					\begin{diag}
						\FMP_{\Spec A/-} \ar[rr, "\T_{A/-}"] \ar[d, "\mathrm{LKE}", hook] & & \Mod_A \ar[d, hook, "i"] \\
						\FMP^!_{\Spec A/-} \ar[rr, "T(A/-)"] & & \ProCoh_A,
					\end{diag}
					where $T(A/-)$ is the pro-coherent relative tangent complex coming from (the Serre dual version of) \cite[Chapter 8]{GR:DAG2}.
				\end{coroll}
				\begin{proof}
					It suffices to note that whenever a prestack admits a quasi-coherent tangent complex then it also admits an ind-coherent tangent complex.
				\end{proof}

			\subsection{Bounded case}
				\label{subsubsection:relation_lie_algebroids_bounded_case}

				In this subsection, we will explain how formal moduli problems and Lie algebroids are related to \cref{defn:beraldo_bounded_case}.
				Let us start with an easy observation.

				\begin{lemma}
					\label{lemma:final_lie_algebroid_de_rham}
					Let $A \in \CAlg^{ \cn }_{ \C }$ be bounded. The final object of $\FMP_{ \Spec A/- }$ is (the restriction to $\CAlg^{ \smal }_{ -/A }$) of the de Rham stack $(\Spec A)_{ \dR }$.
				\end{lemma}
				\begin{proof}
					Observe that, since $\Spec A$ is a derived affine scheme, its de Rham (pre)stack $(\Spec A)_{ \dR }$ admits deformation theory. 
					That is, we have
					\[
						(\Spec A)_{ \dR } \in \FMP^{ ! }_{ \Spec A/- },
					\]
					see \cref{lemma:ind_coh_fmp_from_small_coherent}. 
					It is the final ind-coherent formal moduli problem under $A$ because of the adjunction
					\[
						\adjunction{(-)^{ \red }}{\PreStk_{ \Spec A/- }}{\PreStk_{ (\Spec A)^{ \red }/- }}{(-)_{ \dR }},
					\]
					since each ind-coherent formal moduli problem $\Spec A \to F$ is a nil-isomorphism (i.e.\ $(\Spec A)^{ \red } \simeq F^{ \red }$).
					We conclude by observing that the restriction
					\[
						\FMP^{ ! }_{ \Spec A/- } \to \FMP_{ \Spec A/- }
					\]
					is a right adjoint (to left Kan extension, see \cref{prop:fmp_in_geometry}) and therefore preserves terminal objects.
				\end{proof}

				\begin{coroll}
					\label{coroll:tangent_space_de_rham_stack}
					Let $A$ be as above. There exists an equivalence of tangent complexes
					\[
						\T_{ \Spec A / (\Spec A)_{ \dR } } \simeq \T_{ A } \in \Mod_{ A }.
					\]
				\end{coroll}

				Observe that \cref{coroll:tangent_space_de_rham_stack} holds in a greater generality, by the simple fact that, for any $X$, $X_{ \dR }$ has ``trivial'' deformation theory.	Let us now state the main theorem that we want to prove in this subsection.
				\begin{thm*}
					Let $A \in \CAlg^{ \cn }_{ \C }$ be bounded and let $\g \in \LieAlgbd_{ A }$. There exists a symmetric monoidal fully faithful embedding
					\[
						\LMod_{ U(\g) } \hookrightarrow \ProCoh(\mathrm{MC}(\g)),
					\]
					where $\mathrm{MC}(\g)$ is the formal moduli problem under $A$ associated to $\g$.
					The essential image is given by those pro-coherent sheaves whose $!$-pullback to $\ProCoh_{ A }$ comes from an $A$-module.
				\end{thm*}
				
				Plugging in $\g = \T_{ A }$ we will have a symmetric monoidal equivalence between representations of $\T_{ A }$ and derived D-modules on $\Spec A$ from \cref{defn:beraldo_bounded_case}.
				Some ``historical'' remarks are due.

				\begin{remark}
					\label{remark:representation_lie_algebras}
					The theorem stated above, once tensor products are forgotten, is already known: see \cite[Theorem 7.1]{Nuiten:KoszulLieAlgbd}. We will follow its proof and refine it to obtain a \emph{symmetric monoidal} equivalence.
					This theorem, again without monoidal structure, was already known before for Lie algebras (which can be thought as Lie algebroids with null-homotopic anchor map). Let us explain, as an example, the case for Lie algebras over $\C$, where it can be seen as a corollary of \cite[Theorem 14.3.0.1]{Lurie:SAG}, which relates pointed small formal moduli problems (corresponding to Lie algebras over $\C$) to $\E_{ 1 }^{ \mathrm{aug} }$-formal moduli problems (corresponding to augmented $\E_{ 1 }$-algebras over $\C$). We have the commutative diagram
					\begin{diag}
						\Alg^{\mathrm{aug}}_{\C} \ar[d, "\mathrm{res}"] \ar[r, "\simeq"', "\mathrm{MC}^{(1)}"] & \FMP^{(1)} \ar[d, "\oblv"] \\
						\infcatname{LieAlg}_{\C} \ar[r, "\simeq"', "\mathrm{MC}"] & \FMP_{\C}.
					\end{diag}
					Passing to vertical left adjoints, and recalling that the left adjoint of $\mathrm{res}$ is the universal enveloping algebra, we obtain that
					\[
						\mathrm{MC}^{ (1) }(U(\g)) \simeq \mathrm{LKE}(\mathrm{MC}(\g))
					\]
					where the left Kan extension is along the forgetful functor $\CAlg^{ \smal }_{ \C } \to \Alg^{ \smal }_{ \C }$.
					By \cite[Theorem 14.6.0.1]{Lurie:SAG}, we have
					\[
						\LMod_{ U(\g) } \simeq \Qcoh^{ ! }_{ R }(\mathrm{MC}^{ (1) }(U(\g))) \simeq \Qcoh^{ ! }_{ R }(\mathrm{LKE}(\mathrm{MC}(\g))) \simeq \ProCoh(\mathrm{MC}(\g)),
					\]
					where the last equivalence follows from the fact that both $\Qcoh^{ ! }_{ R }$ and $\ProCoh$ are defined as right Kan extensions.
					Observe that tensor products are not mentioned; moreover, the fact that $\C$ is the base makes the statement simpler.
				\end{remark}
			
				Let us focus first on good dg-Lie algebroids over $A$ and identify their representations with a (full) subcategory of pro-coherent sheaves on the associated formal moduli problem. This has been already studied in \cite[Corollary 7.19]{Nuiten:KoszulLieAlgbd}. We will prove that the equivalence given there can be upgraded to a \emph{symmetric monoidal} one.
				Before doing so, we need to introduce some terminology.

				\begin{defn}
					\label{defn:procoh_with_underlying_qcoh}
					Let $f\colon \Spec A \to \Spec B$ a morphism in $\dAff_{ \C }$. Consider the pullback of $\infty$-categories
					\begin{diag}
						\ProCoh_B \times_{\ProCoh_A} \Mod_A \ar[r] \ar[d] \ar[dr, very near start, phantom, "\lrcorner"] & \ProCoh_B \ar[d, "f^!"] \\
						\Mod_A \ar[r, "i_A"] & \ProCoh_A
					\end{diag}
					and denote by $\ProCoh_{ 0 }(B/A)$ the result.
				\end{defn}

				\begin{remark}
					\label{remark:procoh_with_support_beraldo}
					Let us remark that \cref{defn:procoh_with_underlying_qcoh} can be given for any morphism of prestack $f\colon X \to Y$. Observe that if $f$ is a \emph{nil-isomorphism} (for example if it makes $Y$ a formal moduli problem under $X$) then, using \cref{coroll:serre_duality_prestack}, we obtain a symmetric monoidal equivalence
					\[
						\ProCoh_{ 0 }(Y/X) \simeq \IndCoh_{ 0 }(Y^{ \wedge }_{ X }),
					\]
					where the right-hand side is introduced in \cite[Definition 3.1.1]{Ber:DerivedDMod}.
				\end{remark}

				Let us record here some properties of $\ProCoh_{ 0 }(B/A)$. For a more thorough study of such object $\ProCoh_{ 0 }(Y/X)$ for a nil-isomorphism of stacks $f\colon X \to Y$ with bounded $X$, see \cite[§3]{Ber:DerivedDMod}.

				\begin{prop}
					\label{prop:properties_of_procoh_with_underlying_qcoh}
					Let $f\colon \Spec A \to \Spec B$ and $\ProCoh_{ 0 }(B/A)$ as in \cref{defn:procoh_with_underlying_qcoh}. The following assertions are true:
					\begin{enumerate}[label=(\roman*)]
						\item the pullback diagram can be considered in $\CAlg(\Pr^{ L, st}_{ \C })$, so that $\ProCoh_{ 0 }(B/A)$ is a stable presentable symmetric monoidal $\infty$-category;

						\item if $A$ is bounded then $\ProCoh_{ 0 }(B/A)$ is a colocalization of $\ProCoh_{ B }$;

						\item $!$-pullbacks preserve $\ProCoh_{ 0 }(-/A)$; that is, there exists a natural $\infty$-functor
							\[
								\ProCoh_{0}( -/A) \colon \CAlg^{ \cn, \coh }_{ -/A} \longrightarrow \CAlg(\Pr^{ L, st }_{ \C })_{ -/\Mod_{ A } }.
							\]

					\end{enumerate}
				\end{prop}
				\begin{proof}
					The first claim is true because $i_{ A }$ is fully faithful when $A$ is bounded by \cref{coroll:A_bounded_i_left_adjoint}. The second claim is true because both functors
					\[
						\CAlg(\Pr^{ L, st }_{ \C }) \to \Pr^{ L, st }_{ \C }, \qquad \Pr^{ L, st }_{ \C } \to \CAT
					\]
					preserve limits. 
					The third claim is proved by observing that we have a pullback of $\infty$-functors
					\[
						\ProCoh_{ 0 }(-/A) \simeq \ProCoh \times_{ \ProCoh_{ A } } \Mod_{ A }
					\]
					where we consider $\ProCoh\colon \CAlg^{ \cn, \coh }_{ -/A } \to \CAlg(\Pr^{ L, st }_{ \C })_{ -/\ProCoh_{ A } }$. This concretely means that $!$-pullbacks can be ``restricted'' to $\ProCoh_{ 0 }(-/A)$.
				\end{proof}
	
				The key statement is the following.

				\begin{prop}
					\label{prop:symm_monoidal_rep_good_lie_algebroids}
					Let $A \in \CAlg^{ \cn }_{ \C }$ be bounded and coherent and let $\g \in \LieAlgbd_{ A }$ be a \emph{good} dg-Lie algebroid over $A$.
					There exists a symmetric monoidal left adjoint functor
					\[
						\tilde{G}_{ \g }\colon \ProCoh_{ C^{ * }(\g, A) } \longrightarrow \LMod_{ U(\g) },
					\]
					which is natural in $\g$ ; that is, it is a morphism
					\[
						\ProCoh_{C^{ * }(-, A)} \longrightarrow \Rep_{ * }
					\]
					in $\Fun(\LieAlgbd_{ A }^{ \mathrm{good}, \op }, \CAlg(\Pr^{ L,st }))$ (see \cref{liealgbd:rep_functor} and \cref{prop:functoriality_!_pro_coherent_modules}). Moreover, its restriction to
					\[
						G_{ \g }\colon \ProCoh_{ 0 }(C^{ * }(\g, A)/A) \stackrel{\sim}{\longrightarrow} \Rep_{ \g }
					\]
					is a symmetric monoidal equivalence, again natural in $\g$.
				\end{prop}
				\begin{proof}
					\textbf{Construction of the equivalence for $\g$:}\\
					By \cref{lemma:symm_monoidal_functor_mod_procoh}, it suffices to give a symmetric monoidal exact functor from dually almost perfect $C^{ * }(\g, A)$-modules. We claim that 
					\[
						(A \otimes_{ C^{ * }(\g, A) } -)^{ \vee }\colon \APerf_{ C^{ * }(\g, A) }^{ \op } \longrightarrow \LMod_{ U(\g) }
					\]
					is a good choice (we consider $A$ as a $(U(\g), C^{ * }(\g, A))$-bimodule and we then take the $A$-linear dual from $\LMod_{ U(\g) }^{ \op }$). It is clear that it is exact and let us now prove that it is symmetric monoidal. Consider the augmentation morphism $C^{ * }(\g, A) \to A$: since $\g$ is good then $C^{ * }(\g, A)$ is connective and by \cref{prop:pullback_preserves_almost_perfect} we can factor the above functor as
					\[
						\APerf_{ C^{ * }(\g, A) }^{ \op } \stackrel{(A \otimes_{ C^{ * }(\g, A) } -)^{ \op }}{\longrightarrow} \left( \LMod_{ U(\g) } \times_{ \Mod_{ A } } \APerf_{ A } \right)^{ \op } \stackrel{(-)^{ \vee }}{\longrightarrow} \LMod_{ U(\g) }
					\]
					so it is immediate that it is lax symmetric monoidal (more precisely, the first functor is strong monoidal and the second one is lax symmetric monoidal).
					To check that the restricted duality functor above is also strong monoidal we can then compose it with the conservative and strong monoidal forgetful functor $\LMod_{ U(\g) } \to \Mod_{ A }$ (see \cref{liealgbd:forgetful_conservative_strong_monoidal}). Since the $U(\g)$-module structure is now irrelevant, we can now focus on the $A$-linear duality functor
					\[
						\APerf_{ A }^{ \op } \stackrel{(-)^{ \vee }}{\longrightarrow} \Mod_{ A }.
					\]
					Since $A$ is bounded, it is strong monoidal by \cref{prop:almost_perfect_duality_bounded}.\\
					We claim that $\tilde{G}_{ \g }$ admits a left adjoint $\tilde{\mu}_{ \g }$, which must then necessarily be the $\Ind$-extension of a functor $\LMod_{ U(\g)  }^{ \mathrm{perf} } \to \Coh_{ C^{ * }(\g, A) }^{ \op }$. Let $K \in \Coh_{ C^{ * }(\g, A) }$, $P \in \LMod_{ U(\g) }^{ \mathrm{perf} }$ and consider
					\[
						\Map_{ \g }(P, \tilde{G}_{ \g }(K)) \simeq \Map_{ \g }(P, (A \otimes_{ C^{ * }(\g, A) } K)^{ \vee }) \simeq \Map_{ \g }(A \otimes_{ C^{ * }(\g, A) } K, P^{ \vee }),
					\]
					where we used the fact that $(-)^{ \vee }\colon \LMod_{ U(\g) }^{ \op } \to \LMod_{ U(\g) }$ is the right adjoint of its opposite. Since $A \otimes_{ C^{ * }(\g, A) } -$ is left adjoint to $C^{ * }(\g, -)$, we can continue writing
					\[
						\Map_{ \g }(P, \tilde{G}_{ \g }(K)) \simeq \Map_{ \Mod_{ \C^{ * }(\g, A) } }(K, C^{ * }(\g, P^{ \vee })).
					\]
					Observing that $P^{ \vee } \in \LMod_{ U(\g) }$ is a retract of a finite limit of shifted copied of $U(\g)^{\vee }$, we deduce that also $C^{ * }(\g, P^{ \vee })$ is a retract of a finite limit of shifted copies of $C^{ * }(\g, U(\g)^{ \vee })$. Since the latter is equivalent to $A \in \Coh_{ C^{ * }(\g, A) }$ (since $\g$ is good, see \cite[§7.2]{Nuiten:KoszulLieAlgbd}) we obtain that $C^{ * }(\g, P^{ \vee })$ is a coherent $C^{ * }(\g, A)$-module.
					This proves that we have the adjunction
					\[
						\adjunction{\tilde{\mu}_{ \g }}{\LMod_{ U(\g) }}{\ProCoh_{ C^{ * }(\g, A) }}{\tilde{G}_{ \g }}.
					\]
					where $\tilde{\mu}_{ \g }$ is the $\Ind$-extension of $C^{ * }(\g, (-)^{ \vee })\colon \LMod_{ U(\g) }^{ \mathrm{perf} } \to \Coh_{ C^{ * }(\g, A) }^{ \op }$.
					It is proved in \cite[Proposition 7.16, Lemma 7.18]{Nuiten:KoszulLieAlgbd} that $\tilde{\mu}_{ \g }$ is fully faithful and its essential image is $\ProCoh_{ 0 }(C^{ * }(\g, A)/A)$.
					The assumptions that $A$ is bounded and $\g$ is a good Lie algebroid are again crucial and allows us to use the strong form of PBW theorem as well as making $A$ a coherent $C^{ * }(\g, A)$-module.
					Since $A$ is bounded, $\ProCoh_{ 0 }(C^{ * }(\g, A)/A) \hookrightarrow \ProCoh_{ C^{ * }(\g, A) }$ is a colocalisation and therefore the composition
					\[
						G_{ \g }\colon \ProCoh_{ 0 }(C^{ * }(\g, A)/A) \hookrightarrow \ProCoh_{ C^{ * }(\g, A) } \stackrel{\tilde{G}_{ \g }}{\longrightarrow} \LMod_{ U(\g) }
					\]
					is the sought for (symmetric monoidal and left adjoint) equivalence.\\
					\textbf{Naturality in $\g$:}\\
					It suffices now to observe that the naturality of $G_{ \g }$ descends from the naturality of $\tilde{G}_{ \g }$. This is implied by the naturality of its restriction to dually almost perfect modules. This holds by \cref{lemma:symm_monoidal_functor_mod_procoh}, which is functorial in $A$, observing that $*$-pullback on almost perfect modules is exchanged with $!$-pullback on pro-coherent ones by \cref{prop:pro_coh_!_pullback_compatibilities}.
					Let us now consider the functor
					\[
						(A \otimes_{ C^{ * }(\g, A) } -)^{ \vee }\colon \APerf^{ \op }_{ C^{ * }(\g, A) } \to \LMod_{ U(\g) }
					\]
					Given a morphism of good Lie algebroids over $A$ $f\colon \g \to \h$, it is easy to see that the following diagram commutes
					\begin{diag}
						\APerf^{\op}_{C^*(\h, A)} \ar[r, "(A \otimes_{C^*(\h{,} A)} -)^{\vee}"] \ar[d, "C^*(\g{,}A) \otimes_{C^*(\h{,}A)} -"] & \LMod_{U(\h)} \ar[d, "\oblv_{U(f)}"] \\
						\APerf^{\op}_{C^*(\g, A)} \ar[r, "(A \otimes_{C^*(\g{,} A)} -)^{\vee}"] & \LMod_{U(\g)}.
					\end{diag}
					It is also possible to study the naturality of $\tilde{\mu}_{ \g }$ as in \cite[Proposition 7.16]{Nuiten:KoszulLieAlgbd} and simply pass to right adjoints (and then restrict to almost perfect modules) to obtain this result.
				\end{proof}
				
				We can now extend the result from good Lie algebroids to all Lie algebroids.
				Before doing so, let us recall how sheaves on formal moduli problems are defined.

				\begin{defn}
					\label{defn:sheaves_small_fmp}
					Let $A \in \CAlg^{\cn, \coh}_{ \C }$ be bounded. The right Kan extensions
					\begin{diag}
						\CAlg^{\smal}_{/A} \ar[r, "\Mod"] \ar[d, "\Spf^{\op}"', hook] & \CAlg(\Pr^L) &  \CAlg^{\smal}_{/A} \ar[d, "\Spf^{\op}"', hook] \ar[r, "\Ind(\Coh^{\op})"] & \CAlg(\Pr^L)  \\
							\left(\infcatname{FMP}_{\Spec A/ -}\right)^{\op} \ar[ur, dashed, "\Qcoh^*"'] & &  \left(\infcatname{FMP}_{\Spec A/ -}\right)^{\op} \ar[ur, dashed, "\ProCoh"'] &  
					\end{diag}
					define quasi-coherent sheaves and pro-coherent sheaves on formal moduli problems under $\Spec A$. 
				\end{defn}

				\begin{lemma}
					\label{lemma:sheaves_small_fmp_and_geometry}
					Let $A$ be as above. The following diagrams commute
					\begin{diag}
						(\FMP_{\Spec A/-})^{\op} \ar[r, "\Qcoh^{*}"] \ar[d, "\mathrm{LKE}^{\op}"', hook] & \CAlg(\Pr^L)   & (\FMP_{\Spec A/-})^{\op} \ar[r, "\ProCoh"] \ar[d, "\mathrm{LKE}^{\op}"', hook] & \CAlg(\Pr^L) \\
						(\PreStk_{\Spec A/-})^{\op} \ar[ur, "\Qcoh"'] &  & (\PreStk_{\Spec A/-})^{\op} \ar[ur, "\ProCoh"'] &
					\end{diag}
					where the diagonal arrows are the standard functors of quasicoherent and pro-coherent sheaves on prestacks.
				\end{lemma}
				\begin{proof}
					We will only prove the first assertion, since the second one is similar.
					Consider the commutative diagram
					\begin{diag}
						\CAlg^{\smal}_{-/A} \ar[r, "i", hook] \ar[d, hook] & \CAlg^{\cn}_{-/A} \ar[d, hook] \\
						\Fun(\CAlg^{\smal}_{-/A}, \S)^{\op} \ar[r, "\P(i^{\op})^{\op}", hook] & \Fun(\CAlg_{-/A}^{\cn}, \S)^{\op}
					\end{diag}
					where the vertical arrows are (opposite) Yoneda embeddings. 
					The key observation is that the bottom arrow $\P(i^{ \op })^{ \op }$ is the right adjoint of the (opposite of the) restriction functor along $i$
					\[
						\Fun(\CAlg_{ -/A }^{ \cn }, \S)^{ \op } \longrightarrow \Fun(\CAlg_{ -/A }^{ \smal }, \S)^{ \op },
					\]
					so that it corresponds to the right Kan extension along $i^{ \op }$ and it is fully faithful.
					Observing that, by definition, we can identify the opposite of the left Kan extension along $i$ (present in the statement) with the right Kan extension along $i^{ \op}$, we conclude by using transitivity of Kan extensions.
				\end{proof}

				\begin{thm}
					\label{thm:symm_monoidal_representation_lie_algbd}
					Let $A \in \CAlg_{ \C }^{ \cn }$ be a \emph{bounded} coherent commutative algebra, $\g \in \LieAlgbd_{ A }$. There exists a symmetric monoidal equivalence of $\infty$-categories 
					\[
						\ProCoh_{ 0 }(\MC(\g)/A) \stackrel{\sim}{\longrightarrow} \LMod_{ U(\g) },
					\]
					for $\MC(\g) \in \FMP_{ \Spec A/ - }$ the formal moduli problem corresponding to $\g$.
				\end{thm}
				\begin{proof}
					The idea is that it suffices to right Kan extend the result for good Lie algebroids $\h$, where $X_{ \h } = \Spec C^{ * }(\h, A)$. A crucial point is that $\CAlg(\Pr^{ L }) \to \Pr^{ L }$ commutes with limits. 
					We know that quasi-coherent (and ind-coherent) sheaves on formal moduli problems (seen as prestacks using \cref{prop:fmp_in_geometry}) can be computed as right Kan extensions of their restrictions to small algebras, see \cref{lemma:sheaves_small_fmp_and_geometry}. Then, the same is true for $\ProCoh_{ 0 }(-/A)$ with $!$-pullbacks (see \cref{prop:properties_of_procoh_with_underlying_qcoh}). Consider now the limit-preserving functor 
					\[
							\Rep_{ * }\colon \LieAlgbd_{ A }^{ \op } \longrightarrow \CAlg(\Pr^{ L })_{- /\Mod_{ A } }, \qquad \left( \g \to \T_{ A } \right) \mapsto \left( \left( \LMod_{ U(\g) }, \otimes_{ A } \right) \to (\Mod_{ A }, \otimes_{ A }) \right)
					\]
					sending $\g \to \h$ to the symmetric monoidal forgetful $\LMod_{ U(\h) } \to \LMod_{ U(\g) }$ (see \cref{liealgbd:rep_commutes_limits}). 
					Using the limit-preserving embedding
					\[
						\LieAlgbd_{ A } \hookrightarrow \P(\LieAlgbd_{ A }^{ \good })
					\]
					coming from \cref{prop:lie_algebroids_yoneda_into_good}, we obtain that the functor $\Rep$ is right-Kan extended from its restriction to good Lie algebroids. That is, we have a right Kan extension diagram
					\begin{diag}
						\LieAlgbd_{ A }^{ \good, \op} \ar[r, "\Rep"] \ar[d, hook] & \CAlg(\Pr^L)_{-/\Mod_A} \\
						\Fun(\LieAlgbd_A^{\good, \op}, \S)^{\op} \ar[ur, dashed, "\Rep"] 
					\end{diag}
					The fact that $(C^{ * }(-,A), \mathcal{D})$ is a deformation context (see \cref{prop:classical_CE_deformation_theory}) gives us, by \cite[Proposition 12.3.2.1]{Lurie:SAG}, the commutative diagram
					\begin{diag}
						\LieAlgbd_A^{\mathrm{good}, \op} \ar[r, hook, "C^*(-{,}A)^{\op}"] \ar[d, hook] & \CAlg_{/A}^{\smal} \ar[d, hook, "\Spf^{\op}"] \\
						\Fun(\LieAlgbd_A^{\good, \op}, \S)^{\op} \ar[r, "\MC^{\op}", hook] & \Fun(\CAlg^{\smal}_{/A}, \S)^{\op}
					\end{diag}
					where the bottom arrow is (the opposite of) precomposition with $\mathfrak{D}^{ \op }$, and induces the equivalence between Lie algebroids and formal moduli problems as in \cref{thm:koszul_duality_lie_algbd}.
					By adjunction, we can identify $\MC^{ \op }$ with right Kan extension along $C^{ * }(-, A)^{\op }$. By transitivity of right Kan extension we deduce that the diagrams
					\[
					\begin{tikzcd}
					\LieAlgbd_A^{\good, \op} \ar[r, "\Mod_{C^*(-{,}A)}"] \ar[d, hook] & \CAlg(\Pr^L) &  \LieAlgbd_A^{\good, \op} \ar[d, hook] \ar[r, "\Ind(\Coh^{\op}_{ C^*(-{,}A)})"] & \CAlg(\Pr^L)  \\
						\Fun(\LieAlgbd_A^{\good, \op}, \S)^{\op} \ar[ur, "\Qcoh^*(\MC(-))"'] & &  \Fun(\LieAlgbd_A^{\good, \op}, \S)^{\op} \ar[ur, "\ProCoh(\MC(-))"'] &  
					\end{tikzcd}
					\]
					are also right Kan extension. Thus, the functor
					\[
						\LieAlgbd_{ A }^{ \op } \longrightarrow \CAlg(\Pr^{ L }), \qquad \g \mapsto \ProCoh_{ 0 }(\MC(\g)/A)
					\]
					is also right Kan extended from its restriction to good Lie algebroids. 
					We conclude then that the right Kan extension of \cref{prop:symm_monoidal_rep_good_lie_algebroids} gives the symmetric monoidal equivalence we're after. 
				\end{proof}
				\begin{remark}
					Let us remark that, the equivalence of underlying $\infty$-categories in the statement of \cref{thm:symm_monoidal_representation_lie_algbd}, has been proved in \cite[Theorem 7.1]{Nuiten:KoszulLieAlgbd}. What we did here is to enhance it to a symmetric monoidal one.
				\end{remark}
				
				This theorem tells us that, in particular, in the bounded case, derived D-modules as defined by Beraldo (see \cref{defn:beraldo_bounded_case}) are equivalent, as symmetric monoidal $\infty$-category, to left representations of the tangent dg-Lie algebroid $\T_{ X }$.

				\begin{coroll}
					\label{coroll:nuiten_is_beraldo_monoidal}
					Let $A \in \CAlg^{ \cn, \coh }_{ \C }$ be bounded and define $X = \Spec A$. There is a symmetric monoidal equivalence
					\[
						\infcatname{D}^{ \mathrm{der} }(X) \coloneqq \IndCoh(X_{ \dR }) \times_{ \IndCoh(X) } \Qcoh(X) \simeq \ProCoh_{ 0 }(X_{ \dR }/A) \stackrel{\sim}{\longrightarrow} \LMod_{ U(\T_{ A }) }.
					\]
				\end{coroll}
				\begin{proof}
					Use \cref{coroll:serre_duality_prestack} for $X_{ \dR }$ and \cref{thm:symm_monoidal_representation_lie_algbd} for $\g = \T_{ A }$.
				\end{proof}

				\begin{coroll}
					\label{coroll:representations_liealgbd_self_dual}
					With the same notations as above, there exists an equivalence
					\[
						\LMod_{ U(\g) } \stackrel{\sim}{\longrightarrow} \LMod_{ U(\g) }^{ \vee } \simeq \RMod_{ U(\g) },
					\]
					in $\Pr^{ L }$. That is, the $\infty$-category $\LMod_{ U(\g) }$ is self-dual in $\Pr^{L}$.
				\end{coroll}
				\begin{proof}
					It suffices to observe that $\ProCoh_{ 0 }(\MC(\g)/A) \simeq \IndCoh_{ 0 }( (\MC(\g))_{ \Spec A }^{ \wedge })$ as observed in \cref{remark:procoh_with_support_beraldo} and then conclude, using the boundedness of $A$, by \cite[Proposition 3.2.4]{Ber:DerivedDMod}.
					The second equivalence holds because $\LMod_{ U(\g) }$ is dualizable, being compactly generated, and we have the standard computation of duals in $\Pr^{ L, st }$
					\[
						\Fun^{ L }_{ \ex }(\LMod_{ U(\g) }, \Sp) \simeq \RMod_{ U(\g) }
					\]
					by the Eilenberg-Watts theorem (see \cref{thm:eilenberg_watts} and see also \cite[Proposition 7.2.4.3]{Lurie:HA}).
				\end{proof}

				\begin{remark}
					\label{remark:hopfyness_universal_enveloping_algebra}
					If $\g \simeq \mathrm{diag}(\h)$ where $\h \in \infcatname{LieAlg}_{ A }$, then we know that $U(\h)$ is a cocommutative Hopf algebra and therefore the $\infty$-category $\LMod_{ U(\h) }$ is symmetric monoidal (by the cocommutativity) and self-dual (by the Hopf property).
					The two conclusions are true for any Lie algebroid, even if $U(\g)$ is not an Hopf algebra.
				\end{remark}

				Let us finish this part with a conjecture, that we hope to prove in the future.

				\begin{conj}
					\label{conj:unbounded_version_lie_algbd_rep}
					Let $A \in \CAlg^{ \cn, \lfp }_{ \C }$, possibly unbounded, and let $\g \in \LieAlgbd_{ A }$. There exists a symmetric monoidal equivalence
					\[
						\IndCoh_{ 0 }(\MC(\g)/A) \stackrel{\sim}{\longrightarrow} \LMod_{ U(\g) },
					\]
					where the left-hand side is defined in \cite[Definition 4.1.6]{Ber:DerivedDMod} and $\MC(\g)$ should be the formal moduli problem under $\Spec A$ associated to $\g$ (even though this is not an equivalence). 
				\end{conj}

%% file: construction.tex
	\section{Construction of the functor}
			\label{subsection:chevalley_eilenberg}

			The goal of this section is to give a rigorous construction of a Hodge-filtered version of the Chevalley-Eilenberg functor for dg-Lie algebroids. 
			As a warm up let us review older related constructions.

			\begin{description}
				\item[Classical Lie algebra cohomology:] A reference is \cite[Chapter 7.7]{Weibel:HomologicalAlgebra}.  Let $\g$ be a classical (i.e.\ 1-categorical) Lie algebra over $\C$. Its Chevalley-Eilenberg cohomology is defined to be the right derived functor of the ``invariants functor'' 
				\[
					\LMod^{ \heartsuit }_{ U(\g) } \ni M \mapsto M^{ \g } \simeq \Hom_{ \g }(\C, M) \in \Mod_{ \Z }^{ \heartsuit }.
				\]
				The way we can compute it explicitly is by using the Spencer resolution, as left $U(\g)$-module, of the trivial representation $\C$ given by
				\[
					V^{ * }(\g) \coloneqq U(\g) \otimes_{ \C } \Sym^{ \gr }_{ \C }(\g[1]) \stackrel{\sim}{\longrightarrow} \C
				\]
				equipped with a Koszul-like differential ($\g[1]$ has weight $-1$). This gives us
				\[
					\CE(\g, M) = \mathbb{R}\IntHom_{ \g }(\C, M) \simeq \Sym^{ \gr }(\g^{ \vee }[-1]) \otimes^{ \gr } M
				\]
				(where the last equivalence only holds if $\g$ is finite dimensional over $k$) equipped with the well-known action differential (here $\g^{ \vee }[-1]$ has weight $1$).

				\item[Cohomology of differential-graded Lie algebras:] 	A reference is \cite[Chapter 2.2]{Lurie:DAGX}. We have an analogue history in the differential graded context: let now $\g$ denote a dg-Lie algebra (over $ \C $ as usual). The same formulas as above, interpreted in the dg context, make sense and everything goes through.
					The ``zen'' of differential-graded objects is to make all constructions invariant under quasi-isomorphisms, among which Chevalley-Eilenberg cohomology as well. A slick way of defining its (homology and) cohomology is observing that the Spencer resolution, analogue to $V^{ * }(\g)$, can be obtained as a universal enveloping dg-algebra of another dg-Lie algebra.
					Namely, we consider the chain complex $Cn(\g)$ defined as the mapping cone of $\id_{ \g }$ with Lie bracket given (on homogeneous elements) by 
					\[
						[x+\epsilon y, x' + \epsilon y'] \coloneqq [x, x'] + \epsilon ([y, x'] + (-1)^{ |x| }[x, y'])
					\]
					with the shorthand $x + \epsilon y \in \g \oplus \g[1]$. Let us immediately observe that the unique map $Cn(\g) \to 0$ is a quasi-isomorphism of dg-Lie algebra.
					Passing to universal enveloping dg-algebras we obtain a quasi-isomorphism 
					\[
						U(Cn(\g)) \stackrel{\sim}{\to} U(0) = \C.
					\]
					The canonical map (of dg-Lie algebras) $\g \to Cn(\g)$ makes $U(Cn(\g))$ into a $U(\g)$-dg-bimodule. What is left now to observe is that, using the (dg analogue of) PBW theorem, we have $U(Cn(\g)) \simeq U(\g) \otimes \Sym^{ \gr }(\g[1])$ as a graded vector space (i.e.\ without differential), which is exactly the same formula as $V^{ *}(\g)$ above ($\g[1]$ is of weight $-1$).

					Finally, the Chevalley-Eilenberg complex is then defined as 
					\[
						CE(\g, M) \coloneqq \R\IntHom_{ \g }(\C, M) \simeq \IntHom_{ \g }(U(Cn(\g)), M).
					\]
					If $\g$ is dualizable then we obtain the graded object $\Sym^{ \gr }(\g^{ \vee }[-1]) \otimes^{ \gr }_{ \C } M$ equipped with the usual action differential ($\g^{ \vee }[-1]$ has weight $1$).

				\item[Cohomology of dg-Lie algebroids:] A reference is \cite[§3]{Nuiten:KoszulLieAlgbd}. We can build a similar complex for dg-Lie algebroids, and again everything can be made sense of in the context of dg-Lie algebroids. Consider a cofibrant (connective) cdga $A$ and an $A$-cofibrant dg-Lie algebroid $\g \in \LieAlgbd_{ A }^{ dg }$. We need a resolution of $A$ as a left $U(\g)$-dg-module; as before we can consider the graded vector space $K(\g) \coloneqq U(\g) \otimes_{ A } \Sym_{ A }^{ \gr }(\g[1])$, with $\g[1]$ in weight $-1$, equipped with an appropriate differential encoding the multiplication on $U(\g)$ and the Lie bracket on $\g$, which makes it a left $U(\g)$-dg-module. We then define
			\[
				CE(\g, M) = \IntHom_{ \g }(K(\g), M) \simeq \R\IntHom_{ \g }(A, M).
			\]
			Observe that, when $M=A$, this gives exactly what we denoted $C^{ * }(\g, A)$ in \cref{prop:classical_CE_deformation_theory}.
			If $\g$'s underlying dg-$A$-module is perfect (i.e.\ dualizable) then we obtain $\Sym^{ \gr }_{ A }(\g^{ \vee }[-1]) \otimes_{ A }^{ \gr } M$ with the action differential.

				\item[Complete Hodge-filtered cohomology of dg-Lie algebras:] 	A reference is \cite{Pavia:MixedCE}. A refinement of the standard Chevalley-Eilenberg complex for dg-Lie algebras is its (complete) filtered version, or, equivalently, its mixed-graded version. Let now $\g$ be a dg-Lie algebra.
					The idea is to work internally to mixed-graded objects from the beginning: we can replace $Cn(\g)$ with a mixed graded version. In layman terms, we simply replace $\g \oplus \g[1]$ (with its differential) with $\g(0) \oplus \g[1](-1)$ where we give weight $-1$ to $\g[1]$ and we use identity as mixed differential (i.e.\ we ``split'' the internal cohomological differential of $\g$ from the external mapping cone one). This is then a Lie algebra object in $\Mod^{ \egr }_{ \C }$ and we can consider its (mixed graded) universal enveloping algebra. In this case, the latter won't be equivalent to $\C$ anymore, since we are keeping track of the filtration.
					The definition is now simply, for $M \in \LMod_{ U(\g) }$:
					\[
						CE^{ \egr }(\g, M) \coloneqq \IntHom^{ \egr }_{ U(\g) }(U(Cn^{ \egr }(\g)), M(0)).
					\]
			\end{description}
			Our goal in this section is to generalize the latter example above: we want the construction to work not only on dg-Lie algebras but on all dg-Lie algebroids on arbitrary bases. Let $A$ be our base. One very similar approach is possible: give to $K(\g)$ a non-trivial mixed-graded structure (making it non isomorphic to $A$) and then take mixed-graded hom in $\LMod^{ \egr }_{ U(\g) }$. Given the concrete nature of dg-Lie algebroids, this construction needs to be performed first as $1$-functor and then we can prove that quasi-isomorphisms are preserved (so it goes to the $\infty$-level). This is different and slightly more unpleasant than the Lie algebra case, where dg-operads help us.	What we will do is instead a bit different in this case, and it will provide an alternative equivalent approach to the dg-Lie algebra case as well. 
			
			Let $\g \in \LieAlgbd_{ A }$ be a dg-Lie algebroid and consider the PBW filtration on its universal enveloping dg-algebra $U(\g)^{ \PBW }$, see \cref{subsubsection:universal_enveloping_algebra}. It is an algebra in positively filtered $(A,A)$-bimodules.
			Recall from the proof of \cref{lemma:filtered_modules_degree_0} the adjunction
			\[
				\adjunction{(-)_{ \infty }}{\infcatname{C}}{\infcatname{C}^{ \Fil, \geq 0 }}{\langle 0, - \rangle}
			\]
			where both $\infty$-functors are symmetric monoidal (whenever $\infcatname{C}$ has a symmetric monoidal structure).
			For $\infcatname{C} = (A,A)-\BiMod$, the unit for the (algebra level) adjunction gives a canonical map of (positively) filtered algebras
			\[
				U(\g)^{ \PBW } \to \langle 0, U(\g) \rangle.
			\]
			Therefore, the $(U(\g), A)$-bimodule structure on $A$ upgrades to a $(U(\g)^{ \PBW }, \langle 0, A \rangle)$-bimodule structure on $\langle 0, A \rangle$ (it comes from the $\g$-action on $A$, see \cref{liealgbd:dg_representation_on_A}). Let us start with a preliminary categorical observation, after recalling that 
			\[
				\LMod_{ U(\g)^{ \PBW } }^{ \Fil } \coloneqq \LMod_{ U(\g)^{ \PBW } }(\Mod_{ \C }^{ \Fil })
			\]
			is a symmetric monoidal $\infty$-category when endowed with the relative tensor product over $\langle 0, A \rangle$ (see \cref{liealgbd:injective_structure_filtered_representation}). Complete filtered left $U(\g)^{ \PBW }$-modules, denoted by $\LMod_{ U(\g)^{ \PBW } }^{ \cpl }$, then inherit a symmetric monoidal structure by completion of the previously described one, since $U(\g)^{ \PBW }$ is complete. 
			\begin{prop}
				\label{prop:enrichment_lmod_ug}
				The canonical $\Mod^{ \cpl }_{ \C }$-tensoring of $\LMod_{ U(\g)^{ \PBW } }^{ \cpl }$, denoted by 
				\[
					- \odot - \colon \Mod^{ \cpl }_{ \C } \times \LMod_{ U(\g)^{ \PBW } }^{ \cpl } \to \LMod_{ U(\g)^{ \PBW } }^{ \cpl }
				\]
				is induced by the symmetric monoidal functor 
				\[
					\langle 0, A \rangle \otimes^{ \cpl }_{ \C } -\colon \Mod^{ \cpl }_{ \C } \to \LMod_{ U(\g)^{ \PBW } }^{ \cpl }.
				\]
				That is, $\LMod_{ U(\g) }^{ \cpl }$ is a commutative $\Mod^{ \cpl }_{ \C }$-algebra in $\Pr^{ L, st }$ (meaning that the $\Mod^{ \cpl }_{ \C }$-tensoring and the internal tensor product are naturally compatible).
			\end{prop}
			\begin{proof}
				This enrichment is a particular case of \cref{example:canonical_enrichment_lmod}, since $U(\g)^{ \PBW }$ is an associative algebra in complete filtered $\C$-modules. Since $\LMod_{ U(\g)^{ \PBW } }^{ \cpl }$ can be endowed with a symmetric monoidal structure (see \cref{liealgbd:injective_structure_complete_representation}), what we want to prove here is the compatibility between the enrichment and this monoidal structure. It is immediately seen using the explicit definitions. \qedhere
			\end{proof}

			\begin{remark}
				Since all $\infty$-categories in the above proposition are presentable and stable, by \cref{lemma:internal_hom_categories} there are external hom functors 
				\begin{gather*}
					\IntHom^{ \cpl }_{ U(\g)^{ \PBW } }(-, -)\colon \left( \LMod_{ U(\g)^{ \PBW } }^{ \cpl }\right) ^{ \op } \times \LMod_{ U(\g)^{ \PBW } }^{ \cpl } \to \Mod^{ \cpl }_{ \C },\\
					(M, N) \mapsto \IntHom^{ \cpl }_{ U(\g)^{ \PBW } }(M, N)
				\end{gather*}
				satisfying the well-known universal property
				\[
				\Map_{ \Mod^{ \cpl }_{ \C } }(E, \IntHom^{ \cpl }_{ U(\g)^{ \PBW } }(M, N)) \simeq \Map_{ \LMod^{ \cpl }_{ U(\g)^{ \PBW } } }(E \odot M, N).\]
			\end{remark}
			
			We can finally introduce our main object of interest. See \cref{remark:quick_construction_CE} for an equivalent and possibly faster construction.

			\begin{defn}
				\label{defn:mixed_graded_CE}

				The \emph{complete Hodge-filtered Chevalley-Eilenberg complex of $\g$} is the $\infty$-functor defined as 
				\begin{gather*}
					\CE^{ \cpl }(\g, -)\colon \LMod_{ U(\g) } \hookrightarrow \LMod^{ \cpl }_{ \langle 0, U(\g) \rangle } \to \LMod^{ \cpl }_{ U(\g)^{ \PBW } } \to \Mod^{ \cpl }_{ \C }, \\
					M \longmapsto \IntHom_{ U(\g)^{ \PBW } }^{ \cpl }\left(\langle 0, A \rangle, \langle 0, M \rangle \right),
				\end{gather*}
				where we used the $\Mod_{ \C }^{ \cpl }$-enrichment of $\LMod^{ \cpl }_{ U(\g)^{ \PBW } }$.
			\end{defn}

			\begin{remark}
				Let us observe that a sketch of a similar construction (of the Chevalley-Eilenberg mixed graded \emph{algebra}, i.e.\ with coefficients $A$) appeared already in \cite[Proposition A.3]{CalaqueGrivaux:FormalModuli}. It is defined in a slightly different way but it gives an equivalent object.
			\end{remark}

			\begin{prop}
				\label{prop:CE_lax_monoidal_functor}
				There exists a canonical (dashed) lift
				\begin{diag}
					& \Mod^{\cpl}_{\CE^{\cpl}(\g, A)} \ar[d] \\
					\LMod_{U(\g)} \ar[ur, dashed] \ar[r, "\CE^{\cpl}(\g{,} -)"] & \Mod^{\cpl}_{\C}
				\end{diag}
				which, from now on, we will denote $\CE^{ \cpl }(\g, -)\colon \LMod_{ U(\g) } \to \Mod^{ \cpl }_{ \CE^{ \cpl }(\g, A) }$.
			\end{prop}
			\begin{proof}
				The first observation, required just to make sense of the fact that $\CE^{ \cpl }(\g, A)$ is a commutative algebra (in complete filtered modules), is that the bottom arrow is lax symmetric monoidal.
				Let us first observe that the composition 
				\[
					\LMod_{ U(\g) } \hookrightarrow \LMod_{ \langle 0, U(\g) \rangle }^{ \cpl } \to \LMod_{ U(\g)^{ \PBW } }^{ \cpl }
				\]
				is symmetric monoidal by explicit definition of the tensor product structures (see \cref{lemma:filtered_modules_degree_0}). By definition of internal hom functor, $\IntHom^{ \cpl }_{ U(\g)^{ \PBW }}(\langle 0, A \rangle, -)$ is the right adjoint of 
				\[
					\langle 0, A \rangle \otimes^{ \cpl }_{ \C } -\colon \Mod_{ \C }^{ \cpl } \to \LMod_{ U(\g)^{ \PBW } }^{ \cpl }.
				\]
				We saw in \cref{prop:enrichment_lmod_ug} that it is symmetric monoidal; we then deduce that the functor $\IntHom^{ \cpl }_{ U(\g)^{ \PBW } }(\langle 0, A \rangle, -)$ is lax symmetric monoidal, by \cite[Corollary 7.3.2.7]{Lurie:HA}. Therefore, $\CE^{ \cpl }(\g, -)$ is lax symmetric monoidal.
				We finally conclude by observing that $A$ is the monoidal unit in $\LMod_{ U(\g) }$ and in $\LMod_{ U(\g)^{ \PBW } }^{ \cpl }$, so 
				\[
					\CE^{ \cpl }(\g, -)\colon \LMod_{ U(\g) } \simeq \Mod_{ A }(\LMod_{ U(\g) }) \to \Mod_{ \langle 0, A \rangle }(\LMod_{ U(\g)^{ \PBW } }^{ \cpl }) \to \Mod_{ \CE^{ \cpl }(\g, A) }^{ \cpl }.\qedhere
				\]	
			\end{proof}

			\begin{remark}
				Let us observe that the proof of \cref{prop:CE_lax_monoidal_functor} can be read entirely when $M$ and $N$ have a non-trivial filtration, that is, it implies that 
				\[
					\IntHom^{ \cpl }_{ U(\g)^{ \PBW } }(\langle 0, A \rangle, -)\colon \LMod^{ \cpl }_{ U(\g)^{ \PBW } } \to \Mod^{ \cpl }_{ \C }
				\]
				is lax symmetric monoidal.
				The inclusion $\LMod_{ U(\g) } \hookrightarrow \LMod^{ \cpl }_{ \langle 0, U(\g) \rangle }$ is easily seen to be strong monoidal, as well as the forgetful $ \LMod^{ \cpl }_{ \langle 0, U(\g) \rangle } \to \LMod^{ \cpl }_{ U(\g)^{ \PBW } } $.
			\end{remark}

			\begin{remark}
				\label{remark:adjunction_CE}
				Observe that, by \cref{lemma:standard_bimodule_structure}, $\langle 0, A \rangle$ is canonically a (complete filtered) bimodule over $(U(\g)^{ \PBW }, \IntHom^{ \cpl }_{ U(\g)^{ \PBW } }(\langle 0, A\rangle,\langle 0, A \rangle )^{ \rev })$.
				Denoting 
				\[
					\CE^{ \cpl }(\g, A) \simeq \IntHom^{ \cpl }_{ U(\g)^{ \PBW } }(\langle 0, A \rangle, \langle 0, A \rangle) \in \CAlg^{ \cpl }_{ \C }
				\]
				(so there is no need to reverse the multiplication as above), it gives rise, by \cref{defn:morita_hom_objects}, to an adjunction
				\[
					\adjunction{A \otimes_{ \CE^{ \cpl }(\g, A) }^{ \cpl } -}{\Mod^{ \cpl }_{ \CE^{ \cpl }(\g, A) }}{\LMod_{ U(\g)^{ \PBW } }^{ \cpl }}{\IntHom^{ \cpl }_{ U(\g)^{ \PBW } }(\langle 0, A \rangle, -)}.
				\]
			\end{remark}

			\begin{remark}[Quick construction of the complete Hodge-filtered Chevalley-Eilenberg functor]
				\label{remark:quick_construction_CE}
				Let us start with $\langle 0, A\rangle$ seen as a complete filtered left $U(\g)^{ \PBW }$-module. By formal reasons, it can be seen as a complete filtered $(U(\g)^{ \PBW }, \IntHom^{ \cpl }_{ U(\g)^{ \PBW } }(\langle 0, A \rangle, \langle 0, A \rangle)^{ \rev})$-bimodule. It then induces the Morita adjunction in \cref{remark:adjunction_CE}, and we define the complete Hodge-filtered Chevalley-Eilenberg functor as such right adjoint.
				One only needs to justify why $\IntHom^{ \cpl }_{ U(\g)^{ \PBW } }(\langle 0, A \rangle, \langle 0, A \rangle)$ is a \emph{commutative} algebra, where the final argument of the proof of \cref{prop:CE_lax_monoidal_functor} works.
			\end{remark}
			
			Let us now link our filtered version of the Chevalley-Eilenberg cohomology complex with the classical one. Recall the \emph{realization} functor for filtered objects (see \cref{filtered:constant_functor}), defined as
			\[
				FX \mapsto (FX)_{ \infty } \coloneqq \colim_{ n } F^{ n }X.
			\]
			Recall moreover how filtered hom objects are defined from \cref{filtered:filtered_internal_hom}, and how they are the same both for filtered and for filtered complete objects.
			We can then prove the following proposition.

			\begin{prop}
				\label{prop:CE_lifts_classical}
				The realization of the complete Hodge-filtered Chevalley-Eilenberg coincides with the standard Chevalley-Eilenberg complex, as defined in \cite[§3]{Nuiten:KoszulLieAlgbd}.
			\end{prop}
			\begin{proof}
				We need to prove that
				\[
					\left( \IntHom^{ \cpl }_{ U(\g)^{ \PBW } }(\langle 0, A \rangle, \langle 0, M \rangle) \right)_{ \infty } \simeq \IntHom_{ U(\g) }(A, M) \simeq \CE(\g, M)
				\]
				where we used the enrichment over $\Mod_{ \C }^{ \cpl }$ of $\LMod_{ U(\g)^{ \PBW } }^{ \cpl }$ and the enrichment over $\Mod_{ \C }$ of $\LMod_{ U(\g) }$, respectively.
				The natural map of positively filtered algebras 
				\[
					U(\g)^{ \PBW } \to \langle 0, U(\g) \rangle
				\]
				gives a tensor-forgetful adjunction 
				\[
					\adjunction{\langle 0, U(\g) \rangle \otimes^{ \cpl }_{ U(\g)^{ \PBW } } -}{\LMod_{ U(\g)^{ \PBW } }^{ \cpl }}{\LMod_{ \langle 0, U(\g) \rangle }^{ \cpl }}{\oblv}.
				\]
				It is easily seen that the left adjoint is $\Mod^{ \cpl }_{ \C }$-linear (although only op-lax symmetric monoidal), so that this is a $\Mod^{ \cpl }_{ \C }$-enriched adjunction and we can apply \cref{lemma:enriched_adjunctions} to obtain an equivalence of complete filtered $\C$-modules
				\[
					 \IntHom^{ \cpl }_{ \langle 0, U(\g) \rangle }\left( \langle 0, U(\g) \rangle \otimes^{ \cpl }_{ U(\g)^{ \PBW } } \langle 0, A \rangle, \langle 0, M \rangle \right) \simeq \IntHom^{ \cpl }_{ U(\g)^{ \PBW } }(\langle 0, A \rangle, \langle 0, M \rangle).				
				 \]
				Let us now observe that all filtered objects involved (\emph{inside} the hom, not the hom object itself) are positively-filtered and that (relative) tensor products between them are therefore positively-filtered as well; this is immediate in the absolute case since $\Mod_{ \C }^{ \Fil, \geq 0}$ is a symmetric monoidal subcategory of $\Mod_{ \C }^{ \Fil }$ (and of $\Mod_{ \C }^{ \cpl }$ as well, since we don't need to complete the tensor product in this case). For the relative case, we can just write the bar resolution (whose colimit is the relative tensor product) and observe that all terms in the simplicial diagram are positively filtered and that $\Mod_{ \C }^{ \Fil, \geq 0 } \hookrightarrow \Mod_{ \C }^{ \Fil }$ preserves colimits.

				Using that $\LMod_{ \langle 0, U(\g) \rangle }(\Mod_{ \C }^{ \Fil, \geq 0 }) \simeq (\LMod_{ U(\g) })^{ \Fil, \geq 0 }$ and the formulas for filtered-hom in \cref{filtered:internal_hom_complete}, we finally can write (for positive $n$):
				\begin{gather*}
					F^{ n }\IntHom^{ \Fil }_{ \langle 0, U(\g) \rangle }\left( \langle 0, U(\g) \rangle \otimes^{ \Fil }_{ U(\g)^{ \PBW } } \langle 0, A \rangle, \langle 0, M \rangle \right) \simeq \\
					\simeq \int_{ m \geq 0 } \IntHom_{ U(\g) }(F^{ m }(\langle 0, U(\g) \rangle \otimes^{ \Fil }_{ U(\g)^{ \PBW } } \langle 0, A \rangle), M),
				\end{gather*}
				where we considered only positive $m$'s in the end formula because the negative filtration terms of the tensor product are zero. 
				We now pull inside the end as a coend (hence colimit) to obtain 
				\[
					\IntHom_{ U(\g) }\left( (\langle 0, U(\g) \rangle \otimes^{ \Fil }_{ U(\g)^{ \PBW }  } \langle 0, A \rangle )_{ \infty }, M \right) \simeq \IntHom_{ U(\g) }(A, M),
				\]
				where we used that the realization functor is strong monoidal for the filtered tensor product (and that $U(\g)^{ \PBW }_{ \infty} \simeq U(\g)$). This implies that the same equivalence holds at $\infty$, i.e.\ passing to colimit in $n$. This concludes the proof.
			\end{proof}

			Finally, let us relate the definition of derived left crystals given in \cref{defn:toen_vezzosi_derived_d_modules} with constant modules over this filtered version of Chevalley-Eilenberg (as in \cref{defn:constant_graded_mixed_modules}).
			Given our base $X = \Spec A$ we can associate with it two complete filtered commutative algebras:
			\[
				\CE^{ \cpl }(\T_{ A }, A), \qquad \DR(A)
			\]
			and the question we are interested in now is when they are equivalent.	

			\begin{lemma}
				\label{lemma:CE_foliation_algebra}

				Let $\g \in \LieAlgbd_{ A }$ be a dg-Lie algebroid over $A$ whose underlying $A$-module is perfect. Then 
				\[
					\gr(\CE^{ \cpl }(\g, A)) \simeq \Sym_{ A }^{ \gr }(\g^{ \vee }[-1]),
				\]
				where $\g^{ \vee }[-1]$ has weight $-1$.
			\end{lemma}
			\begin{proof}
				This is proved in \cref{prop:CE_factors_constant_modules} (which does not depend on this statement, from a logical standpoint).
			\end{proof}

			In this case, we have an essentially unique filtered algebra map
			\[
				\DR(A) \longrightarrow \CE^{ \cpl }(\g, A),
			\]
			which comes by adjunction from $A \simeq \gr(\CE^{ \cpl }(\g, A))(0)$, see \cref{prop:deRham_graded_commutative_algebra_structure}.

			The following lemma, proved in \cite[Lemma A.4]{CalaqueGrivaux:FormalModuli}, provides an answer to the initial question (which was to identify the de Rham algebra as an Hodge-filtered Chevalley-Eilenberg complex).
			\begin{lemma}
				\label{lemma:relation_CE_DR}

				The morphism of complete filtered algebras $\DR(A) \to \CE^{ \cpl }(\T_{ A }, A)$ is an equivalence\footnote{A graded equivalence, meaning an equivalence in $\CAlg^{ \cpl }_{ \C }$.} whenever $\L_{ A }$ is a perfect $A$-module.
			\end{lemma}

			Recall that whenever $A$ is locally of finite presentation then $\L_{ A } \in \Perf(A)$ by \cref{thm:locally_finite_presentation_perfect_cotangent_complex}, so that the lemma above applies.
			This means that derived left crystals (introduced in \cref{defn:toen_vezzosi_derived_d_modules}) are then identified with constant modules over $\CE^{ \cpl }(\T_{ A }, A)$ (as defined in \cref{defn:constant_graded_mixed_modules}).

%% file: proof.tex
	 \section{Proof of the theorem}
			\label{subsection:main_proof}
			
			The theorem we want to prove is the following.
			\begin{thm}
				\label{thm:main_thm}
				Let $A$ be a connective commutative algebra and $\g \in \LieAlgbd_{ A }$ with an underlying perfect $A$-module. Then, the complete Hodge-filtered Chevalley-Eilenberg cohomology functor induces an equivalence of symmetric monoidal $\infty$-categories
					\[
						\CE^{\cpl}(\g, -)\colon \LMod_{ U(\g) } \stackrel{\sim}{\longrightarrow} \Mod^{ \const }_{ \CE^{ \cpl }(\g, A) }.
					\]
			\end{thm}
				
			Plugging in $\infcatname{C} = \Mod^{ \cpl }_{ \C }$, $B = U(\g)^{ \PBW }$ and $M = \langle 0, A \rangle \in \LMod_{ U(\g)^{ \PBW } }^{ \cpl }$, \cref{lemma:standard_bimodule_structure} gives us a canonical $(U(\g)^{ \PBW }, \CE^{ \cpl }(\g, A))$-bimodule structure on $\langle 0, A \rangle$, where we observed that
			\[
				\IntHom_{ U(\g)^{ \PBW } }^{ \cpl }(\langle 0, A \rangle, \langle 0, A \rangle) \simeq \IntHom_{ U(\g)^{ \PBW } }^{ \cpl }(\langle 0, A \rangle, \langle 0, A \rangle)^{ \mathrm{rev} } \simeq \CE^{ \cpl }(\g, A)
			\]
			in $\CAlg^{ \cpl }$.

			The first step is to prove that, when $\g \in \LieAlgbd_{ A }$ has an underlying perfect $A$-module, then the complete filtered Chevalley-Eilenberg complexes are all \emph{constant} modules over $\CE^{ \cpl }(\g, A)$.
			\begin{prop}
				\label{prop:CE_factors_constant_modules}
				Let $\g \in \LieAlgbd_{ A }$ with an underlying perfect $A$-module. Then there exists a canonical factorization
				\begin{diag}
					& \Mod^{\const}_{\CE^{\cpl}(\g, A)} \ar[d, hook] \\
					\LMod_{U(\g)} \ar[ur, dashed] \ar[r, "\CE^{\cpl}(\g{,}-)"] & \Mod^{\cpl}_{\CE^{\cpl}(\g, A)},
				\end{diag}
				which we will denote $\CE^{ \cpl }(\g, -)$ by abuse of notation (context should clarify which functor we mean exactly).
			\end{prop}
			\begin{proof}
				Since constant complete filtered modules (see \cref{defn:constant_graded_mixed_modules}) are a full subcategory of complete filtered modules, it suffices to verify that $\CE^{ \cpl }(\g, M) \in \Mod^{ \const}_{ \CE^{ \cpl }(\g, A) }$ for any $\g$-representation $M$. Using \cref{coroll:pbw_cofibrant_lie_algebroids} and that $\gr\colon \LMod_{ U(\g)^{ \PBW } }^{ \cpl } \to \LMod_{ \gr(U(\g)^{ \PBW }) }^{ \gr }$ is strong closed by \cref{prop:gr_closed} we obtain
				\[
					\gr(\CE^{ \cpl }(\g, M)) \simeq \IntHom^{ \gr }_{ \gr(U(\g)^{ \PBW }) }(A, M) \simeq \IntHom^{ \gr }_{ \Sym^{ \gr }_{ A }(\g) }(A, M).
				\]
				Observe now that $M_{ 0 }$ is a (left) graded module over $\Sym^{ \gr }_{ A }(\g)$ whose underlying object is concentrated in weight $0$ (remember that $\g$ has weight $1$).
				By \cref{lemma:graded_modules_weight_0}, we deduce that $M$ is a $\Sym^{ \gr }_{ A }(\g)(0) \simeq A$-module; that is, the $\Sym^{ \gr }_{ A }(\g)$-module structure on $M$ is obtained by restricting along the augmentation map $\Sym^{ \gr }_{ A }(\g) \to A$.
				In other words, $M$ is in the essential image of
				\[
					\Mod^{ \gr }_{ A } \to \Mod^{ \gr }_{ \Sym^{ \gr }_{ A }(\g) }
				\]
				which is just the ``forgetful'' along the natural augmentation of graded algebras $\Sym^{ \gr}_{ A }(\g) \to A$.
				This latter functor has a left adjoint by extension of scalars which is $\Mod^{ \gr }_{ \C }$-linear. By \cref{lemma:enriched_adjunctions} we can write 
				\[
					\gr(\CE^{ \cpl }(\g, M))  \simeq \IntHom^{ \gr }_{ A }(A \otimes^{ \gr }_{ \Sym^{ \gr }_{ A }(\g) } A, M).						
				\]
				Finally, observing that $A \simeq \Sym^{ \gr }_{ A }(0)$, that $\Sym^{ \gr }_{ A }(-)\colon \Mod^{ \gr }_{ A } \to \CAlg^{ \gr }_{ A/- }$ is a left adjoint and that the relative tensor products of $\Sym^{ \gr }_{ A }(\g)$-commutative algebras is their categorical coproduct, we get 
				\[
					\gr(\CE^{ \cpl }(\g, M)) \simeq \IntHom_{ A }^{ \gr }(\Sym_{ A }^{ \gr }(\g[1]), M) \simeq \Sym^{ \gr }_{ A }(\g^{ \vee }[-1]) \otimes^{ \gr }_{ A } M,
				\]
				where the last equivalence is a simple computation of internal hom in $(\Mod_{ A })^{ \gr }$, where we used that we are in char $0$ (for the Sym) and that $\g$ is dualizable (i.e.\ perfect) as an $A$-module with $A$-dual $\g^{ \vee }$ (which is in weight $-1$).
			\end{proof}

			Let us now record here an important property, which will be of use later.
			\begin{prop}
				\label{prop:CE_const_is_conservative}
				If $\g \in \LieAlgbd_{ A }$ has an underlying perfect $A$-module, then the functor 
				\[
					\CE^{ \cpl }(\g, -)\colon \LMod_{ U(\g) } \to \Mod^{ \cpl }_{ \CE^{ \cpl }(\g, A) }
				\]
				is conservative.
			\end{prop}
			\begin{proof}
				By \cref{prop:CE_factors_constant_modules}, the above functor is the composite of 
				\[
					\LMod_{ U(\g) } \stackrel{\CE^{ \cpl }(\g, -)}{\longrightarrow} \Mod^{ \const }_{ \CE^{ \cpl }(\g, A) } \hookrightarrow \Mod^{ \cpl }_{ \CE^{ \cpl }(\g, A) },
				\]
				so we can equivalently prove that the first functor (to constant modules) is conservative.
				By \cref{prop:constant_modules_forgetful_conservative} we can prove, equivalently, that the composite
				\[
					\LMod_{ U(\g) } \stackrel{\CE^{ \cpl }(\g, -)}{\longrightarrow} \Mod_{\CE^{ \cpl }(\g, A)}^{ \const } \stackrel{\gr(-)(0)}{\longrightarrow} \Mod_{ A }
				\]
				is conservative (since $\gr(\CE^{ \cpl }(\g, A))(0) \simeq A$ as commutative algebras).
				It is then easily seen, again by \cref{prop:CE_factors_constant_modules}, that this is exactly the forgetful functor along the algebra map $A \to U(\g)$ \[
					\LMod_{ U(\g) } \to \Mod_{ A }.
				\]
				This is conservative, therefore we conclude.
			\end{proof}

			\begin{remark}
				Let us observe how crucial the ``constant'' hypothesis is for the result of \cref{prop:CE_const_is_conservative}: we will never be able to find the appropriate ``section'' for general non-constant $\CE^{ \cpl }(\g, A)$-modules, since they are not recovered by simply evaluating in $0$.
			\end{remark}

			Let us observe that, in particular, \cref{prop:CE_factors_constant_modules} implies 
			\[
				\gr(\CE^{ \cpl }(\g, A)) \simeq \Sym^{ \gr }_{ A }(\g^{ \vee }[-1]),
			\]
			with $\g[-1]$ in weight $-1$.
			That is, if $\g$ is perfect and of Tor-amplitude $[0, +\infty[$, then $\CE^{ \cpl }(\g, A)$ is a foliation-like algebra (compare with \cref{lemma:CE_foliation_algebra}).

			Let us now link what we do with \cite[Theorem 3.2.1]{TV:AlgFoliations2}, which identifies derived left crystals (which, thanks to \cref{lemma:relation_CE_DR}, can be thought as constant modules over $\CE^{ \cpl }(\T_{ A }, A)$) with (unfiltered) left modules over a certain algebra.
			Since our final theorem (\cref{thm:main_thm}) also identifies derived left crystals with unfiltered left modules over a specific algebra (namely the universal enveloping algebra of the final Lie algebroid), we will get for free the corollary that these two algebras are Morita equivalent.
			We believe, though, that they are already equivalent as algebras (which is a stronger statement). We will explain it here, although it is not essential for our goal.
			The first thing we need to do is to identify
			\[
				S(\g) \coloneqq \IntHom^{ \cpl }_{ \CE^{ \cpl }(\g, A) }(\langle 0, A \rangle, \langle 0, A \rangle) \in \Alg(\Mod^{ \cpl }_{ \C }),
			\]
			where we used the canonical $\Mod^{ \cpl }_{ \C}$-enrichment of $\Mod^{ \cpl }_{ \CE^{ \cpl }(\g, A) }$ (see \cref{lemma:enrichment_mod_A}).
			\begin{prop}
				\label{prop:Ug_as_endomorphism_CE}
				Assume $\g \in \LieAlgbd_{ A }$ has an underlying perfect $A$-module. Then there exists a canonical equivalence 
				\[
					U(\g)^{ \PBW } \simeq \IntHom^{ \cpl }_{ \CE^{ \cpl }(\g, A) }(\langle 0, A \rangle, \langle 0, A \rangle)
				\]
				in $\Alg(\Mod^{ \cpl }_{ \C })$.

			\end{prop}
			\begin{proof}
				First, let us build a morphism
				\[
					U(\g)^{ \PBW } \dashrightarrow \IntHom^{ \cpl }_{\CE^{ \cpl }(\g, A)}(\langle 0, A \rangle, \langle 0, A \rangle)
				\]
				in $\Alg^{ \cpl }_{ \C }$.
				This is equivalent to giving a lift (along $\Lambda^{ 2 }_{ 2 } \hookrightarrow \Delta^{ 2 }$) of the diagram
				\begin{diag}
					& (U(\g)^{\PBW}, \CE^{\cpl}(\g, A))-\BiMod^{\cpl} \ar[d] \\
					\Delta^0 \ar[r, "\langle 0{,} A \rangle"] \ar[ur, dashed] & \Mod^{\cpl}_{\CE^{\cpl}(\g, A)}.
				\end{diag}
				By construction, see \cref{lemma:standard_bimodule_structure}, the $\CE^{ \cpl }(\g, A)$-structure on $\langle 0, A \rangle$ upgrades to
				\[
					\langle 0, A \rangle \in (U(\g), \CE^{ \cpl }(\g, A))-\BiMod^{ \cpl }
				\]
				so we already have a compatible lift, and therefore a morphism of algebras as in the beginning.
			    We have now a map of complete filtered algebras 
				\[
					U(\g)^{ \PBW } \to \IntHom^{ \cpl }_{ \CE^{ \cpl }(\g, A) }(\langle 0, A \rangle, \langle 0, A \rangle)
				\]
				and we want to prove it is an isomorphism.
				The forgetful functor $\oblv_{ \Alg }\colon \Alg^{ \cpl }_{ \C } \to \Mod^{ \cpl }_{ \C }$ is conservative, so we can just prove that it is an isomorphism of complete filtered modules.
				We also know that $\gr\colon \Mod^{ \cpl}_{ \CE^{ \cpl }(\g, A) } \to \Mod^{ \gr }_{ \CE^{ \gr }(\g, A) }$ is strong closed (by \cref{prop:gr_closed}) and conservative (since forgetting module structure is conservative too), therefore it suffices to prove that 
				\[
					\gr(U(\g)^{ \PBW }) \to \gr\left(\IntHom^{ \cpl }_{ \CE^{ \cpl }(\g, A) }(\langle 0, A \rangle, \langle 0, A \rangle)\right) \simeq \IntHom^{ \gr }_{ \CE^{ \gr }(\g, A) }(A, A)
				\]
				is an equivalence of graded modules.
				For the left-hand side, thanks to the PBW theorem (\cref{coroll:pbw_cofibrant_lie_algebroids}), we have 
				\[
					\gr(U(\g)^{ \PBW }) \simeq \Sym^{ \gr }_{ A }(\g) \in \CAlg^{ \gr }_{ A },
				\]
				with $\g$ of weight $1$.
				For the right-hand side, we can just use the computation of \cref{prop:CE_factors_constant_modules} to say
				\[
					\gr\left(\IntHom^{ \cpl }_{ \CE^{ \cpl }(\g, A) }(\langle 0, A \rangle, \langle 0, A \rangle)\right) \simeq \IntHom^{ \gr }_{ \Sym^{ \gr }_{ A }(\g^{ \vee }[-1]) }(A, A) \in \Mod^{ \gr }_{ \C },
				\]
				where $\g^{ \vee }[-1]$ has weight $-1$.
				As before, by \cref{lemma:graded_modules_weight_0}, the action of $\Sym^{ \leq -1 }_{ A }(\g^{ \vee }[-1])$ on $A$ is zero so that $A$ is in the essential image of 
				\[
					\Mod^{ \gr }_{ A } \to \Mod^{ \gr }_{ \Sym^{ \gr }_{ A }(\g^{ \vee }[-1]) }
				\]
				induced by the ``forgetful'' along the natural augmentation map of $\Sym^{ \gr }_{ A }(\g^{ \vee }[-1])$ towards $A$.
				Using the tensor-forgetful adjunction (which is enriched, thanks to \cref{lemma:enriched_adjunctions}) and similar computation as in the proof of \cref{prop:CE_factors_constant_modules}, we obtain
				\[
					\IntHom_{ A }^{ \gr }(A \otimes^{ \gr }_{ \Sym^{ \gr }_{ A }(\g^{ \vee }[-1] ) } A, A ) \simeq \IntHom_{ A }^{ \gr }(\Sym^{ \gr }_{ A }(\g^{ \vee }), A)
				\]
				and being in characteristic $0$ and using the dualizability of $\g$ as an $A$-module we finally get to 
				\[
					\Sym_{ A }^{ \gr }(\g).
				\]
				Therefore the map $\gr(U(\g)^{ \PBW }) \to \gr(\IntHom^{ \cpl }_{ \CE^{ \cpl }(\g, A) }(\langle 0, A \rangle, \langle 0, A \rangle))$ is an equivalence of graded objects (observe that it is the correct map since it satisfies the same universal property in the graded context, accounting for the $(\gr(U(\g)^{ \PBW }), \gr(\CE^{ \cpl }(\g, A)))$-bimodule structure of $A$, using \cref{lemma:graded_modules_weight_0}), and so we can conclude. \qedhere
			\end{proof}

			Let us now get back to work to prove our main result (meaning that the statements from now on will be useful in our final proof).

			\begin{prop}
				\label{prop:left_adjoint_CE}
				The (extended) complete filtered Chevalley-Eilenberg fits in the adjunction
				\[
					\adjunction{A \otimes_{ \CE^{ \cpl }(\g, A) }^{ \cpl } -}{\Mod^{ \cpl }_{ \CE^{ \cpl }(\g, A) }}{\LMod_{ U(\g)^{ \PBW } }^{ \cpl }}{\CE^{ \cpl }(\g, -)}.
				\]
				If $\g$ has an underlying perfect $A$-module, then this adjunction can be restricted to
				\[
					\adjunction{A \otimes_{ \CE^{ \cpl }(\g, A) }^{ \cpl } -}{\Mod^{ \const }_{ \CE^{ \cpl }(\g, A) }}{\LMod_{ U(\g) }}{\CE^{ \cpl }(\g, \langle 0, - \rangle)}.
				\]
			\end{prop}
			\begin{proof}
			    The first adjunction comes from \cref{remark:adjunction_CE}. Let us focus then on the second part of the statement.
				By \cref{prop:CE_factors_constant_modules}, if $\g$ has an underlying perfect $A$-module then we know that 
				\[
					\CE^{ \cpl }(\g, -)\colon \LMod_{ U(\g) } \to \Mod^{ \const }_{ \CE^{ \cpl }(\g, A) }
				\]
				is well-defined, meaning that we have a commutative diagram
				\begin{diag}
						\LMod^{\cpl}_{U(\g)^{\PBW}} \ar[r, "\CE^{\cpl}(\g{,} -)"] & \Mod^{\cpl}_{\CE^{\cpl}(\g, A)} \\
						\LMod_{U(\g)} \ar[u, hook] \ar[r, "\CE^{\cpl}(\g{,} -)"] & \Mod^{\const}_{\CE^{\cpl}(\g, A)} \ar[u, hook].
				\end{diag}
				Observe that the left vertical arrow is fully faithful thanks to the first part of \cref{lemma:filtered_modules_degree_0}, since $(U(\g)^{ \PBW })_{ \infty } \simeq U(\g)$. To prove that the adjunction in the statement restricts to an adjunction between unfiltered $\g$-representations and constant $\CE^{ \cpl }(\g, A)$-modules, it suffices to prove that the functor $\langle 0, A \rangle \otimes^{ \cpl}_{ \CE^{ \cpl }(\g, A) } -\colon \Mod^{ \cpl }_{ \CE^{ \cpl }(\g, A) } \to \LMod^{ \cpl }_{ U(\g)^{ \PBW } }$ restricts to
				\[
					\langle 0, A \rangle \otimes^{\cpl}_{ \CE^{ \cpl }(\g, A) } - \colon \Mod^{ \const }_{ \CE^{ \cpl }(\g, A) } \to \LMod_{ U(\g) }.
				\]
				That is, it needs to send constant modules over $\CE^{ \cpl }(\g, A)$ to unfiltered $\g$-representations. Thanks to \cref{coroll:filtered_complete_modules_gr_0}, we can identify the latter $\infty$-category with the subcategory of $\LMod_{ U(\g)^{ \PBW } }^{ \cpl }$ spanned by those objects whose associated graded is concentrated in weight $0$.
				Let now $M$ be a constant complete filtered $\CE^{ \cpl }(\g, A)$-module; this means that
				\[
					\gr(M) \simeq \Sym^{ \gr }_{ A }(\g^{ \vee }[-1]) \otimes^{ \gr }_{ A } \gr(M)(0),
				\]
				where $\g^{ \vee }[-1]$ has weight $-1$. Let us observe that
				\[
					\gr\left(\langle 0, A \rangle \otimes^{ \cpl }_{ \CE^{ \cpl }(\g, A) } M\right) \simeq A \otimes^{ \gr }_{ \Sym^{ \gr }_{ A }(\g^{ \vee }[-1]) } \left(\Sym^{ \gr }_{ A }(\g^{ \vee }[-1]) \otimes^{ \gr }_{ A } \gr(M)(0)\right) \simeq \gr(M)(0).
				\]
				That is, the left $U(\g)^{ \PBW }$-module $\langle 0, A \rangle \otimes^{ \cpl }_{ \CE^{ \cpl }(\g, A) } M$ has associated graded only in weight $0$ and therefore it comes from a standard left $U(\g)$-module.
				We proved that we have an analogous diagram
				\begin{diag}
					\LMod^{\cpl}_{U(\g)^{\PBW}} & & \ar[ll, "\langle 0{,} A \rangle \otimes^{\cpl}_{\CE^{\cpl}(\g{,} A)} -"]  \Mod^{\cpl}_{\CE^{\cpl}(\g, A)} \\
					\LMod_{U(\g)} \ar[u, hook] & & \ar[ll, "\langle 0{,} A \rangle \otimes^{\cpl}_{\CE^{\cpl}(\g{,} A)} -"]  \Mod^{\const}_{\CE^{\cpl}(\g, A)} \ar[u, hook].
				\end{diag}
					
				We can conclude by observing that the general adjunction at the top level then restricts to the bottom full subcategories (the unit, counit and all the other homotopy coherent data easily restrict), which allows us to finally say that the left adjoint is then simply tensoring with $\langle 0, A \rangle$ seen as a bimodule.
			\end{proof}

			\begin{remark}
				\label{remark:A_tens_CE_forgetful}
				Observe that, using the identification of $\LMod_{ U(\g) } \hookrightarrow \LMod^{ \cpl }_{ U(\g)^{ \PBW } }$ coming from \cref{lemma:filtered_modules_degree_0}, we see that the composition
				\[
					\Mod^{ \const }_{ \CE^{ \cpl }(\g, A) } \stackrel{A \otimes^{ \cpl }_{ \CE^{ \cpl }(\g, A) } -}{\longrightarrow} \LMod_{ U(\g) } \stackrel{\oblv}{\longrightarrow} \Mod_{ A }
				\]
				is simply taking associated graded in weight $0$, $M \mapsto \gr(M)(0)$. 
				The left action of $U(\g)$, though, is highly non-trivial (that is, in the model of mixed graded modules, we only know that $\gr(M)(0)$ is quasi-isomorphic to another mixed graded module equipped with a left action of $U(\g)^{ \PBW }$).
				We will think of this composition as taking the quasi-coherent sheaf underlying the corresponding derived left crystal.
			\end{remark}

			We can now finally conclude the proof of \cref{thm:main_thm}.

			\begin{prop}
				\label{prop:final_proof_main_thm}

				Let $\g \in \LieAlgbd_{ A }$ be a dg-Lie algebroid with a perfect $A$-module underlying. Then the adjunction
				\[
					\adjunction{A \otimes_{ \CE^{ \cpl }(\g, A) }^{ \cpl } -}{\Mod^{ \const }_{ \CE^{ \cpl }(\g, A) }}{\LMod_{ U(\g) }}{\CE^{ \cpl }(\g, \langle 0, -\rangle)}.
				\]
				is an adjoint equivalence.
				Moreover, both functors have natural symmetric monoidal structures, so that they induce equivalences of symmetric monoidal $\infty$-categories.
			\end{prop}
			\begin{proof}
				It suffices to prove that the functor 
					\[
						\langle 0, A \rangle \otimes^{ \cpl }_{ \CE^{ \cpl }(\g, A) } -\colon \Mod^{\const }_{ \CE^{ \cpl }(\g, A) } \to \LMod_{ U(\g) }
					\]
				is an equivalence. Let us start with fully faithfulness, for which it suffices to prove that the unit map
				\[
					M \to \CE^{ \cpl }(\g, \langle 0, A \rangle \otimes^{ \cpl }_{ \CE^{ \cpl }(\g, A) } M)
				\]
				is an equivalence for each $M \in \Mod^{ \const }_{ \CE^{ \cpl }(\g, A) }$.
				Firstly let us recall that the proof of \cref{prop:left_adjoint_CE} gives us 
				\[
					\gr\left( \langle 0, A \rangle \otimes^{ \cpl }_{ \CE^{ \cpl }(\g, A) } M \right) \simeq \gr(M)(0)_{ 0 }.
				\]
				The ``associated graded in weight $0$'' functor 
				\[
					\gr(-)(0)\colon \Mod^{ \const }_{ \CE^{ \cpl }(\g, A) } \to \Mod_{ A }
				\]
				is conservative by \cref{prop:constant_modules_forgetful_conservative}.
				Applying it to the map above, by \cref{prop:CE_factors_constant_modules}, we obtain
				\begin{gather*}
					\gr(M)(0) \to \gr\left(\CE^{ \cpl }(\g, \langle 0, A \rangle \otimes^{ \cpl }_{ \CE^{ \cpl }(\g, A) } M)\right)(0) \simeq \left(\Sym^{ \gr }_{ A }(\g^{ \vee }[-1]) \otimes^{ \gr }_{ A } \gr(M)(0) \right)(0) \\
					\simeq \gr(M)(0)
				\end{gather*}
				which is the identity of $\gr(M)(0)$ as an $A$-module (recall that $\g^{ \vee }[-1]$ has weight $-1$). This proves the fully faithfullness.

				To conclude, it suffices to prove that the right adjoint functor
				\[
					\CE^{ \cpl }(\g, \langle 0, - \rangle)\colon \LMod_{ U(\g) } \to \Mod^{ \const }_{ \CE^{ \cpl }(\g, A) }
				\]
				is conservative. This is the content of \cref{prop:CE_const_is_conservative}.
				Observe, finally, that the (both the full and the restricted) functor $A \otimes^{ \cpl }_{ \CE^{ \cpl }(\g, A) }-$ is symmetric monoidal. This implies that its inverse 
				\[
					\CE^{ \cpl }(\g, \langle 0, - \rangle)\colon \LMod_{ U(\g) } \stackrel{\sim}{\longrightarrow} \Mod^{ \const }_{ \CE^{ \cpl }(\g, A) }
				\]
				has a canonical symmetric monoidal structure as well; that is, it is an equivalence in $\CAlg(\Pr^{ L, st })$.
			\end{proof}
			\begin{remark}
				\label{remark:full_CE_only_lax_monoidal}
				Observe that the ``full'' functor 
				\[
					\CE^{ \cpl }(\g, -)\colon \LMod_{ U(\g)^{ \PBW } }^{ \cpl } \to \Mod^{ \cpl }_{ \CE^{ \cpl }(\g, A) }
				\]
				is only lax symmetric monoidal.
			\end{remark}
			
			We proved \cref{thm:main_thm}. We are now able to answer the question asked in.
			\begin{coroll}
				\label{coroll:left_adjoint_gr_0}
				Let $\g \in \LieAlgbd_{ A }$ with an underlying perfect $A$-module. There is a commutative diagram
				\begin{diag}
					\LMod_{U(\g)} \ar[rr, "\simeq"] \ar[dr, "\oblv_{\g}"'] & & \Mod^{ \const }_{ \CE^{ \cpl }(\g, A) } \ar[dl, "\gr(-)(0)"]  \\
					 & \Mod_{ A } &,
				\end{diag}
				where the top arrow is $\CE^{ \cpl }(\g, -)$. The left adjoint of $\gr(-)(0)$ is then given by the composite
				\[
					\Mod_{ A } \stackrel{U(\g) \otimes_{ A } -}{\longrightarrow} \LMod_{ U(\g) } \stackrel{\CE^{ \cpl }(\g, -)}{\longrightarrow} \Mod^{ \const }_{ \CE^{ \cpl }(\g, A) },
				\]
				sending $A$ to $\CE^{ \cpl }(\g, U(\g))$.
			\end{coroll}
			
		\section{Equivalence of the different definitions of derived D-Modules}
			\label{subsection:derived_D_mod}

			We can finally deduce that on a bounded derived affine scheme $X = \Spec A$, all the definitions of derived D-modules mentioned in this article are equivalent.
			Let us state it clearly.

			\begin{coroll}
				\label{coroll:equivalence_derived_d_mod}

				Let $A$ be a connective commutative algebra which is bounded (or ``eventually coconnective'' using the terminology in \cite{GR:DAG2}) and locally of finite presentation (over $\C$).
				We then have (symmetric monoidal) equivalences of $\infty$-categories
				\[
					\infcatname{D}^{ \der }(\Spec A) \simeq \LMod_{ U(\T_{ A }) } \simeq \Mod^{ \const }_{ \CE^{ \cpl }(\T_{ A }, A) } \simeq \Mod^{ \const }_{ \DR(A) }.
				\]
			\end{coroll}
			\begin{proof}
				From left to right: the first arrow is the equivalence of \cref{coroll:nuiten_is_beraldo_monoidal}, then we have \cref{thm:main_thm} applied to $\g = \T_{ A }$ and, finally, the last arrow comes from the equivalence $\CE^{ \cpl }(\T_{ A }, A) \simeq \DR(A)$ given by \cref{lemma:relation_CE_DR}.
			\end{proof}

			This allows us to use results on derived D-modules in any incarnation.
			For example, the following corollary holds. 

			\begin{coroll}
				\label{coroll:derived_D_modules_bounded_affine}
				Let $A \in \CAlg^{ \afp }_{ \C }$ be bounded. There exists a natural equivalence
				\[
					\infcatname{D}^{ \der }(\Spec A) \stackrel{\sim}{\longrightarrow} \Crys^{ l }(\Spec A) \coloneqq \Qcoh( (\Spec A)_{ \dR } ).
				\]
				That is, derived (left) D-modules on $\Spec A$ are actually equivalent to classical left D-modules on $\Spec A$.
			\end{coroll}
			\begin{proof}
				This comes from \cite[Corollary 4.3.13]{Ber:DerivedDMod}.
			\end{proof}	

			This implies, for example, that $\LMod_{ U(\T_{ A }) }$ is not sensible to the derived structure on $A$ (since $\Qcoh( (\Spec A)_{ \dR } )$ is not); we wouldn't know how to prove such a claim without using the two above corollaries.